\def\JulinPaperOneFlag{1}
\ifnum\pdfstrcmp{\jobname}{Julin_Harnack_Paper_I_Smooth}=0
  \def\JulinPaperOneFlag{1}
\fi
\ifnum\pdfstrcmp{\jobname}{Julin_Harnack_Paper_II_Extensions}=0
  \def\JulinPaperTwoFlag{1}
\fi
\newif\ifJulinPaperOne
\newif\ifJulinPaperTwo
\ifdefined\JulinPaperOneFlag\JulinPaperOnetrue\fi
\ifdefined\JulinPaperTwoFlag\JulinPaperTwotrue\fi
\ifJulinPaperOne\else
  \ifJulinPaperTwo\else\JulinPaperOnetrue\fi
\fi
\def\JulinEndDocument{\end{document}}
\documentclass[11pt,reqno]{amsart}
\usepackage[titletoc]{appendix}
\usepackage{amssymb}
\numberwithin{equation}{section}
\usepackage{epsfig}
\usepackage{color}
\usepackage{bigints}
\usepackage{mathrsfs}
\usepackage{url}
\AtBeginDocument{\color{black}}
\usepackage{bm}
\usepackage{enumerate}
\usepackage[euler]{textgreek}
\usepackage{soul}
\usepackage{pdfcomment}
\usepackage{environ}
\usepackage{xstring}
\usepackage{tikz}

\newtheorem{thm}{Theorem}[section]
\newtheorem{cor}[thm]{Corollary}
\newtheorem{lem}[thm]{Lemma}
\newtheorem{prop}[thm]{Proposition}
\theoremstyle{definition}
\newtheorem{defn}[thm]{Definition}
\newtheorem{rem}[thm]{Remark}

\newcommand{\diver}{\operatorname{div}}
\newcommand{\Ha}{\mathcal{H}}

\newif\ifShowLabels
\newif\ifShowLabelTooltips
\ShowLabelTooltipstrue
\newcommand{\BoxedTagTooltip}[2]{%
     \ifShowLabelTooltips
       \pdftooltip{\fbox{{\ttfamily\small\bfseries\detokenize{#1}}}}{#2}%
     \else
       \fbox{{\ttfamily\small\bfseries\detokenize{#1}}}%
     \fi}
\newcommand{\LabelTagTooltip}[1]{\BoxedTagTooltip{#1}{#1}}
\newif\ifShowReferenceTooltips
\ShowReferenceTooltipsfalse
\newif\ifShowFormulaTooltips
\ShowFormulaTooltipsfalse
\makeatletter
\renewcommand{\pc@annot@tooltip}{%
  /TU (\pc@pdfenc@contents)\space%
  /Contents (\pc@pdfenc@contents)\space%
  /T (tooltip \thezref@unique)\space%
  /C [ ]\space%
  /FT/Btn\space%
  /F 768\space%
  /Ff 65536\space%
  /H/N\space%
  /BS << /W 0 >>\space%
}%
\let\PrismFormulaTooltipText\@empty
\let\PrismOriginalTagform\tagform@
\newcommand{\PrismSetFormulaTooltip}[1]{%
     \begingroup
       \edef\PrismRawFormulaTooltip{\detokenize{#1}}%
       \StrSubstitute{\PrismRawFormulaTooltip}{\@backslashchar}{/}[\PrismCleanFormulaTooltip]%
       \global\let\PrismFormulaTooltipText\PrismCleanFormulaTooltip
     \endgroup}
\renewcommand{\tagform@}[1]{%
     \ifShowFormulaTooltips
       \ifx\PrismFormulaTooltipText\@empty
         \ifx\df@label\@empty
           \PrismOriginalTagform{#1}%
         \else
           \ifcsname PrismFormulaTooltip@\df@label\endcsname
             \pdftooltip{\PrismOriginalTagform{#1}}{\csname PrismFormulaTooltip@\df@label\endcsname}%
           \else
             \pdftooltip{\PrismOriginalTagform{#1}}{label: \detokenize\expandafter{\df@label}}%
           \fi
         \fi
       \else
         \pdftooltip{\PrismOriginalTagform{#1}}{\PrismFormulaTooltipText}%
       \fi
     \else
       \PrismOriginalTagform{#1}%
     \fi}
\let\PrismOriginalEquation\equation
\let\PrismOriginalEndEquation\endequation
\RenewEnviron{equation}{%
     \expandafter\PrismSetFormulaTooltip\expandafter{\BODY}%
     \PrismOriginalEquation
     \BODY
     \PrismOriginalEndEquation
     \global\let\PrismFormulaTooltipText\@empty}
\newcommand{\PrismInstallReferenceTooltips}{%
     \let\PrismOriginalRef\ref
     \renewcommand{\ref}[1]{%
       \ifShowReferenceTooltips
         \pdftooltip{\PrismOriginalRef*{##1}}{label: \detokenize{##1}}%
       \else
         \PrismOriginalRef{##1}%
       \fi}}
\AtBeginDocument{\PrismInstallReferenceTooltips}
\makeatother
\InputIfFileExists{prism_formula_tooltips.tex}{}{}
\newdimen\theight
\def\TeXref#1{
     \leavevmode\vadjust{\setbox0=\hbox{\quad\quad\LabelTagTooltip{#1}}%
     \theight=\ht0
     \advance\theight  by  \dp0
     \advance\theight  by  \lineskip
     \kern -\theight \vbox  to
     \theight{\rightline{\rlap{\box0}}%
      \vss}%
      }}%

    {\begin{thm}\label{#1} \ifShowLabels \TeXref{#1} \fi}%
    {\end{thm}}

    {\begin{def}\label{#1} \ifShowLabels \TeXref{#1} \fi}%
    {\end{def}}

    {\begin{lem}\label{#1} \ifShowLabels \TeXref{#1} \fi}%
    {\end{lem}}

    {\begin{cor}\label{#1} \ifShowLabels \TeXref{#1} \fi}%
    {\end{cor}}

\newcommand{\eqRef}[1]%
     {\ifShowLabels \TeXref{#1} \fi
      \begin{align}\label{#1} }

\ShowLabelsfalse

\newcommand{\cristiancomment}[1]{\textcolor{black}{#1}}

\newcommand{\ahmedcomment}[1]{\textcolor{black}{#1}}

\newcommand{\R}{{\mathbb R}} 

\newcommand{\G}{\mathbb G}

\ifJulinPaperOne
\title[Harnack Inequalities on Carnot groups]{Generalized and Classical Harnack Inequalities for semilinear equations on Carnot Groups}
\else\ifJulinPaperTwo
\title[Nonsmooth Harnack extensions]{Geometric Foundations and Nonsmooth Extensions for Harnack Inequalities on Carnot Groups}
\else
\title[Harnack Inequality \today]{A General Harnack Inequality for Uniformly Elliptic Equations in Stratified Groups\\\today}
\fi\fi

\ifJulinPaperOne
\author[C. E. Guti\'errez and A. Mohammed]{Cristian E. Guti\'errez and Ahmed Mohammed}
\thanks{The authors acknowledge the use of Prism, OpenAI's AI-powered \LaTeX{} editor, in preparing this paper.}
\else\ifJulinPaperTwo
\author[C. E. Guti\'errez and A. Mohammed]{Cristian E. Guti\'errez and Ahmed Mohammed}
\else
\author[C. E. Guti\'errez and A.Mohammed \today]{Cristian E. Guti\'errez and Ahmed Mohammed}
\fi\fi

\begin{document}

\ifJulinPaperOne
\begin{abstract}
We prove generalized and classical Harnack inequalities for non-negative
smooth solutions of divergence-form semilinear subelliptic equations on Carnot groups.
The proof combines smooth horizontal Gauss--Green and coarea formulas, energy
estimates for level sets, a Vitali covering argument, and a horizontal
Keller--Osserman barrier.  The companion paper
\cite{GutierrezMohammedCompanion} develops the geometric foundations and
extends the estimates to weak solutions with nonsmooth data.
\end{abstract}
\else\ifJulinPaperTwo
\begin{abstract}
This paper supplies the geometric and approximation results needed to extend
the smooth Harnack inequalities proved in the companion paper
\cite{GutierrezMohammedSmooth} to weak solutions of divergence-form
subelliptic equations on Carnot groups.  We collect horizontal perimeter,
Gauss--Green, coarea, and regular-level-set tools, then approximate measurable
coefficients and continuous monotone nonlinearities while preserving the
structural bounds and pass the Harnack estimates to the weak limit.
\end{abstract}
\else
\begin{abstract}
We prove a generalized Harnack inequality for non-negative weak solutions of divergence-form subelliptic equations on stratified groups,
\[
-X^*(A(x)Xu)=f(u),
\]
where $A(x)$ is uniformly elliptic and $f:[0,\infty)\to[0,\infty)$ is continuous and non-decreasing.
The argument combines energy estimates for sublevel sets, the Gauss--Green theorem for divergence-measure horizontal fields, and the coarea formula in Carnot groups.
Under an additional growth assumption on $f$, we then deduce the classical Harnack inequality.
\end{abstract}
\fi\fi

\maketitle

\tableofcontents

\ifJulinPaperOne
\section{Introduction}\label{sec:introduction}

Julin proved in \cite{VJU} a generalized Harnack inequality for non-negative
solutions of semilinear uniformly elliptic equations in Euclidean space.  The
purpose of this paper is to establish the corresponding result for smooth
divergence-form equations on Carnot groups and, under an additional growth
assumption on the nonlinearity, to recover the classical Harnack inequality.
The companion paper \cite{GutierrezMohammedCompanion} develops the geometric
foundations used here and extends the estimates to measurable coefficients,
continuous nonlinearities, and weak solutions.

Let \(\mathbb G\) be a Carnot group with stratified Lie algebra
\[
\mathfrak g=\mathcal V_1\oplus\cdots\oplus\mathcal V_\ell
\]
and homogeneous dimension \(Q\), and let \(n=\dim\mathfrak g\) be its
topological dimension.  Throughout the paper we use exponential coordinates
of the first kind.  Since \(\mathbb G\) is connected, simply connected, and
nilpotent, the exponential map is a global diffeomorphism, and we make the
identification
\begin{equation}\label{eq:euclidean-realization-of-carnot-group}
\mathbb R^n\simeq\mathfrak g
\xrightarrow{\ \exp\ }\mathbb G.
\end{equation}
Thus the underlying smooth manifold of \(\mathbb G\) is represented by
\(\mathbb R^n\).  The Carnot group law, inverse, left translations, and
dilations are defined as smooth maps on this Euclidean space.  In general,
the group product is not Euclidean vector addition, and the group dilations
are anisotropic rather than ordinary scalar multiplication.  When we write
\(gh\), \(g^{-1}\), or \(\delta_r g\), these operations are the Carnot-group
operations on the underlying space \(\mathbb R^n\).

Under the identification \eqref{eq:euclidean-realization-of-carnot-group}, a
left-invariant vector field is an ordinary smooth vector field on
\(\mathbb R^n\) whose coefficients are determined by the group law.  In
particular, each horizontal field has the coordinate form
\begin{equation}\label{eq:left-invariant-fields-euclidean-coordinates}
X_j
=
\sum_{k=1}^n b_{jk}(x)\frac{\partial}{\partial x_k},
\qquad j=1,\dots,m,
\end{equation}
and left invariance means
\[
X_j(gh)=d(L_g)_h\bigl(X_j(h)\bigr),
\qquad L_g(h)=gh.
\]
Accordingly, all horizontal derivatives are taken using these differential
operators on \(\mathbb R^n\), although their coefficients and commutators
encode the non-Euclidean group structure.

In exponential coordinates, Haar measure is a constant multiple of
Euclidean Lebesgue measure.  We normalize it so that
\begin{equation}\label{eq:haar-equals-lebesgue}
dg=dx=dx_1\cdots dx_n.
\end{equation}
Thus the notations \(dg\), \(dx\), and \(|E|\) used below refer to this same
measure.  This convention also makes distributional derivatives and
integration by parts compatible with the Euclidean realization.  Further
details about exponential coordinates, left translations, left-invariant
fields, and Haar measure are given in
\cite[Section~2]{GutierrezMohammedCompanion}.

We fix an inner product on the horizontal layer \(\mathcal V_1\), an
orthonormal basis of left-invariant vector fields \(X_1,\dots,X_m\), and the
associated Carnot--Carath\'eodory control distance (see
\cite{GutierrezMohammedCompanion}), which we denote by \(d\).
More generally, all arguments below remain valid for any fixed equivalent
left-invariant homogeneous distance.  The two defining invariances are
\begin{equation}\label{eq:left-invariant-homogeneous-distance}
d(ag,ah)=d(g,h),
\qquad
d(\delta_r g,\delta_r h)=r\,d(g,h).
\end{equation}
Let \(e\) denote the identity element of \(\mathbb G\).  The homogeneous
norm, or gauge, associated with \(d\) is
\begin{equation}\label{eq:distance-norm-compatibility}
N(g):=d(e,g),
\qquad
d(g,h)=N(g^{-1}h).
\end{equation}
The second identity follows from left invariance by translating both points by
\(g^{-1}\).  Since \(d\) is a distance, \(N\) satisfies
\begin{equation}\label{eq:homogeneous-norm-properties}
\begin{gathered}
N(e)=0,
\qquad
N(g^{-1})=N(g),
\qquad
N(gh)\le N(g)+N(h),\\
N(\delta_r g)=rN(g).
\end{gathered}
\end{equation}
Thus \(N\) is not the Euclidean norm on \(\mathbb R^n\); it is the norm
determined by the Carnot-group law and the anisotropic dilations.  The balls
may equivalently be written as
\begin{equation}\label{eq:balls-from-homogeneous-norm}
B(g,r)
=\{h\in\mathbb G:d(g,h)<r\}
=\{h\in\mathbb G:N(g^{-1}h)<r\}
=gB_r,
\qquad
B_r=\{h:N(h)<r\}.
\end{equation}
Whenever a different equivalent homogeneous distance is used locally, its
associated norm is always defined by the same rule \(N(g)=d(e,g)\).  We write
\[
Xu=(X_1u,\dots,X_mu)
\]
for the horizontal gradient and \(\operatorname{div}_H\) for the horizontal
divergence, with the convention in \cite[(2.5)]{GL}.  For horizontal Sobolev spaces, we use
the notation from \cite[Section 2, Preliminaries]{GL}.

In this paper the coefficient field
\(A\in C^\infty(\Omega;\operatorname{Sym}_m)\), the nonlinearity
\(f:[0,\infty)\to[0,\infty)\), and the solution \(u\) are smooth.  We assume
that \(f\) is non-decreasing and that there are constants
\(0<\lambda\le\Lambda\) such that
\[
\lambda|\zeta|^2
\le
\langle A(x)\zeta,\zeta\rangle
\le
\Lambda|\zeta|^2,
\qquad \zeta\in\mathbb R^m,
\]
for every \(x\in\Omega\).  Throughout the paper we set
\[
F(t):=\int_0^t f(s)\,ds.
\]

We study non-negative smooth solutions of
\begin{equation}\label{meq}
\operatorname{div}_H\bigl(A(x)Xu\bigr)=f(u)
\qquad\text{in }\Omega\subset\mathbb G.
\end{equation}
The sign convention for subsolutions and supersolutions is the following.  A
smooth function \(u\) is a subsolution, respectively a supersolution, if
\begin{equation}\label{sod}
-\int_\Omega\langle A(x)Xu,X\phi\rangle\,dx
\ge
\bigl(\text{respectively }\le\bigr)
\int_\Omega f(u)\phi\,dx
\end{equation}
for every non-negative \(\phi\in C_c^1(\Omega)\).  A solution is both a
subsolution and a supersolution, and for smooth functions this formulation is
equivalent to \eqref{meq} pointwise.

Our first main result is the following generalized Harnack inequality.

\begin{thm}[Generalized Harnack inequality]\label{gha}
There is a constant \(C=C(\lambda,\Lambda,Q)>0\) with the following property.
For every \(\xi\in\mathbb G\), \(R>0\), and every non-negative smooth solution $u$
of \eqref{meq} in \(B(\xi,2R)\),
\begin{equation}\label{di}
\int_{m_{\xi,R}}^{M_{\xi,R}}
\frac{dt}{R\sqrt{F(t)}+t}\le C,
\end{equation}
where
\[
M_{\xi,R}:=\sup_{B(\xi,R)}u,
\qquad
m_{\xi,R}:=\inf_{B(\xi,R)}u.
\]
The constant $C$ is independent of \(\xi\), \(R\), \(u\), and $f$.
\end{thm}

The proof is carried out first in the normalized case \(\xi=e\) and \(R=1\),
where \(e\) is the identity of \(\mathbb G\).  The general statement follows
from left translation and homogeneous dilation; this reduction is recorded in
Remark~\ref{rem:general-ball-rescaling}.

The Carnot-group argument requires several changes to Julin's Euclidean
proof.  The horizontal gradient controls only directions in the first layer,
and the relevant level surfaces and normal vectors must therefore be treated
with the horizontal geometry of \(\mathbb G\).  In the present smooth setting
we use the horizontal Gauss--Green formula on regular sublevel sets together
with the horizontal coarea formula.  The companion paper
\cite{GutierrezMohammedCompanion} gives a self-contained treatment of these
facts in the finite-horizontal-perimeter setting, including the Gauss--Green
theorem for divergence-measure horizontal fields, and uses them to justify the
nonsmooth extension.

A second difficulty concerns the covering argument.  Julin's proof uses the
Besicovitch covering lemma, whose bounded-overlap conclusion is not available
in the general metric setting relevant to Carnot groups.  We replace that step
with the \(5r\)-Vitali covering lemma \cite[Theorem~1.2]{HEI} and the resulting
disjoint-ball and density estimates.  The generalized estimate measures the
oscillation of \(u\) in the nonlinear scale determined by \(F\).  Its proof
combines estimates for the distribution functions of the lower and upper
level sets: the supersolution estimates give quantitative information for
\(E_t=\{u<t\}\), while the subsolution estimates control the superlevel sets.
The two estimates are joined in the proof of Theorem~\ref{gha}.

Our second main result is the classical Harnack inequality.  Set
\[
h(t):=\left(\frac{f(t)}{t}\right)^{1/2},
\qquad t>0,
\]
and impose the following growth condition.
\begin{description}
\item[\({\bf(f-I)}\)] The function \(h\) is non-decreasing on
\((0,\infty)\), and there are constants \(1<\varrho<\theta\) such that
\begin{equation}\label{tv}
h(\theta t)\ge\varrho h(t),
\qquad t>0.
\end{equation}
\end{description}

\begin{thm}[Classical Harnack inequality]\label{gha-c}
Suppose, in addition, that \(f\) satisfies \({\bf(f-I)}\).  There is a constant
\(C=C(\lambda,\Lambda,Q,\theta,\varrho)>0\) such that, for every
\(\xi\in\mathbb G\), \(R>0\), and every non-negative smooth solution $u$ of
\eqref{meq} in \(B(\xi,2R)\),
\begin{equation}\label{chi}
\sup_{B(\xi,R)}u\le C\inf_{B(\xi,R)}u.
\end{equation}
The constant is independent of \(\xi\), \(R\), and \(u\).
\end{thm}

The passage from the generalized estimate \eqref{di} to \eqref{chi} is based
on a horizontal Keller--Osserman construction.  A boundary blow-up
supersolution gives a pointwise bound for non-negative subsolutions, and
condition \({\bf(f-I)}\) converts that bound into the local scale estimate
\[
\frac{\sqrt{F(M_r(w)/2)}}{M_r(w)}
\le
\frac{C}{r},
\qquad
M_r(w):=\sup_{B(y_0,r)}w.
\]
Applied to \(u\), this estimate controls \(\sqrt{F(t)}/t\) below the maximum.
The generalized Harnack integral then bounds \(\log(M/m)\), which yields the
classical estimate.

The paper is organized as follows.  Section~\ref{sec:preliminaries} fixes the
notation and records the elementary iteration lemma.  Sections
\ref{sec:supersolution-estimates} and \ref{sec:subsolution-estimates} establish
the estimates for supersolutions and subsolutions, respectively, and
Section~\ref{sec:generalized-harnack-proof} proves Theorem~\ref{gha}.
The remaining sections develop the Keller--Osserman reduction, and
Section~\ref{sec:blue-proof-theorem-22} completes the horizontal barrier
construction and proves Theorem~\ref{gha-c}.  The geometric and
measure-theoretic results underlying these arguments, together with the
extension to weak solutions, are developed in the companion paper
\cite{GutierrezMohammedCompanion}.

\else\ifJulinPaperTwo
\section{Introduction}\label{sec:introduction}
Let \(\mathbb G\) be a Carnot group with homogeneous dimension \(Q\), fixed
homogeneous distance, and horizontal frame \(X_1,\dots,X_m\).  Let
\(A\in L^\infty(\Omega;\operatorname{Sym}_m)\) satisfy
\(\lambda I\le A\le\Lambda I\) almost everywhere, and let
\(f:[0,\infty)\to[0,\infty)\) be continuous and non-decreasing.  Put
\[
Xu=(X_1u,\dots,X_mu),\qquad F(t)=\int_0^t f(s)\,ds,
\]
and consider weak solutions of
\begin{equation}\label{meq}
\operatorname{div}_H(A(x)Xu)=f(u),
\end{equation}
that is,
\begin{equation}\label{sod}
-\int_\Omega\langle AXu,X\phi\rangle=\int_\Omega f(u)\phi
\end{equation}
for all \(\phi\in C_c^1(\Omega)\).  The companion paper
\cite{GutierrezMohammedSmooth} proves the following estimates for smooth data;
the purpose here is to establish the geometric tools and pass them to the
nonsmooth setting.

\begin{thm}[Generalized Harnack inequality]\label{gha}
Every non-negative weak solution in \(B(\xi,2R)\) satisfies
\begin{equation}\label{di}
\int_{m_{\xi,R}}^{M_{\xi,R}}\frac{dt}{R\sqrt{F(t)}+t}\le C,
\end{equation}
where \(C=C(\lambda,\Lambda,Q)\).
\end{thm}
Let \(h(t)=\sqrt{f(t)/t}\), and assume
\begin{equation}\label{tv}
h\text{ is non-decreasing and }h(\theta t)\ge\varrho h(t),
\quad 1<\varrho<\theta.
\end{equation}
This condition is denoted by \({\bf(f-I)}\).
\begin{thm}[Classical Harnack inequality]\label{gha-c}
Under \({\bf(f-I)}\),
\begin{equation}\label{chi}
\sup_{B(\xi,R)}u\le C\inf_{B(\xi,R)}u.
\end{equation}
\end{thm}
\begin{thm}[Horizontal local Keller--Osserman scale estimate]
\label{thm:blue-horizontal-local-scale}
If \(w\) is a non-negative locally bounded weak subsolution in
\(B(y_0,2r)\), then
\begin{equation}\label{eq:blue-carnot-subsolution-scale}
\frac{\sqrt{F(M_r(w)/2)}}{M_r(w)}\le\frac Cr,
\qquad M_r(w)=\sup_{B(y_0,r)}w,
\end{equation}
whenever \(M_r(w)>0\).
\end{thm}

\else
\section{Introduction}\label{sec:introduction}

Julin proved in \cite{VJU} a generalized Harnack inequality for non-negative
solutions of semilinear uniformly elliptic equations in Euclidean space.  The
purpose of this paper is to establish the corresponding result for
divergence-form equations on Carnot groups and, under an additional growth
assumption on the nonlinearity, to recover the classical Harnack inequality.

Let \(\mathbb G\) be a Carnot group with stratified Lie algebra
\[
\mathfrak g=\mathcal V_1\oplus\cdots\oplus\mathcal V_\ell,
\]
and homogeneous dimension \(Q\).  We fix an inner product on the horizontal
layer \(\mathcal V_1\), an orthonormal basis of left-invariant vector fields
\(X_1,\dots,X_m\), and a homogeneous distance on \(\mathbb G\).  We write
\[
Xu=(X_1u,\dots,X_mu)
\]
for the horizontal gradient and \(\operatorname{div}_H\) for the horizontal
divergence.  The coefficient matrix \(A(x)=(a_{ij}(x))\) is assumed measurable,
symmetric, and uniformly elliptic in the horizontal directions: there are
constants \(0<\lambda\le\Lambda\) such that
\[
\lambda|\xi|^2
\le
\langle A(x)\xi,\xi\rangle
\le
\Lambda|\xi|^2,
\qquad \xi\in\mathbb R^m,
\]
for almost every \(x\).  Throughout the paper,
\(f:[0,\infty)\to[0,\infty)\) is continuous and non-decreasing, and
\[
F(t):=\int_0^t f(s)\,ds.
\]

We study non-negative weak solutions of
\begin{equation}\label{meq}
\operatorname{div}_H\bigl(A(x)Xu\bigr)=f(u)
\qquad\text{in }\Omega\subset\mathbb G.
\end{equation}
For an open set \(\Omega\), let
\[
W^{1,2}(\Omega,X)
:=
\{u\in L^2(\Omega):X_ju\in L^2(\Omega),\ j=1,\dots,m\}.
\]
A function \(u\in W^{1,2}(\Omega,X)\) is a weak subsolution, respectively a
weak supersolution, of \eqref{meq} if \(f(u)\in L^1_{\mathrm{loc}}(\Omega)\) and
\begin{equation}\label{sod}
-\int_\Omega\langle A(x)Xu,X\phi\rangle\,dx
\ge
\bigl(\text{respectively }\le\bigr)
\int_\Omega f(u)\phi\,dx
\end{equation}
for every non-negative \(\phi\in C_c^1(\Omega)\).  A weak solution is both a
subsolution and a supersolution.  We use the locally bounded representative
provided by the standard local theory for homogeneous \(X\)-elliptic
subsolutions.

Our first main result is the following generalized Harnack inequality.

\begin{thm}[Generalized Harnack inequality]\label{gha}
There is a constant \(C=C(\lambda,\Lambda,Q)>0\) with the following property.
For every \(\xi\in\mathbb G\), \(R>0\), and every non-negative weak solution
\(u\in W^{1,2}(B(\xi,2R),X)\) of \eqref{meq} in \(B(\xi,2R)\),
\begin{equation}\label{di}
\int_{m_{\xi,R}}^{M_{\xi,R}}
\frac{dt}{R\sqrt{F(t)}+t}\le C,
\end{equation}
where
\[
M_{\xi,R}:=\sup_{B(\xi,R)}u,
\qquad
m_{\xi,R}:=\inf_{B(\xi,R)}u.
\]
Here the balls are defined using the fixed left-invariant homogeneous distance.
In particular, the normalized notation \(B_r=B(e,r)\) corresponds to
\(\xi=e\) and \(R=1\).  See
Section~\ref{sec:background-homogeneous-carnot} for the precise Carnot-group
and distance conventions, and Remark~\ref{rem:general-ball-rescaling} for the
reduction to the normalized case.
\end{thm}

The proof is first carried out in the normalized case \(\xi=e\) and \(R=1\),
that is, on \(B(e,1)\subset B(e,2)\).  The general case follows by left
translation and homogeneous dilation, as detailed in
Remark~\ref{rem:general-ball-rescaling}.

The proof in the Carnot-group setting is substantially more involved than
Julin's Euclidean argument.  It requires geometric-measure-theoretic tools
adapted to Carnot groups, including the Gauss--Green formula for sets of finite
horizontal perimeter \cite{CVA}, as well as a replacement for the Euclidean
covering argument.  In particular, Julin uses the Besicovitch covering lemma,
whose bounded-overlap conclusion is not available in the general metric-space
setting relevant here.  We replace that step with the \(5r\)-Vitali covering
lemma \cite[Theorem~1.2]{HEI} and the corresponding disjoint-ball and density
estimates.

The estimate measures the oscillation of \(u\) in the nonlinear scale
determined by \(F\).  Its proof combines estimates for the distribution
functions of the sublevel and superlevel sets of \(u\).  The supersolution
estimates provide quantitative volume and level-surface information for
\(E_t=\{u<t\}\), while the subsolution estimates control the upper levels.
These two parts are then joined in the proof of Theorem~\ref{gha}.

The second main result is the classical Harnack inequality.  Set
\[
h(t):=\left(\frac{f(t)}{t}\right)^{1/2},
\qquad t>0,
\]
and assume the following growth condition.
\begin{description}
\item[\({\bf(f-I)}\)] The function \(h\) is non-decreasing on
\((0,\infty)\), and there exist constants \(1<\varrho<\theta\) such that
\begin{equation}\label{tv}
h(\theta t)\ge\varrho h(t),
\qquad t>0.
\end{equation}
\end{description}

\begin{thm}[Classical Harnack inequality]\label{gha-c}
Suppose, in addition, that \(f\) satisfies \({\bf(f-I)}\).  There is a constant
\(C=C(\lambda,\Lambda,Q,\theta,\varrho)>0\) such that, for every
\(\xi\in\mathbb G\), \(R>0\), and every non-negative weak solution
\(u\in W^{1,2}(B(\xi,2R),X)\) of \eqref{meq} in \(B(\xi,2R)\),
\begin{equation}\label{chi}
\sup_{B(\xi,R)}u\le C\inf_{B(\xi,R)}u.
\end{equation}
The constant is independent of \(\xi\) and \(R\).
\end{thm}

The passage from the generalized estimate \eqref{di} to \eqref{chi} is based on
a horizontal Keller--Osserman construction.  A boundary blow-up supersolution
gives a pointwise bound for non-negative subsolutions, and condition
\({\bf(f-I)}\) converts that bound into the local scale estimate
\[
\frac{\sqrt{F(M_r(w)/2)}}{M_r(w)}
\le
\frac{C}{r},
\qquad
M_r(w):=\sup_{B(y_0,r)}w.
\]
Applied to \(u\), this estimate controls \(\sqrt{F(t)}/t\) below the maximum.
The generalized Harnack integral then bounds \(\log(M/m)\), which yields
\eqref{chi}.

The proofs use geometric measure theory adapted to the horizontal structure.
In particular, Gauss--Green is applied to sublevel sets of finite horizontal perimeter,
and horizontal coarea is used to integrate the resulting level-surface
estimates.  These tools, together with the required facts about reduced
boundaries, regular level sets, homogeneous balls, and almost-everywhere finite
perimeter, are collected in
Section~\ref{sec:background-homogeneous-carnot}; we do not repeat their detailed
statements here.

The paper is organized as follows.  Section~\ref{sec:preliminaries} fixes the
notation and records the elementary iteration lemma.  Sections
\ref{sec:supersolution-estimates} and \ref{sec:subsolution-estimates} establish
the estimates for supersolutions and subsolutions, respectively, and
Section~\ref{sec:generalized-harnack-proof} proves Theorem~\ref{gha}.
Section~\ref{sec:blue-proof-theorem-22} develops the horizontal barrier
construction and proves Theorem~\ref{gha-c}.  The brief auxiliary
Section~\ref{sec:infinitesimal-generator-dilations} records the infinitesimal
generator of the group dilations used in the cutoff construction.  Finally,
Section~\ref{sec:background-homogeneous-carnot} contains the Carnot-group and
horizontal geometric-measure-theoretic background used throughout.

\fi\fi

\ifJulinPaperTwo
\else

\ifJulinPaperOne
\else

{\color{black}
\section{\bf Preliminaries}\label{sec:preliminaries}

We recall the basic notation used in the two Harnack statements.  A Carnot group
\(\mathbb G\) is a connected, simply connected nilpotent Lie group whose Lie
algebra has a stratification
\[
\mathfrak g=\mathcal V_1\oplus\cdots\oplus\mathcal V_\ell,
\qquad [\mathcal V_1,\mathcal V_j]=\mathcal V_{j+1}.
\]
The horizontal layer \(\mathcal V_1\) is generated by the left-invariant vector
fields \(X_1,\dots,X_m\), and the homogeneous dimension is denoted by \(Q\).  We
write \(B(x,r)\) for balls associated with a fixed homogeneous distance and
\(B_r=B(e,r)\) when the center is the identity.

For a function \(u\), the horizontal gradient and horizontal divergence are
\[
Xu=(X_1u,\dots,X_mu),
\qquad
\operatorname{div}_H Z=\sum_{i=1}^m X_iZ_i.
\]
The coefficient matrix \(A(x)=(a_{ij}(x))\) is assumed measurable, symmetric, and
uniformly elliptic in the horizontal directions:
\[
\lambda |\xi|^2\le \langle A(x)\xi,\xi\rangle\le \Lambda |\xi|^2,
\qquad \xi\in\mathbb R^m.
\]
Weak subsolutions and supersolutions of
\(\operatorname{div}_H(A(x)Xu)=f(u)\) are understood with the sign convention
stated in \eqref{sod}.  The two quantities attached to the nonlinearity are
\[
F(t)=\int_0^t f(s)\,ds,
\qquad
h(t)=\left(\frac{f(t)}{t}\right)^{1/2}\quad (t>0),
\]
and the classical Harnack estimate uses the growth hypothesis {\bf(f-I)} stated
above in \eqref{tv}.

We shall also use the following elementary algebraic iteration lemma, which is
taken from \cite{EGU}.

\begin{lem}\label{eg}
Let $\alpha>0$, and let $(x_j)$ be a sequence of positive numbers such that
\[
x_{j+1}\le CB^j x_j^{1+\alpha}
\]
for some constants $C>0$ and $B>1$. If
\[
x_0\le C^{-1/\alpha}B^{-1/\alpha^2},
\]
then
\[
x_j\le B^{-j/\alpha}x_0,
\]
and in particular
\[
\lim_{j\to\infty}x_j=0.
\]
\end{lem}

The geometric measure theory used in the proof is developed in the companion
paper \cite{GutierrezMohammedCompanion}.  The precise items used later are
Definitions~3.1--3.4 and Propositions~3.3, 3.6, and 3.9--3.13.  They provide
horizontal perimeter, the relative isoperimetric inequality, the reduced
boundary, Gauss--Green, regular sublevel normals, coarea, almost-everywhere
finite perimeter, gauge-ball formulas, and the admissible-level properties
needed in the proof of Proposition~\ref{prop1}.

}
\fi

\ifJulinPaperOne
\section{Preliminaries}\label{sec:preliminaries}
We collect the standard Carnot-group formulas used in the smooth proofs.  The
finite-perimeter formulations and weak-solution extensions are proved in
\cite{GutierrezMohammedCompanion}.
We use throughout the exponential-coordinate realization fixed in
\eqref{eq:euclidean-realization-of-carnot-group}.  Thus \(\mathbb G\) has
underlying manifold \(\mathbb R^n\), while its multiplication, inversion, and
dilations are the Carnot-group operations defined on that space.  The fields
\(X_j\) are the Euclidean-coordinate vector fields in
\eqref{eq:left-invariant-fields-euclidean-coordinates}, and Haar integration
is ordinary Lebesgue integration by \eqref{eq:haar-equals-lebesgue}.  In
particular, an open set \(\Omega\subset\mathbb G\) is also an open subset of
\(\mathbb R^n\), and the symbols \(dg\) and \(dx\) are interchangeable.
\else
{\color{black}
\section*{Proposed compact preliminaries for Sections 3--6: the smooth case}

This is an editorial proposal for a shorter version of the paper.  It is
intended to replace, rather than supplement, the present preliminary section
and the long geometric-measure-theory background when the proofs in
Sections~3--6 are first carried out for smooth \(u\), \(A\), and \(f\).  The
reduction from smooth data to weak solutions can then be handled separately by
the approximation argument at the end of the paper.
\fi

Let \(\mathbb G\) be a Carnot group, let \(X_1,\dots,X_m\) be an orthonormal
basis of the first layer, let \(Q\) be the homogeneous dimension, and let \(d\)
be the fixed left-invariant homogeneous distance from
\eqref{eq:left-invariant-homogeneous-distance}.  Its associated homogeneous
norm is \(N(g)=d(e,g)\), and
\[
B(g,r)=\{z:N(g^{-1}z)<r\},
\qquad
B_r=B(e,r)=\{z:N(z)<r\},
\]
as recorded in \eqref{eq:distance-norm-compatibility} and
\eqref{eq:balls-from-homogeneous-norm}.  We also write
\(Xu=(X_1u,\dots,X_mu)\).  
In the remaning part of this section, we mention several useful facts needed
in the proofs, and we introduce more notation.  Constants depend only on the group, the distance (or,
equivalently, its associated norm), and, when an operator is involved, the
ellipticity constants.

For an open set \(\Omega\subset\mathbb G\), \(BV_H(\Omega)\) denotes the
functions \(v\in L^1_{\mathrm{loc}}(\Omega)\) whose distributional horizontal
gradient \(D_Hv=(D_{X_1}v,\ldots,D_{X_m}v)\) is an \(\mathbb R^m\)-valued
Radon measure; its total variation is denoted by \(|D_Hv|\).  A measurable
set \(E\subset\Omega\) has finite horizontal perimeter in \(\Omega\) if
\(\chi_E\in BV_H(\Omega)\), and its horizontal perimeter measure is
\(P_H(E;\cdot):=|D_H\chi_E|\).

\paragraph{\normalfont\itshape Translations, dilations, balls, and cutoffs.}
Haar measure and the group dilations satisfy
\begin{equation}\label{eq:proposed-prelim-ball-scaling}
|gE|=|E|,
\qquad
|\delta_rE|=r^Q|E|,
\qquad
|B(g,r)|=\omega_Q r^Q; \quad \omega_Q:=|B_1|.
\end{equation}
If \(0<r<R\), one may choose
\(\zeta\in C_c^\infty(B(g,R))\) with
\begin{equation}\label{eq:proposed-prelim-cutoff}
0\le\zeta\le1,
\qquad
\zeta=1\ \hbox{on }B(g,r),
\qquad
|X\zeta|\le \frac{C}{R-r}.
\end{equation}
The \(5r\)-covering lemma will be used in the following form: from every family
of balls with uniformly bounded radii one can select pairwise disjoint balls
\(\{B_i\}\) such that the original union is contained in \(\bigcup_i5B_i\).
These elementary facts follow from left invariance, homogeneity, and the usual
metric Vitali argument; see, for example,
\cite[Chapter~1]{FollandSteinHardy}.

For later normalization, if
\(\widetilde u(z)=u(g\delta_rz)\), then
\begin{equation}\label{eq:proposed-prelim-rescaling}
X_i\widetilde u(z)=r(X_iu)(g\delta_rz),
\qquad
\operatorname{div}_H(\widetilde A X\widetilde u)(z)
=r^2\operatorname{div}_H(AXu)(g\delta_rz),
\end{equation}
where \(\widetilde A(z)=A(g\delta_rz)\).  Thus translations and dilations
preserve the ellipticity constants and reduce arbitrary balls to normalized
balls.

\paragraph{\normalfont\itshape Horizontal calculus and analytic inequalities.}
The horizontal fields are divergence-free with respect to Haar measure, so
\begin{equation}\label{eq:proposed-prelim-integration-by-parts}
\int_{\mathbb G}\phi X_i\psi\,dx
=-
\int_{\mathbb G}\psi X_i\phi\,dx,
\qquad
\operatorname{div}_HZ=\sum_{i=1}^mX_iZ_i,
\end{equation}
whenever one factor is compactly supported.  The usual product and chain rules
hold for smooth functions; in particular,
\(X\Phi(u)=\Phi'(u)Xu\).

The only Sobolev input needed in Sections~3--6 is the horizontal
\(W^{1,1}\) Sobolev--Poincar\'e inequality
\begin{equation}\label{eq:proposed-prelim-sobolev-poincare}
\left(\int_B|v-v_B|^{\frac Q{Q-1}}\,dx\right)^{\frac{Q-1}{Q}}
\le
C\frac{r}{|B|^{1/Q}}\int_B|Xv|\,dx,
\qquad B=B(g,r),
\end{equation}
and its zero-boundary or large-zero-set consequence
\begin{equation}\label{eq:proposed-prelim-zero-set-sobolev}
\left(\int_B|v|^{\frac Q{Q-1}}\,dx\right)^{\frac{Q-1}{Q}}
\le C_\theta\int_B|Xv|\,dx
\end{equation}
when \(v=0\) on a set of measure at least \(\theta|B|\).  We use
\cite[Theorem~1.5 and Remark~1.2]{GAN}.  Applied to \(v=w^2\),
\eqref{eq:proposed-prelim-zero-set-sobolev} gives exactly the exponent
\(2Q/(Q-1)\) used in the subsolution Caccioppoli estimate.

For a set \(E\) of locally finite horizontal perimeter, write
\(\partial_H^*E\) for its horizontal reduced boundary.  The measure
\(P_H(E;\cdot)\) is concentrated on \(\partial_H^*E\), so for a ball \(B\),
\[
P_H(E;B)=\int_{\partial_H^*E\cap B}dP_H(E).
\]
See \cite[Definitions~3.2 and 3.4, Proposition~3.5]{GutierrezMohammedCompanion}
for these definitions and the concentration result.  The relative
isoperimetric inequality in a Carnot ball is
\label{subsec:app-horizontal-bv}
\label{def:app-horizontal-perimeter}
\begin{equation}\label{eq:proposed-prelim-relative-isoperimetric}
P_H(E;B)
\ge
c\min\{|E\cap B|,|B\setminus E|\}^{\frac{Q-1}{Q}}.
\end{equation}
\label{prop:app-relative-isoperimetric-carnot-balls}
\label{prop:app-reduced-boundary-concentration}
For \eqref{eq:proposed-prelim-relative-isoperimetric}, see also
\cite[Theorem~1.18]{GAN}.  In particular, if both phases occupy a fixed
positive proportion of \(B(g,r)\), then
\(P_H(E;B(g,r))\ge c r^{Q-1}\).

Finally, the homogeneous De Giorgi--Moser estimate used for subsolutions is
\begin{equation}\label{eq:proposed-prelim-local-boundedness}
\sup_{B(g,r/2)}v
\le
C\left(\frac1{|B(g,r)|}\int_{B(g,r)}v^2\,dx\right)^{1/2}
\end{equation}
whenever $v\geq 0$ solves
\(\operatorname{div}_H(A(x)Xv)\ge0\) weakly in \(B(g,r)\).  This is the
local boundedness inequality \cite[(4.12)]{GL}.

\paragraph{\normalfont\itshape Smooth level surfaces and coarea.}
Let \(u\in C^\infty(\Omega)\), set \(E_t=\{u<t\}\), and let
\(\Sigma_t=\{u=t\}\).  For almost every \(t\), Sard's theorem and coarea permit
all the following identities.  On the noncharacteristic part of
\(\Sigma_t\),
\label{def:app-horizontal-reduced-boundary}
\label{def:app-euclidean-regular-value}
\label{prop:app-regular-sublevel-normal}
\label{prop:app-prop31-admissible-levels}
\begin{equation}\label{eq:proposed-prelim-level-normal}
\nu_{E_t}=\frac{Xu}{|Xu|},
\qquad
dP_H(E_t)=\frac{|Xu|}{|\nabla u|}\,d\mathcal H^{N-1}.
\end{equation}
\label{eq:app-levelsurface-boundary}
\label{eq:app-DHchi-Et}
\label{eq:app-PH-Et-smooth-level}
\label{eq:app-horizontal-normal-Et}
\label{Pmes}
This is the smooth-boundary representation formula in
\cite[Proposition~1.19]{MontefalconePerimeter}, applied to the defining
function \(u-t\), after adjusting the inward/outward normal convention.
The characteristic part has zero horizontal perimeter measure.  Hence, for
every compactly supported smooth horizontal field \(Z\),
\label{prop:app-gauss-green-finite-perimeter}
\label{subsec:app-gauss-green-theorem}
\begin{equation}\label{eq:proposed-prelim-smooth-gauss-green}
\int_{E_t}\operatorname{div}_HZ\,dx
=
\int_{\Sigma_t}
\left\langle Z,\frac{Xu}{|Xu|}\right\rangle dP_H(E_t).
\end{equation}
\label{eq:app-gauss-green-theorem}
This is the smooth-level form of the Gauss--Green theorem; compare
\cite[Theorem~1.5]{CVA}.

For every nonnegative Borel function \(g\), horizontal coarea gives
\label{prop:app-horizontal-coarea}
\label{prop:app-levelsets-finite-perimeter}
\begin{equation}\label{eq:proposed-prelim-horizontal-coarea}
\int_\Omega g|Xu|\,dx
=
\int_{-\infty}^{\infty}
\int_{\Sigma_t}g\,dP_H(E_t)\,dt.
\end{equation}
\label{eq:app-horizontal-coarea-sobolev}
\label{eq:app-levelset-finite-perimeter-integral}
More generally, the same formula holds for locally Lipschitz functions and
their superlevel sets; this is the version used for the localized density
\(\sigma_A\), defined in \eqref{eq:sigma-A-density-definition}.  See
\cite{VAM}.

If \(B\Subset\Omega\) and
\begin{equation*}
\eta_B(t)=|\{u<t\}\cap B|,
\qquad
\mu_B(t)=|\{u\ge t\}\cap B|,
\end{equation*}
then for almost every admissible level \(t\),
\begin{equation}\label{eq:proposed-prelim-distribution-derivative}
\eta_B'(t)=-\mu_B'(t)
=
\int_{\Sigma_t\cap B}\frac1{|Xu|}\,dP_H(E_t)
=
\int_{\Sigma_t\cap B}\frac1{|\nabla u|}\,d\mathcal H^{N-1}.
\end{equation}
Equations
\eqref{eq:proposed-prelim-level-normal}--
\eqref{eq:proposed-prelim-distribution-derivative} are precisely the
level-set formulas used in the surface estimate and in the differential
inequalities for \(\mu\) and \(\eta\).

\ifJulinPaperOne
We also use the following elementary iteration lemma from \cite{EGU}.
\begin{lem}\label{eg}
Let \(\alpha>0\), and suppose
\(x_{j+1}\le CB^j x_j^{1+\alpha}\), where \(C>0\) and \(B>1\).  If
\(x_0\le C^{-1/\alpha}B^{-1/\alpha^2}\), then
\(x_j\le B^{-j/\alpha}x_0\), and hence \(x_j\to0\).
\end{lem}
This is the complete toolbox needed for the smooth arguments.
\else
Together with the elementary iteration Lemma~\ref{eg}, this is the complete
standard Carnot-analysis toolbox needed for the smooth arguments in
Sections~3--6.  The ball-density function \(\sigma_A\), the truncated boundary
distance \(d_\delta\), the boundary-concentration estimate, and the two Julin
measure lemmas are specific to the present proof and should remain in
Section~3.

If this proposal is adopted, the long background section on horizontal
\(BV\), reduced boundaries, Gauss--Green, and coarea can be commented out in
the smooth version.
}
\fi

\section{Estimates for Supersolutions}\label{sec:supersolution-estimates}

\vspace{.3cm}
{\color{black}
\noindent\emph{Roadmap for the section.}
The goal of this section is to prove the differential measure estimates in
Proposition~\ref{prop:P35-carnot}, the main input for the generalized Harnack
inequality.  The argument proceeds in layers.
Proposition~\ref{prop1}
first gives the basic Caccioppoli estimate and the corresponding level-surface
estimate; the admissible-level lemma just below it explains why Gauss--Green and
coarea may be applied to the relevant sublevel sets.  The Sobolev--Poincare
lemma and Proposition~\ref{prop1} then yield the De Giorgi small-sublevel
estimate.  \textcolor{black}{The \(5r\)-Vitali covering lemma, together with
the balanced-ball estimate in Lemma~\ref{lem:L33-carnot}, lets us convert
volume bounds into level-surface estimates.}  
The horizontal boundary-concentration estimate
Lemma~\ref{bce} is used to prove the measure
Lemma~\ref{carnot-julin-34}, and this closes the large-density branch of
Proposition~\ref{prop:P35-carnot}.

\smallskip
Throughout this section, for $t>0$, we use the notation
\[
E_t:=\{x\in B_2:u(x)<t\}.
\]
We remind the reader of our standing assumptions that $u$, $A$, and $f$ are $C^\infty$, with $A$ uniformly elliptic in the horizontal direction. We first establish the following estimates.

\begin{prop}\label{prop1}
	Suppose $u\in C^\infty(B_2)$ is a non-negative supersolution of (\ref{meq}).
	If $0<r<R<2$, then 
	\begin{equation}\label{-11}
		\int_{E_t\cap B_r}|Xu|^2\,dx\le C\left(\frac{t^2}{(R-r)^2}+F(2t)\right)|E_t\cap B_R|.
	\end{equation}
	
	\vspace{.2cm}
	\noindent Furthermore, for almost every $t>0$, we have the following
	boundary estimate, where $\partial_H^*E_t$ denotes the horizontal reduced
	boundary of $E_t$; see
	\cite[Definition~3.4]{GutierrezMohammedCompanion}.
	\begin{equation}\label{-22}
		\int_{\partial_H^* E_t\cap B_r}|Xu|\,dP_H(E_t)\le C\left[\frac{t}{(R-r)^2}+\frac{F(2t)}{t}\right]|E_t\cap B_R|.
	\end{equation}
	The constants $C$ in (\ref{-11}) and  (\ref{-22}) depend on $\lambda,\Lambda$ and $Q$ only.
\end{prop}

\noindent Before proving Proposition 3.1, we isolate the following admissibility
lemma for the level sets; after the proof of the lemma, we return to the proof
of Proposition 3.1.

\begin{lem}[Admissible levels for the surface estimate]\label{lem:prop31-admissible-surface-levels}
Assume that \(u\in C^\infty(B_2)\), and set
\[
E_t:=\{x\in B_2:u(x)<t\},
\qquad
\Sigma_t:=\{u=t\}\cap B_2.
\]
Then there is a full-measure set \(\mathcal T\subset(0,\infty)\) such that,
for every \(t\in\mathcal T\), the following hold.
\begin{enumerate}[(i)]
\item \(E_t\) has locally finite horizontal perimeter in \(B_2\).
\item The level surface \(\Sigma_t\) is either empty or a \(C^1\) (in fact,
\(C^\infty\)) Euclidean hypersurface.
\item If \(Z\in C_c^1(B_2;\mathbb R^m)\), then
\begin{equation}\label{eq:prop31-admissible-level-GG}
\int_{E_t}\operatorname{div}_H Z\,dx
=
\int_{\partial_H^*E_t}\langle Z,\nu_{E_t}\rangle\,dP_H(E_t).
\end{equation}
\item On every ball \(B_R\Subset B_2\), the smooth-level representation
\eqref{eq:proposed-prelim-level-normal} applies to \(E_t\). In particular,
\begin{equation}\label{eq:prop31-admissible-level-normal}
\nu_{E_t}=\frac{Xu}{|Xu|}
\qquad
P_H(E_t)\text{-a.e. on }\partial_H^*E_t\cap B_R.
\end{equation}
\end{enumerate}
\end{lem}
\noindent{\bf Proof:} By the horizontal coarea formula
\eqref{eq:proposed-prelim-horizontal-coarea}, applied to
\(u\in C^\infty(B_2)\subset W^{1,1}_{H,\mathrm{loc}}(B_2)\), we have, for every
compact set \(K\Subset B_2\),
\begin{equation}\label{eq:prop31-line688-coarea-local}
\int_0^\infty P_H(\{u<t\};K)\,dt
\le
\int_K |Xu|\,dx
<\infty.
\end{equation}
Thus \(E_t\) has locally finite horizontal perimeter for a.e. \(t>0\): the
exceptional levels for which \(P_H(E_t;K)=\infty\) have Lebesgue measure zero
for each fixed compact \(K\Subset B_2\), and a countable exhaustion of \(B_2\)
gives local finite perimeter in \(B_2\). Thus
\eqref{eq:prop31-line688-coarea-local} itself supplies exactly the
finite-perimeter conclusion needed below.

Independently, Sard's theorem gives the regularity of almost every level. In
the scalar case, if \(v\in C^k(U)\) with \(k>n-1\), then the set of values
attained at points where \(\nabla v=0\) has one-dimensional Lebesgue measure
zero.  Applying this with \(U=B_2\subset\mathbb R^n\), and \(v=u\),
we obtain
\begin{equation}\label{eq:prop31-line688-sard}
\mathcal L^1\big(\{s>0:\exists x\in B_2,\,u(x)=s\text{ and }\nabla u(x)=0\}\big)=0,
\end{equation}
where \(\nabla u\) is the ordinary Euclidean gradient.  Hence, outside this
Sard null set, \(\Sigma_t\) is either empty or a \(C^1\) Euclidean hypersurface.

Let \(\mathcal T\) be the complement in \((0,\infty)\) of the union of the
coarea null set and the Sard null set.  For \(t\in\mathcal T\), (i) and (ii)
hold.  Since \(E_t\) has locally finite horizontal perimeter and \(Z\) is
compactly supported in \(B_2\), the Gauss--Green theorem
\cite[Theorem~1.5]{CVA}, in the form
\eqref{eq:proposed-prelim-smooth-gauss-green}, gives (iii).  Finally, the
smooth-level representation \eqref{eq:proposed-prelim-level-normal}
identifies the horizontal perimeter measure and the outer horizontal normal,
which gives (iv).\qed

\vspace{.2cm}
\noindent{\bf Proof of Proposition 3.1:} We first prove the solid Caccioppoli estimate
(\ref{-11}).  Choose \(\zeta\in C_c^\infty(B_R)\) such that
\[
0\le \zeta\le 1,
\qquad
\zeta\equiv 1 \text{ on } B_r,
\qquad
|X\zeta|\le \frac{C}{R-r},
\]
and use the non-negative test function
\[
\phi=(t-u)_+\zeta^2.
\]
Since \(u\) is a supersolution of (\ref{meq}),
\begin{equation}\label{aa-1}
\int_{B_R}\langle A(x)Xu,X\phi\rangle\,dx
\ge
-\int_{B_R}f(u)\phi\,dx.
\end{equation}
Moreover,
\begin{equation}\label{aa-2}
X\phi=\zeta^2X(t-u)_+ +2(t-u)_+\zeta X\zeta,
\end{equation}
and \(X(t-u)_+=-Xu\) on \(E_t\), while \(X(t-u)_+=0\) outside \(E_t\).
Substituting \((\ref{aa-2})\) into \((\ref{aa-1})\) therefore gives
\begin{align*}
\int_{E_t\cap B_R}\zeta^2\langle A(x)Xu,Xu\rangle\,dx
&\le
2\int_{E_t\cap B_R}(t-u)\zeta
|\langle A(x)Xu,X\zeta\rangle|\,dx \\
&\quad+
\int_{E_t\cap B_R}f(u)(t-u)\zeta^2\,dx.
\end{align*}
By horizontal ellipticity, the left-hand side in the preceding inequality is
bounded below by
\[
\lambda\int_{E_t\cap B_R}\zeta^2|Xu|^2\,dx.
\]
Moreover,
\[
|\langle A(x)Xu,X\zeta\rangle|
\le \Lambda |Xu|\,|X\zeta|.
\]
We then use
Young's inequality in the form
\[
2\Lambda ab\le \frac{\lambda}{2}a^2+C(\lambda,\Lambda)b^2,
\qquad
 a=\zeta|Xu|,\quad b=(t-u)|X\zeta|,
\]
and absorb the \((\lambda/2)a^2\) term into the left-hand side.  This gives,
with \(C\) depending only on \(\lambda\) and \(\Lambda\),
\begin{equation*}
\int_{E_t\cap B_R}\zeta^2 |Xu|^2\,dx
\le
C\int_{E_t\cap B_R}(t-u)^2|X\zeta|^2\,dx
+C\int_{E_t\cap B_R}f(u)(t-u)\,dx.
\end{equation*}
Since \(0\le t-u\le t\) on \(E_t\) and
\(|X\zeta|\le C/(R-r)\), the first term on the right satisfies
\begin{equation*}
\int_{E_t\cap B_R}(t-u)^2|X\zeta|^2\,dx
\le
\frac{Ct^2}{(R-r)^2}|E_t\cap B_R|.
\end{equation*}
For the second term, the monotonicity of \(f\) gives
\[
f(u)(t-u)\le tf(t)\quad\text{on }E_t,
\]
and hence
\[
tf(t)\le \int_t^{2t}f(s)\,ds\le F(2t).
\]
Thus
\begin{equation*}
\int_{E_t\cap B_R}f(u)(t-u)\,dx
\le
F(2t)|E_t\cap B_R|.
\end{equation*}
Finally, since \(\zeta\equiv1\) on \(B_r\), we obtain
\begin{equation*}
\int_{E_t\cap B_r}|Xu|^2\,dx
\le
C\left(\frac{t^2}{(R-r)^2}+F(2t)\right)|E_t\cap B_R|,
\end{equation*}
which proves (\ref{-11}).  We prove the surface estimate (\ref{-22}) in
detail, since this is where the boundary term is produced.  The proof uses
three ingredients:


\begin{enumerate}[(i)]
\item Gauss--Green converts a boundary integral over \(\partial_H^*E_t\) into an integral over the set \(E_t\);
\item the identity \(\nu_{E_t}=Xu/|Xu|\), valid on the non-characteristic part of the level surface, turns the boundary integrand into a positive multiple of \(|Xu|\);
\item the two interior terms that remain are estimated by monotonicity of \(f\) and by the already proved solid Caccioppoli estimate (\ref{-11}).
\end{enumerate}
Let \(\mathcal T\subset(0,\infty)\) be the full-measure set supplied by
Lemma~\ref{lem:prop31-admissible-surface-levels}.  It is enough to prove
\((\ref{-22})\) for \(t\in\mathcal T\).  Fix such a \(t\).  If
\(E_t=\emptyset\), or if \(\Sigma_t=\{u=t\}\cap B_2\) is empty, then the
left-hand side of \((\ref{-22})\) has no interior level-surface contribution in
\(B_r\), and the estimate is trivial.  Thus we may assume that \(\Sigma_t\) is
nonempty.

By Lemma~\ref{lem:prop31-admissible-surface-levels}, the set \(E_t\) has
locally finite horizontal perimeter in \(B_2\), Gauss--Green applies to
compactly supported horizontal vector fields, and on every \(B_R\Subset B_2\)
the reduced boundary and the outer horizontal normal are identified by
\eqref{eq:prop31-admissible-level-normal}.  Set 

$$ \widehat R:=\frac{r+R}{2}, $$ 
and choose a cutoff function $$ \eta\in C_c^\infty(B_{\widehat R}) $$ such that $$ 0\le \eta\le 1,\qquad \eta\equiv 1\text{ on }B_r,\qquad |X\eta|\le \frac{C}{\widehat R-r}\le \frac{C}{R-r}. $$ 

\noindent Define the horizontal vector field
\[
Z:=\eta A(x)Xu.
\]
Then \(Z\in C_c^\infty(B_{\widehat R};\mathbb R^m)
\subset C_c^\infty(B_R;\mathbb R^m)\).  We will use the Gauss--Green identity
\eqref{eq:prop31-admissible-level-GG}, furnished by
Lemma~\ref{lem:prop31-admissible-surface-levels}, with this choice of \(Z\).
Here, \(\nu_{E_t}\) denotes the horizontal measure-theoretic outer unit normal
to \(E_t\).  Since \(Z=\eta A(x)Xu\),
\begin{equation}\label{gg-1}
\text{div}_H Z
=
\eta\,\text{div}_H(A(x)Xu)+\langle A(x)Xu,X\eta\rangle.
\end{equation}

\vspace{.2cm}
\noindent We recall that \(\text{div}_H(A(x)Xu)\) and \(f(u)\) are \(C^\infty\).
Therefore, since \(u\in C^\infty(B_2)\) is a supersolution of (\ref{meq}),
we have the following pointwise inequality in \(B_R\):
\begin{equation}\label{pw}
\text{div}_H\left(A(x)Xu(x)\right)\le f(u(x))
\qquad \forall\,x\in B_R.
\end{equation}

\noindent Using \eqref{gg-1} and \eqref{pw} in
\eqref{eq:prop31-admissible-level-GG} gives
 
\begin{align}
\int_{\partial_H^*E_t}\eta\langle A(x)Xu,\nu_{E_t}\rangle\,dP_H(E_t)&=\int_{E_t}\eta\,\text{div}_H(A(x)Xu)\,dx +\int_{E_t} \langle A(x)Xu,X\eta\rangle\,dx\notag\\[.2cm]
&\le \int_{E_t}\eta f(u)\,dx + \int_{E_t}\langle A(x)Xu,X\eta\rangle\,dx. \label{A-1-1}
\end{align}

%
%
\vspace{.2cm}
\noindent Since \(\operatorname{supp}\eta\subset B_{\widehat R}\Subset B_R\),
the boundary integral in \eqref{A-1-1} is supported in
\(\partial_H^*E_t\cap B_R\).  By the normal identity
\eqref{eq:prop31-admissible-level-normal} and horizontal ellipticity,
\begin{equation}\label{A-4}
\int_{\partial_H^*E_t}\eta\langle A(x)Xu,\nu_{E_t}\rangle\,dP_H(E_t)
\ge
\lambda \int_{\partial_H^*E_t\cap B_R} \eta|Xu|\,dP_H(E_t).
\end{equation}

\vspace{.2cm}
\noindent Combining \eqref{A-4} with \eqref{A-1-1}, and using
\(\eta\equiv1\) on \(B_r\), gives
\[
\lambda \int_{\partial_H^*E_t\cap B_r}|Xu|\,dP_H(E_t)
\le
\int_{E_t\cap B_{\widehat R}} f(u)\,dx
+
\int_{E_t\cap B_{\widehat R}} |\langle A(x)Xu,X\eta\rangle|\,dx.
\]

\vspace{.2cm}
\noindent Using $|A(x)Xu|\le \Lambda |Xu|$ and $|X\eta|\le C/(R-r)$, we get 

\begin{equation}\label{b-0}\int_{\partial_H^*E_t\cap B_r}|Xu|\,dP_H(E_t) \le C\int_{E_t\cap B_{\widehat R}} f(u)\,dx + \frac{C}{R-r} \int_{E_t\cap B_{\widehat R}} |Xu|\,dx. 
\end{equation} 

\vspace{.2cm}
\noindent We now estimate the two terms on the right-hand side of (\ref{b-0}). The first term is lower order and is controlled only by monotonicity of $f$. Since $u<t$ on $E_t$ and $f$ is nondecreasing, 
$$ \int_{E_t\cap B_{\widehat R}} f(u)\,dx \le f(t)|E_t\cap B_{\widehat R}| \le f(t)|E_t\cap B_R|. $$ \noindent Furthermore, $$ tf(t)\le \int_t^{2t}f(s)\,ds\le F(2t), $$ and hence 

\begin{equation}\label{t-0}\int_{E_t\cap B_{\widehat R}} f(u)\,dx \le \frac{F(2t)}{t}|E_t\cap B_R|. 
\end{equation}

\noindent For the second term in (\ref{b-0}), we first use Cauchy--Schwarz to reduce the $L^1$ norm of $Xu$ to the $L^2$ norm already estimated in (\ref{-11}):
\begin{equation}\label{t-1}\int_{E_t\cap B_{\widehat R}} |Xu|\,dx \le |E_t\cap B_{\widehat R}|^{1/2} \left( \int_{E_t\cap B_{\widehat R}} |Xu|^2\,dx \right)^{1/2}. 
\end{equation} 

\noindent We now apply the already proved solid Caccioppoli estimate (\ref{-11}) with inner radius $\widehat R$ and outer radius $R$. This is why the intermediate radius $\widehat R=(r+R)/2$ was introduced: it leaves room for the cutoff while keeping the final estimate in terms of $E_t\cap B_R$. Since 

$$R-\widehat R=\frac{R-r}{2}, $$ we obtain 
\begin{equation}\label{t-2} \int_{E_t\cap B_{\widehat R}} |Xu|^2\,dx \le C\left( \frac{t^2}{(R-r)^2}+F(2t) \right) |E_t\cap B_R|. 
\end{equation}

\vspace{.1cm}
\noindent Using (\ref{t-2}) in (\ref{t-1}) leads to  

\begin{equation}\label{t-3} \int_{E_t\cap B_{\widehat R}} |Xu|\,dx \le C\left( \frac{t}{R-r}+\sqrt{F(2t)} \right) |E_t\cap B_R|. 
\end{equation}

\noindent Substituting (\ref{t-0}) and (\ref{t-3}) into (\ref{b-0}) gives an estimate with one mixed term. Indeed, the factor $1/(R-r)$ in (\ref{b-0}) multiplies the two terms in (\ref{t-3}), and therefore

\begin{equation}\label{t-4}\int_{\partial_H^*E_t\cap B_r}|Xu|\,dP_H(E_t) \le C\left( \frac{F(2t)}{t} + \frac{t}{(R-r)^2} + \frac{\sqrt{F(2t)}}{R-r} \right) |E_t\cap B_R|. 
\end{equation}	

\vspace{.1cm}
\noindent Finally, the mixed term has the same size as the geometric mean of the two desired terms. More explicitly, by Young's inequality applied to
\[
 a=\frac{\sqrt t}{R-r},
 \qquad
 b=\sqrt{\frac{F(2t)}{t}},
\]
we get
$$
\frac{\sqrt{F(2t)}}{R-r}=ab
\le \frac12 a^2+\frac12 b^2
=\frac12\frac{t}{(R-r)^2}+\frac12\frac{F(2t)}{t}.
$$
\noindent Using this inequality in (\ref{t-4}) gives 
$$ \int_{\partial_H^*E_t\cap B_r}|Xu|\,dP_H(E_t) \le C\left( \frac{t}{(R-r)^2} + \frac{F(2t)}{t} \right) |E_t\cap B_R|. $$ 
\noindent This proves the surface integral estimate and completes the proof.\qed
}	
	
	\vspace{.3cm}

\noindent We next record a Sobolev--Poincar\'e consequence: if a function vanishes on most of a ball, then its $L^{Q/(Q-1)}$ norm on the ball is controlled by its horizontal gradient.

\begin{lem}\label{pin} Let $u\in W^{1,1}(B,X)$ and suppose $N:=\{u=0\}\cap B$ is such that $|N|\ge(1-\theta)|B|$ for some $0<\theta<1$. Then 
	$$\left(\int_B|u|^{\frac{Q}{Q-1}}\,dx\right)^{\frac{Q-1}{Q}}\le C(\theta,Q)\int_B|X u|\,dx.$$ 
\end{lem}
\noindent{\bf Proof:} We recall the following Sobolev-Poincar\'e inequality (see \cite[Theorem 1.5 and Remark 1.2]{GAN}).
$$\left(\frac{1}{|B|}\int_B|u-u_B|^{\frac{Q}{Q-1}}\,dx\right)^{\frac{Q-1}{Q}}\le \frac{C(Q)\,{\color{black}\text{diam}(B)}}{|B|}\int_B|X u|\,dx.$$
Note that
\begin{align*}|u_B|&=\frac{1}{|B|}\left|\int_Bu\,dx\right|\le \frac{1}{|B|}\int_N|u|\,dx+\frac{1}{|B|}\int_{B\setminus N}|u|\,dx\\
		&\le \frac{|B\setminus N|^{\frac1Q}}{|B|}\left(\int_{B\setminus N}|u|^{\frac{Q}{Q-1}}\,dx\right)^{\frac{Q-1}{Q}}\\
		&\le \frac{\theta^{\frac1Q}}{|B|^{\frac{Q-1}{Q}}}
		\left(\int_B|u|^{\frac{Q}{Q-1}}\,dx\right)^{\frac{Q-1}{Q}}.
\end{align*}
Therefore
$$|u_B||B|^{\frac{Q-1}{Q}}\le \theta^{\frac1Q}\|u\|_{L^{\frac{Q}{Q-1}}(B)}.$$
By Minkowski's inequality we see that
\begin{align*}\|u\|_{L^{\frac{Q}{Q-1}}(B)}&\le \|u-u_B\|_{L^{\frac{Q}{Q-1}}(B)}+\|u_B\|_{L^{\frac{Q}{Q-1}}(B)}\\
		&\le \frac{C(Q)\,{\color{black}\text{diam}(B)}}{|B|^{\frac 1Q}}\|X u\|_{L^1(B)}+|u_B|\,|B|^{\frac{Q-1}{Q}}\\
		&\le \frac{C(Q)\,{\color{black}\text{diam}(B)}}{|B|^{\frac 1Q}}\|X u\|_{L^1(B)} +\theta^{\frac1Q}\|u\|_{L^{\frac{Q}{Q-1}}(B)}.
\end{align*}
The result follows from this.\qed

\vspace{.2cm}
\noindent We now combine the preceding Sobolev--Poincar\'e estimate with the
Caccioppoli inequality \eqref{-11} to show that sufficiently small sublevel sets force a
positive lower bound in the smaller ball.

{\color{black}
\begin{lem}\label{del}
Let \(u\in W^{1,2}_{\text{loc}}(B_2,X)\) be a non-negative supersolution of
\eqref{meq}.  For \(s>0\), recall that
\[
E_s:=\{x\in B_2:u(x)<s\}.
\]
There exists a constant \(\delta_0=\delta_0(\lambda,\Lambda,Q)>0\), such that
if \(t>0\) satisfies
\[
|E_t|
\le
\delta_0\left(\frac{t}{\sqrt{F(2t)}+t}\right)^Q,
\]
then
\[
u\ge \frac{t}{2}\quad\text{a.e. in }B_1.
\]
Equivalently, \(\operatorname*{ess\,inf}_{B_1}u\ge t/2\); for the usual
representative this gives \(\inf_{B_1}u\ge t/2\).
\end{lem}

\noindent{\bf Proof:}
Let \(t>0\) satisfy the smallness assumption.  Choose \(\delta_0\) so small
that
\[
\delta_0\le \frac{|B_1|}{2}.
\]
Since \(t/(\sqrt{F(2t)}+t)\le 1\), we have \(|E_t|\le |B_1|/2\).  Fix
\(0<h<k\le t\) and define the truncation
\[
v:=\min\{(k-u)_+,k-h\}.
\]
Thus \(v=0\) on \(B_2\setminus E_k\), \(v=k-h\) on \(E_h\), and
\(|Xv|\le |Xu|\chi_{E_k}\).  If \(1\le \rho<R\le 2\), then
\(E_k\subset E_t\), and hence
\[
|E_k\cap B_\rho|\le |E_t|\le \frac{|B_1|}{2}\le \frac{|B_\rho|}{2}.
\]
Consequently \(|\{v=0\}\cap B_\rho|\ge |B_\rho|/2\), so Lemma~\ref{pin}
applied in \(B_\rho\) gives
\[
\left(\int_{B_\rho}v^{\frac{Q}{Q-1}}\,dx\right)^{\frac{Q-1}{Q}}
\le
C\int_{B_\rho}|Xv|\,dx.
\]
Since \(v=k-h\) on \(E_h\cap B_\rho\), it follows that
\[
(k-h)|E_h\cap B_\rho|^{\frac{Q-1}{Q}}
\le
C\int_{E_k\cap B_\rho}|Xu|\,dx.
\]
By Cauchy's inequality and the energy estimate \eqref{-11}, applied with level
\(k\), inner radius \(\rho\), and outer radius \(R\),
\begin{align*}
\int_{E_k\cap B_\rho}|Xu|\,dx
&\le |E_k\cap B_\rho|^{1/2}
\left(\int_{E_k\cap B_\rho}|Xu|^2\,dx\right)^{1/2}\\
&\le C\left(\frac{k}{R-\rho}+\sqrt{F(2k)}\right)|E_k\cap B_R|.
\end{align*}
Using also \(E_h\subset E_k\), we obtain the De Giorgi estimate
\[
|E_h\cap B_\rho|
\le
C\left(\frac{k}{(k-h)(R-\rho)}+\frac{\sqrt{F(2k)}}{k-h}\right)
|E_k\cap B_R|^{1+\frac1Q}.
\]

Now set
\[
r_j:=1+2^{-j},
\qquad
k_j:=\frac{t}{2}(1+2^{-j}),
\qquad
x_j:=|E_{k_j}\cap B_{r_j}|.
\]
Applying the preceding estimate with
\(R=r_j\), \(\rho=r_{j+1}\), \(k=k_j\), and \(h=k_{j+1}\), and using
\(2k_j=t(1+2^{-j})\le 2t\), yields
\[
x_{j+1}
\le
C_0\left(1+\frac{\sqrt{F(2t)}}{t}\right)4^j x_j^{1+\frac1Q}.
\]
Finally choose, in addition,
\[
\delta_0\le C_0^{-Q}4^{-Q^2}.
\]
Then the smallness assumption gives
\[
x_0=|E_t|
\le
\left[C_0\left(1+\frac{\sqrt{F(2t)}}{t}\right)\right]^{-Q}4^{-Q^2}.
\]
Lemma~\ref{eg} therefore implies \(x_j\to0\).  Since
\(E_{t/2}\cap B_1\subset E_{k_j}\cap B_{r_j}\) for every \(j\), we get
\(|E_{t/2}\cap B_1|=0\).  Hence \(u\ge t/2\) a.e. in \(B_1\), as claimed.
\qed
}

\vspace{.2cm}
{\color{black}
\noindent \textcolor{black}{We record the perimeter notation and regular-level
facts used in the density estimates below.}  If $E\subset\mathbb G$ has finite horizontal
perimeter, then, in the notation of
\cite[Definition~3.2]{GutierrezMohammedCompanion},
\(P_H(E)=|D_H\chi_E|\).  This is a Radon measure, and
\cite[Proposition~3.5]{GutierrezMohammedCompanion} shows that it is
concentrated on the horizontal reduced boundary \(\partial_H^*E\).  Hence, for
every ball $B(x,r)$,
\begin{equation}\label{restrictionofmeasure}
P_H(E;B(x,r))
=\int_{\partial_H^*E\cap B(x,r)}dP_H(E).
\end{equation}
Moreover, if $u\in C^\infty(B_2)$, $E_t:=\{x\in B_2:u(x)<t\}$, and $t$ is a
Euclidean regular value of $u$ in the sense of
\cite[Definition~3.8]{GutierrezMohammedCompanion}, then
\cite[Proposition~3.9]{GutierrezMohammedCompanion} gives, up to
$P_H(E_t)$-null sets,
\(\partial_H^*E_t\subset\{u=t\}\) and the characteristic part carries no
horizontal perimeter measure:
\begin{equation}\label{Xu}
P_H(E_t)\left(\partial_H^*E_t\cap\{x\in B_2: u(x)=t,\,|Xu(x)|=0\}\right)=0.
\end{equation}
The admissible-level lemma, Lemma~\ref{lem:prop31-admissible-surface-levels},
ensures that these facts are available for almost every level $t>0$ used in
\textcolor{black}{the density estimates below.}
}

\vspace{.2cm}

\vspace{.2cm}

{\color{black}
We recall the notation used in the density estimates below.  Let \(\mathbb G\)
be a Carnot group with identity \(e\), homogeneous dimension \(Q\), Haar
measure \(|\cdot|\), and Carnot--Carath\'eodory distance \(d=d_{cc}\).  Put
\[
N(g):=d(e,g),
\qquad
B_2:=B(e,2),
\qquad
\rho(g):=\operatorname{dist}_{cc}(g,\mathbb G\setminus B_2).
\]
Fix \(0<\delta<1\), and define the truncated intrinsic boundary distance by
\[
d_\delta(g):=\kappa\min\{\rho(g),\delta\},
\qquad g\in B_2,
\]
where \(\kappa\in(0,1/100]\) is a structural constant, depending only on
\(\mathbb G\), to be fixed once and for all.  Notice that
\(\rho(g)\ge 2-N(g)\) by the triangle inequality.  Thus
\(d_\delta=\kappa\delta\) on \(B_{2-\delta}\), while
\(B(g,2d_\delta(g))\Subset B_2\).  Given a measurable set
\(A\subset B_2\), set
	\begin{equation}\label{eq:sigma-A-density-definition}
	\sigma_A(g):=\frac{|A\cap B(g,d_\delta(g))|}{|B(g,d_\delta(g))|},
	\qquad g\in B_2.
	\end{equation}
}

{\color{black}
Unlike Julin's Euclidean argument, which uses the radial gauge deficit
\(2-N(g)\) as a boundary scale, the Carnot-group argument needs the actual
Carnot--Carath\'eodory distance
\(\rho(g)=\operatorname{dist}_{cc}(g,\mathbb G\setminus B_2)\).  The gauge
deficit is not interchangeable with this intrinsic distance: as shown by
Remark~\ref{rem:heisenberg-cap-gauge-deficit}, for a Heisenberg cap the
gauge-deficit radius \(\kappa\min\{2-N(g),\delta\}\) does not give the
boundary-concentration estimate proved in Lemma~\ref{bce} (the truncated
horizontal boundary-concentration lemma).  Thus the use of \(\rho\) in
\(d_\delta(g)=\kappa\min\{\rho(g),\delta\}\) is essential for the proof of
Lemma~\ref{bce}.

The next lemma consolidates the two metric facts about \(d_\delta\) needed
below: its Lipschitz bound and the resulting comparison of balls centered at
nearby points.
}

{\color{black}
\begin{lem}\label{lem:ddelta-metric-facts}
Let \((X,d)\) be a metric space, let \(\Omega\subset X\) be open, and put
\(\rho_\Omega(g):=\operatorname{dist}(g,X\setminus\Omega)\).  Given
\(\delta>0\) and \(0<\kappa<1\), define
\[
d_\delta(g):=\kappa\min\{\rho_\Omega(g),\delta\},
\qquad g\in\Omega.
\]
Then the following hold.
\begin{enumerate}[(i)]
\item For all \(g,\xi\in\Omega\),
\begin{equation}\label{eq:ddelta-lipschitz}
|d_\delta(g)-d_\delta(\xi)|\le \kappa d(g,\xi).
\end{equation}

\item Fix \(g,\xi\in\Omega\), set
\[
s:=(1+\kappa)d(g,\xi),
\]
and assume that
\[
d(g,\xi)\le \frac{d_\delta(g)}{2(1+\kappa)}.
\]
Then
\begin{equation}\label{eq:ddelta-ball-containment}
B(g,d_\delta(g)-s)
\subset B(\xi,d_\delta(\xi))
\subset B(g,d_\delta(g)+s).
\end{equation}
Consequently,
\begin{equation}\label{eq:ddelta-ball-symmetric-difference}
B(g,d_\delta(g))\Delta B(\xi,d_\delta(\xi))
\subset
B(g,d_\delta(g)+s)\setminus B(g,d_\delta(g)-s),
\end{equation}
where \(\Delta\) denotes symmetric difference.
\end{enumerate}
\end{lem}

\begin{proof}
The distance-to-a-closed-set function \(\rho_\Omega\) is \(1\)-Lipschitz, and
for every
\(a,b\in\mathbb R\),
\[
|\min\{a,\delta\}-\min\{b,\delta\}|\le |a-b|.
\]
Applying this with \(a=\rho_\Omega(g)\) and
\(b=\rho_\Omega(\xi)\) gives
\[
\begin{aligned}
|d_\delta(g)-d_\delta(\xi)|
&=\kappa\left|
\min\{\rho_\Omega(g),\delta\}
-\min\{\rho_\Omega(\xi),\delta\}
\right| \\
&\le \kappa |\rho_\Omega(g)-\rho_\Omega(\xi)|
\le \kappa d(g,\xi),
\end{aligned}
\]
which proves \eqref{eq:ddelta-lipschitz}.

We now prove the ball inclusions.  Write
\[
r:=d_\delta(g),
\qquad
\rho:=d_\delta(\xi),
\qquad
h:=d(g,\xi).
\]
Then \(s=(1+\kappa)h\), and \eqref{eq:ddelta-lipschitz} gives
\[
r-\kappa h\le \rho\le r+\kappa h.
\]
The assumption on \(h\) implies \(s\le r/2\), so \(r-s>0\).  If
\(x\in B(g,r-s)\), then
\[
d(x,\xi)
\le d(x,g)+h
< r-s+h
=r-\kappa h
\le \rho,
\]
and hence \(B(g,r-s)\subset B(\xi,\rho)\).  In addition, if
\(x\in B(\xi,\rho)\), then
\[
d(x,g)
\le d(x,\xi)+h
< \rho+h
\le r+\kappa h+h
=r+s,
\]
and hence \(B(\xi,\rho)\subset B(g,r+s)\).  This proves
\eqref{eq:ddelta-ball-containment}.

Finally, set
\[
A:=B(g,r),
\qquad
C:=B(\xi,\rho),
\qquad
A_-:=B(g,r-s),
\qquad
A_+:=B(g,r+s).
\]
The inclusions just proved give \(A_-\subset C\subset A_+\), and since
\(s>0\), also \(A_-\subset A\subset A_+\).  Therefore
\(A_-\subset A\cap C\) and \(A\cup C\subset A_+\), which implies
\[
A\Delta C\subset A_+\setminus A_-.
\]
Returning to the original notation gives
\eqref{eq:ddelta-ball-symmetric-difference}.
\end{proof}
}

We now apply the preceding metric comparison to the localized density
\(\sigma_A\); the next lemma gives the horizontal gradient bound needed in the
coarea and isoperimetric estimates below.

\begin{lem}\label{gradientofsigmaA}
{\color{black}Let $\mathbb G$ be a Carnot group with homogeneous dimension $Q$, endowed
with its Carnot--Carath\'eodory distance $d=d_{cc}$.  Let
$B_2:=B(e,2)$, put
\[
\rho(g):=\operatorname{dist}_{cc}(g,\mathbb G\setminus B_2),
\]
and define
\[
d_\delta(g)=\kappa \min\{\rho(g),\delta\},
\qquad 0<\kappa<1.
\]}
For a measurable set $A\subset \mathbb G$, define
\[
\sigma_A(g):=
\frac{|A\cap B(g,d_\delta(g))|}{|B(g,d_\delta(g))|},
\qquad g\in B_2.
\]
Then there is a constant $C>0$, {\color{black}depending only on the fixed
Carnot-group structure, the homogeneous distance $d$, and the fixed inner
product on the horizontal layer $\mathcal V_1$}, such that
\[
|X\sigma_A(g)|\leq \frac{C}{d_\delta(g)}
\]
for a.e. $g\in B_2$.
\end{lem}

\begin{proof}
{\color{black}Throughout the proof, $C$ denotes a positive constant depending
only on the fixed Carnot-group structure, the homogeneous distance $d$, and
the fixed inner product on $\mathcal V_1$, possibly changing from line to line.  In
particular, $C$ is independent of $A$, $\delta$, $\kappa$, and $g$. Since
$0<\kappa<1$, all constants below may be chosen independently of $\kappa$.}

We first prove a general estimate comparing the balls
\[
B(g,d_\delta(g))
\qquad \text{and} \qquad
B(g',d_\delta(g')).
\]
{\color{black}Set
\[
r:=d_\delta(g),\qquad r':=d_\delta(g'),
\qquad \ell:=d(g,g'),
\qquad s:=(1+\kappa)\ell.
\]}
By the Lipschitz estimate for $d_\delta$, we have
\[
|d_\delta(g)-d_\delta(g')|\leq \kappa d(g,g').
\]
{\color{black}Assume that
\begin{equation}\label{dgg'}
d(g,g')\leq \frac{d_\delta(g)}{2(1+\kappa)}.
\end{equation}

Then
\[
(1+\kappa)d(g,g')\leq \frac r2.
\]
By Lemma~\ref{lem:ddelta-metric-facts},
\[
B(g,r)\Delta B(g',r')
\subset
B(g,r+s)\setminus B(g,r-s).
\]
}
Taking measures {\color{black} and using \eqref{eq:proposed-prelim-ball-scaling}},
we obtain
\[
\begin{aligned}
|B(g,r)\Delta B(g',r')|
&\leq |B(g,r+s)\setminus B(g,r-s)| \\
&= \omega_Q\bigl((r+s)^Q-(r-s)^Q\bigr).
\end{aligned}
\]
Since $s\le r/2$, the mean value theorem gives
\[
(r+s)^Q-(r-s)^Q
\leq C r^{Q-1}s.
\]
Therefore
\begin{equation}\label{eq:ball-sym-diff-general}
|B(g,r)\Delta B(g',r')|
\leq C r^{Q-1}s
\leq C r^{Q-1} d(g,g').
\end{equation}

{\color{black}Now define
\[
\phi(\xi):=|A\cap B(\xi,d_\delta(\xi))|,
\qquad
\psi(\xi):=|B(\xi,d_\delta(\xi))|
=\omega_Qd_\delta(\xi)^Q.
\]
We claim that
\begin{equation}\label{eq:phi-psi-general}
|\phi(g')-\phi(g)|
+
|\psi(g')-\psi(g)|
\leq C d_\delta(g)^{Q-1} d(g,g').
\end{equation}
whenever \eqref{dgg'} holds.}
Indeed, using the elementary inequality
\[
\bigl||E|-|F|\bigr|\leq |E\Delta F|,
\]
we get
\[
\begin{aligned}
|\phi(g')-\phi(g)|
&=
\bigl|
|A\cap B(g',d_\delta(g'))|
-
|A\cap B(g,d_\delta(g))|
\bigr| \\
&\leq
\bigl|
(A\cap B(g',d_\delta(g')))
\Delta
(A\cap B(g,d_\delta(g)))
\bigr| \\
&\leq
|B(g',d_\delta(g'))\Delta B(g,d_\delta(g))|.
\end{aligned}
\]
Thus, by \eqref{eq:ball-sym-diff-general},
\[
|\phi(g')-\phi(g)|
\leq C d_\delta(g)^{Q-1}d(g,g').
\]
Similarly,
\[
\begin{aligned}
|\psi(g')-\psi(g)|
&=
\bigl|
|B(g',d_\delta(g'))|
-
|B(g,d_\delta(g))|
\bigr| \\
&\leq
|B(g',d_\delta(g'))\Delta B(g,d_\delta(g))| \\
&\leq C d_\delta(g)^{Q-1}d(g,g').
\end{aligned}
\]
This proves \eqref{eq:phi-psi-general}.

We next derive a local Lipschitz estimate for $\sigma_A$. Since
\[
\sigma_A=\frac{\phi}{\psi},
\]
we write
\[
\begin{aligned}
\sigma_A(g')-\sigma_A(g)
&=
\frac{\phi(g')}{\psi(g')}
-
\frac{\phi(g)}{\psi(g)} \\
&=
\frac{\phi(g')-\phi(g)}{\psi(g')}
+
\phi(g)
\left(
\frac{1}{\psi(g')}-\frac{1}{\psi(g)}
\right).
\end{aligned}
\]
We need a lower bound for $\psi(g')$. 
By \eqref{eq:phi-psi-general},
\[
|\psi(g')-\psi(g)|
\leq C r^{Q-1}d(g,g').
\]
{\color{black}If $d(g,g')\leq cr$, then
\[
\begin{aligned}
\psi(g')
&\geq \psi(g)-|\psi(g')-\psi(g)| \\
&\geq \omega_Qr^Q-Cr^{Q-1}d(g,g') \\
&\geq (\omega_Q-Cc)r^Q.
\end{aligned}
\]
We choose
\[
0<c\leq \min\left\{\frac14,\frac{\omega_Q}{2C}\right\}.
\]
The condition $c\leq1/4$ guarantees \textcolor{black}{that \eqref{dgg'} holds}
uniformly for $0<\kappa<1$, while
\(c\leq\omega_Q/(2C)\) gives
\[
\psi(g')\geq \frac12\omega_Qr^Q.
\]
Since \(\omega_Q=|B(e,1)|\) and the constant \(C\) in
\eqref{eq:phi-psi-general} depend only on the fixed Carnot-group structure
and the homogeneous distance, the same is true of $c$.  In particular,
$c$ is independent of $A,\delta,\kappa,g$, and $g'$.}

\noindent Also,
\[
0\leq \phi(g)\leq \psi(g)=\omega_Q r^Q.
\]
Therefore, \textcolor{black}{by \eqref{eq:phi-psi-general}},
\[
\begin{aligned}
|\sigma_A(g')-\sigma_A(g)|
&\leq
\frac{|\phi(g')-\phi(g)|}{\psi(g')}
+
\phi(g)
\frac{|\psi(g')-\psi(g)|}{\psi(g')\psi(g)} \\
&\leq
C\frac{r^{Q-1}d(g,g')}{r^Q}
+
C r^Q
\frac{r^{Q-1}d(g,g')}{r^Q r^Q} \\
&\leq
C\frac{d(g,g')}{r}.
\end{aligned}
\]
In other words, if $g'\in B_2$ and $d(g,g')\leq c d_\delta(g)$, then
\begin{equation}\label{eq:sigma-local-lip}
|\sigma_A(g')-\sigma_A(g)|
\leq
C\frac{d(g,g')}{d_\delta(g)}.
\end{equation}
This also shows that $\sigma_A$ is locally Lipschitz in $B_2$, because
$d_\delta(g)>0$ for every $g\in B_2$.

We now apply this estimate to horizontal curves. Fix
\[
1\leq j\leq m
\]
and consider
\[
g_h:=g\exp(hX_j),
\]
where $|h|$ is small. Since $d$ is left-invariant and homogeneous,
\[
\begin{aligned}
d(g,g_h)
&=
d(e,g^{-1}g_h) \\
&=
d(e,\exp(hX_j)).
\end{aligned}
\]
Since $X_j$ belongs to the first layer {\color{black}$\mathcal V_1$} of the stratified Lie algebra, the group dilations satisfy
\[
\delta_\lambda(\exp X_j)=\exp(\lambda X_j),\qquad \lambda>0.
\]
Therefore, taking $\lambda=|h|$ and replacing $X_j$ by $(\operatorname{sgn}h)X_j$, we obtain
\[
\delta_{|h|}\bigl(\exp((\operatorname{sgn}h)X_j)\bigr)
=
\exp\bigl(|h|(\operatorname{sgn}h)X_j\bigr)
=
\exp(hX_j).
\]

{\color{black}
Let $|\cdot|_{\textcolor{black}{\mathcal V_1}}$ denote the norm induced by the fixed inner product on
$\textcolor{black}{\mathcal V_1}$.  Since the map
\[
V\longmapsto d(e,\exp V)
\]
is continuous and the unit sphere of \textcolor{black}{$\mathcal V_1$} is compact, the constant
\[
C_d:=\max_{|V|_{\textcolor{black}{\mathcal V_1}}=1}d(e,\exp V)
\]
is finite.  Because the fields $X_1,\dots,X_m$ are orthonormal, both $X_j$
and $-X_j$ belong to this unit sphere.  Hence
\[
\begin{aligned}
d(g,g_h)
&=
d(e,\exp(hX_j)) \\
&=
d\left(\delta_{|h|}(e),
\delta_{|h|}\bigl(\exp((\operatorname{sgn}h)X_j)\bigr)\right) \\
&=
|h|\, d\left(e,\exp((\operatorname{sgn}h)X_j)\right) \\
&\leq C_d |h|.
\end{aligned}
\]
Therefore, for $|h|$ sufficiently small depending on
$r=d_\delta(g)$, namely for
\[
C_d|h|\leq c r,
\]
}
we may apply \eqref{eq:sigma-local-lip} with $g'=g_h$. We obtain
\[
|\sigma_A(g_h)-\sigma_A(g)|
\leq
C\frac{d(g,g_h)}{d_\delta(g)}
\leq
C\frac{|h|}{d_\delta(g)}.
\]
Thus
\[
\left|
\frac{\sigma_A(g\exp(hX_j))-\sigma_A(g)}{h}
\right|
\leq
\frac{C}{d_\delta(g)}
\]
for all sufficiently small nonzero $h$.

{\color{black}Since $\sigma_A$ is locally Lipschitz in $B_2$, the Carnot-group
Rademacher theorem gives the existence of its horizontal derivatives almost
everywhere; see, for example, \cite{VAM}.} Passing to the limit as $h\to 0$ at points where
$X_j\sigma_A(g)$ exists, we obtain
\[
|X_j\sigma_A(g)|
\leq
\frac{C}{d_\delta(g)}
\]
for a.e. $g\in B_2$.

Finally, since this estimate holds for each $j=1,\dots,m$, we have
\[
|X\sigma_A(g)|
=
\left(\sum_{j=1}^m |X_j\sigma_A(g)|^2\right)^{1/2}
\leq
\frac{C}{d_\delta(g)}
\]
for a.e. $g\in B_2$, after absorbing the factor $\sqrt m$ into the
constant $C$. This proves the lemma.
\end{proof}

\vspace{.1cm}
We also need another preparatory lemma.

{\color{black}
\noindent The proof of \cite[Lemma~3.4]{VJU} uses Euclidean coarea and
Gauss--Green formulas, including a radial vector field to control the weighted
perimeter near the boundary.  In a Carnot group, the corresponding formulas
involve horizontal perimeter \(P_H\).  Since in general
\(P_H(E;\cdot)\le P(E;\cdot)\), a lower bound for the larger measure \(P\)
does not give a lower bound for \(P_H\).  For this reason, Lemma~\ref{bce} is
formulated directly in terms of horizontal perimeter; its intrinsic boundary
scale was explained above.
}

\begin{lem}[Truncated horizontal boundary concentration]\label{bce}
Let \(\mathbb G\) be a Carnot group of homogeneous dimension \(Q\), endowed
with its Carnot--Carath\'eodory distance \(d=d_{cc}\) and Haar measure.  Set
\[
B_2:=B(e,2),
\qquad
\rho(g):=\operatorname{dist}_{cc}(g,\mathbb G\setminus B_2),
\qquad
d_\delta(g):=\kappa\min\{\rho(g),\delta\},
\]
where \(0<\delta\le1\) and \(0<\kappa\le1/100\).  If
\(E\subset B_2\) has locally finite horizontal perimeter in \(B_2\), put
\[
m_E:=\min\{|E|,\ |B_2\setminus E|\}.
\]
Then
\begin{equation}\label{eq:horizontal-boundary-concentration}
\int_{\partial_H^*E\cap B_2}d_\delta(g)\,dP_H(E)(g)
\ge
c\min\left\{
\delta m_E^{\frac{Q-1}{Q}},\ m_E
\right\},
\end{equation}
{\color{black}where \(c>0\) depends only on the fixed Carnot group, the chosen
horizontal inner product and Haar measure, and \(\kappa\).}
\end{lem}

\begin{proof}
If \(m_E=0\), there is nothing to prove.
{\color{black}By the standing assumption \(Q>1\), we may set}
\[
q:=\frac{Q}{Q-1},
\qquad
w_a(g):=\min\{\rho(g),a\},
\qquad 0<a\le1.
\]

We prove the following capped estimate:
\begin{equation}\label{eq:rigorous-capped-estimate}
\int_{\partial_H^*E\cap B_2}w_a\,dP_H(E)
\ge
c\,a\,m_E^{(Q-1)/Q}
\qquad\text{whenever }a\le m_E^{1/Q}.
\end{equation}

\medskip
\noindent\emph{Step 1: Boman chains and the pointwise representation.}

Put \(\Omega=B_2\).  The Carnot--Carath\'eodory metric is a complete
geodesic metric.  Thus, if \(x\in B(g,r)\), a minimizing geodesic from
\(g\) to \(x\) remains in \(B(g,r)\) and satisfies the segment identity
required in \cite[(3.8)]{FLWRepresentation}.

{\color{black}
The representation formula used here is precisely the ``Moreover'' part of
\cite[Proposition~2.12]{FLWRepresentation}, together with the potential
operator in \cite[(2.11)]{FLWRepresentation} defined by
\begin{equation*}\tag{FLW}\label{eq:FLW-potential-operator}
(\mathbf T\psi)(g)
:=
\int_{\mathbb G}
\psi(h)
\frac{d(g,h)}{|B(g,d(g,h))|}\,dh.
\end{equation*}
The cited statement gives, for a control ball \(B\) and a locally Lipschitz
function \(u\),
\begin{equation*}\tag{FLW-R}\label{eq:FLW-representation-inequality}
|u(g)-u_B|
\le
\mathbf C
\mathbf T\bigl(|Xu|\chi_B\bigr)(g),
\qquad g\in B,
\end{equation*}
where \(u_B=|B|^{-1}\int_Bu\).  By
\eqref{eq:proposed-prelim-ball-scaling}, with the earlier notation
\(\omega_Q=|B(e,1)|\), and since \(\chi_B\) restricts the integral to \(B\), \eqref{eq:FLW-potential-operator} leads to
\[
\mathbf T\bigl(|Xu|\chi_B\bigr)(g)
 =\frac1{\omega_Q}\int_B
   \frac{|Xu(h)|}{d(g,h)^{Q-1}}\,dh.
\]
The kernel is locally integrable by the same volume scaling.  Thus, with
\(C_{\mathrm{rep}}:=\mathbf C/\omega_Q\),
\eqref{eq:FLW-representation-inequality} becomes
\begin{equation}\label{eq:FLW-local-kernel-representation}
|u(g)-u_B|
\le
C_{\mathrm{rep}}\int_B
\frac{|Xu(h)|}{d(g,h)^{Q-1}}\,dh,
\qquad g\in B.
\end{equation}
Although
\cite[Proposition~2.12]{FLWRepresentation} is formulated locally, left
translations and homogeneous dilations reduce every Carnot-group ball to a
fixed unit ball, so
\(C_{\mathrm{rep}}\) is uniform over all balls.

The Poincar\'e estimate intended below is the local horizontal
\(L^1\)-Poincar\'e inequality: there are structural constants
\(\Lambda>1\) and \(C_{\mathrm P}>0\) such that, for every control ball
\(B=B(g_B,r_B)\) and every locally Lipschitz function \(v\) on
\(\Lambda B\),
\begin{equation}\label{eq:rigorous-local-L1-poincare}
\int_B |v-v_B|\,dg
\le
C_{\mathrm P}r_B\int_{\Lambda B}|Xv|\,dg,
\qquad
v_B:=\frac{1}{|B|}\int_Bv\,dg.
\end{equation}
This is the Carnot-group form of the Poincar\'e inequality of
\cite{JerisonPoincare}; by left translation and homogeneous dilation, its
constants are independent of the center and radius of \(B\).  Fix such a
number \(\Lambda\).  Finally, fix one number
\(\lambda>\Lambda\), for instance \(\lambda=2\Lambda\).
{\color{black}These are geometric constants for the ball dilation and Boman
covering argument described below, unrelated to the ellipticity constants of the coefficient
matrix \(A\).}

Theorem~3.9 of \cite{FLWRepresentation} states that a control ball satisfies
the Boman condition for any prescribed parameter \(\tau>1\).  Its proof follows
the argument of \cite[Theorem~5.4]{FranchiGutierrezWheedenGrushin}.  Apply
Theorem~3.9 with \(\tau=\lambda\).  For this fixed \(\lambda\), the theorem gives a
constant \(M=M(\lambda)\ge1\), a countable covering \(\mathcal W\) of
\(\Omega=B_2\) by balls \(B=B(g_B,r_B)\), and a distinguished ball \(B_*\).
The covering and the constant \(M\) may depend on the prescribed value of
\(\lambda\); the assertion is not that one fixed covering works
simultaneously for every \(\lambda>1\).  For a covering ball
\(B=B(g_B,r_B)\in\mathcal W\) and a number \(a>0\), we use the standard
notation
\[
aB:=B(g_B,ar_B),
\]
so the center is unchanged and only the radius is multiplied by \(a\).
After increasing \(M\), if necessary, these objects have all of the following
properties.  First, the Boman covering $\mathcal W$ satisfies
\begin{equation}\label{eq:rigorous-boman-overlap}
\lambda B\subset\Omega
\quad\text{for every }B\in\mathcal W,
\qquad
\sum_{B\in\mathcal W}\chi_{\lambda B}\le M
\quad\text{a.e. in }\Omega.
\end{equation}
In particular, since \(\Lambda<\lambda\),
\[
\Lambda B\subset\lambda B\subset\Omega,
\qquad
\sum_{B\in\mathcal W}\chi_{\Lambda B}\le M
\quad\text{a.e. in }\Omega.
\]
}
Second, for every ball \(B\in\mathcal W\), there is a finite chain, which
we take without repeated balls,
\[
\mathcal C(B)=(B=B_0,B_1,\ldots,B_{\ell(B)}=B_*).
\]
For every \(j<\ell(B)\), there is a ball
\(R_j\subset B_j\cap B_{j+1}\) such that
\begin{equation}\label{eq:rigorous-boman-bridge}
B_j\cup B_{j+1}\subset M R_j,
\end{equation}
and every ball \(J\in\mathcal C(B)\) satisfies
\begin{equation}\label{eq:rigorous-boman-shadow-property}
B\subset M J.
\end{equation}
The inclusions \(R_j\subset B_j\cap B_{j+1}\) and
\eqref{eq:rigorous-boman-bridge} first show that the radii of
\(R_j,B_j,B_{j+1}\) are comparable, with constants depending only on
\(M\).  Applying the homogeneous ball-volume formula
\(|B(g,r)|=\omega_Qr^Q\), equivalently the Ahlfors \(Q\)-regularity of Haar
measure, then gives
\begin{equation}\label{eq:rigorous-bridge-measures}
|R_j|\asymp|B_j|\asymp|B_{j+1}|.
\end{equation}

{\color{black}
For \(h\in\Lambda B\), let \(z\in\mathbb G\setminus\Omega\).  Since
\(\lambda B=B(g_B,\lambda r_B)\subset\Omega\), the point \(z\) does not
belong to \(\lambda B\), and therefore
\[
d(g_B,z)\ge\lambda r_B.
\]
On the other hand, \(h\in\Lambda B=B(g_B,\Lambda r_B)\), so
\(d(g_B,h)<\Lambda r_B\).  The reverse triangle inequality now gives
\[
d(h,z)
\ge d(g_B,z)-d(g_B,h)
>(\lambda-\Lambda)r_B.
\]
Taking the infimum over all \(z\in\mathbb G\setminus\Omega\) yields
\[
\rho(h)
=\operatorname{dist}(h,\mathbb G\setminus\Omega)
\ge (\lambda-\Lambda)r_B, \quad \forall h\in \Lambda B.
\]
Consequently
\begin{equation}\label{eq:rigorous-radius-below-rho}
r_B\le \frac{1}{\lambda-\Lambda}\rho(h)
\qquad\text{for every }h\in\Lambda B.
\end{equation}
Since \(B\subset\Lambda B\) and the integrand is nonnegative,
\eqref{eq:FLW-local-kernel-representation} gives, for every Lipschitz
function \(u\),
\begin{equation}\label{eq:rigorous-local-representation}
|u(g)-u_B|
\le
C_{\mathrm{rep}}\int_{\Lambda B}
\frac{|Xu(h)|}{d(g,h)^{Q-1}}\,dh,
\qquad g\in B,
\end{equation}
In particular, \(C_{\mathrm{rep}}\) is independent of \(B\), \(g\), and
\(u\).
}

If \(J,K\) are consecutive chain balls and \(R\) is the bridge
ball\footnote{A bridge ball is an auxiliary Carnot--Carath\'eodory ball,
not necessarily a member of \(\mathcal W\), such that
\(R\subset J\cap K\) and \(J\cup K\subset MR\), where \(MR\) denotes the
ball concentric with \(R\) whose radius is \(M\) times that of \(R\).} from
\eqref{eq:rigorous-boman-bridge}, then
\eqref{eq:rigorous-local-L1-poincare}, together with
\eqref{eq:rigorous-bridge-measures} and \(R\subset J\cap K\), gives a
constant \(C_{\mathrm{av}}>0\), depending only on the fixed Carnot-group
structure, \(\Lambda\), and the Boman constant \(M\), such that
\begin{equation}\label{eq:rigorous-average-difference}
|u_J-u_K|
\le
C_{\mathrm{av}} r_J^{1-Q}
\left(
\int_{\Lambda J}|Xu|
+
\int_{\Lambda K}|Xu|
\right).
\end{equation}
Indeed, one first bounds \(|u_J-u_R|\) by
\[
|u_J-u_R|
\le
\frac{1}{|R|}\int_R|u-u_J|
\le
\frac{C}{|J|}\int_J|u-u_J|,
\]
and then applies \eqref{eq:rigorous-local-L1-poincare} with \(B=J\):
\[
\int_J|u-u_J|
\le
C_{\mathrm P}r_J\int_{\Lambda J}|Xu|.
\]
Since \(|J|\simeq r_J^Q\), this gives the \(J\)-term in
\eqref{eq:rigorous-average-difference}; the estimate involving \(K\) is
identical.  In particular, \(C_{\mathrm{av}}\) is independent of the
consecutive balls \(J,K\) and of \(u\).

Enumerate \(\mathcal W=\{B^1,B^2,\ldots\}\) and form the measurable
disjoint partition
\[
A_{B^1}:=B^1,
\qquad
A_{B^k}:=B^k\setminus\bigcup_{i<k}B^i
\quad(k\ge2).
\]
Thus \(A_B\subset B\), the sets \(A_B\) are pairwise disjoint, and their
union is \(\Omega\).  For nonnegative \(\varphi\in L^1(\Omega)\), define
\[
\mathcal R_B\varphi(g)
:=
C_{\mathrm{rep}}\int_{\Lambda B}
\frac{\varphi(h)}{d(g,h)^{Q-1}}\,dh,
\qquad g\in B,
\]
where the integral is initially understood as a nonnegative extended-valued
measurable function.  The kernel is locally integrable: the homogeneous
ball-volume formula and a decomposition into dyadic annuli give
\[
\int_{B(g,R)}\frac{dh}{d(g,h)^{Q-1}}
\le CR,
\qquad R>0.
\]
Equivalently, this is the case \(\alpha=1\) of the local-integrability
criterion for kernels homogeneous of degree \(\alpha-Q\); see
\cite[Chapter~1, Corollary~(1.17)]{FollandSteinHardy}.  Tonelli's theorem
therefore shows that \(\mathcal R_B\varphi(g)<\infty\) for almost every
\(g\in B\) whenever \(\varphi\in L^1(\Omega)\).

For consecutive chain balls \(J,K\), define
\[
\mathcal A_{J,K}\varphi
:=
C_{\mathrm{av}} r_J^{1-Q}
\left(
\int_{\Lambda J}\varphi+
\int_{\Lambda K}\varphi
\right).
\]
Finally put
\begin{equation}\label{eq:rigorous-positive-operator}
\mathcal T\varphi(g)
:=
\sum_{B\in\mathcal W}\chi_{A_B}(g)
\left[
\mathcal R_B\varphi(g)
+
\sum_{j=0}^{\ell(B)-1}
\mathcal A_{B_j,B_{j+1}}\varphi
\right].
\end{equation}
At each point only one outer summand is nonzero.  If \(g\in A_B\), then
\(g\in B\); telescoping the averages along \(\mathcal C(B)\) and using
\eqref{eq:rigorous-local-representation} and
\eqref{eq:rigorous-average-difference} gives
\begin{equation}\label{eq:rigorous-global-representation}
|u(g)-u_{B_*}|
\le
\mathcal T(|Xu|)(g)
\qquad\text{for a.e. }g\in\Omega.
\end{equation}
Indeed, if \(g\in A_B\) and
\(\mathcal C(B)=(B_0=B,B_1,\ldots,B_{\ell(B)}=B_*)\), then
\[
u(g)-u_{B_*}
=
u(g)-u_{B_0}
+
\sum_{j=0}^{\ell(B)-1}
\bigl(u_{B_j}-u_{B_{j+1}}\bigr).
\]
The triangle inequality, \eqref{eq:rigorous-local-representation}, and
\eqref{eq:rigorous-average-difference} therefore give
\[
|u(g)-u_{B_*}|
\le
\mathcal R_B(|Xu|)(g)
+
\sum_{j=0}^{\ell(B)-1}
\mathcal A_{B_j,B_{j+1}}(|Xu|).
\]
Since \(g\in A_B\), the \(B\)-summand is the unique nonzero outer summand in
\eqref{eq:rigorous-positive-operator}; hence the right-hand side is exactly
\(\mathcal T(|Xu|)(g)\), proving \eqref{eq:rigorous-global-representation}.

\medskip
\noindent\emph{Step 2: endpoint estimates for \(\mathcal T\).}

We first prove
\begin{equation}\label{eq:rigorous-weighted-L1}
\|\mathcal T\varphi\|_{L^1(\Omega)}
\le
C\int_\Omega\rho(g)\varphi(g)\,dg.
\end{equation}
Here Ahlfors \(Q\)-regularity means that the measure of a ball is comparable
to the \(Q\)-th power of its radius; in the present Carnot group we have the
exact identity \(|B(x,s)|=\omega_Qs^Q\) from
\eqref{eq:proposed-prelim-ball-scaling}.  Tonelli's theorem gives
\[
\int_{A_B}\mathcal R_B\varphi(g)\,dg
=
C_{\mathrm{rep}}\int_{\Lambda B}\varphi(h)
\left(
\int_{A_B}\frac{dg}{d(g,h)^{Q-1}}
\right)dh.
\]
If \(h\in\Lambda B\) and \(g\in A_B\subset B\), then
\[
d(g,h)
\le d(g,g_B)+d(g_B,h)
<(1+\Lambda)r_B.
\]
Consequently, with \(R=(1+\Lambda)r_B\),
\[
\int_{A_B}\frac{dg}{d(g,h)^{Q-1}}
\le
\int_{B(h,R)}\frac{dg}{d(g,h)^{Q-1}}
\le CR.
\]
For completeness, the last inequality follows by decomposing \(B(h,R)\)
into the annuli
\[
B(h,2^{-j}R)\setminus B(h,2^{-j-1}R),
\qquad j=0,1,2,\ldots .
\]
On the \(j\)-th annulus the kernel is at most
\((2^{-j-1}R)^{1-Q}\), while Ahlfors regularity bounds its measure by
\(\omega_Q(2^{-j}R)^Q\).  Its contribution is therefore at most
\(C2^{-j}R\), and summing the geometric series gives the claimed bound.
Since \(R=(1+\Lambda)r_B\) and \(\Lambda\) is fixed, we obtain
\[
\int_{A_B}\mathcal R_B\varphi
\le
C r_B\int_{\Lambda B}\varphi.
\]
Summing this inequality, using the bounded overlap in
\eqref{eq:rigorous-boman-overlap} and
\eqref{eq:rigorous-radius-below-rho}, gives
\begin{equation}\label{eq:rigorous-local-weighted-L1}
\sum_B\int_{A_B}\mathcal R_B\varphi
\le
C\int_\Omega\rho\varphi.
\end{equation}

For the chain part, recall that for every \(B\in\mathcal W\) we fixed a
chain
\[
\mathcal C(B)=(B_0,B_1,\ldots,B_{\ell(B)}),
\qquad B_0=B,\quad B_{\ell(B)}=B_*.
\]
Let \(\mathcal E\) be the set of all directed edges in these chains.  More
precisely, \(e=(J,K)\) belongs to \(\mathcal E\) if there are
\(B\in\mathcal W\) and \(0\le j<\ell(B)\) such that
\(J=B_j\) and \(K=B_{j+1}\).  The word ``oriented'' records the order toward
the distinguished ball \(B_*\); thus \((J,K)\) and \((K,J)\) are regarded
as different pairs.  If the same directed edge occurs in several chains, it
is included only once in \(\mathcal E\).  For \(e\in\mathcal E\), set
\[
\Phi_e
:=
\sum_{\substack{B\in\mathcal W:\\
e\text{ occurs in }\mathcal C(B)}}\chi_{A_B}.
\]
Thus \(\Phi_e\) is the indicator of the union of those partition pieces
\(A_B\) whose chosen chain uses the directed edge \(e\).  The disjointness
of the \(A_B\)'s gives \(0\le\Phi_e\le1\).  If
\(e=(J,K)\) occurs in \(\mathcal C(B)\), then
\eqref{eq:rigorous-boman-shadow-property} and \(A_B\subset B\) imply
\[
\operatorname{spt}\Phi_e\subset M J,
\qquad
\int_\Omega\Phi_e\le |MJ|\le Cr_J^Q.
\]
Therefore
\begin{equation}\label{eq:rigorous-chain-L1-first}
\begin{aligned}
&\int_\Omega
\sum_B\chi_{A_B}
\sum_{j=0}^{\ell(B)-1}
\mathcal A_{B_j,B_{j+1}}\varphi\\
&\qquad\le
C\sum_{e=(J,K)\in\mathcal E}
r_J
\left(
\int_{\Lambda J}\varphi+
\int_{\Lambda K}\varphi
\right).
\end{aligned}
\end{equation}
There is a uniform bound for the number of chain neighbors of a fixed
ball \(J\).  Indeed, if \(K\) is such a neighbor, then
\eqref{eq:rigorous-boman-bridge} gives \(r_K\asymp r_J\) and
\(K\subset CJ\).  Since the balls \(\lambda K\) have bounded overlap,
Ahlfors regularity gives
\[
\#\{K:(J,K)\in\mathcal E\text{ or }(K,J)\in\mathcal E\}
\le C.
\]
Consequently the right-hand side of
\eqref{eq:rigorous-chain-L1-first} is bounded by
\[
C\sum_{J\in\mathcal W}r_J\int_{\Lambda J}\varphi
\le
C\int_\Omega\rho\varphi,
\]
where the last inequality again follows from
\eqref{eq:rigorous-boman-overlap} and
\eqref{eq:rigorous-radius-below-rho}.  Together with
\eqref{eq:rigorous-local-weighted-L1}, this proves
\eqref{eq:rigorous-weighted-L1}.

We next establish a weak Sobolev endpoint.  Extend \(\varphi\) by zero outside
\(\Omega\) and define
\[
I_1\varphi(g):=
\int_{\mathbb G}\frac{\varphi(h)}{d(g,h)^{Q-1}}\,dh.
\]
We claim that
\begin{equation}\label{eq:rigorous-potential-domination}
\mathcal T\varphi(g)\le CI_1\varphi(g)
\qquad\text{for a.e. }g\in\Omega.
\end{equation}
For \(g\in A_B\), the local contribution in
\eqref{eq:rigorous-positive-operator} is
\[
\mathcal R_B\varphi(g)
=C_{\mathrm{rep}}\int_{\Lambda B}
\frac{\varphi(h)}{d(g,h)^{Q-1}}\,dh
\le C_{\mathrm{rep}} I_1\varphi(g),
\]
because \(\varphi\ge0\) and \(I_1\varphi\) integrates this same kernel over
all of \(\mathbb G\) (with \(\varphi\) extended by zero outside \(\Omega\)).
It remains to bound the chain contributions.  Let \(g\in A_B\) and
\(J\in\mathcal C(B)\).  From \(g\in B\subset MJ\), we have
\(d(g,g_J)\le Cr_J\).  Thus, for \(h\in\Lambda J\),
\[
d(g,h)\le Cr_J,
\qquad
r_J^{1-Q}\le\frac{C}{d(g,h)^{Q-1}}.
\]
Moreover, by \eqref{eq:rigorous-boman-overlap}, for each fixed \(h\) there
are at most \(M\) balls \(J\in\mathcal W\) for which
\(h\in\Lambda J\).  Since the chosen chain has no repeated balls,
\[
\sum_{\substack{J\in\mathcal C(B):\\h\in\Lambda J}}
r_J^{1-Q}
\le
\frac{C}{d(g,h)^{Q-1}}.
\]
The same argument applies to the successor terms in
\(\mathcal A_{J,K}\), because \(r_J\asymp r_K\).  Fubini's theorem now
proves \eqref{eq:rigorous-potential-domination}.

For completeness, we recall the weak estimate for \(I_1\).  Let
\(\mathcal M\varphi\) be the Hardy--Littlewood maximal function.  Splitting the
integral into \(d(g,h)<R\) and \(d(g,h)\ge R\), and decomposing the first
part into dyadic annuli, gives
\[
I_1\varphi(g)
\le
C\left(R\mathcal M\varphi(g)
+R^{1-Q}\|\varphi\|_{L^1}\right).
\]
Optimizing in \(R\) yields
\[
I_1\varphi(g)
\le
C\|\varphi\|_{L^1}^{1/Q}
\mathcal M\varphi(g)^{(Q-1)/Q}.
\]
The weak \((1,1)\) estimate for \(\mathcal M\), valid on the doubling
homogeneous group, therefore gives
\begin{equation}\label{eq:rigorous-weak-endpoint}
\bigl|\{\mathcal T\varphi>t\}\bigr|
\le
C\left(\frac{\|\varphi\|_{L^1(\Omega)}}{t}\right)^q,
\qquad t>0.
\end{equation}
Here
\[
q=\frac{Q}{Q-1},
\]
as defined at the beginning of the proof; equivalently, \(q\) is the
Hölder-conjugate exponent of the homogeneous dimension \(Q\), since
\(1/Q+1/q=1\).

\medskip
\noindent\emph{Step 3: truncation at the scale \(a\).}

Fix \(0<a\le1\) and recall that
\[
w_a(g):=\min\{\rho(g),a\},
\qquad
\rho(g):=\operatorname{dist}_{cc}(g,\mathbb G\setminus\Omega).
\]
Thus \(w_a\) is the boundary-distance function capped at height \(a\): it
equals \(\rho\) in the boundary layer \(\{\rho<a\}\) and equals the
constant \(a\) in the interior region \(\{\rho\ge a\}\).
For a locally Lipschitz function \(u\) for which the quantities below are
finite, let
\[
V_a(u):=\int_\Omega w_a|Xu|,
\qquad
\varphi_0:=\chi_{\{\rho<a\}}|Xu|,
\qquad
\varphi_1:=\chi_{\{\rho\ge a\}}|Xu|.
\]
The operator \(\mathcal T\) is positive and additive on nonnegative
functions.  Hence \eqref{eq:rigorous-global-representation} gives
\[
|u-u_{B_*}|\le\mathcal T\varphi_0+\mathcal T\varphi_1.
\]
Since \(w_a=\rho\) on \(\{\rho<a\}\) and \(w_a=a\) on
\(\{\rho\ge a\}\), the weighted \(L^1\) endpoint estimate
\eqref{eq:rigorous-weighted-L1} for \(\varphi_0\) and the weak endpoint
estimate \eqref{eq:rigorous-weak-endpoint} for \(\varphi_1\) give
\[
\|\mathcal T\varphi_0\|_{L^1}\le CV_a(u),
\qquad
\bigl|\{\mathcal T\varphi_1>t\}\bigr|
\le
C\left(\frac{V_a(u)}{at}\right)^q.
\]
Using Chebyshev's inequality for the first term and replacing \(t\) by
\(t/2\) in the two estimates, we obtain
\begin{equation}\label{eq:rigorous-distribution-estimate}
\bigl|\{|u-u_{B_*}|>t\}\bigr|
\le
C\left[
\frac{V_a(u)}{t}
+
\left(\frac{V_a(u)}{at}\right)^q
\right],
\qquad t>0.
\end{equation}

The auxiliary lemma below (Lemma~\ref{lem:rigorous-bv-extension})
extends \eqref{eq:rigorous-distribution-estimate} to every
\(u\in BV_H(\Omega)\), with the same constant and with
\[
V_a(u):=\int_\Omega w_a\,d|D_Hu|.
\]
We use this extension in Step~4 below.

\medskip
\noindent\emph{Step 4: application to sets of locally finite horizontal perimeter.}

We now apply the estimate to a set that is only assumed to have locally
finite horizontal perimeter.  Write
\[
J_a:=
\int_{\partial_H^*E\cap B_2}w_a\,dP_H(E).
\]
If \(J_a=+\infty\), the desired conclusion is immediate, so assume
\(J_a<\infty\).  For \(\varepsilon>0\), choose a Lipschitz cutoff
\(\eta_\varepsilon\) such that
\[
0\le\eta_\varepsilon\le1,
\qquad
\eta_\varepsilon=0\quad\text{on }\{\rho\le\varepsilon\},
\qquad
\eta_\varepsilon=1\quad\text{on }\{\rho\ge2\varepsilon\},
\qquad
|X\eta_\varepsilon|\le\frac{C}{\varepsilon},
\]
and set
\[
u_\varepsilon:=\eta_\varepsilon\chi_E.
\]
The support of \(u_\varepsilon\) is compactly contained in \(B_2\), so
local finite perimeter of \(E\) and the product rule give
\(u_\varepsilon\in BV_H(B_2)\).  More precisely, the BV product rule gives
the vector-measure identity
\[
D_Hu_\varepsilon
=
\eta_\varepsilon D_H\chi_E
+
\chi_E X\eta_\varepsilon\,dg.
\]
Therefore, by the triangle inequality for total variation measures,
\[
|D_Hu_\varepsilon|
\le
\eta_\varepsilon|D_H\chi_E|
+
\chi_E|X\eta_\varepsilon|\,dg.
\]
Since \(|D_H\chi_E|=P_H(E;\cdot)\) is concentrated on
\(\partial_H^*E\), and since \(X\eta_\varepsilon=0\) almost everywhere
outside \(\{\varepsilon<\rho<2\varepsilon\}\), it follows that
\[
\begin{aligned}
V_a(u_\varepsilon)
&\le
\int_{\partial_H^*E\cap B_2}
w_a\eta_\varepsilon\,dP_H(E)
+
\int_{\{\varepsilon<\rho<2\varepsilon\}}
w_a|X\eta_\varepsilon|\,dg\\
&\le
J_a
+
\frac{C}{\varepsilon}
\int_{\{\varepsilon<\rho<2\varepsilon\}}w_a\,dg\\
&\le
J_a+C|\{\rho<2\varepsilon\}|.
\end{aligned}
\]
Indeed, on the transition layer one has
\(w_a=\min\{\rho,a\}\le\rho<2\varepsilon\), so the factor
\(\varepsilon^{-1}\) in the cutoff-gradient bound is cancelled.  The first
integral is at most \(J_a\) because \(0\le\eta_\varepsilon\le1\).
Since \(|\{\rho<2\varepsilon\}|\to0\), it follows that
\begin{equation}\label{eq:rigorous-limsup-variation}
\limsup_{\varepsilon\downarrow0}V_a(u_\varepsilon)\le J_a.
\end{equation}

Let
\[
c_E:=(\chi_E)_{B_*}
=
\frac{|E\cap B_*|}{|B_*|}
\in[0,1].
\]
The ball \(B_*\) is one of the Boman covering balls and
\(\lambda B_*\subset\Omega\).  In particular,
\(\overline{B_*}\Subset\Omega\), so
\[
\delta_*:=\inf_{g\in B_*}\rho(g)>0.
\]
If \(0<\varepsilon<\delta_*/2\), then
\(B_*\subset\{\rho\ge2\varepsilon\}\).  Since
\(\eta_\varepsilon=1\) on this set and
\(u_\varepsilon=\eta_\varepsilon\chi_E\), we have
\[
u_\varepsilon=\chi_E\quad\text{on }B_*,
\qquad
(u_\varepsilon)_{B_*}=c_E.
\]

Away from the boundary layer, \(\eta_\varepsilon=1\), so
\(u_\varepsilon=\chi_E\) and therefore takes only the values \(0\) and
\(1\).  We select the value that is at least \(1/2\) away from the average
\(c_E\).  If \(c_E\le1/2\) and
\(g\in E\cap\{\rho\ge2\varepsilon\}\), then
\[
\chi_E(g)=1,
\qquad
\eta_\varepsilon(g)=1,
\qquad
u_\varepsilon(g)=\eta_\varepsilon(g)\chi_E(g)=1.
\]
Consequently
\[
|u_\varepsilon(g)-c_E|
=|1-c_E|
=1-c_E
\ge\frac12.
\]
Thus
\[
E\cap\{\rho\ge2\varepsilon\}
\subset\{|u_\varepsilon-c_E|\ge1/2\};
\]
if \(c_E\ge1/2\) and
\(g\in(B_2\setminus E)\cap\{\rho\ge2\varepsilon\}\), then
\(\chi_E(g)=0\), so
\[
u_\varepsilon(g)=\eta_\varepsilon(g)\chi_E(g)=0,
\qquad
|u_\varepsilon(g)-c_E|=c_E\ge\frac12.
\]
Therefore
\[
(B_2\setminus E)\cap\{\rho\ge2\varepsilon\}
\subset\{|u_\varepsilon-c_E|\ge1/2\}.
\]
Recall that
\[
m_E:=\min\{|E|,\,|B_2\setminus E|\}.
\]
Thus both \(|E|\) and \(|B_2\setminus E|\) are at least \(m_E\).  Let
\[
F:=
\begin{cases}
E, & c_E\le1/2,\\
B_2\setminus E, & c_E>1/2.
\end{cases}
\]
In either case \(|F|\ge m_E\).  The portion of \(F\) outside the boundary
layer satisfies
\[
\begin{aligned}
|F\cap\{\rho\ge2\varepsilon\}|
&=
|F|-|F\cap\{\rho<2\varepsilon\}|\\
&\ge
|F|-|\{\rho<2\varepsilon\}|\\
&\ge
m_E-|\{\rho<2\varepsilon\}|.
\end{aligned}
\]
Thus, in either case, the appropriate one of the two sets on the left above
has measure at least
\[
m_E-|\{\rho<2\varepsilon\}|.
\]
The threshold in \eqref{eq:rigorous-distribution-estimate} is strict: its
left-hand side is \(\{|u_\varepsilon-c_E|>t\}\).  Here we have proved only
that the relevant points satisfy
\(|u_\varepsilon-c_E|\ge1/2\), and equality may occur, for example when
\(c_E=1/2\).  Thus one cannot directly take \(t=1/2\) without separately
handling the equality set.  Any fixed \(t<1/2\) is sufficient; for
convenience we choose \(t=1/4\).  Indeed,
\(\{|u_\varepsilon-c_E|\ge1/2\}
\subset\{|u_\varepsilon-c_E|>1/4\}\), so
\[
m_E-|\{\rho<2\varepsilon\}|
\le
\bigl|\{|u_\varepsilon-c_E|>1/4\}\bigr|.
\]
Applying \eqref{eq:rigorous-distribution-estimate} to
\(u_\varepsilon\) with \(t=1/4\), using the limsup bound
\eqref{eq:rigorous-limsup-variation}, and then letting
\(\varepsilon\downarrow0\), gives
\begin{equation}\label{eq:rigorous-set-distribution}
m_E
\le
C\left[
J_a+
\left(\frac{J_a}{a}\right)^q
\right].
\end{equation}
Suppose \(a\le m_E^{1/Q}\), and denote by \(C_0\) the constant in
\eqref{eq:rigorous-set-distribution}.  Dividing that inequality by \(C_0\)
gives
\[
J_a+\left(\frac{J_a}{a}\right)^q
\ge \frac{m_E}{C_0}.
\]
Since both terms on the left are nonnegative, at least one of them is at
least \(m_E/(2C_0)\).  If
\[
J_a\ge\frac{m_E}{2C_0},
\]
then the assumption \(a\le m_E^{1/Q}\) implies
\[
a\,m_E^{(Q-1)/Q}
\le
m_E^{1/Q}m_E^{(Q-1)/Q}
=m_E,
\]
and hence
\[
J_a
\ge
\frac{1}{2C_0}\,a\,m_E^{(Q-1)/Q}.
\]
In the other case,
\[
\left(\frac{J_a}{a}\right)^q
\ge\frac{m_E}{2C_0}.
\]
Taking the positive \(q\)-th root and multiplying by \(a\) yields
\[
J_a
\ge
a\left(\frac{m_E}{2C_0}\right)^{1/q}
=
(2C_0)^{-1/q}a\,m_E^{(Q-1)/Q},
\]
because \(q=Q/(Q-1)\), and therefore \(1/q=(Q-1)/Q\).  Thus in both
cases
\[
J_a\ge c\,a\,m_E^{(Q-1)/Q},
\]
where \(c>0\) depends only on the structural constant \(C_0\).
This proves \eqref{eq:rigorous-capped-estimate}.

Finally, since $0<\delta\leq 1$, choose
\[
a:=\min\{\delta,m_E^{1/Q}\}.
\]
Then
\[
a\,m_E^{(Q-1)/Q}
=
\min\{\delta m_E^{(Q-1)/Q},m_E\},
\qquad
w_a\le\min\{\rho,\delta\}=\frac{d_\delta}{\kappa}.
\]
Therefore \eqref{eq:rigorous-capped-estimate} implies
\eqref{eq:horizontal-boundary-concentration}, after absorbing \(\kappa\)
into the structural constant.
\end{proof}

The preceding proof applies the distribution estimate to the cutoff functions
\(u_\varepsilon\), which belong to \(BV_H(\Omega)\) but need not be locally
Lipschitz; the next lemma supplies the required extension to \(BV_H(\Omega)\).

\begin{lem}[Auxiliary BV extension of the distribution estimate]
\phantomsection\label{lem:rigorous-bv-extension}
Let \(\Omega=B_2\), let \(B_*\Subset\Omega\) be the distinguished ball
chosen in the proof of Lemma~\ref{bce}, and set
\[
q=\frac{Q}{Q-1},\qquad
\rho(g)=\operatorname{dist}_{cc}(g,\mathbb G\setminus\Omega),\qquad
w_a=\min\{\rho,a\},\quad 0<a\le1.
\]
Assume that \eqref{eq:rigorous-distribution-estimate} holds for locally
Lipschitz functions with finite \(L^1\)-norm and finite horizontal
variation.  For \(u\in BV_H(\Omega)\), define
\[
u_{B_*}=\frac{1}{|B_*|}\int_{B_*}u\,dg,
\qquad V_a(u)=\int_\Omega w_a\,d|D_Hu|.
\]
Then, with the same constant \(C\),
\[
\bigl|\{g\in\Omega:|u(g)-u_{B_*}|>t\}\bigr|
\le
C\left[
\frac{V_a(u)}{t}
+\left(\frac{V_a(u)}{at}\right)^q
\right],
\qquad t>0.
\]
\end{lem}

\begin{proof}
\noindent\emph{1. Strict horizontal BV approximation.}

Fix \(u\in BV_H(\Omega)\).  Write
\(\operatorname{lip}_{cc}v\) for the pointwise local Lipschitz constant of
\(v\) with respect to the Carnot--Carath\'eodory distance.  For every locally
Lipschitz \(v\), the Carnot-group Rademacher theorem implies
\[
\operatorname{lip}_{cc}v=|Xv|
\qquad\text{a.e. in }\Omega.
\]
Indeed, at a differentiability point the Pansu differential of \(v\) is the
horizontal covector determined by \(Xv\), and its operator norm is
\(|Xv|\); see, for example, \cite{VAM}.

The metric definition of \(BV\) is the relaxation of the Lipschitz energy,
and in a sub-Riemannian manifold its variation agrees isometrically with the
horizontal variation; see \cite{MirandaBV} and
\cite[Definition~2.2 and Theorem~3.1]{AmbrosioGhezziMagnaniBV}.  Consequently,
{\color{black}
\[
|D_Hu|(\Omega)
=
\inf\left\{
\liminf_{j\to\infty}\int_\Omega |Xv_j|\,dg:
\begin{array}{l}
v_j\in\operatorname{Lip}_{\mathrm{loc}}(\Omega)\cap L^1(\Omega),\\
\|v_j-u\|_{L^1(\Omega)}\to0
\end{array}
\right\}.
\]
Thus the infimum is taken over all sequences of locally Lipschitz functions
converging to \(u\) in \(L^1(\Omega)\), rather than over individual
functions.  For a fixed admissible sequence, the energies need not converge;
the lower limit means
\[
\liminf_{j\to\infty}\int_\Omega |Xv_j|\,dg
=
\lim_{J\to\infty}\inf_{j\ge J}\int_\Omega |Xv_j|\,dg.
\]
}
For each
\(k\), choose an admissible sequence \((v_j^{(k)})_j\) for which the
liminf is at most \(|D_Hu|(\Omega)+1/k\).  Since
\(v_j^{(k)}\to u\) in \(L^1(\Omega)\), we may choose an index \(j(k)\) such
that, with \(u_k:=v_{j(k)}^{(k)}\),
\[
\|u_k-u\|_{L^1(\Omega)}<\frac1k,
\qquad
\int_\Omega|Xu_k|\,dg
\le
|D_Hu|(\Omega)+\frac2k.
\]
Thus the \(L^1\)-convergence below holds and
\[
\limsup_{k\to\infty}\int_\Omega|Xu_k|\,dg
\le |D_Hu|(\Omega).
\]
Conversely, horizontal variation is lower semicontinuous under
\(L^1\)-convergence.  More explicitly, if
\(\Psi\in C_c^1(\Omega;\mathbb R^m)\) and \(|\Psi|\le1\), then
\[
\int_\Omega u\,\operatorname{div}_H\Psi\,dg
=
\lim_{k\to\infty}
\int_\Omega u_k\,\operatorname{div}_H\Psi\,dg
=
-\lim_{k\to\infty}
\int_\Omega\langle Xu_k,\Psi\rangle\,dg.
\]
Taking absolute values, then the supremum over such \(\Psi\), gives
\[
|D_Hu|(\Omega)
\le
\liminf_{k\to\infty}\int_\Omega|Xu_k|\,dg.
\]
The preceding upper and lower bounds prove
\begin{equation}\label{eq:rigorous-strict-bv-approximation}
u_k\to u\quad\text{in }L^1(\Omega),
\qquad
\int_\Omega|Xu_k|\,dg\longrightarrow |D_Hu|(\Omega).
\end{equation}
\medskip
\noindent\emph{2. Convergence of the variation measures and weighted energies.}

Set \(\mu_k:=|Xu_k|\,dg\).  By
\eqref{eq:rigorous-strict-bv-approximation}, the masses \(\mu_k(\Omega)\)
are uniformly bounded.  Since \(\Omega\) is locally compact and second
countable, bounded sets of finite Radon measures are sequentially weak-*
precompact in
\[
\mathcal M(\Omega)=C_0(\Omega)^*.
\]
Let a subsequence, still denoted by \(\mu_k\), converge weak-* to a positive
Radon measure \(\nu\).  We first prove
\[
|D_Hu|\le\nu
\qquad\text{as measures in }\Omega.
\]
Let \(\phi\in C_c(\Omega)\) be nonnegative and let
\(\Psi\in C_c^1(\Omega;\mathbb R^m)\) satisfy \(|\Psi|\le\phi\).  The
\(L^1\)-convergence in \eqref{eq:rigorous-strict-bv-approximation}, integration
by parts, and weak-* convergence give
\[
\begin{aligned}
\left|\int_\Omega u\,\operatorname{div}_H\Psi\,dg\right|
&=
\lim_{k\to\infty}
\left|\int_\Omega u_k\,\operatorname{div}_H\Psi\,dg\right|\\
&=
\lim_{k\to\infty}
\left|\int_\Omega\langle Xu_k,\Psi\rangle\,dg\right|\\
&\le
\lim_{k\to\infty}\int_\Omega\phi\,d\mu_k
=
\int_\Omega\phi\,d\nu.
\end{aligned}
\]
Taking the supremum over all such horizontal fields \(\Psi\), which is the
dual definition of the variation measure, yields
\[
\int_\Omega\phi\,d|D_Hu|
\le
\int_\Omega\phi\,d\nu.
\]

To compare the total masses, choose functions
\(0\le\zeta_j\le1\) in \(C_c(\Omega)\) with \(\zeta_j\uparrow1\) pointwise.
For each fixed \(j\),
\[
\int_\Omega\zeta_j\,d\nu
=
\lim_{k\to\infty}\int_\Omega\zeta_j\,d\mu_k
\le
\lim_{k\to\infty}\mu_k(\Omega)
=
|D_Hu|(\Omega).
\]
Monotone convergence then gives
\(\nu(\Omega)\le|D_Hu|(\Omega)\).  Together with
\(|D_Hu|\le\nu\), this forces \(\nu=|D_Hu|\).  Every weak-* cluster point
has therefore been identified with the same measure.  Weak-* precompactness
then implies that the full sequence satisfies
\[
\mu_k\stackrel{*}{\rightharpoonup}|D_Hu|
\qquad\text{in }C_0(\Omega)^*.
\]

We now justify that \(w_a\) is an admissible test function for this
convergence.  The function
\[
\rho(g)=\operatorname{dist}_{cc}(g,\mathbb G\setminus\Omega),
\qquad g\in\mathbb G,
\]
is \(1\)-Lipschitz and is identically zero on
\(\mathbb G\setminus\Omega\).  Hence
\(w_a=\min\{\rho,a\}\), extended by zero outside \(\Omega\), is continuous
on \(\mathbb G\).  Its support is contained in
\(\overline\Omega=\overline{B_2}\), which is compact because the
Carnot--Carath\'eodory metric on a Carnot group is proper.  Thus the extension
belongs to \(C_c(\mathbb G)\), and its restriction belongs to
\(C_0(\Omega)\).  Testing the preceding measure convergence with \(w_a\)
gives
\begin{equation}\label{eq:rigorous-weighted-strict-convergence}
V_a(u_k)
=
\int_\Omega w_a\,d\mu_k
\longrightarrow
\int_\Omega w_a\,d|D_Hu|.
\end{equation}
\medskip
\noindent\emph{3. Passage to the distribution estimate.}

The averages satisfy
\[
\left|(u_k)_{B_*}-u_{B_*}\right|
\le
\frac{1}{|B_*|}\int_{B_*}|u_k-u|\,dg
\longrightarrow0.
\]
After passing to a subsequence, the \(L^1\)-convergence also gives
\(u_k\to u\) almost everywhere.  Fix \(0<s<t\).  At every point where this
pointwise convergence holds, if \(|u-u_{B_*}|>t\), then
\(|u_k-(u_k)_{B_*}|>s\) for all sufficiently large \(k\).  Consequently,
up to a null set,
\[
\{|u-u_{B_*}|>t\}
\subset
\liminf_{k\to\infty}
\{|u_k-(u_k)_{B_*}|>s\}.
\]
Fatou's lemma for the corresponding indicator functions, followed by
\eqref{eq:rigorous-distribution-estimate} applied to \(u_k\), gives
\[
\begin{aligned}
\bigl|\{|u-u_{B_*}|>t\}\bigr|
&\le
\liminf_{k\to\infty}
\bigl|\{|u_k-(u_k)_{B_*}|>s\}\bigr|\\
&\le
C\left[
\frac{V_a(u)}{s}
+
\left(\frac{V_a(u)}{as}\right)^q
\right],
\end{aligned}
\]
where \eqref{eq:rigorous-weighted-strict-convergence} was used in the last
step.  Letting \(s\uparrow t\) proves
\eqref{eq:rigorous-distribution-estimate} for every
\(u\in BV_H(\Omega)\), with
\(V_a(u)=\int_\Omega w_a\,d|D_Hu|\).
\end{proof}

We now use the preceding boundary concentration estimate to prove the Carnot
analogue of Julin's measure-theoretic lemma.
\begin{lem}[Carnot analogue of {\cite[Lemma 3.4]{VJU}}]\label{carnot-julin-34}
Assume that \(A\subset B_2\) is measurable and
\[
        |A|\le \frac12 |B_2|.
\]
Let \(\sigma_A\) be the localized density defined in
\eqref{eq:sigma-A-density-definition}.
Then there exist constants \(c_1>0\) and \(\delta_1\in(0,1)\), depending only
on \(\mathbb G\) and on the fixed parameter \(\kappa\), such that for every
\(0<\delta\le \delta_1\),
\[
\left|
\left\{
g\in B_2:\frac15<\sigma_A(g)\le \frac45
\right\}
\right|
\ge
c_1
\min\left\{
\delta
\left|
\left\{\sigma_A>\frac45\right\}
\right|^{\frac{Q-1}{Q}},
\;
\left|
\left\{\sigma_A>\frac45\right\}
\right|
\right\}.
\]
If \(\kappa\) has been fixed once and for all in the construction, then the
dependence on \(\kappa\) may be suppressed.
\end{lem}

{\color{black}
\begin{proof}
For \(s\in[0,1]\), set
\[
A^s:=\{g\in B_2:\sigma_A(g)>s\},
\qquad
T:=\left\{g\in B_2:\frac15<\sigma_A(g)\le\frac45\right\}.
\]
We first record a uniform upper bound for the relevant superlevel sets.  If
\(g\in B_{2-\delta}\), then the triangle inequality gives
\[
\rho(g)=\operatorname{dist}_{cc}(g,\mathbb G\setminus B_2)
\ge 2-d(e,g)\ge\delta.
\]
Hence \(d_\delta(g)=\kappa\delta\) for every
\(g\in B_{2-\delta}\).  By left invariance of Haar measure,
\(|B(g,\kappa\delta)|=|B_{\kappa\delta}|\), and therefore the definition of
\(\sigma_A\) gives
\[
\begin{aligned}
\int_{B_{2-\delta}}\sigma_A(g)\,dg
&=
\frac1{|B_{\kappa\delta}|}
\int_{B_{2-\delta}}
|A\cap B(g,\kappa\delta)|\,dg\\
&=
\frac1{|B_{\kappa\delta}|}
\int_{B_{2-\delta}}
\left(\int_A
\chi_{B(g,\kappa\delta)}(y)\,dy\right)dg\\
&=
\frac1{|B_{\kappa\delta}|}
\int_A
\left(\int_{B_{2-\delta}}
\chi_{B(g,\kappa\delta)}(y)\,dg\right)dy\\
&\le |A|.
\end{aligned}
\]
The second equality uses
\[
|A\cap B(g,\kappa\delta)|
=
\int_A\chi_{B(g,\kappa\delta)}(y)\,dy,
\]
while the third equality is Tonelli's theorem, applicable because the
indicator function is nonnegative.
Indeed, fix \(y\in A\).  The inner integral is the measure of the set of
centers \(g\in B_{2-\delta}\) whose ball of radius \(\kappa\delta\)
contains \(y\):
\[
\begin{aligned}
\int_{B_{2-\delta}}
\chi_{B(g,\kappa\delta)}(y)\,dg
&=
\left|
\left\{g\in B_{2-\delta}:y\in B(g,\kappa\delta)\right\}
\right|\\
&=
\left|B_{2-\delta}\cap B(y,\kappa\delta)\right|\\
&\le
|B(y,\kappa\delta)|
=|B_{\kappa\delta}|.
\end{aligned}
\]
The second equality follows because symmetry of the control distance gives
\[
y\in B(g,\kappa\delta)
\quad\Longleftrightarrow\quad
d(g,y)<\kappa\delta
\quad\Longleftrightarrow\quad
g\in B(y,\kappa\delta),
\]
and the last equality follows from left invariance of the distance and Haar
measure.  Thus the inner integral is bounded by the denominator
\(|B_{\kappa\delta}|\), uniformly in \(y\), and integration over \(y\in A\)
gives the asserted bound by \(|A|\).  Since
\(\sigma_A>3/5\) on \(A^{3/5}\), we have
\[
\frac35|A^{3/5}\cap B_{2-\delta}|
\le
\int_{A^{3/5}\cap B_{2-\delta}}\sigma_A(g)\,dg
\le
\int_{B_{2-\delta}}\sigma_A(g)\,dg
\le |A|.
\]
Since, by assumption,  
\(|A|\le |B_2|/2\) we get
\[
|A^{3/5}\cap B_{2-\delta}|
\le\frac53|A|
\le\frac56|B_2|.
\]
Choose \(\delta_1\in(0,1)\), depending only on \(Q\), so that
\[
|B_2\setminus B_{2-\delta}|\le\frac1{24}|B_2|
\qquad(0<\delta\le\delta_1).
\]
It follows that
\begin{equation}\label{eq:corrected-superlevel-upper-bound}
|A^s|\le|A^{3/5}|\le\frac78|B_2|
\qquad\text{for }3/5\le s\le1.
\end{equation}

By Lemma~\ref{gradientofsigmaA},
\[
d_\delta(g)|X\sigma_A(g)|\le C
\qquad\text{for a.e. }g\in B_2.
\]
Moreover,
\[
A^{7/10}\setminus A^{4/5}
=\left\{\frac7{10}<\sigma_A\le\frac45\right\}
\subset T.
\]
Therefore the horizontal coarea formula, applied to the locally Lipschitz
function \(\sigma_A\), yields
\begin{equation}\label{eq:corrected-horizontal-coarea-step}
\begin{aligned}
|T|
&\ge
c\int_{A^{7/10}\setminus A^{4/5}}
d_\delta(g)|X\sigma_A(g)|\,dg\\
&=
c\int_{7/10}^{4/5}
\left(
\int_{\partial_H^*A^s\cap B_2}
d_\delta(g)\,dP_H(A^s)(g)
\right)ds.
\end{aligned}
\end{equation}
For almost every \(s\), the coarea theorem also guarantees that \(A^s\) has
finite horizontal perimeter in \(B_2\), so Lemma~\ref{bce} applies.  Put
\[
m_s:=\min\{|A^s|,|B_2\setminus A^s|\},
\qquad
h:=|A^{4/5}|.
\]
{\color{black}
If \(h=0\), then the right-hand side of the asserted inequality is zero,
and the conclusion is immediate.  We may therefore assume \(h>0\).
}
Since \(A^{4/5}\subset A^s\) for \(s<4/5\), while
\eqref{eq:corrected-superlevel-upper-bound} gives
\(|B_2\setminus A^s|\ge|B_2|/8\), we have
\[
m_s\ge\min\left\{h,\frac18|B_2|\right\}.
\]
The function
\[
\Phi_\delta(r):=\min\left\{
\delta r^{\frac{Q-1}{Q}},r
\right\}
\]
is increasing.  On the fixed interval \(0\le r\le|B_2|\), the preceding
lower bound implies
\begin{equation}\label{eq:corrected-ms-comparison}
\Phi_\delta(m_s)\ge c\Phi_\delta(h),
\end{equation}
with a structural constant independent of \(s\) and \(\delta\).  To see this
explicitly, the assertion is immediate when \(h\le|B_2|/8\).  If
\(h>|B_2|/8\), then
\[
\Phi_\delta(m_s)
\ge
\min\left\{
\delta(|B_2|/8)^{\frac{Q-1}{Q}},|B_2|/8
\right\}
\ge c\delta,
\]
whereas \(\Phi_\delta(h)\le\delta|B_2|^{(Q-1)/Q}\); this proves
\eqref{eq:corrected-ms-comparison}.

Combining Lemma~\ref{bce} with
\eqref{eq:corrected-ms-comparison}, we obtain for a.e.
\(s\in(7/10,4/5)\)
\[
\int_{\partial_H^*A^s\cap B_2}
d_\delta(g)\,dP_H(A^s)(g)
\ge
c\min\left\{
\delta h^{\frac{Q-1}{Q}},h
\right\}.
\]
Inserting this estimate into
\eqref{eq:corrected-horizontal-coarea-step} 
gives
\[
|T|
\ge
c_1\min\left\{
\delta|A^{4/5}|^{\frac{Q-1}{Q}},|A^{4/5}|
\right\}.
\]
Since \(A^{4/5}=\{\sigma_A>4/5\}\), this is exactly the asserted estimate.
\end{proof}
}

\begin{rem}[Heisenberg-cap obstruction]
\label{rem:heisenberg-cap-gauge-deficit}
The use of the intrinsic boundary distance \(\rho\) above is essential in the Carnot setting.
Indeed, 
Julin's Euclidean boundary scale is the gauge-deficit scale 
\[
\widehat d_\delta(g)
:=
\kappa\min\{2-N(g),\delta\},
\]
where \(N(g)=d(e,g)\).  On the Euclidean ball
\(B_2=B_{\mathbb R^n}(0,2)\), with the Euclidean distance, we have
\(\rho(g)=2-N(g)\), and hence \(\widehat d_\delta(g)=d_\delta(g)\).
The corresponding density \(\widehat\sigma_A\) in \cite{VJU} is defined by averaging
on \(B(g,\widehat d_\delta(g))\).  On the Kor\'anyi ball below, however,
these scales need not agree.  Using the gauge-deficit scale, the Carnot
analogue of \cite[Lemma 3.4]{VJU} is false; this does
not contradict Lemma~\ref{carnot-julin-34}, whose definition uses the
intrinsic distance \(\rho\).

To see this, consider the first Heisenberg group
\(\mathbb H^1=\mathbb R^2\times\mathbb R\), equipped with the Kor\'anyi gauge
\[
N(z,\tau)=\bigl(|z|^4+\tau^2\bigr)^{1/4},
\qquad
B_2=\{N<2\},
\]
and, for \(0<\varepsilon\ll\delta\), the cap
\[
A_\varepsilon
:=
\{(z,\tau)\in B_2:4-\varepsilon<\tau<4\}.
\]
Direct integration of the horizontal sections gives
\[
|A_\varepsilon|\asymp\varepsilon^{3/2}.
\]
For the gauge-deficit averaging radius, balls meeting this cap have radius
\(O(\varepsilon)\).  The points where such a ball meets both the cap and its
complement lie in a box of horizontal size \(O(\varepsilon^{1/4})\) and
vertical thickness \(O(\varepsilon^{5/4})\).  Consequently,
\begin{equation}\label{eq:heisenberg-cap-old-density}
\left|\left\{\widehat\sigma_{A_\varepsilon}>\frac45\right\}\right|
\ge c\varepsilon^{3/2},
\qquad
\left|\left\{\frac15<\widehat\sigma_{A_\varepsilon}\le\frac45\right\}\right|
\le C\varepsilon^{7/4}.
\end{equation}
Since \(Q=4\), the right-hand side predicted by the old measure lemma is at
least a constant multiple of \(\varepsilon^{3/2}\), whereas the transition
region in \eqref{eq:heisenberg-cap-old-density} is
\(O(\varepsilon^{7/4})\).  This is a contradiction.

Thus the gauge-deficit version of Lemma~3.4, and the large-density branch of
the earlier attempted Proposition~3.5 based on it, cannot be used.  The
corrected statements above replace \(2-N(g)\) by the genuine
Carnot--Carath\'eodory distance \(\rho(g)\) to the complement and use
horizontal perimeter throughout.  The cap calculation does not contradict
the corrected lemma because its estimate
\(\widehat d_\delta=O(\varepsilon)\) is specific to the discarded
gauge-deficit scale.
\end{rem}

We next combine the preceding surface estimate with the relative isoperimetric inequality to obtain the Carnot-group analogue of \cite[Lemma 3.3]{VJU}.


\begin{lem}\label{lem:L33-carnot}
Let \(\mathbb G\) be a Carnot group of homogeneous dimension \(Q\), endowed
{\color{black}with the Carnot--Carath\'eodory distance \(d=d_{cc}\)}, Haar measure
\(|\cdot|\), horizontal gradient \(X\), and horizontal perimeter \(P_H\).
Let \(B(g,r)\) denote metric balls and set
\[
        B_2:=B(e,2).
\]
Assume that \(u\in C^\infty(B_2)\) is nonnegative and satisfies
\[
        \operatorname{div}_H(A(g)Xu)\le f(u)
        \qquad \text{in }B_2,
\]
where \(A\) is uniformly elliptic and bounded in the horizontal
directions. {
\color{black}
We use the standing assumptions that \(A\) and \(f\) are smooth and
that \(f\colon[0,\infty)\to[0,\infty)\) is nondecreasing.}
Let
\[
        F(s):=\int_0^s f(\tau)\,d\tau.
\]
Define
\[
        A_t:=\{u\ge t\},
        \qquad
        E_t:=\{u<t\}.
\]
{\color{black}
Then, for almost every \(t>0\), the following assertion holds.  More
precisely, \(t\) is taken in the full-measure admissible set supplied by
Lemma~\ref{lem:prop31-admissible-surface-levels}, the boundary estimate in
Proposition~\ref{prop1}, and the distribution-function coarea formula
\eqref{eq:proposed-prelim-distribution-derivative}.}
Suppose that
\[
        B(g_0,2r)\subset B_2,
        \qquad
        r\le \frac{t}{\sqrt{F(2t)}},
\]
with the usual convention that the condition is void if \(F(2t)=0\). Assume
also that, for some \(\theta\in(0,1/2]\),
\[
        \theta |B(g_0,r)|
        \le
        |A_t\cap B(g_0,r)|
        \le
        (1-\theta)|B(g_0,r)|.
\]
Then
\[
        {\color{black}
        \int_{\partial_H^*E_t\cap B(g_0,r)}
        \frac1{|Xu|}\,dP_H(E_t)}
        \ge
        \frac{c(\theta)}{t}\,|B(g_0,r)|.
\]
{\color{black}
The quotient \(|Xu|^{-1}\) is needed only
\(P_H(E_t)\)-almost everywhere: by the smooth-level representation,
the characteristic part \(\{Xu=0\}\cap\partial_H^*E_t\) is
\(P_H(E_t)\)-null.  
The constant \(c(\theta)>0\) depends only on
\(\theta\), on the group, and on the ellipticity and boundedness constants of
\(\mathcal A\).}
\end{lem}

\begin{proof}
{\color{black}Fix a level \(t\) in the full-measure admissible class specified
in the statement.}
We organize the proof as follows.  First, the surface estimate from
Proposition~\ref{prop1} gives an upper bound for
{\color{black}
\(\int_{\partial_H^*E_t\cap B(g_0,r)}|Xu|\,dP_H(E_t)\).}
The radius assumption absorbs the
\(F(2t)/t\) term, and volume growth turns this into the scale-invariant bound
\(C t r^{Q-2}\).  Second, the density hypothesis and the relative
isoperimetric inequality give a lower bound for the perimeter of
\(E_t=\{u<t\}\) in \(B(g_0,r)\).  {\color{black}Finally,
Cauchy--Schwarz on the horizontal reduced boundary combines the upper bound
for the weighted perimeter with the lower bound for the perimeter itself.}

\smallskip
\noindent\emph{Step 1: Upper bound for the weighted level measure.}
{\color{black}
We justify explicitly both the localization to \(B(g_0,2r)\) and the use of
the estimate at the particular level \(t\).  Proposition~\ref{prop1} states
its boundary estimate for almost every level and for balls centered at \(e\),
but its proof is local and is invariant under left translations.  For the
admissible level fixed in the statement, Lemma~\ref{lem:prop31-admissible-surface-levels}
gives Gauss--Green and
\[
\nu_{E_t}=\frac{Xu}{|Xu|}
\qquad P_H(E_t)\text{-a.e. on }\partial_H^*E_t.
\]
The characteristic part is \(P_H(E_t)\)-null, so no global
noncharacteristic hypothesis is required.

For completeness, set \(\widehat r=3r/2\) and choose
\(\eta\in C_c^\infty(B(g_0,\widehat r))\) such that
\[
0\le\eta\le1,
\qquad
\eta\equiv1\text{ on }B(g_0,r),
\qquad
|X\eta|\le\frac Cr.
\]
Applying Gauss--Green to the horizontal field
\(Z=\eta\mathcal A(g)Xu\), using the supersolution inequality and horizontal
ellipticity, gives
\[
\int_{\partial_H^*E_t\cap B(g_0,r)}|Xu|\,dP_H(E_t)
\le
C\int_{E_t\cap B(g_0,\widehat r)}f(u)
+\frac Cr\int_{E_t\cap B(g_0,\widehat r)}|Xu|.
\]
Since \(f\) is nondecreasing,
\[
\int_{E_t\cap B(g_0,\widehat r)}f(u)
\le
\frac{F(2t)}{t}|E_t\cap B(g_0,2r)|.
\]
The solid Caccioppoli argument from Proposition~\ref{prop1}, with inner ball
\(B(g_0,\widehat r)\) and outer ball \(B(g_0,2r)\), gives
\[
\int_{E_t\cap B(g_0,\widehat r)}|Xu|^2
\le
C\left(\frac{t^2}{r^2}+F(2t)\right)
|E_t\cap B(g_0,2r)|.
\]
Consequently, Cauchy--Schwarz yields
\[
\int_{E_t\cap B(g_0,\widehat r)}|Xu|
\le
C\left(\frac tr+\sqrt{F(2t)}\right)
|E_t\cap B(g_0,2r)|.
\]
Substitution in the Gauss--Green estimate, followed by
\[
\frac{\sqrt{F(2t)}}r
\le
\frac12\frac{t}{r^2}
+\frac12\frac{F(2t)}t,
\]
proves the localized regular-level estimate
\begin{equation}\label{eq:L33-prop31}
        \int_{\partial_H^*E_t\cap B(g_0,r)}
        |Xu|\,dP_H(E_t)
        \le
        C\left(
        \frac{t}{r^2}
        +
        \frac{F(2t)}{t}
        \right)
        |E_t\cap B(g_0,2r)|.
\end{equation}
Here \(C\) is a structural constant.
}

If \(F(2t)>0\), the radius hypothesis gives
\[
        r\le \frac{t}{\sqrt{F(2t)}},
\]
and hence
\[
        \frac{F(2t)}{t}\le \frac{t}{r^2}.
\]
If \(F(2t)=0\), the same inequality is immediate.
Therefore \eqref{eq:L33-prop31} simplifies to
\begin{equation}\label{eq:L33-grad-bound-1}
        {\color{black}
        \int_{\partial_H^*E_t\cap B(g_0,r)}
        |Xu|\,dP_H(E_t)}
        \le
        C\frac{t}{r^2}
        |E_t\cap B(g_0,2r)|.
\end{equation}
Since
\[
        |E_t\cap B(g_0,2r)|
        \le
        |B(g_0,2r)|
        =
        C r^Q,
\]
we obtain
\begin{equation}\label{eq:L33-grad-bound-2}
        {\color{black}
        \int_{\partial_H^*E_t\cap B(g_0,r)}
        |Xu|\,dP_H(E_t)}
        \le
        C t r^{Q-2}.
\end{equation}

\smallskip
\noindent\emph{Step 2: Lower bound for the perimeter from density.}
We next use the density assumption. Since
\[
        \theta |B(g_0,r)|
        \le
        |A_t\cap B(g_0,r)|
        \le
        (1-\theta)|B(g_0,r)|,
\]
we also have
\[
        |E_t\cap B(g_0,r)|
        =
        |B(g_0,r)|-|A_t\cap B(g_0,r)|
        \ge
        \theta |B(g_0,r)|.
\]
Thus both phases \(A_t\) and \(E_t\) occupy at least a fixed proportion
\(\theta\) of \(B(g_0,r)\).

{\color{black}
Apply the same-ball relative isoperimetric inequality
\eqref{eq:proposed-prelim-relative-isoperimetric} directly to \(E_t\) in
\(B=B(g_0,r)\).  The two density bounds give
\[
P_H(E_t;B)
\ge
c\min\bigl\{|E_t\cap B|,|B\setminus E_t|\bigr\}^{\frac{Q-1}{Q}}
\ge
c(\theta)|B|^{\frac{Q-1}{Q}}.
\]
Using the homogeneity of Haar measure,
\[
        |B(g_0,r)|=|B(e,1)|r^Q,
\]
this becomes
\begin{equation}\label{eq:L33-per-lower}
        P_H(E_t;B(g_0,r))
        \ge
        c(\theta)r^{Q-1}.
\end{equation}
}

\smallskip
\noindent\emph{Step 3: Perimeter lower bound on the reduced boundary.}
{\color{black}
By the definition of horizontal perimeter and its concentration on the
horizontal reduced boundary,
\begin{equation}\label{eq:L33-per-equals-surface}
        P_H(E_t;B(g_0,r))
        =
        \int_{\partial_H^*E_t\cap B(g_0,r)}dP_H(E_t).
\end{equation}
Combining \eqref{eq:L33-per-lower} and
\eqref{eq:L33-per-equals-surface}, we obtain
\begin{equation}\label{eq:L33-surface-lower}
        \int_{\partial_H^*E_t\cap B(g_0,r)}dP_H(E_t)
        \ge
        c(\theta)r^{Q-1}.
\end{equation}
}

\smallskip
\noindent\emph{Step 4: Cauchy--Schwarz and conclusion.}
{\color{black}
We apply Cauchy--Schwarz with respect to \(P_H(E_t)\) on the reduced
boundary.  Since \(|Xu|>0\) there for \(P_H(E_t)\)-almost every point,
\[
\begin{aligned}
        \left(
        \int_{\partial_H^*E_t\cap B(g_0,r)}dP_H(E_t)
        \right)^2
        &=
        \left(
        \int_{\partial_H^*E_t\cap B(g_0,r)}
        |Xu|^{1/2}|Xu|^{-1/2}\,dP_H(E_t)
        \right)^2                                      \\
        &\le
        \left(
        \int_{\partial_H^*E_t\cap B(g_0,r)}
        |Xu|\,dP_H(E_t)
        \right)
        \left(
        \int_{\partial_H^*E_t\cap B(g_0,r)}
        \frac1{|Xu|}\,dP_H(E_t)
        \right).
\end{aligned}
\]
Therefore,
\[
        \int_{\partial_H^*E_t\cap B(g_0,r)}
        \frac1{|Xu|}\,dP_H(E_t)
        \ge
        \frac{
        \left(
        \int_{\partial_H^*E_t\cap B(g_0,r)}dP_H(E_t)
        \right)^2
        }{
        \int_{\partial_H^*E_t\cap B(g_0,r)}
        |Xu|\,dP_H(E_t)
        }.
\]
Using \eqref{eq:L33-surface-lower} in the numerator and
\eqref{eq:L33-grad-bound-2} in the denominator gives
\[
\begin{aligned}
        \int_{\partial_H^*E_t\cap B(g_0,r)}
        \frac1{|Xu|}\,dP_H(E_t)
        &\ge
        \frac{c(\theta)r^{2Q-2}}{C t r^{Q-2}}       \\
        &=
        \frac{c(\theta)}{t}r^Q.
\end{aligned}
\]
Finally, since
\[
        |B(g_0,r)|=|B(e,1)|r^Q,
\]
we conclude that
\[
        \int_{\partial_H^*E_t\cap B(g_0,r)}
        \frac1{|Xu|}\,dP_H(E_t)
        \ge
        \frac{c(\theta)}{t}|B(g_0,r)|.
\]
}
This proves the claim.
\end{proof}

\textcolor{black}{We now prove the main result of this section: the
Carnot-group analogue of \cite[Proposition 3.5]{VJU}.}

\begin{prop}\label{prop:P35-carnot}
Let \(\mathbb G\) be a Carnot group of homogeneous dimension \(Q\), endowed
{\color{black}with the Carnot--Carath\'eodory distance \(d=d_{cc}\)}, Haar measure
\(|\cdot|\), and horizontal gradient \(X\). Let
\[
        B_2:=B(e,2).
\]
Assume that \(u\in C^\infty(B_2)\) is nonnegative and satisfies
\[
        \operatorname{div}_H(A(g)Xu)\le f(u)
        \qquad\text{in }B_2,
\]
where \(A\) is uniformly elliptic and bounded in the horizontal
directions. {\color{black}We use the standing assumptions that
\(A\) and \(f\) are smooth and that
\(f\colon[0,\infty)\to[0,\infty)\) is nondecreasing.}
Define
\[
        A_t:=\{u\ge t\},
        \qquad
        E_t:=\{u<t\},
\]
and
\[
        \mu(t):=|A_t\cap B_2|,
        \qquad
        \eta(t):=|E_t\cap B_2|.
\]
Finally, set
\[
        F(s):=\int_0^s f(\tau)\,d\tau.
\]

Then there exists a constant \(c>0\), depending only on the group and on the
structural constants of \(\mathcal A\), such that for almost every \(t>0\)
satisfying
\[
        \mu(t)\le \frac12 |B_2|,
\]
one has
\[
        -\mu'(t)
        \ge
        c\min\left\{
        \frac1{\sqrt{F(2t)}}\,\mu(t)^{\frac{Q-1}{Q}},
        \frac1t\,\mu(t)
        \right\}.
\]
Similarly, for almost every \(t>0\) satisfying
\[
        \eta(t)\le \frac12 |B_2|,
\]
one has
\[
        \eta'(t)
        \ge
        c\min\left\{
        \frac1{\sqrt{F(2t)}}\,\eta(t)^{\frac{Q-1}{Q}},
        \frac1t\,\eta(t)
        \right\}.
\]
If \(F(2t)=0\), the expression \(1/\sqrt{F(2t)}\) is interpreted as
\(+\infty\), so the corresponding minimum reduces to the second term.
\end{prop}

\begin{proof}
\emph{Map of the proof.}
We prove the estimate for \(\mu(t)\), the distribution function of the
superlevel sets \(A_t=\{u\ge t\}\) introduced in the statement.  The proof has
six stages.  First, we restrict to the full-measure set of levels for which the
{\color{black}finite-perimeter, Gauss--Green, boundary-estimate, and coarea
identities hold.}
Second, we introduce the scale \(\delta\), the truncated boundary distance
\(d_\delta\), and the localized density \(\sigma\); this is the density
construction \eqref{eq:sigma-A-density-definition} applied to \(A=A_t\).  Third,
we treat the case where the high-density set \(\{\sigma>4/5\}\) is small: a
Vitali covering by balanced balls and Lemma~\ref{lem:L33-carnot} give the
\(t^{-1}\mu(t)\) lower bound.  Fourth, if the high-density set is large, Lemma
\ref{carnot-julin-34} forces a large transition region
\(\{1/5<\sigma\le4/5\}\); another Vitali covering and
Lemma~\ref{lem:L33-carnot} estimate \(-\mu'(t)\) from this transition region.
Fifth, we convert the scale \(\delta\) into the two terms appearing in the
minimum.  Finally, we explain why the same argument gives the estimate for
\(\eta(t)=|E_t\cap B_2|\).

\smallskip
\noindent\emph{Stage 1: regular levels and the coarea identity.}
We prove the estimate for \(\mu(t)\). The proof for \(\eta(t)\) is the same,
after replacing the superlevel set \(A_t\) by the sublevel set \(E_t\).

{\color{black}
Let \(\mathcal T_*\subset(0,\infty)\) be the intersection of the full-measure
sets for which Lemma~\ref{lem:prop31-admissible-surface-levels}, the boundary
estimate in Proposition~\ref{prop1}, Lemma~\ref{lem:L33-carnot}, and the
distribution-function identity
\eqref{eq:proposed-prelim-distribution-derivative} are valid.  Fix
\(t\in\mathcal T_*\) such that \(\mu(t)\le |B_2|/2\).  The set \(E_t\) then
has locally finite horizontal perimeter, and the characteristic part of its
reduced boundary is \(P_H(E_t)\)-null.  For a Borel set \(U\subset B_2\), put
\[
\lambda_t(U)
:=
\int_{\{u=t\}\cap U}\frac1{|\nabla u|}\,d\mathcal H^{N-1}.
\]
The Euclidean coarea formula gives
\begin{equation}\label{eq:P35-rigorous-coarea-measure}
-\mu'(t)=\eta'(t)=\lambda_t(B_2).
\end{equation}
On the other hand, the smooth-level representation gives, for every Borel
\(U\subset B_2\),
\[
\lambda_t(U)
\ge
\int_{\partial_H^*E_t\cap U}\frac1{|Xu|}\,dP_H(E_t),
\]
Indeed, on the noncharacteristic part of the regular level surface, the
smooth-level representation gives
\[
dP_H(E_t)
=
\frac{|Xu|}{|\nabla u|}\,d\mathcal H^{N-1}.
\]
Consequently,
\[
\begin{aligned}
\int_{\partial_H^*E_t\cap U}\frac1{|Xu|}\,dP_H(E_t)
&=
\int_{\partial_H^*E_t\cap U\cap\{|Xu|>0\}}
\frac1{|\nabla u|}\,d\mathcal H^{N-1}\\
&\le
\int_{\{u=t\}\cap U}
\frac1{|\nabla u|}\,d\mathcal H^{N-1}\\
&=
\lambda_t(U).
\end{aligned}
\]
The inequality, rather than equality, is used because the Euclidean coarea
measure is taken over the entire level surface, whereas the horizontal
perimeter integral is supported on its noncharacteristic reduced boundary.
Moreover, the value assigned to \(|Xu|^{-1}\) on the
\(P_H(E_t)\)-null characteristic part is irrelevant.
This is the rigorous coarea formulation even when the level surface contains
characteristic points; no assertion that \(|Xu|>0\) everywhere on
\(\{u=t\}\) is needed.}

\smallskip
\noindent\emph{Stage 2: the intrinsic scale and the local density function.}
Let \(\delta_1>0\) be the small constant appearing in the Carnot-group
analogue of Lemma 3.4, namely Lemma~\ref{carnot-julin-34}. Define
\[
        \delta:=
        \min\left\{
        \frac{t}{\sqrt{F(2t)}},\delta_1
        \right\}.
\]
If \(F(2t)=0\), we interpret \(t/\sqrt{F(2t)}=+\infty\), so that
\(\delta=\delta_1\).

{\color{black}
Let
\[
        \rho(g):=\operatorname{dist}_{cc}(g,\mathbb G\setminus B_2)
\]
and define
\[
        d_\delta(g)
        :=
        \frac1{100}\min\{\rho(g),\delta\},
        \qquad g\in B_2.
\]
This is the truncated intrinsic boundary distance used in the density estimates
of this section, with the fixed structural choice \(\kappa=1/100\).}
For the
fixed level set \(A_t\), define the local density function
\[
        \sigma(g)
        :=
        \frac{|A_t\cap B(g,d_\delta(g))|}
        {|B(g,d_\delta(g))|}.
\]
Thus \(\sigma=\sigma_{A_t}\) in the notation of
\eqref{eq:sigma-A-density-definition}. The role of \(d_\delta\) is twofold: its
upper bound by \(\delta\) enforces the radius hypothesis in
Lemma~\ref{lem:L33-carnot}, while its distance-to-boundary factor guarantees
that the doubled balls used below remain inside \(B_2\).
We split the argument into two alternatives, depending on the size of the set
where this local density is large.

\medskip

\noindent
\emph{Stage 3: the small high-density alternative.}
Assume first that the high-density set is small, namely
\[
        |\{\sigma>4/5\}|\le \frac12 \mu(t).
\]
Since
\[
        \mu(t)=|A_t\cap B_2|,
\]
we have
\[
\begin{aligned}
        |A_t\cap \{\sigma\le 4/5\}|
        &=
        |A_t\cap B_2|-|A_t\cap\{\sigma>4/5\}|        \\
        &\ge
        \mu(t)-|\{\sigma>4/5\}|              \\
        &\ge
        \frac12\mu(t).
\end{aligned}
\]

{\color{black}
Let \(G_t\) be the set of density points of \(A_t\) that belong to
\(A_t\cap\{\sigma\le4/5\}\).  By the Lebesgue differentiation theorem,
\[
|G_t|=|A_t\cap\{\sigma\le4/5\}|\ge\frac12\mu(t).
\]
For every \(g\in G_t\),
\[
        \frac{|A_t\cap B(g,\rho)|}{|B(g,\rho)|}\longrightarrow 1
        \qquad\text{as }\rho\downarrow 0.
\]
On the other hand, by the definition of \(\sigma\),
\[
        \frac{|A_t\cap B(g,d_\delta(g))|}
        {|B(g,d_\delta(g))|}
        \le \frac45.
\]
Therefore, by continuity in the radius, there exists
\[
        \rho_g\in (0,d_\delta(g)]
\]
such that
\[
        \frac15
        \le
        \frac{|A_t\cap B(g,\rho_g)|}{|B(g,\rho_g)|}
        \le
        \frac45.
\]
Here the asserted continuity is explicit: if \(r_j\to r>0\), then
\[
B(g,r_j)\mathbin{\triangle}B(g,r)
\subset
B(g,r+\varepsilon_j)\setminus B(g,r-\varepsilon_j),
\qquad \varepsilon_j:=|r_j-r|,
\]
for all sufficiently large \(j\), and the measure of this annulus tends to
zero because \(|B(g,s)|=\omega_Qs^Q\).  Thus both
\(r\mapsto|A_t\cap B(g,r)|\) and the density ratio are continuous for
\(r>0\).}
Set
\[
        B_g:=B(g,\rho_g).
\]
Then
\[
        \rho_g\le d_\delta(g)\le \delta
        \le
        \frac{t}{\sqrt{F(2t)}}.
\]
Moreover, because
{\color{black}
\[
        d_\delta(g)\le \frac1{100}\rho(g),
\]}
we have
\[
        B(g,2\rho_g)\subset B_2.
\]
{\color{black}Indeed,
\(2\rho_g\le2d_\delta(g)<\rho(g)\), and the definition of
\(\rho(g)\) therefore gives \(B(g,2\rho_g)\Subset B_2\).}

We now apply the \(5r\)-Vitali covering lemma to the family
\[
        {\color{black}\{B_g:g\in G_t\}}.
\]
We obtain a pairwise disjoint subcollection \(\{B_i\}\), say
\[
        B_i=B(g_i,\rho_i),
\]
such that
\[
        {\color{black}G_t}
        \subset
        \bigcup_i 5B_i.
\]
Hence
\[
        \frac12\mu(t)
        \le
        {\color{black}|G_t|}
        \le
        \sum_i |5B_i|.
\]
By homogeneity of Haar measure,
\[
        |5B_i|=5^Q |B_i|.
\]
Therefore
\[
        \sum_i |B_i|\ge c\,\mu(t).
\]

Each ball \(B_i=B(g_i,\rho_i)\) satisfies the hypotheses of
Lemma~\ref{lem:L33-carnot}. Indeed,
\[
        B(g_i,2\rho_i)\subset B_2,
        \qquad
        \rho_i\le \frac{t}{\sqrt{F(2t)}},
\]
and
\[
        \frac15 |B_i|
        \le
        |A_t\cap B_i|
        \le
        \frac45 |B_i|.
\]
Thus the Carnot-group analogue of Julin's Lemma 3.3 applies to each \(B_i\),
with \(\theta=1/5\); in the notation above this is
Lemma~\ref{lem:L33-carnot}. It gives
\[
        {\color{black}
        \lambda_t(B_i)
        \ge
        \int_{\partial_H^*E_t\cap B_i}
        \frac1{|Xu|}\,dP_H(E_t)}
        \ge
        \frac{c}{t}|B_i|.
\]
Since the balls \(B_i\) are pairwise disjoint, summing over \(i\) gives
\[
\begin{aligned}
        -\mu'(t)
        &={\color{black}\lambda_t(B_2)}                  \\
        &\ge
        {\color{black}\sum_i\lambda_t(B_i)}                  \\
        &\ge
        \frac{c}{t}\sum_i |B_i|              \\
        &\ge
        \frac{c}{t}\mu(t).
\end{aligned}
\]
This proves the desired estimate in the first alternative, since
\[
        \min\left\{
        \frac1{\sqrt{F(2t)}}\mu(t)^{\frac{Q-1}{Q}},
        \frac1t\mu(t)
        \right\}
        \le
        \frac1t\mu(t).
\]

\medskip

\noindent
\emph{Stage 4: the large high-density alternative.}
Assume now that
\[
        |\{\sigma>4/5\}|\ge \frac12 \mu(t).
\]
Set
\[
        T:=\{g\in B_2:1/5<\sigma(g)\le 4/5\}.
\]
The set \(T\) is precisely the transition region in Lemma~\ref{carnot-julin-34}
for the density \(\sigma=\sigma_{A_t}\).  For each \(g\in T\), set
\[
        B_g:=B(g,d_\delta(g)).
\]
Then, by definition of \(T\),
\[
        \frac15 |B_g|
        <
        |A_t\cap B_g|
        \le
        \frac45 |B_g|.
\]
Also,
\[
        d_\delta(g)\le \delta
        \le
        \frac{t}{\sqrt{F(2t)}},
\]
and, as above,
\[
        B(g,2d_\delta(g))\subset B_2.
\]

Apply again the \(5r\)-Vitali covering lemma to the family
\[
        \{B_g:g\in T\}.
\]
We obtain a pairwise disjoint subcollection \(\{B_i\}\) such that
\[
        T\subset \bigcup_i 5B_i
\]
up to a null set. Therefore
\[
        |T|\le \sum_i |5B_i|=5^Q\sum_i |B_i|,
\]
and hence
\[
        \sum_i |B_i|\ge c |T|.
\]

Applying Lemma~\ref{lem:L33-carnot} to each \(B_i\), again with
\(\theta=1/5\), yields
\[
        {\color{black}
        \lambda_t(B_i)
        \ge
        \int_{\partial_H^*E_t\cap B_i}
        \frac1{|Xu|}\,dP_H(E_t)}
        \ge
        \frac{c}{t}|B_i|.
\]
Summing over the disjoint balls, we obtain
\[
\begin{aligned}
        -\mu'(t)
        &\ge
        {\color{black}\sum_i\lambda_t(B_i)}                      \\
        &\ge
        \frac{c}{t}\sum_i |B_i|                  \\
        &\ge
        \frac{c}{t}|T|.
\end{aligned}
\]

We now estimate \(|T|\) from below. Since \(\mu(t)=|A_t\cap B_2|\le |B_2|/2\),
we may apply Lemma~\ref{carnot-julin-34} to the set \(A_t\), with the density
\(\sigma=\sigma_{A_t}\) and the scale \(\delta\le\delta_1\). It gives
\[
        |T|
        \ge
        c
        \min\left\{
        \delta |\{\sigma>4/5\}|^{\frac{Q-1}{Q}},
        |\{\sigma>4/5\}|
        \right\}.
\]
Using the Case 2 assumption
\[
        |\{\sigma>4/5\}|\ge \frac12\mu(t),
\]
we get
\[
        |T|
        \ge
        c
        \min\left\{
        \delta \mu(t)^{\frac{Q-1}{Q}},
        \mu(t)
        \right\}.
\]
Consequently,
\begin{equation}\label{eq:P35-case2-intermediate}
        -\mu'(t)
        \ge
        \frac{c}{t}
        \min\left\{
        \delta \mu(t)^{\frac{Q-1}{Q}},
        \mu(t)
        \right\}.
\end{equation}

\smallskip
\noindent\emph{Stage 5: conversion of the scale \(\delta\).}
It remains to replace \(\delta/t\) by the desired expression. Recall that
\[
        \delta=
        \min\left\{
        \frac{t}{\sqrt{F(2t)}},\delta_1
        \right\}.
\]

If
\[
        \delta=\frac{t}{\sqrt{F(2t)}},
\]
then \eqref{eq:P35-case2-intermediate} gives directly
\[
        -\mu'(t)
        \ge
        c
        \min\left\{
        \frac1{\sqrt{F(2t)}}\mu(t)^{\frac{Q-1}{Q}},
        \frac1t\mu(t)
        \right\}.
\]

It remains to consider the case
\[
        \delta=\delta_1
        <
        \frac{t}{\sqrt{F(2t)}}.
\]
Then \eqref{eq:P35-case2-intermediate} gives
\[
        -\mu'(t)
        \ge
        \frac{c}{t}
        \min\left\{
        \delta_1 \mu(t)^{\frac{Q-1}{Q}},
        \mu(t)
        \right\}.
\]
Since
\[
        \mu(t)\le \frac12 |B_2|,
\]
we have
\[
        \mu(t)^{\frac{Q-1}{Q}}
        =
        \frac{\mu(t)}{\mu(t)^{1/Q}}
        \ge
        \frac{\mu(t)}{|B_2|^{1/Q}}.
\]
Therefore
\[
        \delta_1\mu(t)^{\frac{Q-1}{Q}}
        \ge
        c\,\mu(t).
\]
Thus
\[
        \min\left\{
        \delta_1 \mu(t)^{\frac{Q-1}{Q}},
        \mu(t)
        \right\}
        \ge
        c\,\mu(t).
\]
Hence
\[
        -\mu'(t)\ge \frac{c}{t}\mu(t).
\]
Since
\[
        \min\left\{
        \frac1{\sqrt{F(2t)}}\mu(t)^{\frac{Q-1}{Q}},
        \frac1t\mu(t)
        \right\}
        \le
        \frac1t\mu(t),
\]
we again obtain the desired estimate.

Combining Stages 3--5 proves
\[
        -\mu'(t)
        \ge
        c\min\left\{
        \frac1{\sqrt{F(2t)}}\mu(t)^{\frac{Q-1}{Q}},
        \frac1t\mu(t)
        \right\}.
\]

\smallskip
\noindent\emph{Stage 6: the estimate for \(\eta\).}
Fix an admissible level \(t\) for which
\(\eta(t)=|E_t\cap B_2|\le |B_2|/2\), and repeat Stages 2--5 with
\(A=E_t\) in the density definition.  The coarea measure is unchanged:
\[
        \eta'(t)=\lambda_t(B_2).
\]
If a selected ball contains between \(1/5\) and \(4/5\) of \(E_t\),
then, since the admissible level set \(\{u=t\}\) has Haar measure zero, the
complementary set \(A_t=\{u\ge t\}\) also occupies between \(1/5\) and
\(4/5\) of that ball.  Hence Lemma~\ref{lem:L33-carnot} gives
\[
\lambda_t(B_i)
\ge
\int_{\partial_H^*E_t\cap B_i}\frac1{|Xu|}\,dP_H(E_t)
\ge
\frac{c}{t}|B_i|
\]
on every balanced ball selected in either density alternative.

In the large-high-density alternative, apply
Lemma~\ref{carnot-julin-34} directly to \(E_t\); its measure hypothesis is
exactly \(\eta(t)\le |B_2|/2\).  The resulting transition-region estimate
is therefore
\[
\eta'(t)
\ge
\frac{c}{t}
\min\left\{
\delta\eta(t)^{\frac{Q-1}{Q}},\eta(t)
\right\}.
\]
The scale conversion in Stage 5, with \(\mu(t)\) replaced by \(\eta(t)\),
then yields
\[
\eta'(t)
\ge
c\min\left\{
\frac{\eta(t)^{\frac{Q-1}{Q}}}{\sqrt{F(2t)}},
\frac{\eta(t)}{t}
\right\}.
\]
This proves the second assertion and completes the proof.
\end{proof}

\section{Estimates for subsolutions}\label{sec:subsolution-estimates}
We now record the subsolution estimates needed later.  They are the horizontal
analogues of the Euclidean estimates, with the Euclidean dimension replaced by
the homogeneous dimension \(Q\), the gradient replaced by the horizontal
gradient \(X\), and Euclidean balls replaced by balls for the fixed homogeneous
distance.  Constants below depend only on the Carnot group and on the
ellipticity constants of \(A\).  We use locally bounded representatives of
subsolutions, which are available by the local boundedness theory recalled in
the introduction.

We begin with the homogeneous linear equation
\begin{equation}\label{eq:sub-homogeneous-linear}
\operatorname{div}_H(A(x)Xu)=0.
\end{equation}

\begin{lem}[Small superlevel set for homogeneous subsolutions]\label{lem:sub-linear-small-superlevel}
Let \(u\in W^{1,2}(B(g_0,2r),X)\) be a weak subsolution of
\eqref{eq:sub-homogeneous-linear} in \(B(g_0,2r)\), and set
\[
M_r:=\sup_{B(g_0,r)}u.
\]
There exists \(\varepsilon_0>0\), depending only on the structural constants,
such that, if \(u(g_0)>0\) and
\[
\big|\{u\ge u(g_0)/2\}\cap B(g_0,r)\big|
\le
\varepsilon_0 |B(g_0,r)|,
\]
then
\[
M_r\ge 8u(g_0).
\]
\end{lem}

\begin{proof}
The De Giorgi local boundedness estimate for homogeneous \(X\)-elliptic
subsolutions gives, for every non-negative subsolution \(v\) of
\eqref{eq:sub-homogeneous-linear},
\begin{equation}\label{eq:sub-linear-local-boundedness}
\sup_{B(g_0,r/2)}v
\le
C |B(g_0,r)|^{-1/2}\|v\|_{L^2(B(g_0,r))}
\le
C r^{-Q/2}\|v\|_{L^2(B(g_0,r))}.
\end{equation}
The truncation
\[
v:=(u-u(g_0)/2)_+
\]
is again a weak subsolution of the homogeneous linear equation.  This follows
from the usual truncation argument: test the subsolution inequality for \(u\)
with approximations of \(\chi_{\{u>u(g_0)/2\}}\phi\), where
\(0\le\phi\in W^{1,2}_0(B(g_0,2r),X)\).

Since \(v(g_0)=u(g_0)/2\) and \(g_0\in B(g_0,r/2)\), \eqref{eq:sub-linear-local-boundedness} gives
\[
\frac{u(g_0)}2
\le
C |B(g_0,r)|^{-1/2}\|v\|_{L^2(B(g_0,r))}.
\]
Moreover, \(0\le v\le M_r\) in \(B(g_0,r)\), and \(v\) is supported in
\(\{u\ge u(g_0)/2\}\cap B(g_0,r)\).  Hence the assumed measure
smallness implies
\[
\|v\|_{L^2(B(g_0,r))}
\le
M_r\,\varepsilon_0^{1/2}|B(g_0,r)|^{1/2}.
\]
Combining the last two estimates yields
\[
u(g_0)\le C\varepsilon_0^{1/2}M_r.
\]
Choosing \(\varepsilon_0\) so small that \(C\varepsilon_0^{1/2}\le 1/8\), we get
\(M_r\ge 8u(g_0)\), as claimed.
\end{proof}

For \(t>0\), write
\[
A_t:=\{u\ge t\}.
\]
The next estimate is the subsolution counterpart of the Caccioppoli estimates
used in Section 3 for supersolutions.

\begin{lem}[Caccioppoli inequality for subsolutions]\label{lem:sub-caccioppoli}
Let \(u\in W^{1,2}(B(g_0,2r),X)\) be a non-negative weak subsolution of
\eqref{meq} in \(B(g_0,2r)\), equivalently
\[
\operatorname{div}_H(A(x)Xu)\ge f(u)
\qquad\text{weakly in }B(g_0,2r),
\]
and set
\[
M_r:=\sup_{B(g_0,r)}u.
\]
Then, for every \(0<\rho<r\) and every \(0<t<M_r\),
\begin{equation}\label{eq:sub-caccioppoli}
\int_{B(g_0,\rho)}(u-t)_+^2\,dx
\le
\frac{C}{(r-\rho)^2}
|A_t\cap B(g_0,r)|^{1/Q}
\frac{M_r}{\sqrt{F(t)}}
\int_{B(g_0,r)}(u-t)_+^2\,dx,
\end{equation}
with the convention that the right-hand side is \(+\infty\) if \(F(t)=0\).
\end{lem}

\begin{proof}
If \(F(t)=0\), there is nothing to prove, so assume \(F(t)>0\).  Let
\(v:=(u-t)_+\), and choose \(\zeta\in C_0^1(B(g_0,r))\) such that
\(0\le\zeta\le1\), \(\zeta=1\) in \(B(g_0,\rho)\), and
\(|X\zeta|\le C/(r-\rho)\).  Using \(\phi=v\zeta^2\) as a test function in the
subsolution inequality \eqref{sod}, and observing that \(Xu=Xv\) on \(\{u>t\}\),
we obtain
\[
-\int_{B(g_0,r)}\zeta^2\langle A(x)Xv,Xv\rangle\,dx
-2\int_{B(g_0,r)}v\zeta\langle A(x)Xv,X\zeta\rangle\,dx
\ge
\int_{B(g_0,r)}f(u)v\zeta^2\,dx.
\]
By ellipticity and Young's inequality, this gives the following estimate in
detail.  If \(\lambda,\Lambda>0\) are the ellipticity constants, then
\[
\lambda\int_{B(g_0,r)}\zeta^2|Xv|^2\,dx
+
\int_{B(g_0,r)}f(u)v\zeta^2\,dx
\le
2\Lambda\int_{B(g_0,r)}v\zeta |Xv|\,|X\zeta|\,dx.
\]
Indeed, the preceding inequality implies that the sum of the elliptic energy
and the reaction term is controlled by the absolute value of the mixed term.
Young's inequality then yields
\[
2\Lambda v\zeta |Xv|\,|X\zeta|
\le
\frac{\lambda}{2}\zeta^2|Xv|^2+C v^2|X\zeta|^2,
\]
and after absorbing the first term on the right into the left-hand side,
\[
\int_{B(g_0,r)}\zeta^2|Xv|^2\,dx
+
\int_{B(g_0,r)}f(u)v\zeta^2\,dx
\le
C\int_{B(g_0,r)}v^2|X\zeta|^2\,dx.
\]
Since \(|X(v\zeta)|^2\le 2\zeta^2|Xv|^2+2v^2|X\zeta|^2\) and
\(|X\zeta|\le C/(r-\rho)\), we obtain
\begin{equation}\label{eq:sub-basic-energy}
\int_{B(g_0,r)}|X(v\zeta)|^2\,dx
+
\int_{B(g_0,r)}f(u)v\zeta^2\,dx
\le
\frac{C}{(r-\rho)^2}
\int_{B(g_0,r)}v^2\,dx.
\end{equation}

On \(A_t\cap B(g_0,r)\), monotonicity of \(f\), the identity
\(F(t)=\int_0^t f(s)\,ds\), and \(0<t<M_r\) give
\[
f(u)\ge f(t)\ge \frac{F(t)}{t}
\ge
\frac{F(t)}{M_r^2}(u-t).
\]
Therefore
\begin{equation}\label{eq:sub-reaction-lower}
\int_{B(g_0,r)}f(u)v\zeta^2\,dx
\ge
\frac{F(t)}{M_r^2}
\int_{B(g_0,r)}v^2\zeta^2\,dx.
\end{equation}
Set \(w:=v\zeta\in W^{1,2}_0(B(g_0,r),X)\).  Combining
\eqref{eq:sub-basic-energy} and \eqref{eq:sub-reaction-lower}, and then using
the horizontal Sobolev inequality applied to \(w^2\) (the same Sobolev input as
in Lemma~\ref{pin}), we get
\[
\begin{aligned}
\frac{C}{(r-\rho)^2}\int_{B(g_0,r)}v^2\,dx
&\ge
\int_{B(g_0,r)}|Xw|^2\,dx
+
\frac{F(t)}{M_r^2}\int_{B(g_0,r)}w^2\,dx \\
&\ge
c\frac{\sqrt{F(t)}}{M_r}
\int_{B(g_0,r)}|w|\,|Xw|\,dx \\
&=
c\frac{\sqrt{F(t)}}{M_r}
\int_{B(g_0,r)}|X(w^2)|\,dx \\
&\ge
c\frac{\sqrt{F(t)}}{M_r}
\left(
\int_{B(g_0,r)}w^{\frac{2Q}{Q-1}}\,dx
\right)^{\frac{Q-1}{Q}}.
\end{aligned}
\]
Since \(w=v\) in \(B(g_0,\rho)\), H\"older's inequality on the support
\(A_t\cap B(g_0,\rho)\) gives
\[
\int_{B(g_0,\rho)}v^2\,dx
\le
|A_t\cap B(g_0,\rho)|^{1/Q}
\left(
\int_{B(g_0,\rho)}v^{\frac{2Q}{Q-1}}\,dx
\right)^{\frac{Q-1}{Q}}.
\]
The last two estimates, together with
\(|A_t\cap B(g_0,\rho)|\le |A_t\cap B(g_0,r)|\), imply
\eqref{eq:sub-caccioppoli}.
\end{proof}

The preceding Caccioppoli inequality yields the scale estimate needed for
subsolutions by a standard De Giorgi iteration.

\begin{lem}[De Giorgi scale estimate for subsolutions]\label{lem:sub-degiorgi-scale}
Let \(u\in W^{1,2}(B(g_0,2r),X)\) be a non-negative weak subsolution of
\eqref{meq} in \(B(g_0,2r)\), and set
\[
M_r:=\sup_{B(g_0,r)}u.
\]
If \(0<M_r\le 2u(g_0)\), then
\begin{equation}\label{eq:sub-degiorgi-scale}
\frac{u(g_0)}{\sqrt{F(u(g_0)/2)}}
\ge
c r,
\end{equation}
where \eqref{eq:sub-degiorgi-scale} is understood as automatic if
\(F(u(g_0)/2)=0\).
\end{lem}

\begin{proof}
If \(F(u(g_0)/2)=0\), the left-hand side in \eqref{eq:sub-degiorgi-scale} is
interpreted as \(+\infty\), and there is nothing to prove.  Assume therefore
that \(F(u(g_0)/2)>0\).

By left translation and the group dilation \(\delta_r\), define
\[
\widetilde u(g):=u(g_0\delta_r g),
\qquad
\widetilde A(g):=A(g_0\delta_r g),
\qquad
\widetilde f(s):=r^2f(s).
\]
The homogeneity of the horizontal vector fields implies that \(\widetilde u\) is
a weak subsolution of
\[
\operatorname{div}_H(\widetilde A(g)X\widetilde u)\ge \widetilde f(\widetilde u)
\qquad\text{in }B_2,
\]
and \(\widetilde A\) has the same ellipticity constants as \(A\).  Also
\(M_1(\widetilde u)=M_r(u)\le2u(g_0)=2\widetilde u(e)\).  Thus it is enough to
prove, in the normalized case \(g_0=e\) and \(r=1\), that
\[
\Lambda_0:=
\frac{u(e)}{\sqrt{F_0(u(e)/2)}}
\ge c,
\qquad
F_0(s):=\int_0^s\widetilde f(\tau)\,d\tau=r^2F(s).
\]
For notational simplicity we drop the tildes and write \(F\) for this rescaled
primitive during the normalized argument.

Let \(1/2<\rho<\tau<1\) and let \(u(e)/2\le h<k\le u(e)\).  Applying
Lemma~\ref{lem:sub-caccioppoli} in the ball \(B_\tau\), using
\(M_\tau\le M_1\le2u(e)\) and \(F(k)\ge F(u(e)/2)\), gives
\[
\int_{A_k\cap B_\rho}(u-k)_+^2\,dx
\le
C\Lambda_0
|A_k\cap B_\tau|^{1/Q}
\frac1{(\tau-\rho)^2}
\int_{A_h\cap B_\tau}(u-h)_+^2\,dx.
\]
Since
\[
|A_k\cap B_\tau|
\le
\frac1{(k-h)^2}
\int_{A_h\cap B_\tau}(u-h)_+^2\,dx,
\]
we obtain
\begin{equation}\label{eq:sub-degiorgi-iteration-ineq}
\int_{A_k\cap B_\rho}(u-k)_+^2\,dx
\le
C\Lambda_0
\frac1{(k-h)^{2/Q}}
\frac1{(\tau-\rho)^2}
\left(
\int_{A_h\cap B_\tau}(u-h)_+^2\,dx
\right)^{1+1/Q}.
\end{equation}

Set
\[
r_i:=\frac12(1+2^{-i}),
\qquad
k_i:=u(e)\left(\frac34-4^{-1-i}\right),
\qquad
x_i:=M_1^{-2}
\int_{A_{k_i}\cap B_{r_i}}(u-k_i)_+^2\,dx.
\]
Using \eqref{eq:sub-degiorgi-iteration-ineq} with
\(h=k_i\), \(k=k_{i+1}\), \(\tau=r_i\), and \(\rho=r_{i+1}\), and using
\(M_1\le2u(e)\), we get constants \(C_0>1\) and \(C_1>1\) such that
\begin{equation}\label{eq:sub-degiorgi-recursion}
x_{i+1}
\le
C_0C_1^i\Lambda_0 x_i^{1+1/Q}.
\end{equation}
We claim that
\begin{equation}\label{eq:sub-degiorgi-lambda-lower}
\Lambda_0
\ge
|B_1|^{-1/Q}C_0^{-1}C_1^{-Q}.
\end{equation}
Indeed, if \eqref{eq:sub-degiorgi-lambda-lower} failed, then
\[
x_0
\le
M_1^{-2}\int_{B_1}u^2\,dx
\le
|B_1|
\le
C_0^{-Q}C_1^{-Q^2}\Lambda_0^{-Q}.
\]
The algebraic iteration lemma, Lemma~\ref{eg}, applied to
\eqref{eq:sub-degiorgi-recursion}, would imply \(x_i\to0\).  Hence
\((u-3u(e)/4)_+=0\) a.e. in \(B_{1/2}\), so
\[
\sup_{B_{1/2}}u\le \frac34u(e),
\]
contradicting \(e\in B_{1/2}\) and the chosen representative value \(u(e)>0\).
Thus \eqref{eq:sub-degiorgi-lambda-lower} holds, so \(\Lambda_0\ge c\).  Undoing
the scaling gives \eqref{eq:sub-degiorgi-scale}.
\end{proof}

\begin{thm}[Integral lower bound for subsolutions]
\label{thm:sub-integral-lower-bound}
Let \(u\) be a non-negative weak subsolution of \eqref{meq} in \(B(g_0,2R)\),
with
\[
u\in W^{1,2}(B(g_0,2R),X)\cap C(B(g_0,2R)).
\]
Equivalently,
\[
\operatorname{div}_H(A(x)Xu)\ge f(u)
\qquad\text{weakly in }B(g_0,2R).
\]
Set
\[
M_R:=\sup_{B(g_0,R)}u,
\qquad
m_R:=\inf_{B(g_0,R)}u.
\]
If \(u(g_0)>0\), then
\[
\int_{m_R/4}^{M_R}\frac{dt}{\sqrt{F(t)}}\ge cR,
\]
where the integral is interpreted as \(+\infty\) if \(F(t)=0\) on a nontrivial
subinterval, and \(c>0\) depends only on the Carnot group, the fixed homogeneous
distance, and the ellipticity constants of \(A\).
\end{thm}

\begin{proof}
Write \(\mu_0:=u(g_0)>0\).  If, in the construction below, one of the intervals
of integration contains a nontrivial interval on which \(F=0\), then the desired
inequality is immediate.  Thus we may carry out the estimates with positive
denominators; the stated convention covers the remaining cases.

Let \(R_0:=0\).  We choose radii
\[
0=R_0<R_1<\cdots<R_K=R
\]
so that, with
\[
\mu_k:=\sup_{B(g_0,R_k)}u,
\]
one has
\[
\mu_k=2\mu_{k-1}\quad\text{for }k=1,\dots,K-1,
\qquad
\mu_K=M_R\le 2\mu_{K-1}.
\]
Indeed, if \(M_R\le2\mu_0\), we take \(K=1\) and \(R_1=R\).  Otherwise, by
continuity of \(u\) on compact balls and monotonicity of
\(s\mapsto\sup_{B(g_0,s)}u\), we stop at the first dyadic level above \(M_R\)
and choose the intermediate radii at which the successive dyadic levels are
attained.  For notational uniformity set
\[
\mu_{-1}:=\frac{\mu_0}{2},
\qquad
\mu_{-2}:=\frac{\mu_0}{4}.
\]

Fix \(k\in\{1,
\dots,K\}\) and put \(r_k:=R_k-R_{k-1}\).  Choose a point
\(g_{k-1}\in\overline{B(g_0,R_{k-1})}\) such that
\(u(g_{k-1})=\mu_{k-1}\).  This is possible because \(u\) is continuous and the
closed balls for the homogeneous distance are compact.  By the triangle
inequality and left invariance of the distance,
\[
B(g_{k-1},r_k)\subset B(g_0,R_k),
\qquad
B(g_{k-1},2r_k)\subset B(g_0,2R).
\]
Therefore \(u\) is a weak subsolution in \(B(g_{k-1},2r_k)\), and
\[
\sup_{B(g_{k-1},r_k)}u
\le
\sup_{B(g_0,R_k)}u
=\mu_k
\le 2\mu_{k-1}
=2u(g_{k-1}).
\]
Lemma~\ref{lem:sub-degiorgi-scale}, applied with center \(g_{k-1}\) and radius
\(r_k\), gives
\begin{equation}\label{eq:sub-integral-step-scale}
\frac{\mu_{k-1}}{\sqrt{F(\mu_{k-1}/2)}}
\ge
c r_k.
\end{equation}

Since \(\mu_{k-1}=4\mu_{k-3}\) and \(\mu_{k-1}/2=\mu_{k-2}\), the monotonicity of
\(F\) gives
\[
\begin{aligned}
\int_{\mu_{k-3}}^{\mu_{k-2}}\frac{dt}{\sqrt{F(t)}}
&\ge
\frac{\mu_{k-2}-\mu_{k-3}}{\sqrt{F(\mu_{k-2})}}  \\
&=
\frac{\mu_{k-3}}{\sqrt{F(\mu_{k-1}/2)}}
=
\frac14\frac{\mu_{k-1}}{\sqrt{F(\mu_{k-1}/2)}}
\ge
c r_k,
\end{aligned}
\]
where in the last inequality we used \eqref{eq:sub-integral-step-scale}, after
absorbing the factor \(1/4\) into \(c\).

The intervals \([\mu_{k-3},\mu_{k-2}]\), \(k=1,\dots,K\), are adjacent.  Hence,
after summing the preceding estimate,
\[
\int_{\mu_{-2}}^{\mu_{K-2}}\frac{dt}{\sqrt{F(t)}}
\ge
c\sum_{k=1}^K r_k
=cR.
\]
Finally, \(m_R\le u(g_0)=\mu_0\), so \(m_R/4\le\mu_{-2}\), and
\(\mu_{K-2}\le M_R\).  Since the integrand is non-negative,
\[
\int_{m_R/4}^{M_R}\frac{dt}{\sqrt{F(t)}}
\ge
\int_{\mu_{-2}}^{\mu_{K-2}}\frac{dt}{\sqrt{F(t)}}
\ge cR,
\]
which proves the theorem.
\end{proof}

\begin{rem}
The positive-base-point assumption is part of the scale argument: the first
dyadic interval used in the proof is
\([u(g_0)/4,u(g_0)/2]\).  The lower limit \(m_R/4\) is also natural for this
reason.  Without additional assumptions on the behavior of \(f\) near the
minimum level, one cannot generally replace the interval \([m_R/4,M_R]\) by
\([m_R,M_R]\).
\end{rem}

\begin{cor}[A necessary Keller--Osserman condition]
\label{cor:sub-keller-osserman-necessary}
Suppose that there exists a nonconstant, continuous, non-negative weak
subsolution \(u\in W^{1,2}_{\mathrm{loc}}(\mathbb G,X)\) of
\[
\operatorname{div}_H(A(x)Xu)\ge f(u)
\]
in all of \(\mathbb G\).  Then necessarily
\[
\int_0^\infty\frac{dt}{\sqrt{F(t)}}=\infty.
\]
\end{cor}

\begin{proof}
Since \(u\ge0\) and \(u\) is not identically zero, there is a point
\(g_0\in\mathbb G\) such that \(u(g_0)>0\).  For each \(R>0\), the global
subsolution hypothesis allows us to apply Theorem~\ref{thm:sub-integral-lower-bound}
in \(B(g_0,2R)\).  If
\[
M_R:=\sup_{B(g_0,R)}u,
\qquad
m_R:=\inf_{B(g_0,R)}u,
\]
then
\[
\int_{m_R/4}^{M_R}\frac{dt}{\sqrt{F(t)}}\ge cR.
\]
Because \(u\ge0\), we have \(0\le m_R/4\le M_R\), and therefore
\[
\int_0^\infty\frac{dt}{\sqrt{F(t)}}
\ge
\int_{m_R/4}^{M_R}\frac{dt}{\sqrt{F(t)}}
\ge cR
\]
for every \(R>0\).  Letting \(R\to\infty\) gives the claimed divergence.
\end{proof}

\section{Proof of the generalized Harnack inequality}\label{sec:generalized-harnack-proof}

This section proves Theorem~\ref{gha}.  We write
\[
B_R:=B(e,R)
\]
for balls with respect to the fixed left-invariant homogeneous distance.

\begin{proof}[Proof of Theorem~\ref{gha}]
We use the following four inputs proved above:
\begin{enumerate}[(i)]
\item Lemma~\ref{del};
\item Proposition~\ref{prop:P35-carnot};
\item Lemma~\ref{lem:sub-linear-small-superlevel};
\item Lemma~\ref{lem:sub-degiorgi-scale}.
\end{enumerate}
The horizontal
coarea formula \eqref{eq:proposed-prelim-horizontal-coarea} and the
admissible-level reduction in Lemma~\ref{lem:prop31-admissible-surface-levels}
justify applying the distribution-function estimates for a.e. level.  As in
the proof of Proposition~\ref{prop:P35-carnot}, one may first argue for smooth representatives
and then pass to weak solutions by the standard horizontal Sobolev approximation
and lower semicontinuity of the relevant distribution-function inequalities.
Thus the estimates below hold for the locally bounded representative of the
weak solution; to simplify notation we write pointwise suprema and infima.

Set
\[
E_t:=\{u<t\},
\qquad
A_t:=\{u\ge t\},
\qquad
\eta(t):=|E_t\cap B_2|,
\qquad
\mu(t):=|A_t\cap B_2|.
\]
Also set
\[
m:=\inf_{B_1}u,
\qquad
M:=\sup_{B_1}u,
\qquad
G(t):=\sqrt{F(2t)}+t.
\]
If \(M=m\), there is nothing to prove.  Hence assume \(M>m\).  We first prove
\begin{equation}\label{eq:gha-target-F2}
\int_m^M\frac{dt}{G(t)}\le C.
\end{equation}
This implies \eqref{di}.  Indeed, with the change of variables \(s=2t\),
\begin{equation}\label{eq:gha-change-variable}
\begin{aligned}
\int_m^M\frac{ds}{\sqrt{F(s)}+s}
&\le
2\int_{m/2}^{M/2}\frac{dt}{\sqrt{F(2t)}+t}  \\
&\le
2\int_{m/2}^{m}\frac{dt}{t}
+2\int_m^M\frac{dt}{G(t)}
\le 2\log2+C,
\end{aligned}
\end{equation}
where the first integral over \([m/2,m]\) is understood as zero when \(m=0\).

We choose a median level for the sublevel distribution.  Let
\begin{equation}\label{eq:gha-median-level}
t_0:=\sup\left\{t\ge0:\eta(t)\le\frac{|B_2|}{2}\right\}.
\end{equation}
The set in \eqref{eq:gha-median-level} is nonempty because \(u\ge0\) gives
\(\eta(0)=|\{u<0\}\cap B_2|=0\).  Moreover,
\(\eta(t)\uparrow |B_2|\) as \(t\to\infty\), so this set is bounded above
and \(t_0<\infty\).
Then \(\eta(t)\le |B_2|/2\) for every \(t<t_0\), while \(\mu(t)\le |B_2|/2\)
for every \(t>t_0\), apart from irrelevant level-set jumps.  All differential
inequalities below are used only for almost every level, as in
Proposition~\ref{prop:P35-carnot}.

Let \(\delta_0>0\) be the constant in Lemma~\ref{del}.  We claim first that
either
\begin{equation}\label{eq:gha-lower-alternative-bounded}
\int_0^{t_0}\frac{dt}{G(t)}\le C,
\end{equation}
or there is a level \(t_\delta\in(0,t_0)\) such that
\begin{equation}\label{eq:gha-lower-alternative-small}
\eta(t_\delta)
\le
\delta_0\left(\frac{t_\delta}{G(t_\delta)}\right)^Q,
\qquad
\int_{t_\delta}^{t_0}\frac{dt}{G(t)}\le C.
\end{equation}
To prove the claim, suppose that the first inequality in
\eqref{eq:gha-lower-alternative-small} fails at a level \(t<t_0\).  Since
\(\eta(t)\le |B_2|/2\), the \(\eta\)-estimate in
Proposition~\ref{prop:P35-carnot} gives
\[
\eta'(t)
\ge
c\min\left\{
\frac{\eta(t)^{(Q-1)/Q}}{\sqrt{F(2t)}},
\frac{\eta(t)}{t}
\right\}.
\]
The failure of the smallness condition means
\(\eta(t)^{1/Q}\ge \delta_0^{1/Q}t/G(t)\).  Therefore both terms in the minimum
are bounded from below by a constant multiple of
\(\eta(t)^{(Q-1)/Q}/G(t)\), and hence
\begin{equation}\label{eq:gha-eta-differential}
\frac{d}{dt}\eta(t)^{1/Q}
\ge
\frac{c}{G(t)}
\end{equation}
for a.e. such \(t\).  Integrating \eqref{eq:gha-eta-differential} over
\((t,t_0)\), as long as the smallness condition fails on that interval, gives
\begin{equation}\label{eq:gha-eta-integral-control}
\int_t^{t_0}\frac{ds}{G(s)}
\le
C\left(\eta(t_0)^{1/Q}-\eta(t)^{1/Q}\right)
\le C|B_2|^{1/Q}.
\end{equation}
Consequently, if \eqref{eq:gha-lower-alternative-bounded} is false, then the
smallness condition must hold for some \(t_\delta\in(0,t_0)\), and
\eqref{eq:gha-eta-integral-control} gives the second estimate in
\eqref{eq:gha-lower-alternative-small}.  This proves the claim.

We now control the part of the integral below \(t_0\).  If
\eqref{eq:gha-lower-alternative-bounded} holds, then this part is already
bounded.  Otherwise \eqref{eq:gha-lower-alternative-small} and Lemma~\ref{del}
imply
\begin{equation}\label{eq:gha-lower-min-bound}
m\ge \frac{t_\delta}{2}.
\end{equation}
Thus, if \(m<t_0\),
\begin{equation}\label{eq:gha-lower-part}
\int_m^{t_0}\frac{dt}{G(t)}
\le
\int_{t_\delta/2}^{t_\delta}\frac{dt}{t}
+
\int_{t_\delta}^{t_0}\frac{dt}{G(t)}
\le C.
\end{equation}
If \(m\ge t_0\), there is no lower part to estimate.  Thus, in all cases, the
portion of the integral lying below \(t_0\) is bounded by a structural constant.

It remains to control the part above \(t_0\).  If \(M\le t_0\), there is nothing
to prove, so assume \(M>t_0\).  Choose \(\varepsilon>0\), to be fixed below, and
define
\begin{equation}\label{eq:gha-teps-definition}
t_\varepsilon:=\sup\big(\{t\ge t_0:\mu(t)>\varepsilon\}\cup\{t_0\}\big).
\end{equation}
For a.e. \(t\in(t_0,t_\varepsilon)\) we have
\(\varepsilon<\mu(t)\le |B_2|/2\).  Proposition~\ref{prop:P35-carnot} gives
\[
-\mu'(t)
\ge
c\min\left\{
\frac{\mu(t)^{(Q-1)/Q}}{\sqrt{F(2t)}},
\frac{\mu(t)}{t}
\right\}.
\]
Since \(\mu(t)^{1/Q}\ge\varepsilon^{1/Q}\), the second term in the minimum is at
least \(c(\varepsilon)\mu(t)^{(Q-1)/Q}/G(t)\), while the first term is larger
than \(c\mu(t)^{(Q-1)/Q}/G(t)\).  Hence
\begin{equation}\label{eq:gha-mu-pre-tail}
-\frac{d}{dt}\mu(t)^{1/Q}
\ge
\frac{c(\varepsilon)}{G(t)}
\qquad\text{for a.e. }t\in(t_0,t_\varepsilon).
\end{equation}
After integration,
\begin{equation}\label{eq:gha-upper-pre-tail}
\int_{t_0}^{t_\varepsilon}\frac{dt}{G(t)}
\le C(\varepsilon).
\end{equation}

We next prove that, if \(\varepsilon\) is chosen depending only on the structural
constants, then
\begin{equation}\label{eq:gha-tail-doubling-bound}
M\le2t_\varepsilon.
\end{equation}
Assume by contradiction that \(M>2t_\varepsilon\).  Then, by the definition of
\(t_\varepsilon\),
\begin{equation}\label{eq:gha-M-half-small}
|\{u\ge M/2\}\cap B_2|=\mu(M/2)\le\varepsilon.
\end{equation}
Since \(u\) is a solution of \eqref{meq} and \(f\ge0\), it is also a subsolution
of the homogeneous equation \eqref{eq:sub-homogeneous-linear}.  By the weak
maximum principle for homogeneous \(X\)-elliptic subsolutions, see \cite{GL}, choose
\(p\in\partial B_1\) with \(u(p)=M\).  If
\(\varepsilon\le\varepsilon_0|B_{1/4}|\), where \(\varepsilon_0\) is the constant
from Lemma~\ref{lem:sub-linear-small-superlevel}, then
\eqref{eq:gha-M-half-small} permits us to apply
Lemma~\ref{lem:sub-linear-small-superlevel} in the ball \(B(p,1/4)\).  Since
\(B(p,1/2)\subset B_2\), we obtain
\begin{equation}\label{eq:gha-first-growth-step}
\sup_{B_{5/4}}u\ge \sup_{B(p,1/4)}u\ge 8M.
\end{equation}
Let \(N:=\sup_{B_{3/2}}u\).  Again by the maximum principle, choose
\(q\in\partial B_{3/2}\) with \(u(q)=N\).  From
\eqref{eq:gha-first-growth-step}, \(N\ge8M\), and hence
\(N/2\ge4M>t_\varepsilon\).  Therefore \(\mu(N/2)\le\varepsilon\), and a second
application of Lemma~\ref{lem:sub-linear-small-superlevel}, now in
\(B(q,1/4)\), gives
\begin{equation}\label{eq:gha-second-growth-step}
\sup_{B_{7/4}}u\ge \sup_{B(q,1/4)}u\ge 8N=8\sup_{B_{3/2}}u.
\end{equation}

We now construct dyadic radii.  Put \(R_0:=1\) and \(M_0:=M\).  Recursively
choose radii \(R_k>R_{k-1}\) such that
\begin{equation}\label{eq:gha-dyadic-maxima}
M_k:=\sup_{B_{R_k}}u=2M_{k-1}.
\end{equation}
We stop at the first index \(K\) for which \(R_K\ge3/2\).  The growth estimates
\eqref{eq:gha-first-growth-step} and \eqref{eq:gha-second-growth-step} imply
\begin{equation}\label{eq:gha-radii-location}
K\ge4,
\qquad
R_3\le\frac54,
\qquad
R_K\le\frac74<2.
\end{equation}
Indeed, \eqref{eq:gha-first-growth-step} gives the first three dyadic levels
before radius \(5/4\), while \eqref{eq:gha-second-growth-step} ensures that the
first dyadic level crossing radius \(3/2\) is reached before radius \(7/4\).
Set
\begin{equation}\label{eq:gha-rk-definition}
r_k:=R_k-R_{k-1}.
\end{equation}
Then
\begin{equation}\label{eq:gha-first-three-radii}
r_1+r_2+r_3=R_3-R_0\le\frac14.
\end{equation}

Fix \(k=4,\dots,K\).  We claim that
\begin{equation}\label{eq:gha-mu-rk-lower}
\mu(t)^{1/Q}\ge c r_k
\qquad\text{for every }t\in(M_{k-4},M_{k-2}].
\end{equation}
Since \(\mu\) is non-increasing, it suffices to prove this at
\(t=M_{k-2}\).  Suppose, to the contrary, that
\(\mu(M_{k-2})^{1/Q}<c r_k\), where \(c>0\) is chosen small enough that
\[
\mu(M_{k-2})
\le
\varepsilon_0 |B_{r_k/2}|.
\]
Choose \(z_{k-1}\in\partial B_{R_{k-1}}\) with
\(u(z_{k-1})=M_{k-1}\).  Since \(M_{k-2}=M_{k-1}/2\), the preceding inequality
implies
\[
|\{u\ge M_{k-1}/2\}\cap B(z_{k-1},r_k/2)|
\le
\varepsilon_0 |B_{r_k/2}|.
\]
Moreover \(B(z_{k-1},r_k)\subset B_{R_k}\subset B_2\).  Applying
Lemma~\ref{lem:sub-linear-small-superlevel} with radius \(r_k/2\) gives
\[
\sup_{B(z_{k-1},r_k/2)}u
\ge
8M_{k-1}.
\]
But \(B(z_{k-1},r_k/2)\subset B_{R_k}\), whereas
\(\sup_{B_{R_k}}u=M_k=2M_{k-1}\), a contradiction.  This proves
\eqref{eq:gha-mu-rk-lower}.

The same point \(z_{k-1}\) and radius \(r_k/2\) may be used in the De Giorgi
scale estimate.  Since
\[
\sup_{B(z_{k-1},r_k/2)}u\le M_k=2M_{k-1}=2u(z_{k-1}),
\]
Lemma~\ref{lem:sub-degiorgi-scale} gives
\begin{equation}\label{eq:gha-subsolution-scale-step}
\frac{M_{k-1}}{\sqrt{F(M_{k-1}/2)}}
\ge
c r_k.
\end{equation}
For \(t\in(M_{k-4},M_{k-3}]\), we have
\(M_{k-1}=8M_{k-4}<8t\) and
\(2t\le2M_{k-3}=M_{k-1}/2\).  Since \(F\) is non-decreasing,
\eqref{eq:gha-subsolution-scale-step} implies
\begin{equation}\label{eq:gha-t-over-F-step}
\frac{t}{\sqrt{F(2t)}}
\ge c r_k
\qquad\text{for }t\in(M_{k-4},M_{k-3}],
\end{equation}
with the usual convention that the left-hand side is \(+\infty\) if
\(F(2t)=0\).

For a.e. \(t\in(M_{k-4},M_{k-3}]\), we have \(t>M>t_\varepsilon\), hence
\(\mu(t)\le\varepsilon\le |B_2|/2\), and
Proposition~\ref{prop:P35-carnot} applies.  Combining its \(\mu\)-estimate with
\eqref{eq:gha-mu-rk-lower} and \eqref{eq:gha-t-over-F-step} gives
\begin{equation}\label{eq:gha-mu-tail-differential}
-\mu'(t)
\ge
\frac{c r_k}{t}\mu(t)^{(Q-1)/Q}.
\end{equation}
Equivalently,
\[
-\frac{d}{dt}\mu(t)^{1/Q}
\ge
\frac{c r_k}{t}.
\]
Integrating from \(M_{k-4}\) to \(M_{k-3}=2M_{k-4}\), we get
\begin{equation}\label{eq:gha-one-dyadic-drop}
\mu(M_{k-4})^{1/Q}-\mu(M_{k-3})^{1/Q}
\ge c r_k.
\end{equation}
Summing \eqref{eq:gha-one-dyadic-drop} over \(k=4,\dots,K\), and using the
 telescoping of the adjacent dyadic intervals, yields
\begin{equation}\label{eq:gha-telescoping-tail}
\mu(M)^{1/Q}
\ge
c\sum_{k=4}^K r_k.
\end{equation}
Since \(R_K\ge3/2\), \(R_0=1\), and \eqref{eq:gha-first-three-radii} holds,
\[
\sum_{k=4}^K r_k
=R_K-R_0-(r_1+r_2+r_3)
\ge \frac14.
\]
Thus \eqref{eq:gha-telescoping-tail} gives \(\mu(M)\ge c\).  On the other hand,
\(M>2t_\varepsilon\) implies \(M>t_\varepsilon\), so the definition of
\(t_\varepsilon\) gives \(\mu(M)\le\varepsilon\).  Choosing \(\varepsilon<c\),
after also satisfying the smallness requirements above, gives a contradiction.
Therefore \eqref{eq:gha-tail-doubling-bound} is proved.

Combining \eqref{eq:gha-upper-pre-tail} with
\eqref{eq:gha-tail-doubling-bound}, we find
\begin{equation}\label{eq:gha-upper-part}
\int_{t_0}^M\frac{dt}{G(t)}
\le
\int_{t_0}^{t_\varepsilon}\frac{dt}{G(t)}
+
\int_{t_\varepsilon}^{2t_\varepsilon}\frac{dt}{t}
\le C.
\end{equation}
Together, the lower estimate, namely either
\eqref{eq:gha-lower-alternative-bounded} or \eqref{eq:gha-lower-part}, and the
upper estimate \eqref{eq:gha-upper-part} prove \eqref{eq:gha-target-F2}.  The
change-of-variables estimate
\eqref{eq:gha-change-variable} then proves \eqref{di}, and the theorem follows.
\end{proof}

\begin{rem}[Reduction from general balls to the normalized ball]
\label{rem:general-ball-rescaling}
The proof above was written for the concentric balls \(B(e,1)\subset B(e,2)\).
We now explain why it gives the scale-invariant statement in
Theorem~\ref{gha}.  Fix \(\xi\in\mathbb G\) and \(R>0\), and define
\[
T_{\xi,R}(y):=\xi\circ\delta_R(y).
\]
Since the distance is left invariant and homogeneous,
\[
T_{\xi,R}(B(e,s))=B(\xi,Rs),
\qquad s>0.
\]
Suppose that \(u\) is a non-negative weak solution of \eqref{meq} in
\(B(\xi,2R)\), and set, for \(y\in B(e,2)\),
\[
v(y):=u(T_{\xi,R}(y)),
\qquad
A_{\xi,R}(y):=A(T_{\xi,R}(y)),
\qquad
f_R(t):=R^2f(t).
\]
The horizontal vector fields are left invariant and homogeneous of degree one,
while the Haar measure scales by the homogeneous dimension \(Q\).  Thus
\[
Xv(y)=R(Xu)(T_{\xi,R}(y)),
\qquad
dx=R^Q\,dy
\quad\text{under }x=T_{\xi,R}(y).
\]
More explicitly, if \(\psi\in C_c^1(B(e,2))\) and
\(\phi=\psi\circ T_{\xi,R}^{-1}\), then
\[
(X\phi)(T_{\xi,R}(y))=R^{-1}X\psi(y).
\]
Substitution in the weak formulation for \(u\), followed by division by
\(R^{Q-2}\), gives
\[
-\int_{B(e,2)}
\langle A_{\xi,R}(y)Xv,X\psi\rangle\,dy
=
\int_{B(e,2)}R^2f(v)\psi\,dy.
\]
Consequently, \(v\) solves
\[
\operatorname{div}_H(A_{\xi,R}(y)Xv)=f_R(v)
\qquad\text{in }B(e,2),
\]
and \(A_{\xi,R}\) has the same ellipticity constants as \(A\).  The primitive
of the rescaled nonlinearity is
\[
F_R(t):=\int_0^t f_R(s)\,ds=R^2F(t).
\]
Moreover,
\[
\sup_{B(e,1)}v=M_{\xi,R},
\qquad
\inf_{B(e,1)}v=m_{\xi,R}.
\]
Applying the normalized estimate proved above to \(v\) therefore yields
\[
\int_{m_{\xi,R}}^{M_{\xi,R}}
\frac{dt}{\sqrt{F_R(t)}+t}
=
\int_{m_{\xi,R}}^{M_{\xi,R}}
\frac{dt}{R\sqrt{F(t)}+t}
\le C,
\]
which is exactly \eqref{di} on the general ball.

The same rescaling also preserves the hypothesis used in
Theorem~\ref{gha-c}.  Indeed, if
\(h(t)=\sqrt{f(t)/t}\), then the corresponding function for \(f_R\) is
\[
h_R(t):=\sqrt{\frac{f_R(t)}{t}}=R h(t).
\]
Hence \(h_R\) is non-decreasing and
\(h_R(\theta t)\ge\varrho h_R(t)\) with the same constants
\(\theta\) and \(\varrho\).  Applying the normalized classical Harnack
inequality to \(v\) gives
\[
\sup_{B(\xi,R)}u
=\sup_{B(e,1)}v
\le C\inf_{B(e,1)}v
=C\inf_{B(\xi,R)}u,
\]
which proves the general-ball form \eqref{chi}, with a constant independent of
\(\xi\) and \(R\).
\end{rem}

{\color{black}
\section{Proof of Theorem~\ref{gha-c}}
\label{sec:blue-proof-theorem-22}

We now prove the classical Harnack inequality in the setting of the paper, that
is, for non-negative weak solutions of the horizontal divergence equation
\eqref{meq} on a Carnot group.  The proof has two ingredients: the generalized
Harnack inequality, Theorem~\ref{gha}, and a local upper estimate for
subsolutions under \({\bf(f-I)}\).

\medskip
\noindent\textbf{Map of the section.}
The argument is divided into five parts.  First, we establish boundary
attainment for finite obstacle problems at regular boundary points.  Second,
we use those finite problems to construct a horizontal Keller--Osserman
supersolution which blows up on the boundary of a suitable regular domain.
Third, comparison with this barrier gives a pointwise Keller--Osserman bound
and then the local scale estimate needed for \(M=\sup_{B_1}u\).  Fourth, that
scale estimate controls \(\sqrt{F(t)}/t\) for levels below \(M/4\).  Finally,
we compare the resulting logarithmic lower bound with the generalized Harnack
inequality and conclude that \(\sup_{B_1}u\le C\inf_{B_1}u\).

\subsection{The horizontal equation and finite obstacle boundary values}

This subsection fixes the horizontal equation and proves the boundary fact
needed in the later barrier construction.  The main point is that a finite
obstacle solution assumes its constant boundary datum at every boundary point
which is regular for the linear principal part of the operator.

We work with the variable-coefficient horizontal equation
\begin{equation}\label{eq:blue-horizontal-model-f}
\operatorname{div}_H(A(y)Xw)=f(w)
\end{equation}
in a domain \(\Omega\subset\mathbb G\).  Throughout this subsection
\(A=(a_{ij})\) is measurable, symmetric, and uniformly \(X\)-elliptic, namely
\[
\lambda |\xi|^2\le \langle A(y)\xi,\xi\rangle\le \Lambda |\xi|^2,
\qquad \xi\in\mathbb R^m,
\]
for a.e. \(y\).  A non-negative weak subsolution of
\eqref{eq:blue-horizontal-model-f} means that
\[
-\int_\Omega\langle A(y)Xw,X\phi\rangle\,dy
\ge
\int_\Omega f(w)\phi\,dy
\]
for every non-negative \(\phi\in C_0^1(\Omega)\).  As throughout the paper,
\[
F(t):=\int_0^t f(\tau)\,d\tau,
\]
and in this subsection we also write
\[
\mathcal F(s,a):=\int_a^s f(\tau)\,d\tau,
\qquad 0<a<s.
\]

\medskip
\noindent\textbf{Finite-obstacle terminology.}
Fix constants \(a<M\).  The number \(a\) is the lower obstacle: every
admissible function \(v\) must satisfy \(v\ge a\) almost everywhere in
\(\Omega\).  The number \(M\) is the constant boundary datum in the Sobolev
sense.  More precisely, the condition
\[
v-M\in W_0^{1,2}(\Omega,X)
\]
means that \(v\) has the same Sobolev trace as the constant function \(M\);
it does not require the pointwise boundary values of \(v\) to be known in
advance.  Here \(W_0^{1,2}(\Omega,X)\) is the closure of
\(C_c^\infty(\Omega)\) in the horizontal Sobolev norm.

The admissible functions form a closed convex set.  The variational
inequality below says that the first variation of the energy is non-negative
in every admissible direction \(V-U\).  The inequality is one-sided because
the obstacle constraint may prevent perturbations in the opposite direction.

The \emph{linear principal operator} is the second-order part obtained by
removing the term \(f(u)\) from the equation.  With the sign convention used
in the variational formulation, it is
\[
\mathcal L_Au
:=
X^*(A(y)Xu)
=
\sum_{i,j=1}^m X_i^*\bigl(a_{ij}(y)X_ju\bigr)
=
-\operatorname{div}_H(A(y)Xu).
\]
Thus the equation \(\operatorname{div}_H(A(y)Xu)=f(u)\) is equivalently
\(\mathcal L_Au+f(u)=0\).  In weak form, \(\mathcal L_A\) is represented by
the bilinear form
\[
\mathscr L(u,\phi)
:=
\int_\Omega\langle A(y)Xu,X\phi\rangle\,dy,
\qquad
\langle\mathcal L_Au,\phi\rangle=\mathscr L(u,\phi).
\]

A point \(z\in\partial\Omega\) is \emph{regular for \(\mathcal L_A\)} if
every solution of the homogeneous linear problem \(\mathcal L_Au=0\) with
continuous boundary datum \(g\) satisfies
\[
\lim_{\Omega\ni y\to z}u(y)=g(z).
\]
Equivalently, the linear operator admits a barrier at \(z\).  Finally, a
boundary inequality is understood in the Sobolev sense: \(v\le w\) on
\(\partial G\) means
\[
(v-w)_+:=\max\{v-w,0\}\in W_0^{1,2}(G,X).
\]

\color{black}

\begingroup
\color{black}

\begin{lem}[Scale-invariant nonhomogeneous maximum principle]
\label{lem:scale-invariant-nonhomogeneous-maximum-principle}
Let \(p>Q/2\), let \(\mathcal O\subset B(z,r)\) be an open set, and let
\[
L=-X^*(A X),
\]
where \(A\) is measurable and uniformly horizontally elliptic.  Suppose that
\(h\in L^p(\mathcal O)\) and that \(v\) is a weak subsolution of \(Lv=h\) in
\(\mathcal O\).
Then
\begin{equation*}\tag{SI-MP}\label{eq:GL-scaled-inhomogeneous-maximum-principle}
\sup_{\mathcal O} v^+
\le
\sup_{\partial\mathcal O}v^+
+C r^{2-Q/p}\|h\|_{L^p(\mathcal O)},
\qquad p>Q/2,
\end{equation*}
where the boundary supremum is understood in the Sobolev sense.  The constant
may be written as
\[
C=C(p,Q,\lambda,\Lambda,C_{\mathrm{SP}},C_{\mathrm{vol}}),
\]
where \(\lambda,\Lambda\) are the horizontal ellipticity constants,
\(C_{\mathrm{SP}}\) is the constant in the Sobolev--Poincar\'e inequality
recalled in the proof of Lemma~\ref{pin}, and
\(C_{\mathrm{vol}}:=|B(e,1)|\), so that
\(|B(z,r)|=C_{\mathrm{vol}}r^Q\) \textcolor{black}{by the ball-scaling identity
\eqref{eq:proposed-prelim-ball-scaling}.}  In particular, \(C\) is independent
of \(z\), \(r\), the shape of \(\mathcal O\subset B(z,r)\), and the functions
\(v,h\).
The estimate is the scale-invariant form of the nonhomogeneous maximum
principle stated at the beginning of the proof below.
\end{lem}

\begin{proof}
The only analytic input is Theorem~3.1 of Guti\'errez--Lanconelli \cite{GL},
the nonhomogeneous maximum principle for equations with an \(L^p\) right-hand
side.  In the present notation, if \(p>Q/2\),
\(h\in L^p(\mathcal O)\), and \(v\) is a weak subsolution of \(Lv=h\) in a
bounded domain \(\mathcal O\), then
\begin{equation}\label{eq:GL-inhomogeneous-maximum-principle}
\sup_{\mathcal O} v^+
\le
\sup_{\partial\mathcal O}v^+
+C\|h\|_{L^p(\mathcal O)},
\end{equation}
where the boundary supremum is understood in the Sobolev sense.  The constant
in this unscaled estimate depends on \(p\), the domain \(\mathcal O\), and the
structural constants of the operator and vector fields, but it is independent
of \(v\) and \(h\).  Theorem~5.5 of \cite{GL}, which is the interior
nonhomogeneous Harnack inequality, is not needed here.

We now give the complete rescaling argument.  An estimate is called
\emph{scale invariant} when, after the change of variables
\[
y=T_{z,r}(\xi):=z\delta_r\xi,
\qquad
\xi=T_{z,r}^{-1}(y)=\delta_{1/r}(z^{-1}y),
\]
it becomes the corresponding estimate on a unit-scale domain.  Here
\(z\delta_r\xi\) denotes the Carnot-group product of \(z\) with the dilated
point \(\delta_r\xi\).  Thus \(T_{z,r}\) first dilates by the factor \(r\)
about the identity \(e\), and then left-translates by \(z\).  Homogeneity and
left invariance of the distance give
\[
T_{z,r}(B(e,1))=B(z,r).
\]
Consequently,
\[
\widetilde{\mathcal O}:=T_{z,r}^{-1}(\mathcal O)\subset B(e,1),
\]
so the change of variables moves the center \(z\) to \(e\) and normalizes the
radius \(r\) to one.

Define
\[
\widetilde v(\xi):=v(T_{z,r}(\xi)),
\qquad
\widetilde A(\xi):=A(T_{z,r}(\xi)),
\qquad
\widetilde h(\xi):=r^2h(T_{z,r}(\xi)).
\]
The horizontal vector fields are homogeneous of degree one,
\[
X_j(\varphi\circ\delta_r)
=r(X_j\varphi)\circ\delta_r,
\]
and are invariant under left translations.  Therefore the weak inequality
\(Lv\le h\), with the sign convention for subsolutions used in
\eqref{eq:GL-inhomogeneous-maximum-principle}, transforms into
\[
\widetilde L\widetilde v\le\widetilde h
\qquad\text{in }\widetilde{\mathcal O},
\qquad
\widetilde L=-X^*(\widetilde A X).
\]
The matrix \(\widetilde A\) has exactly the same ellipticity constants as
\(A\).

Applying \eqref{eq:GL-inhomogeneous-maximum-principle} on the unit-scale
domain \(\widetilde{\mathcal O}\subset B(e,1)\) gives
\[
\sup_{\widetilde{\mathcal O}}\widetilde v^+
\le
\sup_{\partial\widetilde{\mathcal O}}\widetilde v^+
+C\|\widetilde h\|_{L^p(\widetilde{\mathcal O})}.
\]
The constant is uniform for \(\widetilde{\mathcal O}\subset B(e,1)\): in the
proof of the nonhomogeneous maximum principle, the domain dependence is
controlled by the unit-scale Sobolev and volume constants, and functions in
\(W_0^{1,2}(\widetilde{\mathcal O},X)\) are extended by zero before the global
Sobolev inequality is applied.

Haar measure has homogeneous degree \(Q\), namely
\(d(T_{z,r}\xi)=r^Qd\xi\).  Hence
\[
\begin{aligned}
\|\widetilde h\|_{L^p(\widetilde{\mathcal O})}
&=
r^2\left(\int_{\widetilde{\mathcal O}}
       |h(T_{z,r}(\xi))|^p\,d\xi\right)^{1/p} \\
&=
r^{2-Q/p}\left(\int_{\mathcal O}|h(y)|^p\,dy\right)^{1/p}
=r^{2-Q/p}\|h\|_{L^p(\mathcal O)}.
\end{aligned}
\]
The interior and Sobolev boundary suprema are unchanged by the bijection
\(T_{z,r}:\widetilde{\mathcal O}\to\mathcal O\).  Substitution in the
unit-scale estimate gives
\eqref{eq:GL-scaled-inhomogeneous-maximum-principle}, proving explicitly that
\textup{(SI-MP)} follows from \eqref{eq:GL-inhomogeneous-maximum-principle}.

More precisely, one may take
\[
C=C(p,Q,\lambda,\Lambda,C_{\mathrm{SP}},C_{\mathrm{vol}}),
\]
with the constants defined in the statement of the lemma.  If the Sobolev and
volume estimates are only locally uniform, then \(z\) must be restricted to a
fixed compact set and \(r\le r_0\), and \(C\) may additionally depend on that
compact set and on \(r_0\).
\end{proof}

\begin{lem}[Boundary attainment for the finite obstacle problem]
\label{lem:finite-obstacle-boundary-attainment}
Let \(\Omega\Subset\mathbb G\) be bounded, let \(a<M\), and put
\[
\mathcal K_{M,a}(\Omega)
:=
\{v\in W^{1,2}(\Omega,X):v-M\in W_0^{1,2}(\Omega,X),\ v\ge a\}.
\]
Suppose that \(U\in\mathcal K_{M,a}(\Omega)\), \(a\le U\le M\), and
\begin{equation*}\tag{VI}\label{eq:finite-obstacle-vi-general}
\int_\Omega\langle A(y)XU,X(V-U)\rangle\,dy
+\int_\Omega f(U)(V-U)\,dy\ge0
\end{equation*}
for every \(V\in\mathcal K_{M,a}(\Omega)\).  If
\(z\in\partial\Omega\) is regular for \(\mathcal L_A\), then
{\color{black}\(U\) has a representative, still denoted by \(U\), for which}
\begin{equation}\label{eq:finite-obstacle-boundary-attainment}
\lim_{\substack{y\to z\\y\in\Omega}}U(y)=M.
\end{equation}
Thus the boundary-attainment assertion is a consequence of linear boundary
regularity; it is not an additional hypothesis on the obstacle solution.
\end{lem}
\endgroup

{\color{black}

\begin{proof}
\medskip
\noindent\textbf{Step 1: Introduce the torsion function and formulate the
boundary reduction.}
Let \(\tau\in W_0^{1,2}(\Omega,X)\) be the torsion function
\[
\int_\Omega\langle AX\tau,X\phi\rangle\,dy
=
\int_\Omega\phi\,dy,
\qquad
\phi\in W_0^{1,2}(\Omega,X).
\]
{\color{black}
Since the negative part \(\tau_-:=\max\{-\tau,0\}\) belongs to
\(W_0^{1,2}(\Omega,X)\) and
\(X\tau_-=-\chi_{\{\tau<0\}}X\tau\), it is an admissible test function.
Taking \(\phi=\tau_-\) gives
\[
-\int_{\{\tau<0\}}\langle AX\tau_-,X\tau_-\rangle\,dy
=
\int_\Omega\tau_-\,dy.
\]
The left-hand side is non-positive by ellipticity, whereas the right-hand side
is non-negative.  Hence both sides vanish, so \(\tau_-=0\) almost everywhere
and therefore \(\tau\ge0\) almost everywhere.  Moreover, after enclosing
\(\Omega\) in a fixed ball, Lemma~\ref{lem:scale-invariant-nonhomogeneous-maximum-principle},
applied with constant right-hand side, gives \(\tau\in L^\infty(\Omega)\).
We use the standard interior H\"older-continuous representative of \(\tau\).
Its existence follows from the nonhomogeneous Harnack inequality for
\(X\)-elliptic equations \cite[Theorem~5.5]{GL}, by the standard
oscillation-decay argument.  The same convention will be used below for the
localized torsion functions \(q_r\) and for the bounded harmonic functions
\(h_r\) and \(H_r\).  We shall prove the following bounded-source consequence
of linear boundary regularity:}
\begin{equation}\label{eq:short-proof-torsion-limit}
\tau(y)\longrightarrow0
\qquad\text{as }\Omega\ni y\to z.
\end{equation}

\medskip
\noindent\textbf{Step 2: Estimate the localized torsion function.}
For completeness, this follows by localization.  {\color{black}Choose a bounded
open set \(D\Subset\mathbb G\) with \(\overline\Omega\subset D\), and fix
\(r>0\) sufficiently small that \(\overline{B(z,r)}\subset D\).}  Put
\(\Omega_r=\Omega\cap B(z,r)\), and let \(q_r\in W_0^{1,2}(\Omega_r,X)\)
solve
\[
\int_{\Omega_r}\langle AXq_r,X\phi\rangle\,dy
=
\int_{\Omega_r}\phi\,dy.
\]
The same negative-part test used above for \(\tau\) gives \(q_r\ge0\) almost
everywhere.  We now derive the upper bound directly from
Lemma~\ref{lem:scale-invariant-nonhomogeneous-maximum-principle}.  Recall that
the operator in
Lemma~\ref{lem:scale-invariant-nonhomogeneous-maximum-principle} is
\[
L=-X^*(A X).
\]
The variational identity defining \(q_r\) says
\[
\langle X^*(A Xq_r),\phi\rangle
=
\int_{\Omega_r}\phi\,dy,
\]
and consequently
\[
Lq_r=-1
\qquad\text{weakly in }\Omega_r.
\]
Thus Lemma~\ref{lem:scale-invariant-nonhomogeneous-maximum-principle}
applies with \(\mathcal O=\Omega_r\subset B(z,r)\), \(v=q_r\), and
\(h=-1\).  Moreover, \(q_r\in W_0^{1,2}(\Omega_r,X)\), and hence $q_r^+\in W^{1,2}_0(\Omega_r,X)$. Therefore,
\(\sup_{\partial\Omega_r}q_r^+=0\) in the sense used in that lemma.  Hence, for
any fixed \(p>Q/2\),
\[
\sup_{\Omega_r}q_r
=
\sup_{\Omega_r}q_r^+
\le
C r^{2-Q/p}\|1\|_{L^p(\Omega_r)}
=
C r^{2-Q/p}|\Omega_r|^{1/p}.
\]
Since \(\Omega_r\subset B(z,r)\) and Haar measure has homogeneous growth,
\[
|\Omega_r|^{1/p}
\le |B(z,r)|^{1/p}
=|B(e,1)|^{1/p}r^{Q/p}.
\]
Absorbing the fixed factor \(|B(e,1)|^{1/p}\) into \(C\) gives
\[
0\le q_r\le Cr^2.
\]
{\color{black}The inequalities just obtained hold almost everywhere.  Since the
representative of \(q_r\) fixed in Step~1 is locally continuous, they hold at
every point of \(\Omega_r\).}

{\color{black}
\medskip
\noindent\textbf{Step 3: Use boundary regularity to control the torsion
function.}
We compare the torsion function directly with a harmonic barrier localized
near \(z\).  By
Lemma~3.2 of \cite{TU}, the regularity of \(z\) for \(\Omega\) implies its
regularity for \(\Omega_r=\Omega\cap B(z,r)\).  Choose a Lipschitz function
\(\eta_r\) satisfying
\[
0\le\eta_r\le1,
\qquad
\eta_r=0\text{ on }B(z,r/3),
\qquad
\eta_r=1\text{ outside }B(z,2r/3),
\]
and let \(H_r\) be the weakly \(\mathscr L\)-harmonic function in \(\Omega_r\)
with Sobolev boundary datum \(\eta_r\).  Thus
\[
H_r-\eta_r\in W_0^{1,2}(\Omega_r,X).
\]
The weak maximum principle gives \(0\le H_r\le1\), while the regularity of
\(z\) gives
\[
H_r(y)\longrightarrow0
\qquad\text{as }\Omega_r\ni y\to z.
\]

Set
\[
K:=\|\tau\|_{L^\infty(\Omega)}.
\]
We claim that
\begin{equation}\label{eq:step-3-direct-torsion-comparison}
0\le\tau\le q_r+KH_r
\qquad\text{in }\Omega_r.
\end{equation}
The lower bound \(0\le\tau\) was proved in Step~1.  To obtain the upper bound
from the weak maximum principle applied to \(\tau-q_r-KH_r\), it remains to
verify that
\[
\tau\le q_r+KH_r
\qquad\text{on }\partial\Omega_r
\]
in the Sobolev sense; equivalently, we must show that
\((\tau-q_r-KH_r)_+\in W_0^{1,2}(\Omega_r,X)\).  Choose
\(\zeta\in C_c^\infty(B(z,r))\) with \(0\le\zeta\le1\) and
\(\zeta=1\) on \(B(z,2r/3)\).  Then
\(\zeta\tau\in W_0^{1,2}(\Omega_r,X)\).  Define
\[
\psi:=(\tau-q_r-K\eta_r)_+.
\]
On \(\Omega_r\setminus B(z,2r/3)\), one has \(\eta_r=1\),
\(q_r\ge0\), and \(\tau\le K\) almost everywhere, so \(\psi=0\) there.
On \(B(z,2r/3)\), where \(\zeta=1\),
\[
0\le\psi
=(\tau-q_r-K\eta_r)_+
\le\tau
=(\zeta\tau)_+,
\]
because \(q_r\ge0\), \(K\eta_r\ge0\), and \(\tau\ge0\).  Together with
\(\psi=0\) outside \(B(z,2r/3)\), this proves the global domination
\[
0\le\psi\le(\zeta\tau)_+
\qquad\text{almost everywhere in }\Omega_r.
\]
Here \(\psi\in W^{1,2}(\Omega_r,X)\), since it is the positive part of a
Sobolev function, while \((\zeta\tau)_+\in W_0^{1,2}(\Omega_r,X)\).
Recall that the lattice-ideal property says that if
\(g\in W_0^{1,2}(\Omega_r,X)\), \(f\in W^{1,2}(\Omega_r,X)\), and
\(|f|\le g\) almost everywhere with \(g\ge0\), then
\(f\in W_0^{1,2}(\Omega_r,X)\).  Applying this property with
\(f=\psi\) and \(g=(\zeta\tau)_+\), using the global domination above, gives
\[
(\tau-q_r-K\eta_r)_+\in W_0^{1,2}(\Omega_r,X).
\]
Furthermore,
\[
(\tau-q_r-KH_r)_+
\le
(\tau-q_r-K\eta_r)_++K|H_r-\eta_r|.
\]
Since both terms on the right belong to \(W_0^{1,2}(\Omega_r,X)\), the same
lattice-ideal property gives
\[
w:=(\tau-q_r-KH_r)_+\in W_0^{1,2}(\Omega_r,X).
\]

Indeed, for every \(\phi\in W_0^{1,2}(\Omega_r,X)\), extended by zero to
\(\Omega\), the defining equations for \(\tau\), \(q_r\), and \(H_r\) give
\[
\int_{\Omega_r}\langle AX\tau,X\phi\rangle\,dy
=\int_{\Omega_r}\phi\,dy,
\qquad
\int_{\Omega_r}\langle AXq_r,X\phi\rangle\,dy
=\int_{\Omega_r}\phi\,dy,
\]
and
\[
\int_{\Omega_r}\langle AXH_r,X\phi\rangle\,dy=0.
\]
Subtracting the second identity and \(K\) times the third from the first gives
\[
\int_{\Omega_r}
\langle AX(\tau-q_r-KH_r),X\phi\rangle\,dy=0.
\]
Thus \(\tau-q_r-KH_r\) is weakly \(\mathscr L\)-harmonic in \(\Omega_r\).
Testing this identity with \(\phi=w\) gives
\[
0
=
\int_{\Omega_r}
\langle AX(\tau-q_r-KH_r),Xw\rangle\,dy
\ge
\lambda\int_{\Omega_r}|Xw|^2\,dy.
\]
The Poincar\'e inequality implies \(w=0\), proving
\eqref{eq:step-3-direct-torsion-comparison}.  The inequality, initially
almost everywhere, holds pointwise for the locally continuous representatives.
Using \(0\le q_r\le Cr^2\) and letting \(y\to z\), we obtain
\[
0\le
\limsup_{\substack{y\to z\\y\in\Omega}}\tau(y)
\le Cr^2.
\]
Finally, letting \(r\downarrow0\) proves
\eqref{eq:short-proof-torsion-limit}.
}

\medskip
\noindent\textbf{Step 4: Derive an inequality for the boundary deficit.}
Now set
\[
v:=M-U.
\]
Since \(U\le M\) and \(U-M\in W_0^{1,2}(\Omega,X)\),
\[
v\ge0,
\qquad
v\in W_0^{1,2}(\Omega,X).
\]
If \(0\le\phi\in W_0^{1,2}(\Omega,X)\), then
\(V=U+\phi\) is an admissible competitor in
\eqref{eq:finite-obstacle-vi-general}.  Therefore
\[
\int_\Omega\langle AXU,X\phi\rangle\,dy
+\int_\Omega f(U)\phi\,dy
\ge0.
\]
Using \(Xv=-XU\), \(U\le M\), and the monotonicity of \(f\), we obtain
\begin{equation}\label{eq:short-proof-v-subsolution}
\int_\Omega\langle AXv,X\phi\rangle\,dy
\le
\int_\Omega f(U)\phi\,dy
\le
f(M)\int_\Omega\phi\,dy.
\end{equation}

\medskip
\noindent\textbf{Step 5: Compare the boundary deficit with the torsion
function.}
Put \(C_M=f(M)\) and
\[
w:=(v-C_M\tau)_+\in W_0^{1,2}(\Omega,X).
\]
Subtracting \(C_M\) times the torsion equation from
\eqref{eq:short-proof-v-subsolution} and taking \(\phi=w\) gives
\[
\begin{aligned}
0
&\ge
\int_\Omega
\langle AX(v-C_M\tau),Xw\rangle\,dy\\
&=
\int_{\{v>C_M\tau\}}
\langle AX(v-C_M\tau),X(v-C_M\tau)\rangle\,dy\\
&\ge
\lambda\int_\Omega|Xw|^2\,dy.
\end{aligned}
\]
Thus \(w=0\), and hence
\begin{equation}\label{eq:short-proof-boundary-deficit-ae}
0\le M-U=v\le f(M)\tau.
\end{equation}

\medskip
\noindent\textbf{Step 6: Choose a pointwise representative and pass to the
boundary.}
{\color{black}This step fixes the representative whose existence is asserted in
the statement of the lemma.}
At this stage, inequality~\eqref{eq:short-proof-boundary-deficit-ae} holds
only almost everywhere, because it was obtained in the Sobolev space.  As usual, the weak
solution is an element of a Sobolev equivalence class and the notation \(U\)
refers to a chosen representative of that class.  Let \(N\subset\Omega\) be a
null set outside of which \eqref{eq:short-proof-boundary-deficit-ae} holds.
Redefine the chosen
representative by setting
\[
U(y):=M
\qquad\text{for }y\in N,
\]
and leave it unchanged on \(\Omega\setminus N\).  This modification does not
change the Sobolev element, its weak derivatives, or the variational
inequality, since all of these are insensitive to changes on null sets.  We
retain the same notation \(U\) for the modified representative.  Because
\(\tau\ge0\) pointwise for the representative fixed above, the choice
\(U=M\) on \(N\) extends the almost-everywhere squeeze to the pointwise
inequality
\[
0\le M-U(y)\le f(M)\tau(y)
\qquad\text{for every }y\in\Omega.
\]
Letting \(y\to z\) and using
\eqref{eq:short-proof-torsion-limit} yields
\[
0\le M-U(y)\longrightarrow0.
\]
Therefore
\[
U(y)\longrightarrow M
\qquad\text{as }\Omega\ni y\to z,
\]
which proves the lemma.
\end{proof}

}

\subsection{Construction of a horizontal blow-up barrier}

The next lemma constructs the supersolution used for comparison.  The
construction starts from finite obstacle minimizers with boundary value
\(M\), chooses the domain radius so that the obstacle is attained at the
center, and then lets \(M\to\infty\).  Lemma~\ref{lem:finite-obstacle-boundary-attainment}
ensures that the limiting supersolution blows up at the boundary.

\medskip
\noindent\textbf{Barrier terminology.}
A bounded domain \(D\subset\mathbb G\) is called \emph{Wiener regular} for
\(\mathcal L_A\) if every point \(z\in\partial D\) is regular; equivalently,
for every continuous boundary datum \(g\), the solution \(H_g\) of
\[
\mathcal L_AH_g=0\quad\text{in }D,
\qquad
H_g=g\quad\text{on }\partial D,
\]
satisfies
\[
\lim_{\substack{y\to z\\y\in D}}H_g(y)=g(z)
\qquad\text{for every }z\in\partial D.
\]

A \emph{horizontal Keller--Osserman barrier} centered at \(y_0\in D\) with
height \(a>0\) is a function
\[
W\in W_{\mathrm{loc}}^{1,2}(D,X)\cap C(D),
\qquad W\ge0,
\]
which is a weak supersolution,
\[
-\int_D\langle A(y)XW,X\phi\rangle\,dy
\le
\int_D f(W)\phi\,dy
\quad
\text{for every }0\le\phi\in C_c^1(D),
\]
has the prescribed center value \(W(y_0)=a\), and blows up at the boundary in
the precise sense
\[
\forall L>0\ \exists\delta_L>0:
\quad
\operatorname{dist}(y,\partial D)<\delta_L
\Longrightarrow W(y)>L.
\]
For each prescribed center height \(a>0\), the construction produces a domain
\(D=D_{a,y_0}\) and a barrier \(W=W_{a,y_0}\).  The number \(R_f(a)>0\) is a
chosen \emph{outer-radius bound} for this domain, meaning
\[
y_0\in D\subset B(y_0,R_f(a)).
\]
Equivalently,
\[
\sup_{y\in D}d(y_0,y)\le R_f(a).
\]
It is not claimed that \(D\) is a ball, that \(R_f(a)\) is the smallest
possible radius, or that \(R_f(a)\) describes the rate at which \(W\) blows up.
We use the following terminology.  The \emph{Keller--Osserman tail integral}
is the improper integral
\[
\mathcal K_f(t):=\int_t^\infty\frac{ds}{\sqrt{F(s)}},
\qquad t>0.
\]
We say that \(f\) satisfies the \emph{Keller--Osserman condition} if this
integral is finite for some, equivalently every, \(t>0\).  In this paper we
record the condition at \(t=1\):
\begin{equation}\label{eq:Keller-Osserman-condition}
\int_1^\infty\frac{ds}{\sqrt{F(s)}}<\infty.
\end{equation}
For a prescribed height \(a>0\), the corresponding \emph{shifted
Keller--Osserman integral} is
\begin{equation}\label{eq:shifted-Keller-Osserman-integral}
\mathcal I_f(a):=
\int_a^\infty\frac{ds}{\sqrt{F(s)-F(a)}}.
\end{equation}
The integral is understood as improper both at \(s=a\) and at infinity.  Since
\(f(a)>0\), its singularity at \(a\) is of order \((s-a)^{-1/2}\) and is
integrable; condition~\eqref{eq:Keller-Osserman-condition} guarantees
integrability at infinity.  The dependence of the outer-radius bound on the
nonlinearity is quantified by this shifted integral: Lemma~\ref{blup} proves
that one may choose
\[
R_f(a)
\le
C_0\mathcal I_f(a),
\]
where \(C_0\) depends only on the ellipticity constants, the Carnot group, and
the fixed homogeneous distance.

\phantomsection
\label{lem:blue-horizontal-barrier-existence}
\begin{lem}[Existence of the horizontal
Keller--Osserman barrier]\label{blup}
Assume that \(f:[0,\infty)\to[0,\infty)\) is continuous and non-decreasing,
that \(f(t)>0\) for \(t>0\), and that it satisfies the Keller--Osserman
condition~\eqref{eq:Keller-Osserman-condition}.
There is a constant
\(C_0=C_0(\lambda,\Lambda,\mathbb G,d)>0\) with the following
property.  For every \(a>0\) and every \(y_0\in\mathbb G\), there are a radius
\(R_f(a)>0\), a bounded Wiener-regular domain
\(\mathcal D_{a,y_0}\), and a function
\[
W_{a,y_0}\in
W_{\mathrm{loc}}^{1,2}(\mathcal D_{a,y_0},X)
\cap C(\mathcal D_{a,y_0})
\]
such that
\[
y_0\in\mathcal D_{a,y_0}
\subset B(y_0,R_f(a)),
\qquad W_{a,y_0}\ge0,
\]
\[
-\int_{\mathcal D_{a,y_0}}
\langle A(y)XW_{a,y_0},X\phi\rangle\,dy
\le
\int_{\mathcal D_{a,y_0}}f(W_{a,y_0})\phi\,dy
\]
for every non-negative \(\phi\in C_c^1(\mathcal D_{a,y_0})\),
\[
W_{a,y_0}(y_0)=a,
\qquad
W_{a,y_0}(y)\longrightarrow\infty
\quad\text{as }\operatorname{dist}(y,\partial\mathcal D_{a,y_0})\downarrow0,
\]
and
\[
R_f(a)
\le
C_0\int_a^\infty
\frac{ds}{\sqrt{F(s)-F(a)}}.
\]
\end{lem}

\begingroup
\color{black}

\begin{proof}
\medskip
\noindent\textbf{Step 1: the Keller--Osserman integral.}
Set
\[
I(t):=\int_t^\infty\frac{ds}{\sqrt{F(s)-F(t)}}
=\int_0^\infty
\frac{ds}{\sqrt{\int_0^s f(t+\xi)\,d\xi}},
\qquad t>0.
\]
The first integral is finite: near its lower endpoint this follows from
\(f(t)>0\), and at infinity it follows from the Keller--Osserman
condition.  The second representation shows that \(I\) is decreasing.
Moreover, the Keller--Osserman condition forces \(f(t)\to\infty\);
hence dominated convergence gives
\begin{equation}\label{eq:concise-92-I-zero}
I(t)\downarrow0\qquad\text{as }t\to\infty.
\end{equation}

\medskip
\noindent\textbf{Step 2: a uniformly regular family of reference domains.}
Choose the bounded connected reference domain \(D_*\) furnished by
\cite[Lemma~3.8]{UG}.  Thus there are structural constants
\(0<\sigma<\Sigma\), \(s_*>0\), and \(\vartheta>0\) such that
\[
B(e,\sigma)\subset D_*\subset B(e,\Sigma)
\]
and
\[
|B_s(z)\setminus D_*|\ge\vartheta|B_s(z)|
\qquad
(z\in\partial D_*,\ 0<s<s_*).
\]
In particular,
\[
d_*:=\operatorname{dist}(e,\partial D_*)\ge\sigma>0.
\]
For
\[
D_\rho(y_0):=y_0\delta_\rho(D_*),
\]
left invariance and homogeneity give
\[
B(y_0,\sigma\rho)\subset D_\rho(y_0)
\subset B(y_0,\Sigma\rho),
\qquad
\operatorname{dist}(y_0,\partial D_\rho(y_0))=d_*\rho.
\]
The scaled exterior density estimate, together with
\cite[Proposition~4.4]{TU}, implies that every boundary point of every
\(D_\rho(y_0)\) is Wiener regular. {\color{black}Let us justify the assumptions of Proposition 4.4 in \cite{TU}. It is enough to check the exterior-density hypothesis.  If
\(z\in\partial\Omega_\rho\), then
\[
z=y_0\delta_\rho(z_*)
\qquad\text{for some }z_*\in\partial D_*.
\]
For balls of the fixed left-invariant homogeneous distance \(d\),
homogeneity gives
\[
B_d(z,s)=y_0\delta_\rho\bigl(B_d(z_*,s/\rho)\bigr).
\]
Since Haar measure scales by \(\rho^Q\),
\[
\begin{aligned}
|B_d(z,s)\setminus\Omega_\rho|
&=\rho^Q
|B_d(z_*,s/\rho)\setminus D_*|\\
&\ge
\vartheta\rho^Q|B_d(z_*,s/\rho)|
=\vartheta|B_d(z,s)|
\end{aligned}
\]
whenever \(0<s<\rho s_*\).

If the fixed distance \(d\) is chosen to be \(d_X\), this is precisely
the density condition in Proposition~4.4.  For an arbitrary fixed
homogeneous distance, use its equivalence with \(d_X\): there are
constants \(\alpha,\beta>0\), depending only on the two distances, such
that
\[
B_d(z,\alpha r)\subset B_{d_X}(z,r)
\subset B_d(z,\beta r).
\]
It follows, for \(0<r<\rho s_*/\alpha\), that
\[
\begin{aligned}
|B_{d_X}(z,r)\setminus\Omega_\rho|
&\ge |B_d(z,\alpha r)\setminus\Omega_\rho|\\
&\ge \vartheta|B_d(z,\alpha r)|\\
&=\vartheta_X|B_{d_X}(z,r)|,
\end{aligned}
\]
where
\[
\vartheta_X
:=
\vartheta\,
\frac{|B_d(e,1)|\alpha^Q}{|B_{d_X}(e,1)|}>0.
\]
Thus, at every \(z\in\partial\Omega_\rho\), the hypothesis of
\cite[Proposition~4.4]{TU} holds with, for example,
\[
r_1:=\frac{\rho s_*}{2\alpha}
\qquad\text{and}\qquad
\theta:=\vartheta_X.
\]
The proposition therefore shows that every boundary point of
\(D_\rho(y_0)\) is regular for the present operator.  Since \(\rho\) and
\(y_0\) were arbitrary, this proves the asserted Wiener regularity for
the whole family of scaled domains.
}

\medskip
\noindent\textbf{Step 3: the admissible class and coercivity.}
Steps 3 and 4 together implement the direct method of the calculus of
variations.  In this step we verify the properties needed for compactness:
the admissible class is nonempty, convex, and weakly closed, and every
minimizing sequence is bounded.  Step 4 then uses these facts to extract a
weak limit, prove that it minimizes the energy, and establish uniqueness.
For the moment fix an arbitrary \(\rho>0\).  To lighten notation, we write
\(U_M\), \(J_M\), and \(\mathcal K_M\) instead of \(U_{\rho,M}\),
\(J_{\rho,M}\), and \(\mathcal K_{\rho,M}\); the dependence of the domain
\(D_\rho(y_0)\) on \(\rho\) remains explicit.
For \(M>a\), let \(U_M\) minimize
\[
J_M(u):=
\frac12\int_{D_\rho(y_0)}\langle AXu,Xu\rangle
+\int_{D_\rho(y_0)}F(u)
\]
over
\[
\mathcal K_M:=\{u\in M+W_0^{1,2}(D_\rho(y_0),X):u\ge a\}.
\]
We give the details of existence and uniqueness.  The class
\(\mathcal K_M\) is nonempty because it contains the constant function \(M\), and it
is convex.  Write \(u=M+w\), where
\(w\in W_0^{1,2}(D_\rho(y_0),X)\).  Since \(u\ge a>0\), one has
\(F(u)\ge0\), and ellipticity gives
\[
J_M(u)\ge\frac{\lambda}{2}
\int_{D_\rho(y_0)}|Xw|^2.
\]
To verify explicitly that every minimizing sequence is bounded in
\(W^{1,2}(D_\rho(y_0),X)\), let
\((u_j)\subset \mathcal K_M\) be a minimizing sequence and write
\(u_j=M+w_j\), \(w_j\in W_0^{1,2}(D_\rho(y_0),X)\).  By definition,
\emph{minimizing sequence} means that
\[
J_M(u_j)\longrightarrow m_M:=\inf_{u\in \mathcal K_M}J_M(u).
\]
Since \(M\in \mathcal K_M\),
\[
m_M\le J_M(M)<\infty.
\]
Taking \(\varepsilon=1\) in the definition of convergence, there is an
index \(j_0\) such that, for every \(j\ge j_0\),
\[
J_M(u_j)\le m_M+1\le J_M(M)+1.
\]
After discarding the finitely many terms with \(j<j_0\), this estimate
holds for the entire sequence under consideration. It follows from the preceding coercivity estimate that
\[
\int_{D_\rho(y_0)}|Xw_j|^2
\le\frac{2}{\lambda}\bigl(J_M(M)+1\bigr).
\]
Since \(w_j\) has zero Sobolev boundary values, the horizontal Poincar\'e
inequality gives
\[
\int_{D_\rho(y_0)}|w_j|^2
\le C_P\int_{D_\rho(y_0)}|Xw_j|^2.
\]
Consequently
\[
\sup_j\|u_j\|_{W^{1,2}(D_\rho(y_0),X)}^2
\le C\bigl(M^2|D_\rho(y_0)|+J_M(M)+1\bigr)<\infty.
\]
Therefore, we have shown that the minimizing sequence is bounded.

We also verify that \(\mathcal K_M\) is closed.  Suppose that
\(u_j\in \mathcal K_M\) and \(u_j\to u\) strongly in \(W^{1,2}\).  Since
\(W_0^{1,2}(D_\rho(y_0),X)\) is norm closed and
\(u_j-M\in W_0^{1,2}\), it follows that
\(u-M\in W_0^{1,2}\).  Choose measurable representatives of \(u\) and
the countable family \((u_j)\).  Since each \(u_j\ge a\) almost everywhere,
we may discard the union of the corresponding null sets and obtain, for
every \(j\),
\[
(a-u(y))_+\le |u(y)-u_j(y)|
\qquad\text{for almost every }y\in D_\rho(y_0).
\]
Indeed, the left-hand side is zero when \(u(y)\ge a\); when \(u(y)<a\),
the inequality follows from
\[
a-u(y)\le u_j(y)-u(y)=|u_j(y)-u(y)|,
\]
because \(u_j(y)\ge a>u(y)\).
Thus
\[
\|(a-u)_+\|_{L^2}
\le\|u-u_j\|_{L^2}\longrightarrow0,
\]
so \(u\ge a\) almost everywhere.  Hence \(u\in \mathcal K_M\), and \(\mathcal K_M\) is
strongly (that is, norm) closed in \(W^{1,2}\).

Clearly \(\mathcal K_M\) is also convex. A norm-closed convex subset of a
Banach space is weakly closed; this follows, for example, from the
Hahn--Banach separation theorem.  Consequently, \(\mathcal K_M\) is both strongly
closed and weakly closed in \(W^{1,2}(D_\rho(y_0),X)\).  Here convexity is
essential for the last implication: a norm-closed set need not, in
general, be weakly closed.

\medskip
\noindent\textbf{Step 4: existence and uniqueness of the minimizer.}
We now apply the compactness and closedness established in Step 3.
The space \(W^{1,2}(D_\rho(y_0),X)\), endowed with the norm
\[
\|v\|_{W^{1,2}}^2
=
\|v\|_{L^2}^2+\|Xv\|_{L^2}^2,
\]
is a Hilbert space and hence is reflexive.  Therefore the bounded
minimizing sequence has a subsequence, still denoted by \((u_j)\), and a
function \(u\in W^{1,2}(D_\rho(y_0),X)\) such that
\[
u_j\rightharpoonup u
\qquad\text{weakly in }W^{1,2}(D_\rho(y_0),X).
\]
Since \(\mathcal K_M\) is weakly closed and every \(u_j\) belongs to \(\mathcal K_M\), the
weak limit also satisfies
\[
u\in \mathcal K_M.
\]

To obtain stronger convergence of the functions, put
\(w_j:=u_j-M\).  The sequence \((w_j)\) is bounded in
\(W_0^{1,2}(D_\rho(y_0),X)\).  The horizontal Rellich compactness theorem,
applied on the bounded domain \(D_\rho(y_0)\), says that the embedding
\[
W_0^{1,2}(D_\rho(y_0),X)
\hookrightarrow L^2(D_\rho(y_0))
\]
is compact.  Hence, after passing to a further subsequence,
\(w_j\to w\) strongly in \(L^2\) for some \(w\).  On the other hand, the
weak \(W^{1,2}\) convergence above implies
\(w_j=u_j-M\rightharpoonup u-M\) weakly in \(L^2\).  Since a strongly
convergent sequence has the same weak and strong limit, \(w=u-M\), and
therefore
\[
u_j\longrightarrow u
\qquad\text{strongly in }L^2(D_\rho(y_0)).
\]
Finally, strong \(L^2\) convergence yields an almost-everywhere
convergent subsequence.  
We spell out the passage to the limit in the energy.  Weak convergence in
\(W^{1,2}\) gives
\[
Xu_j\rightharpoonup Xu
\qquad\text{in }L^2(D_\rho(y_0);\mathbb R^m).
\]
Since \(A(y)\) is symmetric and uniformly elliptic, its positive square
root \(A^{1/2}(y)\) is measurable and bounded, and multiplication by
\(A^{1/2}\) is a bounded linear operator on \(L^2\).  Hence
\[
A^{1/2}Xu_j\rightharpoonup A^{1/2}Xu
\qquad\text{in }L^2.
\]
The \(L^2\)-norm is weakly lower semicontinuous, and therefore
\[
\int_{D_\rho(y_0)}\langle AXu,Xu\rangle\,dy
\le
\liminf_{j\to\infty}
\int_{D_\rho(y_0)}\langle AXu_j,Xu_j\rangle\,dy.
\]
For the potential term, the almost-everywhere convergence and the
continuity of \(F\) give \(F(u_j)\to F(u)\) almost everywhere.  Moreover,
\(u_j\ge a>0\) and \(f\ge0\), so \(F(u_j)\ge0\).  Fatou's lemma thus yields
\[
\int_{D_\rho(y_0)}F(u)\,dy
\le
\liminf_{j\to\infty}\int_{D_\rho(y_0)}F(u_j)\,dy.
\]
Combining the two estimates gives
\[
J_M(u)\le\liminf_{j\to\infty}J_M(u_j)=m_M.
\]
Since \(u\in \mathcal K_M\), the definition of \(m_M\) gives the reverse inequality
\(m_M\le J_M(u)\).  Thus \(J_M(u)=m_M\), and the limit is a minimizer,
which we denote by \(U_M\).  Since \(f\) is non-decreasing, \(F\) is
convex; the Dirichlet term is strictly convex on the affine space, by
ellipticity and Poincar\'e's inequality.  Consequently the minimizer is
unique.

\medskip
\noindent\textbf{Step 5: bounds, variational inequality, and regularity.}
The obstacle constraint already gives \(U_M\ge a\).  Moreover,
\(\widetilde U_M:=\min\{U_M,M\}\) belongs to \(\mathcal K_M\).  The Sobolev chain
rule, the monotonicity of \(F\), and \(f\ge0\) give
\[
J_M(\widetilde U_M)\le J_M(U_M).
\]
Uniqueness therefore implies \(U_M=\widetilde U_M\), so \(U_M\le M\).

It remains to prove that \(U_M\) satisfies the obstacle variational
inequality with obstacle \(a\) and boundary datum \(M\), namely the integral
inequality in \eqref{eq:concise-92-VI} for every competitor
\(V\in\mathcal K_M\).  First take a bounded
\(V\in \mathcal K_M\), put \(H=V-U_M\), and set
\(U_t=U_M+tH\), \(0\le t\le1\).  Convexity of \(\mathcal K_M\) gives
\(U_t\in \mathcal K_M\), and minimality gives
\[
0\le\frac{J_M(U_t)-J_M(U_M)}{t}.
\]
Because \(U_M\) and \(V\) are bounded, differentiation at \(t=0+\) is
justified by dominated convergence and yields
\[
\int_{D_\rho(y_0)}\langle AXU_M,XH\rangle\,dy
+\int_{D_\rho(y_0)}f(U_M)H\,dy\ge0.
\]
For an arbitrary \(V\in \mathcal K_M\), take
\(V_N:=\min\{V,N\}\), with \(N>M\).  Then \(V_N\in \mathcal K_M\),
and \(V_N\to V\) in \(W^{1,2}\).  Indeed,
\(V-V_N=(V-N)_+\to0\) in \(L^2\), while the Sobolev chain rule gives
\[
X(V-V_N)=\chi_{\{V>N\}}XV\longrightarrow0
\qquad\text{in }L^2,
\]
by dominated convergence.  Moreover, \(f(U_M)\in L^\infty\) because
\(a\le U_M\le M\).  Consequently,
\[
\left|\int_{D_\rho(y_0)}
\langle AXU_M,X(V_N-V)\rangle\,dy\right|
\le
\Lambda\|XU_M\|_{L^2}\|X(V_N-V)\|_{L^2}\longrightarrow0,
\]
and, since \(D_\rho(y_0)\) has finite measure,
\[
\left|\int_{D_\rho(y_0)}f(U_M)(V_N-V)\,dy\right|
\le
\|f(U_M)\|_{L^\infty}|D_\rho(y_0)|^{1/2}
\|V_N-V\|_{L^2}\longrightarrow0.
\]
Passing to the limit in the variational inequality for \(V_N\) therefore gives
\begin{equation}\label{eq:concise-92-VI}
a\le U_M\le M,
\qquad
\int_{D_\rho(y_0)}\langle AXU_M,X(V-U_M)\rangle\,dy
+\int_{D_\rho(y_0)} f(U_M)(V-U_M)\,dy\ge0
\end{equation}
for every \(V\in \mathcal K_M\).  We next justify the continuity needed in
the subsequent pointwise arguments: \(U_M\) has a locally H\"older continuous
representative in \(D_\rho(y_0)\), which extends continuously to the boundary
and assumes there the value \(M\), as recorded in
\eqref{eq:concise-92-boundary-data}.  Interior continuity follows from the full
obstacle variational inequality \eqref{eq:concise-92-VI}.  Fix
\(x_0\in D_\rho(y_0)\) and \(r>0\) such
that
\[
B(x_0,2r)\Subset D_\rho(y_0).
\]
Choose \(\eta\in C_c^1(B(x_0,2r))\) with \(0\le\eta\le1\), and let
\(k\ge a\).  With
\(w_+=(U_M-k)_+\), the function
\[
V=U_M-\eta^2w_+
\]
still belongs to \(\mathcal K_M\): on \(\{U_M>k\}\) one has
\[
V=(1-\eta^2)U_M+\eta^2k\ge k\ge a,
\]
whereas \(V=U_M\ge a\) on \(\{U_M\le k\}\).  Substitution in
\eqref{eq:concise-92-VI}, followed by ellipticity, Cauchy--Schwarz and
Young's inequality, gives
\begin{equation}\label{eq:concise-92-upper-caccioppoli}
\int\eta^2|Xw_+|^2
\le C\int w_+^2|X\eta|^2.
\end{equation}

{\color{black}

\begingroup
\color{black}
We justify this as follows. Here and in the next estimate all integrals are over
\(D_\rho(y_0)\).  We give the details.  The horizontal Sobolev chain
rule for truncations gives
\[
Xw_+=\chi_{\{U_M>k\}}XU_M
\qquad\text{a.e.}
\]
Taking \(V-U_M=-\eta^2w_+\) in
\eqref{eq:concise-92-VI}, multiplying the resulting inequality by
\(-1\), and using that \(w_+=0\) and \(Xw_+=0\) a.e. outside
\(\{U_M>k\}\), we obtain
\[
\begin{aligned}
0
&\ge
\int\langle AXU_M,X(\eta^2w_+)\rangle
+\int f(U_M)\eta^2w_+\\
&=
\int\eta^2\langle AXw_+,Xw_+\rangle
+2\int\eta w_+\langle AXw_+,X\eta\rangle
+\int f(U_M)\eta^2w_+.
\end{aligned}
\]
The last integral is non-negative because
\(f(U_M),\eta,w_+\ge0\), and hence it may be discarded from the
left-hand side.  If
\[
\lambda|\xi|^2\le\langle A(y)\xi,\xi\rangle,
\qquad
|A(y)\xi|\le\Lambda|\xi|,
\]
then
\[
\lambda\int\eta^2|Xw_+|^2
\le
2\Lambda\int\eta w_+|Xw_+|\,|X\eta|.
\]
Young's inequality, in the form
\[
2\Lambda\eta w_+|Xw_+|\,|X\eta|
\le
\frac{\lambda}{2}\eta^2|Xw_+|^2
+\frac{2\Lambda^2}{\lambda}w_+^2|X\eta|^2,
\]
allows the first term on the right to be absorbed.  This proves
\eqref{eq:concise-92-upper-caccioppoli}, for example with
\(C=4\Lambda^2/\lambda^2\).
\endgroup
Similarly, for \(w_-=(k-U_M)_+\), the upward variation
\[
V=U_M+\eta^2w_-
\]
is admissible and yields
\begin{equation}\label{eq:concise-92-lower-caccioppoli}
\int\eta^2|Xw_-|^2
\le C\int w_-^2|X\eta|^2
+C\int f(U_M)\eta^2w_-.
\end{equation}
\begingroup
\color{black}
For this estimate the truncation chain rule gives
\[
Xw_-=-\chi_{\{U_M<k\}}XU_M
\qquad\text{a.e.}
\]
Substituting \(V-U_M=\eta^2w_-\) in
\eqref{eq:concise-92-VI} and restricting the gradient terms to
\(\{U_M<k\}\), where \(XU_M=-Xw_-\), gives
\[
\begin{aligned}
0
&\le
\int\langle AXU_M,X(\eta^2w_-)\rangle
+\int f(U_M)\eta^2w_-\\
&=
-\int\eta^2\langle AXw_-,Xw_-\rangle
-2\int\eta w_-\langle AXw_-,X\eta\rangle
+\int f(U_M)\eta^2w_-.
\end{aligned}
\]
Consequently, ellipticity and the bound on \(A\) imply
\[
\lambda\int\eta^2|Xw_-|^2
\le
2\Lambda\int\eta w_-|Xw_-|\,|X\eta|
+\int f(U_M)\eta^2w_-.
\]
Applying the same Young inequality as above and absorbing
\((\lambda/2)\int\eta^2|Xw_-|^2\) yields
\[
\int\eta^2|Xw_-|^2
\le
\frac{4\Lambda^2}{\lambda^2}
\int w_-^2|X\eta|^2
+\frac{2}{\lambda}
\int f(U_M)\eta^2w_-.
\]
Thus \eqref{eq:concise-92-lower-caccioppoli} holds with a constant
depending only on the ellipticity constants.  Notice that, unlike in
the upper estimate, the non-negative zeroth-order term cannot be
discarded because it occurs on the right after the energy term is
isolated.
\endgroup

}

Since \(a\le U_M\le M\), one has \(0\le f(U_M)\le f(M)\).  Hence
\eqref{eq:concise-92-upper-caccioppoli} and
\eqref{eq:concise-92-lower-caccioppoli} are the local obstacle De Giorgi
estimates with a bounded inhomogeneous term.  Thus \(U_M\) is a
variational solution of an obstacle problem with constant obstacle and
bounded datum, which is precisely the setting of the interior obstacle
regularity theorem in \cite{GMobs}; equivalently, one may apply the
horizontal De Giorgi oscillation theorem of \cite{MarchiDG}.  Consequently
\(U_M\) has a locally H\"older continuous representative.  On
\(\{U_M>a\}\), both signs of compactly supported variations are
admissible, and therefore
\begin{equation}\label{eq:concise-92-noncontact}
\operatorname{div}_H(AXU_M)=f(U_M)
\qquad\text{weakly in }\{U_M>a\}.
\end{equation}

\begingroup
\color{black}
We spell out this conclusion.  Since the representative of \(U_M\) is
continuous, the non-contact set
\[
\mathcal O_M:=\{y\in D_\rho(y_0):U_M(y)>a\}
\]
is open.  Let \(\phi\in C_c^1(\mathcal O_M)\), and put
\(K:=\operatorname{supp}\phi\).  Compactness of \(K\), continuity of
\(U_M\), and the strict inequality \(U_M>a\) on \(K\) give
\[
\delta_K:=\min_{y\in K}\bigl(U_M(y)-a\bigr)>0.
\]
If \(\phi\not\equiv0\), choose
\[
0<t<\frac{\delta_K}{\|\phi\|_{L^\infty}};
\]
when \(\phi\equiv0\), the desired identity is immediate.  Because
\(t\phi\in W_0^{1,2}(D_\rho(y_0),X)\), both functions
\[
V_+:=U_M+t\phi,
\qquad
V_-:=U_M-t\phi
\]
have the same Sobolev boundary values as \(U_M\).  Moreover, on \(K\),
\[
V_\pm
\ge U_M-t|\phi|
\ge a+\delta_K-t\|\phi\|_{L^\infty}
\ge a,
\]
whereas outside \(K\) one has \(V_\pm=U_M\ge a\).  Thus
\(V_+,V_-\in\mathcal K_M\).

Applying \eqref{eq:concise-92-VI} first with \(V=V_+\) and dividing by
\(t>0\) yields
\[
\int_{D_\rho(y_0)}\langle AXU_M,X\phi\rangle\,dy
+\int_{D_\rho(y_0)}f(U_M)\phi\,dy\ge0.
\]
Applying the same inequality with \(V=V_-\) and dividing by \(-t<0\)
reverses the inequality and gives the opposite bound.  Consequently,
\[
\int_{D_\rho(y_0)}\langle AXU_M,X\phi\rangle\,dy
+\int_{D_\rho(y_0)}f(U_M)\phi\,dy=0
\qquad
\text{for every }\phi\in C_c^1(\mathcal O_M).
\]
Equivalently,
\[
-\int_{D_\rho(y_0)}\langle AXU_M,X\phi\rangle\,dy
=\int_{D_\rho(y_0)}f(U_M)\phi\,dy,
\]
which is exactly the weak formulation of
\(\operatorname{div}_H(AXU_M)=f(U_M)\) in \(\mathcal O_M\), and proves
\eqref{eq:concise-92-noncontact}.
\endgroup
Since the domain is regular, Lemma
\ref{lem:finite-obstacle-boundary-attainment} also gives
\begin{equation}\label{eq:concise-92-boundary-data}
U_M\in C(\overline{D_\rho(y_0)}),
\qquad U_M=M\quad\text{on }\partial D_\rho(y_0).
\end{equation}

\medskip
\noindent\textbf{Step 6: distance from a non-contact point to the boundary.}
We next record the only quantitative estimate needed.  Suppose
\(x\in D_\rho(y_0)\) and \(b:=U_M(x)>a\).  Put
\[
v:=(U_M-b)_+,\qquad
g_b(0):=0,\qquad g_b(t):=f(b+t)\quad(t>0).
\]
A smooth convex approximation of the positive part, tested in
\eqref{eq:concise-92-noncontact}, gives
\[
\operatorname{div}_H(AXv)\ge g_b(v)
\qquad\text{weakly in }D_\rho(y_0).
\]
Indeed, choose smooth convex functions \(\theta_\epsilon\) such that
\[
\theta_\epsilon=0\text{ on }(-\infty,b],
\qquad 0\le\theta_\epsilon'\le1,
\qquad \theta_\epsilon''\ge0,
\qquad
\theta_\epsilon'(t)\longrightarrow1\quad\text{for every }t>b.
\]
These properties imply
\(0\le\theta_\epsilon(t)\le(t-b)_+\) and
\(\theta_\epsilon(t)\to(t-b)_+\) for every \(t\in\mathbb R\).
If \(0\le\phi\in C_c^1(D_\rho(y_0))\), then
\(\theta_\epsilon'(U_M)\phi\) vanishes outside \(\{U_M>b\}\), and its
support is contained in the compact set
\[
\operatorname{supp}\phi\cap\{U_M\ge b\}\Subset\{U_M>a\}.
\]
It is therefore an admissible test function in
\eqref{eq:concise-92-noncontact}: it belongs to the relevant local
zero-boundary Sobolev space, and its use is justified by approximation
with smooth compactly supported test functions.  The Sobolev chain rule
then gives
\begin{align*}
-\int\langle AX\theta_\epsilon(U_M),X\phi\rangle
&=
\int f(U_M)\theta_\epsilon'(U_M)\phi \\
&\quad+
\int\theta_\epsilon''(U_M)
\langle AXU_M,XU_M\rangle\phi \\
&\ge
\int f(U_M)\theta_\epsilon'(U_M)\phi.
\end{align*}
We give the limit passage explicitly.  Let
\(K\Subset D_\rho(y_0)\).  Since
\(0\le\theta_\epsilon(U_M)\le (U_M-b)_+=v\) and
\(\theta_\epsilon(U_M)\to v\) pointwise, dominated convergence gives
\[
\theta_\epsilon(U_M)\longrightarrow v
\qquad\text{in }L^2(K).
\]
The Sobolev chain rule gives
\[
X\theta_\epsilon(U_M)
=\theta_\epsilon'(U_M)XU_M,
\qquad
Xv=\chi_{\{U_M>b\}}XU_M
\quad\text{a.e.}
\]
The factors \(\theta_\epsilon'(U_M)\) converge pointwise to
\(\chi_{\{U_M>b\}}\), including on \(\{U_M=b\}\), where both sides
vanish.  Moreover, their difference is bounded by \(1\).  Hence
\[
\int_K
\left|X\theta_\epsilon(U_M)-Xv\right|^2
\le
\int_K |XU_M|^2
\left|\theta_\epsilon'(U_M)
-\chi_{\{U_M>b\}}\right|^2
\longrightarrow0
\]
by dominated convergence.  Thus
\[
\theta_\epsilon(U_M)\longrightarrow v
\quad\text{in }W_{\mathrm{loc}}^{1,2}(D_\rho(y_0),X),
\]
as claimed.  Since \(a\le U_M\le M\), one also has
\[
0\le f(U_M)\theta_\epsilon'(U_M)\le f(M).
\]
Therefore dominated convergence on \(K\) gives
\[
f(U_M)\theta_\epsilon'(U_M)
\longrightarrow
f(U_M)\chi_{\{U_M>b\}}=g_b(v)
\quad\text{in }L^1_{\mathrm{loc}}.
\]
Indeed, the last identity holds because \(v=U_M-b>0\) on
\(\{U_M>b\}\), whereas \(v=0\) and \(g_b(v)=g_b(0)=0\) on
\(\{U_M\le b\}\).  Finally, the strong \(L^2\)-convergence of the
horizontal gradients and the boundedness of \(A\) allow passage to the
limit on the left-hand side, while the \(L^1\)-convergence above allows
passage on the right-hand side.  Thus, for every non-negative
\(\phi\in C_c^1(D_\rho(y_0))\),
\[
-\int\langle AXv,X\phi\rangle
\ge
\int g_b(v)\phi,
\]
which is precisely the asserted weak subsolution inequality.

We also record explicitly why the jump of \(g_b\) at zero does not
invalidate the next application.  The proof of Theorem
\ref{thm:sub-integral-lower-bound} remains valid for a non-negative,
non-decreasing function \(g\) that is continuous on \((0,\infty)\), even
if \(g(0)<g(0+)\).  Indeed, every center value \(\mu_{k-1}\) occurring in
that proof is strictly positive.  Lemma
\ref{lem:sub-degiorgi-scale} is consequently used only at the positive
level \(\mu_{k-1}/2\), and the subsequent dyadic summation uses only the
monotonicity of the primitive
\[
G(t):=\int_0^t g(s)\,ds.
\]
Neither argument uses continuity of \(g\) at zero.  Hence Theorem
\ref{thm:sub-integral-lower-bound}, with \(f,F\) replaced by \(g_b,G_b\),
applies to \(v\).

Choose \(x_j\to x\) with \(U_M(x_j)>b\).  Such points exist by the
strong maximum principle applied to
\eqref{eq:concise-92-noncontact}; otherwise \(x\) would be a local
maximum of \(U_M\) in the non-contact set, which is impossible since
\(f(U_M)>0\) there; see, for example, \cite{GL}.  Put
\[
r_j:=\frac14\operatorname{dist}(x_j,\partial D_\rho(y_0)).
\]
Then \(B(x_j,2r_j)\Subset D_\rho(y_0)\).
Applying Theorem \ref{thm:sub-integral-lower-bound} to \(v\) in
\(B(x_j,2r_j)\), write
\[
m_j:=\inf_{B(x_j,r_j)}v,
\qquad M_j:=\sup_{B(x_j,r_j)}v.
\]
Since the primitive of \(g_b\) is
\(G_b(t)=F(b+t)-F(b)\) for \(t\ge0\), that theorem gives
\[
cr_j
\le\int_{m_j/4}^{M_j}\frac{dt}{\sqrt{G_b(t)}}
\le\int_0^{M-b}\frac{dt}{\sqrt{F(b+t)-F(b)}}.
\]
Letting \(j\to\infty\) therefore gives
\begin{equation}\label{eq:concise-92-distance}
\operatorname{dist}(x,\partial D_\rho(y_0))
\le
C\int_b^M\frac{ds}{\sqrt{F(s)-F(b)}}
\le CI(b),
\end{equation}
where \(C\) is structural.

\medskip
\noindent\textbf{Step 7: fixing the domain and forcing contact at the center.}
Now fix
\begin{equation}\label{eq:concise-92-rho}
\rho:=\frac{2C}{d_*}I(a),
\qquad D:=D_\rho(y_0).
\end{equation}
For this choice every minimizer touches the obstacle at the center.
Indeed, if \(U_M(y_0)>a\), then
\eqref{eq:concise-92-distance} and the monotonicity of \(I\) imply
\[
2CI(a)=d_*\rho
=\operatorname{dist}(y_0,\partial D)
\le CI(U_M(y_0))
\le CI(a),
\]
which is impossible.  Hence
\begin{equation}\label{eq:concise-92-center}
U_M(y_0)=a\qquad(M>a).
\end{equation}

\medskip
\noindent\textbf{Step 8: construction of \(W_{a,y_0}\) by monotone limit
and boundary blow-up.}
The minimizers increase with \(M\).  In fact, if
\(a<M_1\le M_2\), the lattice competitors
\(\min\{U_{M_1},U_{M_2}\}\) and
\(\max\{U_{M_1},U_{M_2}\}\) in
\eqref{eq:concise-92-VI} give, with
\(w=(U_{M_1}-U_{M_2})_+\),
\[
\lambda\int_D|Xw|^2
+\int_{\{w>0\}}
\bigl(f(U_{M_1})-f(U_{M_2})\bigr)w
\le0.
\]

\begingroup
\color{black}
Here are the details of this lattice argument.  Write
\[
U_1:=U_{M_1},\qquad U_2:=U_{M_2},
\qquad w:=(U_1-U_2)_+.
\]
Put \(z_i:=U_i-M_i\in W_0^{1,2}(D,X)\).  Since \(M_1\le M_2\),
\[
w
=\bigl(z_1-z_2-(M_2-M_1)\bigr)_+.
\]
The function \(T(s):=(s-(M_2-M_1))_+\) is Lipschitz and satisfies
\(T(0)=0\).  The Sobolev truncation theorem, applied to
\(z_1-z_2\in W_0^{1,2}(D,X)\), therefore gives
\[
w=(U_1-U_2)_+\in W_0^{1,2}(D,X),
\qquad
Xw=\chi_{\{U_1>U_2\}}(XU_1-XU_2)
\quad\text{a.e.}
\]
Moreover,
\[
V_1:=\min\{U_1,U_2\}=U_1-w
\]
has boundary value \(M_1\), while
\[
V_2:=\max\{U_1,U_2\}=U_2+w
\]
has boundary value \(M_2\).  Both functions are at least \(a\), because
\(U_1,U_2\ge a\).  Hence
\[
V_1\in\mathcal K_{M_1},
\qquad
V_2\in\mathcal K_{M_2}.
\]

Use \(V_1\) in the variational inequality for \(U_1\).  Since
\(V_1-U_1=-w\), one obtains
\[
\int_D\langle AXU_1,Xw\rangle\,dy
+\int_D f(U_1)w\,dy\le0.
\]
Use \(V_2\) in the variational inequality for \(U_2\).  Since
\(V_2-U_2=w\), this gives
\[
\int_D\langle AXU_2,Xw\rangle\,dy
+\int_D f(U_2)w\,dy\ge0.
\]
Subtracting the second left-hand side from the first yields
\[
\int_D\langle A(XU_1-XU_2),Xw\rangle\,dy
+\int_D\bigl(f(U_1)-f(U_2)\bigr)w\,dy
\le0.
\]
By the truncation identity for \(Xw\), the first integral equals
\[
\int_{\{w>0\}}\langle AXw,Xw\rangle\,dy
\ge\lambda\int_D|Xw|^2\,dy.
\]
The second integral is supported on \(\{w>0\}=\{U_1>U_2\}\), and on
this set the monotonicity of \(f\) gives
\[
\bigl(f(U_1)-f(U_2)\bigr)w\ge0.
\]
Combining these observations proves the displayed inequality above.
In fact, both terms on its left are non-negative, so
\(\int_D|Xw|^2=0\).  Since \(w\in W_0^{1,2}(D,X)\), the horizontal
Poincar\'e inequality implies \(w=0\) a.e. in \(D\).  Thus
\(U_{M_1}\le U_{M_2}\) a.e.; for their continuous representatives the
inequality holds everywhere in \(D\).
\endgroup

Thus \(U_{M_1}\le U_{M_2}\). We now prove that the family of minimizers \(\{U_M:M>a\}\) is locally
bounded independently of \(M\).  If
\(K\Subset D\) and
\(\delta=\operatorname{dist}(K,\partial D)>0\), choose \(T_K>a\) so
large that
\[
CI(T_K)<\delta,
\]
which is possible by \eqref{eq:concise-92-I-zero}.  If \(x\in K\) and
\(U_M(x)>a\), then \eqref{eq:concise-92-distance} gives
\[
\delta\le CI(U_M(x)).
\]
Since \(I\) is decreasing, \(U_M(x)>T_K\) would imply
\(CI(U_M(x))\le CI(T_K)<\delta\), a contradiction.  Points at which
\(U_M(x)=a\) satisfy the same bound because \(a<T_K\).  Hence
\begin{equation}\label{eq:concise-92-local-Linfty}
\sup_{M>a}\|U_M\|_{L^\infty(K)}\le T_K.
\end{equation}

We next make the compactness argument explicit.  Let
\(\Omega_0\Subset\Omega_1\Subset D\) be relatively compact open sets,
and choose \(\eta\in C_c^1(\Omega_1)\) with \(0\le\eta\le1\) and
\(\eta=1\) on \(\Omega_0\).
Taking \(k=a\) in \eqref{eq:concise-92-upper-caccioppoli} gives
\[
\int_{\Omega_0}|XU_M|^2
\le
C\int_{\Omega_1}(U_M-a)^2|X\eta|^2
\le C_{\Omega_0,\Omega_1},
\]
where the last constant is independent of \(M\) by
\eqref{eq:concise-92-local-Linfty}, applied to
\(\operatorname{supp}\eta\Subset D\).  Thus \((U_M)\) is uniformly
bounded in \(W^{1,2}(\Omega_0,X)\).  The same local \(L^\infty\) bound
controls \(f(U_M)\), so the local obstacle estimates
\cite{GMobs,MarchiDG} give a common modulus of continuity on compact
subsets of \(D\).

We now construct the desired function \(W_{a,y_0}\) as the monotone
limit of the minimizers.  Let \(M_n\uparrow\infty\).  Monotonicity and
the local bound show that
the pointwise limit
\[
W_{a,y_0}:=\lim_{n\to\infty}U_{M_n}
=\sup_{M>a}U_M
\]
is finite in \(D\).  Arzel\`a--Ascoli, applied on each compact subset,
shows that every subsequence of \((U_{M_n})\) has a further subsequence
converging locally uniformly.  Its limit must equal the pointwise limit
\(W_{a,y_0}\).  It follows that the whole sequence converges locally
uniformly and that \(W_{a,y_0}\in C(D)\).

We now make the diagonal weak-compactness argument explicit.  Choose an
exhaustion
\[
\Omega_1\Subset\Omega_2\Subset\cdots\Subset D,
\qquad
\bigcup_{k=1}^{\infty}\Omega_k=D.
\]
For each fixed \(k\), the preceding estimate shows that
\((U_{M_n})\) is bounded in the Hilbert space
\(W^{1,2}(\Omega_k,X)\).  By weak compactness, first extract a subsequence
converging weakly on \(\Omega_1\); from that subsequence extract one
converging weakly on \(\Omega_2\), and continue inductively.  The diagonal
subsequence converges weakly in \(W^{1,2}(\Omega_k,X)\) for every \(k\).
On the other hand, the already established local uniform convergence
implies strong \(L^2(\Omega_k)\) convergence to \(W_{a,y_0}\).  Hence the
function component of each weak \(W^{1,2}\) limit is \(W_{a,y_0}\), so
\(W_{a,y_0}\in W^{1,2}(\Omega_k,X)\) and the weak limit is
\(W_{a,y_0}\) on every \(\Omega_k\).  Relabeling the diagonal subsequence,
we therefore have, for every relatively compact open set
\(\Omega_0\Subset D\),
\[
U_{M_n}\rightharpoonup W_{a,y_0}
\qquad\text{weakly in }W^{1,2}(\Omega_0,X).
\]
Indeed, \(\overline{\Omega_0}\subset\Omega_k\) for some \(k\), and the
convergence on \(\Omega_0\) follows by restriction.  Thus
\(W_{a,y_0}\in W_{\mathrm{loc}}^{1,2}(D,X)\), and
\eqref{eq:concise-92-center} gives \(W_{a,y_0}(y_0)=a\).

Finally, fix \(0\le\phi\in C_c^1(D)\).  The function
\(V=U_{M_n}+\phi\) belongs to \(\mathcal K_{M_n}\), so
\eqref{eq:concise-92-VI} gives
\[
\int_D\langle AXU_{M_n},X\phi\rangle
+\int_D f(U_{M_n})\phi\ge0.
\]
On \(\operatorname{supp}\phi\), the first integral converges by weak
\(W^{1,2}\) convergence.  The local uniform convergence, the local
bound, and the continuity of \(f\) give
\[
f(U_{M_n})\longrightarrow f(W_{a,y_0})
\qquad\text{uniformly on }\operatorname{supp}\phi.
\]
Passing to the limit yields
\[
-\int_D\langle AXW_{a,y_0},X\phi\rangle
\le
\int_D f(W_{a,y_0})\phi,
\]
so \(W_{a,y_0}\) is a weak supersolution.

Finally, let \(L>a\) and choose \(M>L\).  By
\eqref{eq:concise-92-boundary-data}, \(U_M=M\) on the compact boundary
of \(D\).  Since \(U_M\in C(\overline D)\) and \(M>L\), compactness of
\(\partial D\) yields a number \(\delta_L>0\) such that
\[
U_M(y)>L
\qquad\text{whenever }y\in D
\text{ and }\operatorname{dist}(y,\partial D)<\delta_L.
\]
The set
\[
\{y\in D:\operatorname{dist}(y,\partial D)<\delta_L\}
\]
is the boundary collar referred to here.  Since
\(W_{a,y_0}\ge U_M\), the same collar satisfies \(W_{a,y_0}>L\).
Therefore
\[
W_{a,y_0}(y)\longrightarrow\infty
\qquad\text{as }\operatorname{dist}(y,\partial D)\downarrow0.
\]
Taking
\[
\mathcal D_{a,y_0}:=D,
\qquad
R_f(a):=\Sigma\rho
=\frac{2C\Sigma}{d_*}
\int_a^\infty\frac{ds}{\sqrt{F(s)-F(a)}}
\]
proves the lemma after absorbing \(2C\Sigma/d_*\) into \(C_0\).
\end{proof}

\endgroup

\subsection{Pointwise and local Keller--Osserman estimates}

The barrier just constructed yields two estimates.  The first bounds a
subsolution at a point in terms of its distance from the boundary.  The
second combines that pointwise bound with the growth condition \({\bf(f-I)}\)
to control the nonlinear scale \(\sqrt{F(M_r/2)}/M_r\) on a ball.  Here
\(\Psi_f\) and its inverse \(\Phi_f\) are defined in
\eqref{eq:blue-horizontal-Psi-Phi-definition}, and \(M_r(w)\) denotes the
supremum of \(w\) on the concentric ball of radius \(r\).

\begin{thm}[Horizontal Keller--Osserman pointwise bound]\label{thm:blue-horizontal-KO-bound}
Assume that \(f:[0,\infty)\to[0,\infty)\) is continuous, non-decreasing,
\(f(t)>0\) for \(t>0\), and satisfies the Keller--Osserman condition
\[
\int_1^\infty\frac{ds}{\sqrt{F(s)}}<\infty.
\]
There is a structural constant
\(C_0=C_0(\lambda,\Lambda,\mathbb G,d)\) such that, if
\begin{equation}\label{eq:blue-horizontal-Psi-Phi-definition}
\begin{aligned}
\Psi_f(a)
&:=C_0\int_a^\infty\frac{ds}{\sqrt{\mathcal F(s,a)}},
&& a>0,\\
\Phi_f
&:=\Psi_f^{-1}
&& \text{on }(0,\Psi_f(0+)),\\
\Phi_f(\Psi_f(0+))
&:=0
&& \text{if }\Psi_f(0+)<\infty,
\end{aligned}
\end{equation}
then
any non-negative locally bounded weak subsolution \(w\) of
\eqref{eq:blue-horizontal-model-f} in \(\Omega\) satisfies
\begin{equation}\label{eq:blue-horizontal-global-bound}
w(y)\le \mathcal Q_f(d_\Omega(y)),
\qquad
\mathcal Q_f(r):=\Phi_f\bigl(\min\{r,\Psi_f(0+)\}\bigr),
\end{equation}
where \(d_\Omega(y):=\operatorname{dist}(y,\partial\Omega)\) is computed with the
fixed homogeneous distance.
\end{thm}

\begin{proof}
By Lemma \ref{blup} above, for every \(a>0\) and \(y_0\in\mathbb G\) there are
\(R_f(a)\le\Psi_f(a)\), a regular domain
\(\mathcal D_{a,y_0}\subset B(y_0,R_f(a))\), and a non-negative weak
supersolution \(W_{a,y_0}\) in \(\mathcal D_{a,y_0}\) such that
\(W_{a,y_0}(y_0)=a\) and
\(W_{a,y_0}(y)\to\infty\) as
\(\operatorname{dist}(y,\partial\mathcal D_{a,y_0})\downarrow0\).

Fix \(y_0\in\Omega\).  First suppose that
\(0<d_\Omega(y_0)<\Psi_f(0+)\), and choose
\(a>\Phi_f(d_\Omega(y_0))\).  Since \(\Psi_f\) is decreasing,
\[
R_f(a)\le \Psi_f(a)<d_\Omega(y_0),
\]
so \(\mathcal D_{a,y_0}\Subset\Omega\).  Let \(W=W_{a,y_0}\) be the
corresponding barrier and write \(\mathcal D=\mathcal D_{a,y_0}\).  Since
\(\overline{\mathcal D}\Subset\Omega\), local boundedness gives
\(\sup_{\mathcal D}w<\infty\).  The boundary blow-up of \(W\) therefore yields
\(\delta>0\) such that
\[
w\le W
\qquad\text{whenever }y\in\mathcal D
\text{ and }\operatorname{dist}(y,\partial\mathcal D)<2\delta.
\]
Let \(\mathcal D_\delta\) be the connected component containing \(y_0\) of
\[
\{y\in\mathcal D:\operatorname{dist}(y,\partial\mathcal D)>\delta\}.
\]
After decreasing \(\delta\), this component contains \(y_0\), and the preceding
collar inequality gives the weak boundary order \(w\le W\) on
\(\partial\mathcal D_\delta\).  Applying the weak comparison principle in
\(\mathcal D_\delta\) gives \(w\le W\) there.  In particular,
\[
w(y_0)\le W(y_0)=a.
\]
Letting \(a\downarrow\Phi_f(d_\Omega(y_0))\) yields
\[
w(y_0)\le \Phi_f(d_\Omega(y_0)).
\]

It remains to consider the case in which \(\Psi_f(0+)<\infty\) and
\(d_\Omega(y_0)\ge\Psi_f(0+)\).  By definition,
\[
\Psi_f(0+)=\lim_{a\downarrow0}\Psi_f(a).
\]
Since \(\Psi_f\) is strictly decreasing, every \(a>0\) satisfies
\[
\Psi_f(a)<\Psi_f(0+)\le d_\Omega(y_0).
\]
Thus the barrier corresponding to this arbitrary height \(a\) is contained
compactly in \(\Omega\), and the same comparison argument used above gives
\[
w(y_0)\le a.
\]
Because this holds for every \(a>0\) and \(w\ge0\), letting \(a\downarrow0\)
gives \(w(y_0)=0\).  The endpoint convention in
\eqref{eq:blue-horizontal-Psi-Phi-definition} gives
\[
\mathcal Q_f(d_\Omega(y_0))
=\Phi_f\bigl(\Psi_f(0+)\bigr)=0,
\]
because
\(\min\{d_\Omega(y_0),\Psi_f(0+)\}=\Psi_f(0+)\).  Thus the asserted estimate
also holds in this second distance regime.  Since \(y_0\) was arbitrary,
\eqref{eq:blue-horizontal-global-bound} follows.
\end{proof}

Under the additional growth assumption \({\bf(f-I)}\), the preceding
pointwise bound yields the following quantitative control of the nonlinear
scale on interior balls.

\begin{thm}[Horizontal local Keller--Osserman scale estimate]\label{thm:blue-horizontal-local-scale}
Assume in addition that \(f\) satisfies \({\bf(f-I)}\), that is,
\[
h(t):=\left(\frac{f(t)}{t}\right)^{1/2}
\]
is non-decreasing on \((0,\infty)\), and there are \(1<\varrho<\theta\) such
that \(h(\theta t)\ge\varrho h(t)\) for all \(t>0\).  If \(w\) is a
non-negative locally bounded weak subsolution of \eqref{eq:blue-horizontal-model-f}
in \(B(y_0,2r)\), and
\[
M_r(w):=\sup_{B(y_0,r)}w,
\]
then
\begin{equation}\label{eq:blue-carnot-subsolution-scale}
\frac{\sqrt{F(M_r(w)/2)}}{M_r(w)}
\le
\frac{C}{r},
\end{equation}
where \(C=C(\lambda,\Lambda,Q,\theta,\varrho)\).
\end{thm}

\begin{proof}
If \(M_r(w)=0\), there is nothing to prove.  We assume \(M_r(w)>0\).
We first record the one-dimensional consequences of \({\bf(f-I)}\).  Put
\[
\nu:=\frac{\log\varrho}{\log\theta}>0.
\]
If \(s\ge t>0\), choose \(k\ge0\) with \(\theta^k t\le s<\theta^{k+1}t\).  Then
\[
h(s)
\ge h(\theta^k t)
\ge \varrho^k h(t)
\ge c_{\theta,\varrho}\left(\frac{s}{t}\right)^\nu h(t).
\]
Consequently,
\begin{equation}\label{eq:blue-horizontal-tail-h}
\int_t^\infty\frac{ds}{s h(s)}
\le
\frac{C_{\theta,\varrho}}{h(t)},
\qquad t>0.
\end{equation}
Since \(s h(s)=\sqrt{s f(s)}\), this implies the Keller--Osserman condition.
The remaining estimates are purely one-dimensional.  Namely, if \(s\ge2t\),
then
\(F(s)\ge F(2t)\ge2F(t)\), because
\[
F(2t)-F(t)=\int_t^{2t}f(\tau)\,d\tau
\ge t f(t)
\ge F(t).
\]
Thus
\begin{equation}\label{eq:blue-horizontal-F-lower-far}
\mathcal F(s,t)=F(s)-F(t)
\ge \frac12F(s),
\qquad s\ge2t.
\end{equation}
For \(t<s<2t\), monotonicity of \(f\) gives
\begin{equation}\label{eq:blue-horizontal-F-lower-near}
\mathcal F(s,t)
\ge (s-t)f(t)
=t f(t)\left(\frac{s}{t}-1\right).
\end{equation}
We now split the defining integral for \(\Psi_f\).  From
\eqref{eq:blue-horizontal-F-lower-near},
\[
\int_t^{2t}\frac{ds}{\sqrt{\mathcal F(s,t)}}
\le
\frac{1}{\sqrt{t f(t)}}\int_t^{2t}\left(\frac{s}{t}-1\right)^{-1/2}ds
\le
\frac{C}{h(t)}.
\]
From \eqref{eq:blue-horizontal-F-lower-far}, the change of variables
\(s=2\sigma\), and the bound
\(F(2\sigma)\ge\int_\sigma^{2\sigma}f(\tau)\,d\tau\ge\sigma f(\sigma)\),
\[
\begin{aligned}
\int_{2t}^\infty\frac{ds}{\sqrt{\mathcal F(s,t)}}
&\le C\int_{2t}^\infty\frac{ds}{\sqrt{F(s)}}  \\
&\le C\int_t^\infty\frac{d\sigma}{\sqrt{\sigma f(\sigma)}}
= C\int_t^\infty\frac{d\sigma}{\sigma h(\sigma)}
\le \frac{C}{h(t)},
\end{aligned}
\]
where the last inequality is \eqref{eq:blue-horizontal-tail-h}.  Hence
\begin{equation}\label{eq:blue-horizontal-Psi-h}
\Psi_f(t)
\le
\frac{C}{h(t)},
\qquad t>0.
\end{equation}
Also \(\Psi_f(0+)=\infty\): indeed, since \(h\) is non-decreasing,
\(F(s)\le C s^2\) for \(0<s\le1\), and hence
\[
\Psi_f(t)
\ge c\int_t^1\frac{ds}{\sqrt{F(s)}}
\ge c\int_t^1\frac{ds}{s}\to\infty
\quad\text{as }t\downarrow0.
\]
Thus \(\Phi_f\) is defined on all of \((0,\infty)\).  Substituting
\(t=\Phi_f(\sigma)\) in \eqref{eq:blue-horizontal-Psi-h} gives
\[
h(\Phi_f(\sigma))\le\frac{C}{\sigma},
\qquad \sigma>0.
\]
Since \(F(t)\le t f(t)=t^2h(t)^2\), we also have
\begin{equation}\label{eq:blue-horizontal-Phi-F}
\frac{\sqrt{F(\Phi_f(\sigma))}}{\Phi_f(\sigma)}
\le
\frac{C}{\sigma},
\qquad \sigma>0.
\end{equation}

We now apply Theorem~\ref{thm:blue-horizontal-KO-bound} in the domain
\(\Omega=B(y_0,2r)\).  For \(y\in B(y_0,r)\),
\(d_\Omega(y)\ge r\), and therefore
\[
w(y)\le \Phi_f(d_\Omega(y)).
\]
Using monotonicity of \(f\), \eqref{eq:blue-horizontal-Psi-h}, and the monotonicity
of \(\Psi_f\), we obtain, for every \(y\in B(y_0,r)\) with \(w(y)>0\),
\[
\frac{\sqrt{F(w(y)/2)}}{w(y)}
\le C h(w(y))
\le \frac{C}{\Psi_f(w(y))}
\le \frac{C}{d_\Omega(y)}
\le \frac{C}{r}.
\]
Choose a sequence \(y_j\in B(y_0,r)\) such that \(w(y_j)\to M_r(w)\).  Passing
to the limit in the last display gives \eqref{eq:blue-carnot-subsolution-scale}.
\end{proof}

The scale estimate in Theorem~\ref{thm:blue-horizontal-local-scale}, together
with the classical Harnack inequality, implies the following global
Liouville property.
\begin{cor}[Liouville property for global solutions]
\label{cor:blue-liouville-global}
Assume that $f:[0,\infty)\to[0,\infty)$ is continuous and non-decreasing,
$f(0)=0$, $f(t)>0$ for every $t>0$, and ${\bf(f-I)}$ holds.  Then every
non-negative global solution $u$ of
\eqref{eq:blue-horizontal-model-f} to which the classical Harnack inequality
\eqref{chi} applies is identically zero.
\end{cor}

\begin{proof}
Fix $x\in\mathbb G$.  If $u(x)=0$, there is nothing to prove at this point.
Suppose instead that $u(x)>0$.  For each $R>0$, set
\[
M_R:=\sup_{B(x,R)}u.
\]
Theorem~\ref{thm:blue-horizontal-local-scale}, applied with center $x$ and
radius $R$, gives
\[
\frac{\sqrt{F(M_R/2)}}{M_R}\le\frac{C}{R}.
\]
The classical Harnack inequality \eqref{chi} (Theorem~\ref{gha-c}), applied
on $B(x,2R)$, yields
\[
M_R\le C_H\inf_{B(x,R)}u\le C_Hu(x),
\]
where $C_H$ is independent of $x$, $R$, and $u$.  Since $M_R\ge u(x)$ and
$F$ is non-decreasing, it follows that
\[
\sqrt{F(u(x)/2)}
\le \sqrt{F(M_R/2)}
\le \frac{CM_R}{R}
\le \frac{CC_Hu(x)}{R}.
\]
Letting $R\to\infty$ gives $F(u(x)/2)=0$.  But $f(t)>0$ for every $t>0$
implies $F(t)>0$ for every $t>0$, contradicting $u(x)>0$.  Thus $u(x)=0$;
since $x$ was arbitrary, $u\equiv0$ in $\mathbb G$.
\end{proof}

Applying Theorem~\ref{thm:blue-horizontal-local-scale} to the
solution \(u\) in \(B_2=B(e,2)\), with center \(e\) and radius \(1\), gives
\begin{equation}\label{eq:blue-M-scale}
\frac{\sqrt{F(M/2)}}{M}
\le
C,
\qquad
M:=\sup_{B_1}u.
\end{equation}
If \(M=0\), then \(u\equiv0\) on \(B_1\) and there is nothing to prove.  Hence
we assume \(M>0\).

\subsection{Control of the nonlinear integrand below the maximum}

We now specialize the local scale estimate to the solution \(u\).  The goal
of this short subsection is to convert the estimate at the top level
\(M=\sup_{B_1}u\)
into a uniform bound for \(\sqrt{F(t)}/t\) at every lower level
\(0<t\le M/4\).

Let
\[
m:=\inf_{B_1}u.
\]
We claim that there is a structural constant \(C\) such that
\begin{equation}\label{eq:blue-integrand-control}
\frac{\sqrt{F(t)}}{t}\le C,
\qquad 0<t\le \frac{M}{4}.
\end{equation}
Indeed, since \(F(t)\le t f(t)\) and \(h(t)=\sqrt{f(t)/t}\) is non-decreasing,
\[
\frac{\sqrt{F(t)}}{t}
\le h(t)
\le h(M/4).
\]
Moreover, by monotonicity of \(f\),
\[
F(M/2)-F(M/4)
=
\int_{M/4}^{M/2}f(s)\,ds
\ge
\frac{M}{4}f(M/4).
\]
Thus
\[
h(M/4)^2
=
\frac{f(M/4)}{M/4}
\le
\frac{16F(M/2)}{M^2}.
\]
Taking square roots and using \eqref{eq:blue-M-scale} proves
\eqref{eq:blue-integrand-control}.

\subsection{Comparison with the generalized Harnack integral}

The preceding control makes the integrand comparable from below with \(1/t\)
between the levels \(m/4\) and \(M/4\).  This produces
\(\log(M/m)\).  The generalized Harnack inequality gives the matching upper
bound and completes the proof.

From \eqref{eq:blue-integrand-control}, for \(0<t\le M/4\),
\[
\sqrt{F(t)}+t\le (C+1)t.
\]
Hence
\begin{equation}\label{eq:blue-log-lower}
\int_{m/4}^{M/4}\frac{dt}{\sqrt{F(t)}+t}
\ge
c\int_{m/4}^{M/4}\frac{dt}{t}
=
c\log\frac{M}{m},
\end{equation}
provided \(m>0\).  If \(m=0<M\), the left-hand side in
\eqref{eq:blue-log-lower} is infinite, while Theorem~\ref{gha} gives a finite
upper bound below; hence this case is impossible.  Thus either \(M=0\), already
handled, or \(m>0\).

On the other hand, we now give the details for the upper bound.  Since

\[
0<m\le M,
\qquad
\sqrt{F(t)}+t\ge t,
\]
the portion of the integral below the level \(m\) is controlled by the
logarithmic integral.  More precisely, if \(M/4\le m\), then
\([m/4,M/4]\subset[m/4,m]\), while if \(M/4>m\), then we split the interval at
\(m\).  In both cases,
\[
\int_{m/4}^{M/4}\frac{dt}{\sqrt{F(t)}+t}
\le
\int_{m/4}^{m}\frac{dt}{t}
+
\int_m^M\frac{dt}{\sqrt{F(t)}+t}.
\]
The first term is exactly \(\log 4\), and the second term is bounded by the
generalized Harnack estimate \eqref{di}, applied in \(B_2\) with
\(m=\inf_{B_1}u\) and \(M=\sup_{B_1}u\).  Therefore
\begin{equation}\label{eq:blue-log-upper}
\begin{aligned}
\int_{m/4}^{M/4}\frac{dt}{\sqrt{F(t)}+t}
&\le
\int_{m/4}^{m}\frac{dt}{t}
+
\int_m^M\frac{dt}{\sqrt{F(t)}+t}  \\
&\le
\log4+C,
\end{aligned}
\end{equation}
where the last inequality is precisely the generalized Harnack estimate
\eqref{di} from Theorem~\ref{gha}.  Combining \eqref{eq:blue-log-lower} and
\eqref{eq:blue-log-upper} gives
\[
\log\frac{M}{m}\le C.
\]
Exponentiating, we obtain
\[
\sup_{B_1}u=M\le C m=C\inf_{B_1}u.
\]
This is \eqref{chi}, and Theorem~\ref{gha-c} is proved.
}

\fi 

\section*{Author addresses and e-mail}
\noindent
\textbf{Cristian E. Guti\'errez}\\
Department of Mathematics, Temple University\\
Philadelphia, PA 19122, USA\\
E-mail: \texttt{cristian.gutierrez@temple.edu}

\medskip
\noindent
\textbf{Ahmed Mohammed}\\
Department of Mathematical Sciences, Ball State University\\
Muncie, IN 47306, USA\\
E-mail: \texttt{amohammed@bsu.edu}
\JulinEndDocument

\ifJulinPaperTwo\color{black}\fi
\section{Infinitesimal generator of dilations}\label{sec:infinitesimal-generator-dilations}
The infinitesimal generator $Z$ of the group dilations $\delta_\lambda$ is the vector field obtained by differentiating the dilation flow with respect to the parameter. More precisely, for a smooth function $\varphi$,
$$Z\varphi(g)=\left.\frac{d}{d\lambda}\varphi(\delta_\lambda g)\right|_{\lambda=1}.$$
Thus $Z(g)$ is the tangent vector at $g$ to the curve
$$\lambda\mapsto \delta_\lambda(g)$$ at $\lambda=1$.  In coordinates,
suppose 
$$g=(x_1^{(1)},x_2^{(2)},\cdots, x_r^{(r)})$$ are exponential coordinates adapted to the stratification, where $x_j^{(j)}$ corresponds to $V_j$, then
$$\delta_\lambda\left(x_1^{(1)},x_1^{(2)},\cdots, x_1^{(r)}\right)=\left(\lambda x_1^{(1)},\lambda^2x_1^{(2)},\cdots, \lambda ^rx_1^{(r)}\right),$$ and therefore
$$Z=\sum_{\alpha\in V_1}x_\alpha \partial_{x_\alpha}+2\sum_{\beta\in V_2}x_\beta \partial_{x_\beta}+\cdots+r\sum_{\gamma\in V_r}x_\gamma \partial_{x_\gamma}.$$ So $Z$ is the analogue of the Euclidean vector field
$$x\cdot\nabla.$$

\vspace{1cm}

\color{black}
\ifJulinPaperTwo\color{black}\fi
\section{Carnot-group and horizontal geometric measure theory background}
\label{sec:background-homogeneous-carnot}

We give a concise description of the notions needed for the horizontal
calculus on a Carnot group. Throughout, the underlying space is
\(\mathbb R^N\). Thus a point is written
\[
x=(x_1,\dots,x_N)\in \mathbb R^N.
\]

\subsection{Group law}

A \emph{group law} on \(\mathbb R^N\) is an operation
\[
(x,y)\mapsto x\circ y
\]
which assigns to every pair \(x,y\in\mathbb R^N\) another point
\(x\circ y\in\mathbb R^N\), and satisfies the following properties.

\begin{enumerate}
\item There is an identity element, denoted by \(e\), such that
\[
e\circ x=x\circ e=x
\]
for every \(x\in\mathbb R^N\).

\item Every \(x\in\mathbb R^N\) has an inverse, denoted by \(x^{-1}\), such
that
\[
x\circ x^{-1}=x^{-1}\circ x=e.
\]

\item The operation is associative:
\[
(x\circ y)\circ z=x\circ (y\circ z)
\]
for all \(x,y,z\in\mathbb R^N\).
\end{enumerate}

Usually, after a change of coordinates, one may take \(e=0\). In general,
the group law \(x\circ y\) is not the ordinary sum \(x+y\).

\subsection{Lie group structure}

We assume that the group law is smooth. This means that the maps
\[
(x,y)\mapsto x\circ y
\qquad\text{and}\qquad
x\mapsto x^{-1}
\]
are \(C^\infty\) maps on \(\mathbb R^N\). Then \(\mathbb R^N\), together
with the operation \(\circ\), is a \emph{Lie group}.

For our purposes, this simply means that the group operation is compatible
with differentiation.

\subsection{Left translations}

For each fixed \(x\in\mathbb R^N\), the \emph{left translation by \(x\)} is
the map
\[
L_x:\mathbb R^N\to\mathbb R^N
\]
defined by
\[
L_x(y)=x\circ y.
\]

Thus left translation means multiplication on the left by the element \(x\).

\subsection{Tangent vectors and vector fields}

A tangent vector at a point \(x\in\mathbb R^N\) may be identified with an
ordinary vector in \(\mathbb R^N\). Equivalently, a tangent vector is a
direction in which one can differentiate functions.

A \emph{vector field} \(X\) on \(\mathbb R^N\) assigns to each point
\(x\in\mathbb R^N\) a tangent vector \(X(x)\in\mathbb R^N\). In coordinates,
it has the form
\[
X(x)=\sum_{k=1}^N a_k(x)\frac{\partial}{\partial x_k}.
\]

Thus, if \(u\) is a smooth function, then
\[
Xu(x)
=
\sum_{k=1}^N a_k(x)\frac{\partial u}{\partial x_k}(x).
\]

\subsection{Left-invariant vector fields}

For \(x,y\in\mathbb G\), let
\[
(dL_x)_y:T_y\mathbb G\longrightarrow T_{x\circ y}\mathbb G
\]
denote the differential of the left translation \(L_x(z)=x\circ z\), evaluated
at \(y\).

A vector field \(X\) is called \emph{left invariant} if it is unchanged by
left translations; in symbols,
\[
(dL_x)_yX(y)=X(L_x(y))=X(x\circ y),
\qquad x,y\in\mathbb G.
\]
In particular, taking \(y=e\) gives
\[
X(x)=(dL_x)_eX(e).
\]

Equivalently, \(X\) is completely determined by its value at the identity.

If \(v\in\mathbb R^N\) is a tangent vector at the identity, the associated
left-invariant vector field \(X_v\) is defined by
\[
X_v(x)=dL_x(e)v.
\]
Here \(dL_x(e)\) denotes the derivative of the map \(L_x\) at the identity.

Equivalently, if \(\gamma(t)\) is a smooth curve such that
\[
\gamma(0)=e,
\qquad
\gamma'(0)=v,
\]
then
\[
X_v u(x)
=
\left.\frac{d}{dt}\right|_{t=0}u(x\circ \gamma(t)).
\]

Thus \(X_vu(x)\) differentiates \(u\) along the curve obtained by starting
at \(x\) and multiplying on the right by a curve through the identity.

In coordinates, every left-invariant vector field can be written as
\[
X_v=\sum_{k=1}^N a_k(x)\frac{\partial}{\partial x_k}.
\]

The coefficients \(a_k(x)\) are determined by the group law.

\subsection{The Lie algebra}

The collection of all left-invariant vector fields forms a vector space
called the \emph{Lie algebra} of the group. It is denoted by
\[
\mathfrak g.
\]

The Lie algebra has dimension \(N\). Choosing a basis
\[
v_1,\dots,v_N
\]
of tangent vectors at the identity gives a basis of left-invariant vector
fields
\[
X_1,\dots,X_N.
\]

The Lie algebra is equipped with the \emph{commutator} of vector fields.
For two vector fields \(X\) and \(Y\), their commutator is
\[
[X,Y]=XY-YX.
\]
This means that, when applied to a smooth function \(u\),
\[
[X,Y]u=X(Yu)-Y(Xu).
\]

The commutator measures the failure of \(X\) and \(Y\) to commute.

\subsection{Nilpotent Lie algebra}

Starting from the Lie algebra \(\mathfrak g\), define
\[
\mathfrak g_1=\mathfrak g,
\]
and recursively
\[
\mathfrak g_{k+1}=[\mathfrak g,\mathfrak g_k],
\]
where \([\mathfrak g,\mathfrak g_k]\) denotes the vector space generated by
all commutators
\[
[X,Y],
\qquad
X\in\mathfrak g,\quad Y\in\mathfrak g_k.
\]

The Lie algebra is called \emph{nilpotent} if this sequence eventually
becomes zero. That is, there exists an integer \(s\) such that
\[
\mathfrak g_{s+1}=\{0\}.
\]

The integer \(s\) is called the \emph{step} of the group.

\subsection{Dilations}

A \emph{dilation} is a map
\[
\delta_\lambda:\mathbb R^N\to\mathbb R^N,
\qquad \lambda>0,
\]
of the form
\[
\delta_\lambda(x_1,\dots,x_N)
=
(\lambda^{\alpha_1}x_1,\dots,\lambda^{\alpha_N}x_N),
\]
where
\[
1\leq \alpha_1\leq \alpha_2\leq \cdots \leq \alpha_N
\]
are positive numbers, usually integers. These numbers are called the
\emph{weights} of the dilation.

The family \(\{\delta_\lambda\}_{\lambda>0}\) satisfies
\[
\delta_\lambda\circ\delta_\mu=\delta_{\lambda\mu},
\qquad
\delta_1=\mathrm{Id}.
\]

A dilation is compatible with the group law if it is a group automorphism,
that is,
\[
\delta_\lambda(x\circ y)
=
\delta_\lambda(x)\circ \delta_\lambda(y)
\]
for every \(x,y\in\mathbb R^N\) and every \(\lambda>0\).

\subsection{Homogeneous groups}

A \emph{homogeneous group} is a Lie group on \(\mathbb R^N\) satisfying the
following conditions.

\begin{enumerate}
\item The group is nilpotent.

\item There is a family of dilations \(\{\delta_\lambda\}_{\lambda>0}\).

\item Each dilation \(\delta_\lambda\) is a group automorphism:
\[
\delta_\lambda(x\circ y)
=
\delta_\lambda(x)\circ\delta_\lambda(y).
\]
\end{enumerate}

The \emph{homogeneous dimension} of the group is
\[
Q=\alpha_1+\cdots+\alpha_N.
\]

This number appears because the Jacobian determinant of \(\delta_\lambda\)
is
\[
\det D\delta_\lambda=\lambda^Q.
\]
Consequently, for every measurable set \(E\subset\mathbb R^N\),
\[
|\delta_\lambda(E)|=\lambda^Q |E|,
\]
where \(|E|\) denotes Lebesgue measure.

\subsection{Homogeneous vector fields}

A vector field \(X\) is said to be homogeneous of degree \(\alpha\) if
\[
X(u\circ\delta_\lambda)
=
\lambda^\alpha (Xu)\circ\delta_\lambda
\]
for every smooth function \(u\) and every \(\lambda>0\).

For example, in the Euclidean case,
\[
\delta_\lambda(x)=\lambda x,
\]
and
\[
\frac{\partial}{\partial x_j}
\]
is homogeneous of degree \(1\), because
\[
\frac{\partial}{\partial x_j}(u(\lambda x))
=
\lambda \frac{\partial u}{\partial x_j}(\lambda x).
\]

In a homogeneous group, one can choose a basis of left-invariant vector
fields
\[
X_1,\dots,X_N
\]
adapted to the dilations, so that each \(X_j\) is homogeneous of degree
\(\alpha_j\).

\subsection{Carnot groups}

A \emph{Carnot group} is a special homogeneous group. Its Lie algebra
\(\mathfrak g\) can be decomposed as
\[
\mathfrak g
=
V_1\oplus V_2\oplus\cdots\oplus V_s,
\]
where each \(V_k\) is a vector subspace, and
\[
[V_1,V_k]=V_{k+1},
\qquad k=1,\dots,s-1,
\]
while
\[
[V_1,V_s]=\{0\}.
\]

This decomposition is called a \emph{stratification} of the Lie algebra.

The first layer
\[
V_1
\]
is called the \emph{horizontal layer}. It generates all the other layers by
taking commutators.

If
\[
m=\dim V_1,
\]
we choose a basis of left-invariant vector fields
\[
X_1,\dots,X_m
\]
spanning \(V_1\). These are called the \emph{horizontal vector fields}.

The remaining directions in the group are obtained by taking commutators
of the horizontal vector fields.

\subsection{Horizontal gradient}

Let \(u\) be a smooth function on an open set \(\Omega\subset\mathbb R^N\).
The \emph{horizontal gradient} of \(u\) is
\[
Xu=(X_1u,\dots,X_m u).
\]

This is an \(\mathbb R^m\)-valued object. It records only the derivatives
of \(u\) in the horizontal directions.

It is not, in general, the full Euclidean gradient
\[
\nabla u
=
\left(
\frac{\partial u}{\partial x_1},
\dots,
\frac{\partial u}{\partial x_N}
\right).
\]

\subsection{Horizontal divergence}

Let
\[
\varphi=(\varphi_1,\dots,\varphi_m)
\]
be an \(\mathbb R^m\)-valued function. The \emph{horizontal divergence} of
\(\varphi\) is
\[
\operatorname{div}_H\varphi
=
\sum_{j=1}^m X_j\varphi_j.
\]

Thus \(\operatorname{div}_H\) is the divergence taken only with respect to
the horizontal vector fields.

\subsection{Formal adjoints}\label{subsec:app-formal-adjoints}

The formal adjoint of \(X_j\), denoted by \(X_j^*\), is defined by the
identity
\[
\int_{\mathbb R^N} X_j u\, v\,dx
=
\int_{\mathbb R^N} u\, X_j^*v\,dx
\]
for smooth compactly supported functions \(u\) and \(v\).

For completeness, we give the proof that the left-invariant horizontal
fields are divergence-free.  Let \(Y\) be any left-invariant vector field and
write \(V=Y(e)\in\mathfrak g\).  The flow of \(Y\) is
\[
\Phi_t(x)=x\circ\exp(tV),
\]
because the curve \(t\mapsto x\circ\exp(tV)\) has velocity
\(dL_{x\circ\exp(tV)}V=Y(x\circ\exp(tV))\).  Thus the flow is right translation
by the one-parameter subgroup \(\exp(tV)\).

A homogeneous group, and in particular a Carnot group, is nilpotent and hence
unimodular.  Here \(\operatorname{ad}_V:\mathfrak g\to\mathfrak g\) denotes the
linear map defined by
\(\operatorname{ad}_V(W)=[V,W]\).  Equivalently,
\(\operatorname{tr}(\operatorname{ad}_V)=0\) for every
\(V\in\mathfrak g\), since \(\operatorname{ad}_V\) is nilpotent.  Therefore right
translations preserve Haar measure.  In exponential coordinates Haar measure
is Lebesgue measure, so \(\Phi_t\) preserves \(dx\).  Hence, for every
\(\eta\in C_c^1(\mathbb R^N)\),
\[
\int_{\mathbb R^N}Y\eta\,dx
=
\left.\frac{d}{dt}\right|_{t=0}
\int_{\mathbb R^N}\eta(\Phi_t(x))\,dx
=
\left.\frac{d}{dt}\right|_{t=0}
\int_{\mathbb R^N}\eta(x)\,dx
=0.
\]
Writing \(Y=\sum_{\alpha=1}^N b_\alpha\partial_{x_\alpha}\), the last identity
means
\[
0=\int_{\mathbb R^N}Y\eta\,dx
=-\int_{\mathbb R^N}\eta\,
\sum_{\alpha=1}^N\frac{\partial b_\alpha}{\partial x_\alpha}\,dx
\qquad\forall\,\eta\in C_c^1(\mathbb R^N),
\]
and therefore \(\operatorname{div}Y=0\) in the Euclidean-coordinate sense.  In
particular, each horizontal field \(X_j\) is divergence-free.  Therefore
\[
X_j^*=-X_j.
\]

Hence
\[
\int_{\mathbb R^N} X_j u\,v\,dx
=
-\int_{\mathbb R^N} u\,X_jv\,dx.
\]

This is the basic integration-by-parts formula for the horizontal vector
fields.

\subsection{Summary}

The structure needed for subelliptic PDEs is the following.

\begin{itemize}
\item The space is \(\mathbb R^N\), but with a possibly noncommutative group
law \(x\circ y\).

\item The group law determines left-invariant vector fields
\[
X_1,\dots,X_N.
\]

\item A homogeneous group has dilations
\[
\delta_\lambda(x_1,\dots,x_N)
=
(\lambda^{\alpha_1}x_1,\dots,\lambda^{\alpha_N}x_N)
\]
compatible with the group law.

\item A Carnot group has distinguished first-layer vector fields
\[
X_1,\dots,X_m,
\]
called horizontal vector fields.

\item These define the horizontal gradient
\[
Xu=(X_1u,\dots,X_m u)
\]
and the horizontal divergence
\[
\operatorname{div}_H\varphi=\sum_{j=1}^m X_j\varphi_j.
\]

\item Since the \(X_j\)'s are left invariant and divergence-free by the argument
above,
\[
X_j^*=-X_j.
\]
\end{itemize}
\subsection{Horizontal BV and horizontal perimeter}\label{subsec:app-horizontal-bv}

This subsection fixes the notation used in the level-surface estimates of
the companion smooth paper \cite{GutierrezMohammedSmooth}.  The
main point is that the horizontal derivative of a characteristic function is a
vector-valued Radon measure, and its total variation is the horizontal
perimeter measure.

\begin{defn}[Horizontal BV functions]\label{def:app-horizontal-bv}\itshape
Let \(\Omega\subset\mathbb G\) be open.  A function
\(u\in L^1_{\mathrm{loc}}(\Omega)\) belongs to \(BV_H(\Omega)\) if there exists
an \(\mathbb R^m\)-valued Radon measure on \(\Omega\),
\[
D_Hu=(D_{X_1}u,\ldots,D_{X_m}u)
\]
such that, for every \(\Phi=(\Phi_1,\ldots,
\Phi_m)\in C_c^1(\Omega;\mathbb R^m)\),
\begin{equation}\label{eq:app-bvh-distributional-definition}
\int_\Omega u\,\operatorname{div}_H\Phi\,dx
=
-\int_\Omega \langle \Phi,dD_Hu\rangle .
\end{equation}
Such a measure, if it exists, is uniquely determined by \(u\): any two
measures satisfying \eqref{eq:app-bvh-distributional-definition} agree as
distributions, hence agree as Radon measures.
Here
\[
\operatorname{div}_H\Phi:=\sum_{j=1}^m X_j\Phi_j.
\]
The total variation of this measure is denoted by \(|D_Hu|\). Equivalently,
\begin{equation}\label{eq:app-bvh-total-variation}
|D_Hu|(\Omega)
=
\sup\left\{
\left.
\int_\Omega u\,\operatorname{div}_H\Phi\,dx
\,\right|\,
\Phi\in C_c^1(\Omega;\mathbb R^m),\ |\Phi|\le 1
\right\}.
\end{equation}
\end{defn}

The following comments clarify Definition~\ref{def:app-horizontal-bv} and fix
the conventions used below.

The sign convention in \eqref{eq:app-bvh-distributional-definition} is the one
used throughout the paper.  If \(u\in C^1(\Omega)\), then integration by parts
using the divergence-free property proved in
Subsection~\ref{subsec:app-formal-adjoints} gives
\begin{equation}\label{eq:app-smooth-bv-gradient}
D_Hu=Xu\,dx,
\qquad
|D_Hu|=|Xu|\,dx.
\end{equation}
Indeed,
\[
\int_\Omega u\,X_j\Phi_j\,dx
=
-\int_\Omega (X_ju)\Phi_j\,dx,
\]
and summing in \(j\) gives \eqref{eq:app-smooth-bv-gradient}.

If \(\mu=(\mu_1,\ldots,\mu_m)\) is an \(\mathbb R^m\)-valued Radon measure and
\(\Phi\in C_c^1(\Omega;\mathbb R^m)\), we write
\begin{equation}\label{eq:app-vector-measure-pairing}
\int_\Omega \langle \Phi,d\mu\rangle
:=
\sum_{j=1}^m\int_\Omega \Phi_j\,d\mu_j.
\end{equation}
Thus \eqref{eq:app-bvh-distributional-definition} means that
\(D_Hu\) is the horizontal distributional gradient of \(u\).

\begin{defn}[Horizontal perimeter]\label{def:app-horizontal-perimeter}\itshape
A measurable set \(E\subset\Omega\) has finite horizontal perimeter in
\(\Omega\) if \(\chi_E\in BV_H(\Omega)\).  Its horizontal perimeter measure is
\begin{equation}\label{Pmes}
P_H(E;\cdot):=|D_H\chi_E|.
\end{equation}
We also write \(P_H(E)\) for this measure.  Thus, for a Borel set
\(A\subset\Omega\),
\[
P_H(E;A)=|D_H\chi_E|(A),
\]
and, by the regularity of Radon measures,
\[
P_H(E;A)=\inf\{P_H(E;U):A\subset U\Subset\Omega,\ U\text{ open}\}
\]
whenever the right-hand side is taken over relatively compact open
neighborhoods of \(A\).  We say that \(E\) has locally finite horizontal
perimeter in \(\Omega\) if \(P_H(E;U)<\infty\) for every \(U\Subset\Omega\).
\end{defn}

\begin{prop}[Relative isoperimetric inequality in Carnot balls]\label{prop:app-relative-isoperimetric-carnot-balls}
There exists a constant \(c_{\mathrm{iso}}>0\), depending only on the Carnot
group and on the fixed homogeneous distance, such that the following holds.  If
\(B(x,r)\Subset\Omega\) and \(E\subset\Omega\) has finite horizontal perimeter in
\(\Omega\), then
\begin{equation}\label{eq:app-relative-isoperimetric-carnot-balls}
P_H(E;B(x,r))
\ge
c_{\mathrm{iso}}
\min\left\{|E\cap B(x,r)|,\ |B(x,r)\setminus E|\right\}^{\frac{Q-1}{Q}}.
\end{equation}
In particular, if for some \(\theta\in(0,1/2]\),
\[
|E\cap B(x,r)|\ge \theta |B(x,r)|,
\qquad
|B(x,r)\setminus E|\ge \theta |B(x,r)|,
\]
then
\begin{equation}\label{eq:app-relative-isoperimetric-carnot-balls-density}
P_H(E;B(x,r))\ge c(\theta)r^{Q-1}.
\end{equation}
\end{prop}

\begin{proof}
The relative isoperimetric inequality \eqref{eq:app-relative-isoperimetric-carnot-balls}
is the Carnot--Carath\'eodory isoperimetric inequality of
\cite[Theorem~1.18]{GAN}, written with the horizontal perimeter notation of
Definition~\ref{def:app-horizontal-perimeter}.  The constant is uniform on balls
because the distance is left invariant and homogeneous.  The second assertion
follows immediately from \eqref{eq:app-relative-isoperimetric-carnot-balls} and
\(|B(x,r)|=|B(e,1)|r^Q\).
\end{proof}

For open \(U\Subset\Omega\), Definition~\ref{def:app-horizontal-bv} gives the
useful variational formula
\begin{equation}\label{eq:app-perimeter-variational-formula}
P_H(E;U)=
\sup\left\{
\left.
\int_E \operatorname{div}_H\Phi\,dx
\,\right|\,
\Phi\in C_c^1(U;\mathbb R^m),\ |\Phi|\le 1
\right\}.
\end{equation}
This is the formula used when passing from volume integrals over \(E\) to
surface integrals over its horizontal reduced boundary.

In the simple case where \(E\) has \(C^1\) boundary, \(P_H(E)\) is just
Euclidean surface measure weighted by the length of the horizontal projection
of the outer normal.  More precisely, if \(\nu\) is the Euclidean outer unit
normal to \(\partial E\), set
\[
\nu_H^{\rm raw}(x):=
\big(\langle \nu(x),X_1(x)\rangle,\ldots,
\langle \nu(x),X_m(x)\rangle\big).
\]
Here \(N\) denotes the Euclidean, or topological, dimension of the Carnot
group.  Under the standard identification of the underlying manifold with
\(\mathbb R^N\), points are written as
\(x=(x_1,\ldots,x_N)\in\mathbb R^N\); this \(N\) should not be confused with
the homogeneous dimension \(Q\).  Then, for every Borel set \(A\subset\Omega\),
\begin{equation}\label{eq:app-smooth-horizontal-perimeter}
P_H(E;A)=
\int_{A\cap\partial E}|\nu_H^{\rm raw}(x)|\,d\mathcal H^{N-1}(x).
\end{equation}
Equivalently, if locally \(E=\{g<0\}\) with \(g\in C^1\) and
\(\nabla g\ne0\) on \(\partial E\), then
\(|\nu_H^{\rm raw}|=|Xg|/|\nabla g|\) on \(\partial E\).  Thus characteristic
points, where the Euclidean normal has zero horizontal projection, contribute
zero to \(P_H(E)\).  At non-characteristic points the outer horizontal normal
is \(\nu_E=\nu_H^{\rm raw}/|\nu_H^{\rm raw}|\).

\subsection{Reduced boundary and the Gauss--Green formula}\label{subsec:app-gauss-green-theorem}

Let \(E\subset\Omega\) have finite horizontal perimeter.  The polar
decomposition of the vector measure \(D_H\chi_E\) is
\begin{equation}\label{eq:app-polar-decomposition}
D_H\chi_E=-\nu_E\,P_H(E;\cdot),
\end{equation}
where \(\nu_E:\Omega\to\mathbb R^m\) is \(P_H(E)\)-measurable and
\begin{equation}\label{eq:app-normal-unit-length}
|\nu_E|=1
\qquad P_H(E)\text{-a.e.}
\end{equation}
This vector field is unique \(P_H(E)\)-a.e.  Indeed,
\(P_H(E;\cdot)=|D_H\chi_E|\), and the Radon--Nikodym derivative
\(dD_H\chi_E/d|D_H\chi_E|\) in the polar decomposition of the vector measure
\(D_H\chi_E\) is unique \(|D_H\chi_E|\)-a.e.; with the sign convention in
\eqref{eq:app-polar-decomposition},
\(\nu_E=-dD_H\chi_E/dP_H(E)\).
The minus sign in \eqref{eq:app-polar-decomposition} is chosen so that
\(\nu_E\) is the outer horizontal normal.  With this convention, taking
\(u=\chi_E\) in \eqref{eq:app-bvh-distributional-definition} gives
\begin{equation}\label{eq:app-divergence-polar-identity}
\int_E\operatorname{div}_H\Phi\,dx
=
\int_\Omega\langle \Phi,\nu_E\rangle\,dP_H(E)
\end{equation}
for every \(\Phi\in C_c^1(\Omega;\mathbb R^m)\).

\begin{defn}[Horizontal reduced boundary]\label{def:app-horizontal-reduced-boundary}\itshape
The horizontal reduced boundary of \(E\), denoted by \(\partial_H^*E\), is the
set of points \(x\in\Omega\) at which the following two conditions hold.
First, the horizontal perimeter measure is nontrivial at every scale:
\begin{equation}\label{eq:app-reduced-boundary-positive-density}
P_H(E;B(x,\rho))>0
\qquad\text{for every }\rho>0\text{ such that }B(x,\rho)\subset\Omega.
\end{equation}
Second, the averaged outer horizontal normal has a limiting unit vector as
\(\rho\downarrow0\) through such radii:
\begin{equation}\label{eq:app-reduced-boundary-normal-limit}
\lim_{\substack{\rho\downarrow0\\ B(x,\rho)\subset\Omega}}
\frac{1}{P_H(E;B(x,\rho))}
\int_{B(x,\rho)}\nu_E\,dP_H(E)
\end{equation}
exists in \(\mathbb R^m\) and has norm \(1\).
\end{defn}

The following examples may help fix the meaning of \(\partial_H^*E\).

\medskip
\noindent\textbf{Example 1 (Euclidean model).}
If \(\mathbb G=\mathbb R^n\) and the horizontal fields are the standard
coordinate fields \(X_j=\partial_{x_j}\), \(j=1,\ldots,n\), then
\(D_H\chi_E=D\chi_E\) and \(P_H(E)\) is the classical perimeter measure.  For
instance, if \(E=\{x_n<0\}\), then
\[
D\,\chi_E=-e_n\,\mathcal H^{n-1}\!\lfloor_{\{x_n=0\}},
\qquad
P_H(E;A)=\mathcal H^{n-1}(A\cap\{x_n=0\}),
\]
\[
\partial_H^*E=\partial^*E=\{x_n=0\}\cap\Omega,
\qquad
\nu_E=e_n.
\]
More generally, if \(E\subset\mathbb R^n\) has \(C^1\) boundary, then
\(P_H(E;A)=\mathcal H^{n-1}(A\cap\partial E)\),
\(\partial_H^*E=\partial E\cap\Omega\), and \(\nu_E\) is the usual Euclidean
outer unit normal.  Thus the Gauss--Green formula above reduces exactly to the
classical divergence theorem.

\medskip
\noindent\textbf{Example 2 (A non-characteristic vertical half-space in
\(\mathbb H^1\)).}
Consider the first Heisenberg group \(\mathbb H^1\), with coordinates
\((x,y,\tau)\) and standard horizontal fields
\[
X=\partial_x-\frac y2\partial_\tau,
\qquad
Y=\partial_y+\frac x2\partial_\tau .
\]
Let \(E=\{x<0\}\).  On \(\partial E=\{x=0\}\), the Euclidean outer normal is
\((1,0,0)\), hence the horizontal normal vector is
\((\langle X,(1,0,0)\rangle,\langle Y,(1,0,0)\rangle)=(1,0)\).  Therefore
\[
P_H(E;A)=\mathcal H^2(A\cap\{x=0\}),
\qquad
\partial_H^*E=\partial E\cap\Omega,
\]
and \(\nu_E=(1,0)\) on the reduced boundary.

\medskip
\noindent\textbf{Example 3 (A horizontal plane with one characteristic point in
\(\mathbb H^1\)).}
With the same horizontal fields, let \(E=\{\tau<0\}\).  On
\(\partial E=\{\tau=0\}\), the Euclidean outer normal is \((0,0,1)\), so the
horizontal normal vector is
\[
\bigl(\langle X,(0,0,1)\rangle,\langle Y,(0,0,1)\rangle\bigr)
=\left(-\frac y2,\frac x2\right).
\]
Thus
\[
P_H(E;A)=
\int_{A\cap\{\tau=0\}}\frac{\sqrt{x^2+y^2}}{2}\,d\mathcal H^2.
\]
The point \((0,0,0)\) is characteristic: the horizontal projection of the
Euclidean normal vanishes there.  It contributes no \(P_H(E)\)-mass and is not
part of the horizontal reduced boundary; away from that point,
\[
\partial_H^*E=(\{\tau=0\}\setminus\{(0,0,0)\})\cap\Omega,
\qquad
\nu_E=\frac{(-y,x)}{\sqrt{x^2+y^2}}.
\]

We record the following basic consequence of the definition of the horizontal
reduced boundary, since it is the reason all perimeter integrals may be written
over \(\partial_H^*E\).

\begin{prop}[Concentration on the horizontal reduced boundary]\label{prop:app-reduced-boundary-concentration}
Let \(E\subset\Omega\) have finite horizontal perimeter in \(\Omega\).  Then the
perimeter measure is concentrated on the reduced boundary:
\begin{equation}\label{eq:app-perimeter-concentrated-reduced-boundary}
P_H(E;\Omega\setminus\partial_H^*E)=0.
\end{equation}
Consequently,
\begin{equation}\label{eq:app-integral-reduced-boundary-only}
\int_\Omega \psi\,dP_H(E)
=
\int_{\partial_H^*E}\psi\,dP_H(E)
\end{equation}
for every \(P_H(E)\)-integrable function \(\psi\).
\end{prop}

\begin{proof}
Let \(\mu:=P_H(E;\cdot)\).  The polar decomposition
\eqref{eq:app-polar-decomposition} gives a \(\mu\)-measurable representative of
\(\nu_E\), and \eqref{eq:app-normal-unit-length} gives
\(|\nu_E|=1\) for \(\mu\)-a.e. point.

First, since \(\mu\) is a Radon measure on \(\Omega\), it is concentrated on its
support, where the support is understood relative to \(\Omega\).  Hence, for
\(\mu\)-a.e. \(x\in\Omega\), every relative neighborhood of \(x\) has positive
\(\mu\)-measure.  In particular,
\[
P_H(E;B(x,\rho))=\mu(B(x,\rho))>0
\]
for every \(\rho>0\) such that \(B(x,\rho)\subset\Omega\).

Second, by the Lebesgue differentiation theorem for Radon measures on the
Carnot group, for \(\mu\)-a.e. \(x\in\Omega\),
\[
\lim_{\substack{\rho\downarrow0\\ B(x,\rho)\subset\Omega}}
\frac{1}{\mu(B(x,\rho))}
\int_{B(x,\rho)}|\nu_E(y)-\nu_E(x)|\,d\mu(y)=0.
\]
For such a point \(x\), the triangle inequality gives
\[
\left|
\frac{1}{\mu(B(x,\rho))}\int_{B(x,\rho)}\nu_E\,d\mu
-
\nu_E(x)
\right|
\le
\frac{1}{\mu(B(x,\rho))}
\int_{B(x,\rho)}|\nu_E(y)-\nu_E(x)|\,d\mu(y),
\]
and therefore
\[
\lim_{\substack{\rho\downarrow0\\ B(x,\rho)\subset\Omega}}
\frac{1}{P_H(E;B(x,\rho))}
\int_{B(x,\rho)}\nu_E\,dP_H(E)
=
\nu_E(x).
\]
Because \(|\nu_E(x)|=1\) at these points, this limit has norm \(1\).  Thus
\(\mu\)-a.e. point of \(\Omega\) satisfies both defining conditions of
\(\partial_H^*E\), and so
\[
\mu(\Omega\setminus\partial_H^*E)=0,
\]
which is \eqref{eq:app-perimeter-concentrated-reduced-boundary}.

It remains to derive \eqref{eq:app-integral-reduced-boundary-only}.  Since
\(\mu(\Omega\setminus\partial_H^*E)=0\), the restrictions of \(\mu\) to
\(\Omega\) and to \(\partial_H^*E\) agree.  If \(\psi\ge0\) is measurable, then
\[
\int_\Omega \psi\,d\mu
=
\int_{\partial_H^*E}\psi\,d\mu
+
\int_{\Omega\setminus\partial_H^*E}\psi\,d\mu
=
\int_{\partial_H^*E}\psi\,d\mu,
\]
because the second term is an integral over a \(\mu\)-null set.  Applying this
identity to the positive and negative parts of an integrable real-valued
\(\psi\), or to the real and imaginary parts if \(\psi\) is complex-valued,
gives \eqref{eq:app-integral-reduced-boundary-only} for every
\(P_H(E)\)-integrable \(\psi\).
\end{proof}

\begin{prop}[Gauss--Green for finite horizontal perimeter sets]\label{prop:app-gauss-green-finite-perimeter}
Let \(E\subset\Omega\) have finite horizontal perimeter in \(\Omega\).  Then,
for every \(\Phi\in C_c^1(\Omega;\mathbb R^m)\),
\begin{equation}\label{eq:app-gauss-green-theorem}
\int_E\operatorname{div}_H\Phi\,dx
=
\int_\Omega\langle\Phi,\nu_E\rangle\,dP_H(E)
=
\int_{\partial_H^*E}\langle\Phi,\nu_E\rangle\,dP_H(E).
\end{equation}
\end{prop}

\begin{proof}
We first prove the first identity in \eqref{eq:app-gauss-green-theorem}.  Since
\(E\) has finite horizontal perimeter in \(\Omega\), Definition~\ref{def:app-horizontal-perimeter}
means that \(\chi_E\in BV_H(\Omega)\) and
\(P_H(E;\cdot)=|D_H\chi_E|\).  Applying the distributional definition
\eqref{eq:app-bvh-distributional-definition} with \(u=\chi_E\) gives
\[
\int_\Omega \chi_E\operatorname{div}_H\Phi\,dx
=
-\int_\Omega\langle\Phi,dD_H\chi_E\rangle .
\]
The left-hand side is exactly \(\int_E\operatorname{div}_H\Phi\,dx\), because
\(\chi_E\) is the characteristic function of \(E\) and
\(\operatorname{div}_H\Phi\) is integrable with compact support in \(\Omega\).
For the right-hand side, use the polar decomposition
\eqref{eq:app-polar-decomposition}, namely
\(D_H\chi_E=-\nu_E P_H(E;\cdot)\).  Writing this componentwise,
\(D_{X_j}\chi_E=-(\nu_E)_jP_H(E;\cdot)\), and using the vector-measure pairing
\eqref{eq:app-vector-measure-pairing}, we obtain
\[
-\int_\Omega\langle\Phi,dD_H\chi_E\rangle
=
-\sum_{j=1}^m\int_\Omega \Phi_j\,dD_{X_j}\chi_E
=
\sum_{j=1}^m\int_\Omega \Phi_j(\nu_E)_j\,dP_H(E).
\]
Therefore
\[
-\int_\Omega\langle\Phi,dD_H\chi_E\rangle
=
\int_\Omega\langle\Phi,\nu_E\rangle\,dP_H(E),
\]
and the first equality in \eqref{eq:app-gauss-green-theorem} follows.

The second equality follows from the concentration property
\eqref{eq:app-perimeter-concentrated-reduced-boundary}, applied through
\eqref{eq:app-integral-reduced-boundary-only} with
\(\psi=\langle\Phi,\nu_E\rangle\).  This is the compactly supported version of
the Gauss--Green theorem in stratified groups used in the proof of Proposition
3.1; see \cite[Theorem 1.5]{CVA}.
\end{proof}

\subsection{Smooth boundaries and regular level surfaces}\label{subsec:app-smooth-boundaries}

The next proposition explains why, for smooth level surfaces, the horizontal
outer normal is obtained by projecting the Euclidean normal onto the horizontal
bundle.

\begin{prop}[Smooth Euclidean boundary]\label{prop:app-smooth-boundary-horizontal-normal}
Let \(E\Subset\Omega\subset\mathbb G=\mathbb R^N\) have \(C^1\) Euclidean
boundary, and let \(\nu\) be the Euclidean outward unit normal to \(\partial E\).
Define
\begin{equation}\label{eq:app-smooth-horizontal-normal-vector}
N_H
:=
\bigl(\langle X_1,\nu\rangle,\ldots,\langle X_m,\nu\rangle\bigr).
\end{equation}
Then, as \(\mathbb R^m\)-valued Radon measures in \(\Omega\),
\begin{equation}\label{eq:app-smooth-boundary-DHchi}
D_H\chi_E
=
-N_H\,\mathcal H^{N-1}\!\lfloor_{\partial E},
\end{equation}
and hence
\begin{equation}\label{eq:app-smooth-boundary-perimeter}
P_H(E;\cdot)
=
|N_H|\,\mathcal H^{N-1}\!\lfloor_{\partial E}.
\end{equation}
Here \(\mathcal H^{N-1}\!\lfloor_{\partial E}\) denotes the restriction of
\((N-1)\)-dimensional Hausdorff measure to \(\partial E\).  Thus, for example,
the measure in \eqref{eq:app-smooth-boundary-perimeter} assigns to a Borel set
\(A\subset\Omega\) the value
\[
\int_{A\cap\partial E}|N_H|\,d\mathcal H^{N-1}.
\]
The same restriction notation is used in
Proposition~\ref{prop:app-regular-sublevel-normal}, specifically in
\eqref{eq:app-DHchi-Et} and \eqref{eq:app-PH-Et-smooth-level}, with
\(\Sigma_t\) in place of \(\partial E\).
On the non-characteristic part \(\partial E\cap\{|N_H|>0\}\),
\begin{equation}\label{eq:app-smooth-boundary-unit-normal}
\nu_E=\frac{N_H}{|N_H|}
\qquad P_H(E)\text{-a.e.}
\end{equation}
Moreover,
\begin{equation}\label{eq:app-smooth-boundary-characteristic-zero}
P_H(E;\partial E\cap\{|N_H|=0\})=0.
\end{equation}
\end{prop}

\begin{proof}
Let \(\Phi=(\Phi_1,\ldots,\Phi_m)\in C_c^1(\Omega;\mathbb R^m)\).  By the
definition of horizontal divergence,
\[
\int_E\operatorname{div}_H\Phi\,dx
=
\sum_{j=1}^m\int_E X_j\Phi_j\,dx.
\]
By the divergence-free property proved in
Subsection~\ref{subsec:app-formal-adjoints},
\(\operatorname{div}X_j=0\) in the Euclidean-coordinate sense.  Thus, writing
\[
X_j=\sum_{\alpha=1}^N a_{\alpha j}(x)\frac{\partial}{\partial x_\alpha},
\]
the divergence-free condition is
\[
\operatorname{div}X_j=
\sum_{\alpha=1}^N\frac{\partial a_{\alpha j}}{\partial x_\alpha}=0.
\]
Therefore the Euclidean product rule gives
\[
\operatorname{div}(\Phi_jX_j)
=
\sum_{\alpha=1}^N
\frac{\partial}{\partial x_\alpha}\bigl(\Phi_j a_{\alpha j}\bigr)
=
X_j\Phi_j+
\Phi_j\operatorname{div}X_j
=
X_j\Phi_j.
\]
Equivalently,
\[
X_j\Phi_j=\operatorname{div}(\Phi_jX_j).
\]
The classical Euclidean divergence theorem therefore gives
\[
\int_E X_j\Phi_j\,dx
=
\int_{\partial E}\Phi_j\langle X_j,\nu\rangle\,d\mathcal H^{N-1}.
\]
Summing in \(j\),
\begin{equation}\label{eq:app-smooth-boundary-divergence-computation}
\int_E\operatorname{div}_H\Phi\,dx
=
\int_{\partial E}\langle \Phi,N_H\rangle\,d\mathcal H^{N-1}.
\end{equation}
On the other hand, by \eqref{eq:app-bvh-distributional-definition},
\[
\int_E\operatorname{div}_H\Phi\,dx
=
-\int_\Omega\langle \Phi,dD_H\chi_E\rangle.
\]
Comparing this identity with
\eqref{eq:app-smooth-boundary-divergence-computation} for all test fields
\(\Phi\) gives \eqref{eq:app-smooth-boundary-DHchi}.  Taking total variations
gives \eqref{eq:app-smooth-boundary-perimeter}.  The polar decomposition
\eqref{eq:app-polar-decomposition} then yields
\eqref{eq:app-smooth-boundary-unit-normal} on the set \(|N_H|>0\).  Finally,
\eqref{eq:app-smooth-boundary-perimeter} gives
\eqref{eq:app-smooth-boundary-characteristic-zero}, because the density
\(|N_H|\) vanishes on \(\{|N_H|=0\}\).
\end{proof}

\subsection{Regular sublevel sets}\label{mdetail}

The preceding proposition becomes the exact normal identity used for regular
sublevel sets in \cite{GutierrezMohammedSmooth}.  In
this subsection, \(B_R=B(e,R)=\{x\in\mathbb G:d(e,x)<R\}\) denotes the open ball
for the fixed left-invariant homogeneous distance \(d\), viewed in the
Euclidean coordinates \(\mathbb G\simeq\mathbb R^N\).  All boundaries and
perimeters in Proposition~\ref{prop:app-regular-sublevel-normal} are therefore
relative to the open set \(B_R\), unless explicitly stated otherwise.  Thus
``Euclidean'' critical and regular values refer to the ordinary Euclidean
gradient in these coordinates, not to the metric \(d\).
\begin{defn}[Euclidean critical and regular values]\label{def:app-euclidean-regular-value}\itshape
Let \(u\in C^1(B_R)\). A number \(s\in\mathbb R\) is called a Euclidean
critical value of \(u\) in \(B_R\) if it is the value of \(u\) at some
Euclidean critical point of \(u\), that is, if
\begin{equation}\label{eq:app-euclidean-critical-values}
s\in
\mathcal C_u(B_R)
:=
\{\tau\in\mathbb R:\exists x\in B_R\text{ such that }u(x)=\tau
\text{ and }\nabla u(x)=0\}.
\end{equation}
A number \(t\in\mathbb R\) is called a Euclidean regular value of \(u\) in
\(B_R\) if it is not a Euclidean critical value, namely if
\(t\notin\mathcal C_u(B_R)\). Equivalently,
\begin{equation}\label{eq:app-euclidean-regular-value}
\nabla u(x)\ne0
\qquad\text{for every }x\in B_R\text{ such that }u(x)=t.
\end{equation}
Thus, for each Euclidean regular value \(t\), the level set
\(\{u=t\}\cap B_R\) contains no Euclidean critical point of \(u\).
\end{defn}

\begin{prop}[Regular sublevel sets]\label{prop:app-regular-sublevel-normal}
Let \(u\in C^1(B_R)\), let \(t\) be a Euclidean regular value of \(u\) in the
sense of Definition~\ref{def:app-euclidean-regular-value}, equivalently
\eqref{eq:app-euclidean-regular-value}, and set
\[
E_t:=\{x\in B_R:u(x)<t\},
\qquad
\Sigma_t:=\{x\in B_R:u(x)=t\}.
\]
Then, relative to the open ball \(B_R\),
\begin{equation}\label{eq:app-levelsurface-boundary}
\partial_{B_R} E_t:=\partial E_t\cap B_R=\Sigma_t.
\end{equation}
If \(E_t=\emptyset\), then \(\Sigma_t=\emptyset\).  If \(E_t=B_R\), then the
interior level surface \(\{u=t\}\cap B_R\) is empty, while the ambient boundary
of \(E_t\) as a subset of \(\mathbb G\) is \(\partial B_R\); equivalently, with
the alternative convention \(\Sigma_t:=\partial_{\mathbb G}E_t\), this last case
gives \(\Sigma_t=\partial B_R\).  The identities below are relative to the open
set \(B_R\), so this ball-boundary piece is not included.  If one worked
instead with the closed ball \(\overline{B_R}\) as a set in \(\mathbb G\), then
the boundary portion on \(\partial B_R\) would have to be added separately to
the ambient boundary and to the perimeter.  Whenever the
interior surface \(\Sigma_t\) is nonempty, it is a \(C^1\) hypersurface.
The Euclidean outward unit normal to \(E_t\) on \(\Sigma_t\) is
\begin{equation}\label{eq:app-levelsurface-euclidean-normal}
\nu=\frac{\nabla u}{|\nabla u|}.
\end{equation}
Consequently, as \(\mathbb R^m\)-valued Radon measures in \(B_R\),
\begin{equation}\label{eq:app-DHchi-Et}
D_H\chi_{E_t}
=
-
\frac{Xu}{|\nabla u|}
\,\mathcal H^{N-1}\!\lfloor_{\Sigma_t},
\end{equation}
and
\begin{equation}\label{eq:app-PH-Et-smooth-level}
P_H(E_t;\cdot)
=
|D_H\chi_{E_t}|(\cdot)
=
\frac{|Xu|}{|\nabla u|}
\,\mathcal H^{N-1}\!\lfloor_{\Sigma_t}.
\end{equation}
Equivalently, the last measure sends a Borel set \(A\subset B_R\) to
\[
\int_{A\cap\Sigma_t}\frac{|Xu|}{|\nabla u|}\,d\mathcal H^{N-1}.
\]
In particular,
\begin{equation}\label{eq:app-characteristic-zero-Et}
P_H(E_t;\Sigma_t\cap\{|Xu|=0\})=0.
\end{equation}
On the non-characteristic part of the level surface,
\begin{equation}\label{eq:app-horizontal-normal-Et}
\nu_{E_t}
=
\frac{Xu}{|Xu|}
\qquad
P_H(E_t)\text{-a.e. on }\partial_H^*E_t\cap B_R.
\end{equation}
Equivalently, modulo \(P_H(E_t)\)-null sets,
\begin{equation}\label{eq:app-reduced-boundary-levelset-description}
\partial_H^*E_t\cap B_R
=
\{x\in B_R:u(x)=t,\ |Xu(x)|>0\}.
\end{equation}
\end{prop}

\begin{proof}
If \(E_t=\emptyset\), the asserted identities are immediate.  If \(E_t=B_R\),
then there is no interior level surface in the open ball \(B_R\); the only
ambient boundary is \(\partial B_R\), and that boundary is outside the relative
boundary in \(B_R\) and is not seen by test fields compactly supported in
\(B_R\).  Thus the asserted identities are again immediate as Radon-measure
identities in \(B_R\).  This is exactly the point at which the open-ball
convention matters: the proof records only the interior level-surface
contribution.  With a closed-ball convention, an additional boundary
contribution supported on \(\partial B_R\) would appear and would have to be
tracked separately.
We may therefore assume that \(\Sigma_t\ne\emptyset\).

Since \(t\) is a Euclidean regular value, \(\nabla u\ne0\) on \(\Sigma_t\).  The
implicit function theorem gives that \(\Sigma_t\) is a \(C^1\) hypersurface.
Because \(E_t=\{u<t\}\) is open in \(B_R\), its relative boundary in \(B_R\) is
exactly \(\Sigma_t\), proving \eqref{eq:app-levelsurface-boundary}.  The
exterior of \(E_t\) in \(B_R\) corresponds to larger values of \(u\), so the
Euclidean outer normal is \eqref{eq:app-levelsurface-euclidean-normal}.

For this set, the horizontal projection of the Euclidean normal in
\eqref{eq:app-smooth-horizontal-normal-vector} is
\[
N_H
=
\biggl(\frac{X_1u}{|\nabla u|},\ldots,
\frac{X_mu}{|\nabla u|}\biggr)
=
\frac{Xu}{|\nabla u|}.
\]
Applying Proposition~\ref{prop:app-smooth-boundary-horizontal-normal} gives
\eqref{eq:app-DHchi-Et} and \eqref{eq:app-PH-Et-smooth-level}.  Formula
\eqref{eq:app-characteristic-zero-Et} follows immediately from
\eqref{eq:app-PH-Et-smooth-level}, since the density \(|Xu|/|\nabla u|\) is
zero on \(\Sigma_t\cap\{|Xu|=0\}\).  On \(\{|Xu|>0\}\), the polar decomposition
\eqref{eq:app-polar-decomposition} and \eqref{eq:app-DHchi-Et} give
\[
-\nu_{E_t}\,P_H(E_t;\cdot)
=
D_H\chi_{E_t}
=
-\frac{Xu}{|Xu|}\,P_H(E_t;\cdot),
\]
which proves \eqref{eq:app-horizontal-normal-Et}.  Since the characteristic
part has zero \(P_H(E_t)\)-measure, the reduced boundary may be represented by
\eqref{eq:app-reduced-boundary-levelset-description} for the purposes of all
surface integrals with respect to \(P_H(E_t)\).
\end{proof}

\subsection{Horizontal coarea and almost-everywhere finite perimeter}\label{subsec:app-horizontal-coarea-formula}

\begin{prop}[Horizontal coarea formula]\label{prop:app-horizontal-coarea}
Let \(u\in BV_H(\Omega)\).  Then
\begin{equation}\label{eq:app-horizontal-coarea-bv}
|D_Hu|(\Omega)
=
\int_{-\infty}^{\infty}P_H(\{u>t\};\Omega)\,dt.
\end{equation}
More generally, for every non-negative Borel function \(g\),
\begin{equation}\label{eq:app-horizontal-coarea-weighted-bv}
\int_\Omega g\,d|D_Hu|
=
\int_{-\infty}^{\infty}
\int_{\partial_H^*(\{u>t\})}g\,dP_H(\{u>t\})\,dt.
\end{equation}
If \(u\in W^{1,1}_H(\Omega)\), then \(d|D_Hu|=|Xu|\,dx\), and therefore
\begin{equation}\label{eq:app-horizontal-coarea-sobolev}
\int_\Omega g(x)|Xu(x)|\,dx
=
\int_{-\infty}^{\infty}
\int_{\partial_H^*(\{u>t\})}g\,dP_H(\{u>t\})\,dt.
\end{equation}
The same formulas hold with the sublevel sets
\(\{u<t\}:=\{x\in\Omega:u(x)<t\}\) in place of the superlevel sets
\(\{u>t\}:=\{x\in\Omega:u(x)>t\}\), up to the usual null set of levels.
\end{prop}

\begin{proof}
We first justify the open-set identity used below.  If \(U\subset\Omega\) is
open and \(v\in L^1_{\mathrm{loc}}(U)\), set
\[
\mathcal V_H(v;U):=
\sup\left\{
\left.
\int_U v\,\operatorname{div}_H\Phi\,dx
\,\right|\,
\Phi\in C_c^1(U;\mathbb R^m),\ |\Phi|\le 1
\right\}.
\]
When \(v\in BV_H(\Omega)\), this variational quantity is exactly the total
variation of the restricted horizontal derivative:
\[
\mathcal V_H(v;U)=|D_Hv|(U).
\]
Indeed, the inequality \(\mathcal V_H(v;U)\le |D_Hv|(U)\) follows immediately
from \eqref{eq:app-bvh-distributional-definition}.  Conversely, by the polar
decomposition \(D_Hv=\sigma_v|D_Hv|\), with \(|\sigma_v|=1\)
\(|D_Hv|\)-a.e., and by the regularity of \(|D_Hv|\), one first restricts to
compact sets \(K\Subset U\), then approximates \(-\sigma_v\chi_K\) in
\(L^1(|D_Hv|)\) by smooth compactly supported horizontal vector fields with
norm at most one.  Testing with these fields and letting \(K\uparrow U\) gives
the reverse inequality.

We now prove the coarea identity for \(\mathcal V_H\) on \(U\).  First suppose
that \(w\in C^\infty(U)\).  For a.e. \(t\), Sard's theorem and
Proposition~\ref{prop:app-regular-sublevel-normal} give
\[
\mathcal V_H(\chi_{\{w>t\}};U)=P_H(\{w>t\};U)
=
\int_{U\cap\{w=t\}}\frac{|Xw|}{|\nabla w|}\,d\mathcal H^{N-1},
\]
where the quotient is taken to be zero on \(\{|\nabla w|=0\}\).  The Euclidean
coarea formula applied to \(w\) then yields
\[
\int_{-\infty}^{\infty}\mathcal V_H(\chi_{\{w>t\}};U)\,dt
=
\int_U |Xw|\,dx
=
\mathcal V_H(w;U).
\]

Let now \(u\in BV_H(\Omega)\).  Choose a standard strict horizontal BV
approximation \(u_k\in C^\infty(U)\) such that \(u_k\to u\) in
\(L^1_{\mathrm{loc}}(U)\) and
\[
\int_U |Xu_k|\,dx\longrightarrow \mathcal V_H(u;U).
\]
For a.e. \(t\), the level set \(\{u=t\}\) has Lebesgue measure zero; after
passing to a subsequence if necessary, this gives
\(\chi_{\{u_k>t\}}\to\chi_{\{u>t\}}\) in \(L^1_{\mathrm{loc}}(U)\).  By the lower
semicontinuity of horizontal variation and Fatou's lemma,
\[
\int_{-\infty}^{\infty}\mathcal V_H(\chi_{\{u>t\}};U)\,dt
\le
\liminf_{k\to\infty}
\int_{-\infty}^{\infty}\mathcal V_H(\chi_{\{u_k>t\}};U)\,dt
=
\mathcal V_H(u;U),
\]
where the last equality follows from the smooth case and the strict
convergence of \(u_k\).

For the opposite inequality, let \(\Phi\in C_c^1(U;\mathbb R^m)\) with
\(|\Phi|\le1\).  For \(N>0\), set
\[
u_N:=\max\{-N,\min\{u,N\}\}.
\]
The layer-cake identity
\[
u_N=-N+\int_{-N}^N\chi_{\{u>t\}}\,dt
\]
and the fact that \(\int_U\operatorname{div}_H\Phi\,dx=0\) give
\[
\int_U u_N\operatorname{div}_H\Phi\,dx
=
\int_{-N}^N\int_{\{u>t\}}\operatorname{div}_H\Phi\,dx\,dt.
\]
Hence
\[
\left|\int_U u_N\operatorname{div}_H\Phi\,dx\right|
\le
\int_{-N}^N\mathcal V_H(\chi_{\{u>t\}};U)\,dt.
\]
Letting \(N\to\infty\) and then taking the supremum over all such \(\Phi\) gives
\[
\mathcal V_H(u;U)
\le
\int_{-\infty}^{\infty}\mathcal V_H(\chi_{\{u>t\}};U)\,dt.
\]
Combining the two inequalities proves
\[
\mathcal V_H(u;U)
=
\int_{-\infty}^{\infty}\mathcal V_H(\chi_{\{u>t\}};U)\,dt.
\]
For a.e. \(t\), the set \(\{u>t\}\) has locally finite horizontal perimeter in
\(U\), and by Definition~\ref{def:app-horizontal-perimeter} together with
\eqref{eq:app-perimeter-variational-formula},
\(\mathcal V_H(\chi_{\{u>t\}};U)=P_H(\{u>t\};U)\).  Changing the integrand on the
exceptional null set of levels does not change the integral.  Therefore,
\begin{equation}\label{eq:app-coarea-open-sets-proof}
|D_Hu|(U)
=
\int_{-\infty}^{\infty}P_H(\{u>t\};U)\,dt.
\end{equation}
Taking \(U=\Omega\) gives \eqref{eq:app-horizontal-coarea-bv}.  In particular,
if \(U\Subset\Omega\), then \(P_H(\{u>t\};U)<\infty\) for a.e. \(t\).  A
countable exhaustion \(U_1\Subset U_2\Subset\cdots\Subset\Omega\) therefore
gives a full-measure set of levels \(\mathcal T\subset\mathbb R\) such that
\(\{u>t\}\) has locally finite horizontal perimeter in \(\Omega\) for every
\(t\in\mathcal T\).  For levels outside \(\mathcal T\) we do not need to choose a
perimeter measure, since they do not affect the one-dimensional integral in
\(t\).

We next derive the weighted form.  Define, for Borel sets \(A\subset\Omega\),
\[
\mu(A):=\int_{\mathcal T}P_H(\{u>t\};A)\,dt.
\]
Since \(\mathbb R\setminus\mathcal T\) has one-dimensional measure zero, these
integrals may equivalently be written over all \(t\in\mathbb R\), with arbitrary
values assigned on the exceptional levels.
By Tonelli's theorem, \(\mu\) is a positive Borel measure: for disjoint Borel
sets \(A_k\), countable additivity follows from
\[
P_H(\{u>t\};\cup_k A_k)=\sum_k P_H(\{u>t\};A_k)
\]
for each fixed \(t\), and then integration in \(t\).  Formula
\eqref{eq:app-coarea-open-sets-proof} gives \(\mu(U)=|D_Hu|(U)\) on every open
\(U\subset\Omega\).  Since both \(\mu\) and \(|D_Hu|\) are Radon measures, equality
on open sets implies equality on all Borel sets.  Hence
\begin{equation}\label{eq:app-coarea-measure-identity-proof}
|D_Hu|(A)
=
\int_{-\infty}^{\infty}P_H(\{u>t\};A)\,dt
\qquad\text{for every Borel }A\subset\Omega.
\end{equation}

If \(g=\sum_{k=1}^N a_k\chi_{A_k}\) is a non-negative simple Borel function,
then \eqref{eq:app-coarea-measure-identity-proof} gives
\[
\int_\Omega g\,d|D_Hu|
=
\sum_{k=1}^N a_k |D_Hu|(A_k)
=
\int_{-\infty}^{\infty}
\sum_{k=1}^N a_k P_H(\{u>t\};A_k)\,dt.
\]
The last expression is
\[
\int_{-\infty}^{\infty}
\int_\Omega g\,dP_H(\{u>t\})\,dt.
\]
For a general non-negative Borel function \(g\), choose simple functions
\(g_k\uparrow g\).  Applying the preceding identity to \(g_k\) and using the
monotone convergence theorem on both sides yields
\begin{equation}\label{eq:app-coarea-weighted-before-reduced-boundary}
\int_\Omega g\,d|D_Hu|
=
\int_{-\infty}^{\infty}
\int_\Omega g\,dP_H(\{u>t\})\,dt.
\end{equation}
Finally, for a.e. \(t\) the perimeter measure of \(\{u>t\}\) is concentrated on
its horizontal reduced boundary, by
\eqref{eq:app-perimeter-concentrated-reduced-boundary}.  Therefore the inner
integral in \eqref{eq:app-coarea-weighted-before-reduced-boundary} can be
written as an integral over \(\partial_H^*(\{u>t\})\), which proves
\eqref{eq:app-horizontal-coarea-weighted-bv}.

Assume now that \(u\in W^{1,1}_H(\Omega)\).  Let
\(u_k\in C^\infty(\Omega)\) be a standard horizontal Sobolev approximation on
compact subsets, so that \(u_k\to u\) in \(L^1_{\mathrm{loc}}(\Omega)\) and
\(Xu_k\to Xu\) in \(L^1_{\mathrm{loc}}(\Omega;\mathbb R^m)\).  For every
\(\Phi\in C_c^1(\Omega;\mathbb R^m)\), the smooth identity
\eqref{eq:app-smooth-bv-gradient} gives
\[
\int_\Omega u_k\operatorname{div}_H\Phi\,dx
=
-\int_\Omega\langle\Phi,Xu_k\rangle\,dx.
\]
Passing to the limit gives
\[
\int_\Omega u\operatorname{div}_H\Phi\,dx
=
-\int_\Omega\langle\Phi,Xu\rangle\,dx.
\]
Thus \(D_Hu=Xu\,dx\), and hence \(|D_Hu|=|Xu|\,dx\).  Substituting this identity
into \eqref{eq:app-horizontal-coarea-weighted-bv} gives
\eqref{eq:app-horizontal-coarea-sobolev}.

The formulas with sublevel sets follow by applying the already proved
superlevel formulas to \(-u\).  Indeed,
\(D_H(-u)=-D_Hu\), so \(|D_H(-u)|=|D_Hu|\), and
\(\{-u>s\}=\{u<-s\}\).  The change of variables \(t=-s\) gives the same
identities with \(\{u<t\}\) in place of \(\{u>t\}\).  As usual, all statements
about individual level sets are understood after discarding the null set of
levels for which the relevant perimeter may fail to be locally finite.
\end{proof}

\begin{prop}[Almost-everywhere finite perimeter of sublevel sets]\label{prop:app-levelsets-finite-perimeter}
Let \(u\in W^{1,1}_{H,\mathrm{loc}}(\Omega)\).  Then, for every compact set
\(K\Subset\Omega\),
\begin{equation}\label{eq:app-levelset-finite-perimeter-integral}
\int_{-\infty}^{\infty}P_H(\{u<t\};K)\,dt
\le
\int_K |Xu|\,dx
<\infty.
\end{equation}
In particular, for a.e. \(t\in\mathbb R\), the set \(\{u<t\}\) has locally
finite horizontal perimeter in \(\Omega\).  If \(\mathcal O\Subset\Omega\) has
finite horizontal perimeter in \(\Omega\), then, for a.e. \(t\),
\begin{equation}\label{eq:app-intersection-finite-perimeter}
\{u<t\}\cap\mathcal O
\quad\text{has finite horizontal perimeter in }\Omega.
\end{equation}
\end{prop}

\begin{proof}
Apply Proposition~\ref{prop:app-horizontal-coarea} to the restriction of \(u\)
to an open neighborhood of \(K\), and use \(g=\chi_K\).  This gives
\eqref{eq:app-levelset-finite-perimeter-integral}; hence the set of levels for
which \(P_H(\{u<t\};K)=\infty\) has one-dimensional Lebesgue measure zero.  A
countable exhaustion \(K_1\Subset K_2\Subset\cdots\Subset\Omega\) then gives a
single null set of levels outside of which \(P_H(\{u<t\};K_j)<\infty\) for every
\(j\), which is precisely local finite perimeter in \(\Omega\).

For the intersection statement, use the standard BV product rule for
characteristic functions: the intersection of two locally finite perimeter
sets has locally finite perimeter, and on compact subsets the perimeter is
controlled by the sum of the two perimeter measures.  Since \(\mathcal O\) has
finite perimeter and \(\{u<t\}\) has locally finite perimeter for the levels
just chosen, \eqref{eq:app-intersection-finite-perimeter} follows.
\end{proof}

\subsection{Gauge balls and the auxiliary coarea estimate}\label{subsec:app-gauge-coarea}

Let \(d\) be the fixed homogeneous distance on \(\mathbb G\), and let
\(N(x)=d(e,x)\) be its associated homogeneous norm, as in
\eqref{eq:distance-norm-compatibility}.  Thus, by
\eqref{eq:balls-from-homogeneous-norm},
\[
B_s:=B(e,s)=\{x\in\mathbb G:N(x)<s\}.
\]
The gauge \(N\) is locally Lipschitz, hence
\(N\in W^{1,\infty}_{H,\mathrm{loc}}(\mathbb G)\) and
\(|XN|\in L^\infty_{\mathrm{loc}}(\mathbb G)\).

\begin{prop}[Finite perimeter of gauge balls and coarea on annuli]\label{prop:app-gauge-balls}
For every compact annulus \(\Omega\cap\{a<N<b\}\Subset\mathbb G\),
\begin{equation}\label{eq:app-gauge-ball-finite-perimeter-ae}
\int_a^b P_H(B_s;\Omega)\,ds
=
\int_{\Omega\cap\{a<N<b\}}|XN|\,dx
<\infty.
\end{equation}
Consequently, \(B_s\) has finite horizontal perimeter in \(\Omega\) for a.e.
\(s\in(a,b)\).  Moreover, for every non-negative Borel function \(g\),
\begin{equation}\label{eq:app-gauge-coarea-formula}
\int_{\{a<N<b\}}g(x)|XN(x)|\,dx
=
\int_a^b\int_{\partial_H^*B_s}g\,dP_H(B_s)\,ds.
\end{equation}
In particular, if \(E_t\) is measurable and \(u\in W^{1,1}_{H,\mathrm{loc}}\),
then
\begin{equation}\label{eq:app-gauge-coarea-Et-Xu}
\int_a^b\int_{E_t^{(1)}\cap\partial_H^*B_s}|Xu|\,dP_H(B_s)\,ds
=
\int_{E_t\cap\{a<N<b\}}|Xu|\,|XN|\,dx,
\end{equation}
and therefore, on compact annuli,
\begin{equation}\label{eq:app-gauge-coarea-Et-estimate}
\int_a^b\int_{E_t^{(1)}\cap\partial_H^*B_s}|Xu|\,dP_H(B_s)\,ds
\le
C\int_{E_t\cap B_b}|Xu|\,dx,
\end{equation}
where \(C=\|XN\|_{L^\infty(B_b)}\).
\end{prop}

\begin{proof}
Apply Proposition~\ref{prop:app-horizontal-coarea} to the locally Lipschitz
function \(N\).  With \(g=\chi_\Omega\chi_{\{a<N<b\}}\), the coarea formula gives
\eqref{eq:app-gauge-ball-finite-perimeter-ae}; the finiteness follows from
\(|XN|\in L^\infty_{\mathrm{loc}}\) and the compactness of the annulus.  This
also shows that \(P_H(B_s;\Omega)<\infty\) for a.e. \(s\in(a,b)\).
The weighted formula \eqref{eq:app-gauge-coarea-formula} is the same coarea
identity with a general non-negative Borel weight \(g\).  Taking
\(g=|Xu|\chi_{E_t}\) gives \eqref{eq:app-gauge-coarea-Et-Xu}; the measure
\(P_H(B_s)\) is supported on \(\partial_H^*B_s\), and \(\chi_{E_t}\) is read in
the measure-theoretic sense as \(\chi_{E_t^{(1)}}\) on that boundary.  The
estimate \eqref{eq:app-gauge-coarea-Et-estimate} follows from the local bound
\(|XN|\le C\) on \(B_b\).
\end{proof}

\subsection{Admissible sets for the smooth level-surface estimate}\label{subsec:app-prop31-admissible-sets}

The following statement collects exactly what is needed in the proof of the
level-surface estimate in \cite{GutierrezMohammedSmooth}.

\begin{prop}[Admissible level sets for Gauss--Green]\label{prop:app-prop31-admissible-levels}
Assume \(u\in C^\infty(B_2)\), and set
\[
E_t:=\{x\in B_2:u(x)<t\}.
\]
Then there is a full-measure set of levels \(\mathcal T\subset(0,\infty)\) such
that, for every \(t\in\mathcal T\), the following hold.
\begin{enumerate}[(i)]
\item \(E_t\) has locally finite horizontal perimeter in \(B_2\).
\item \(t\) is a Euclidean regular value of \(u\), so Proposition~\ref{prop:app-regular-sublevel-normal} applies locally in every ball \(B_R\Subset B_2\).
\item If \(Z\in C_c^1(B_2;\mathbb R^m)\), then
\begin{equation}\label{eq:app-prop31-GG-ready}
\int_{E_t}\operatorname{div}_H Z\,dx
=
\int_{\partial_H^*E_t}\langle Z,\nu_{E_t}\rangle\,dP_H(E_t).
\end{equation}
\item If one also intersects with gauge balls, then \(E_t\cap B_s\) has finite horizontal perimeter in \(B_2\) for a.e. \(s\in(0,2)\).
\end{enumerate}
\end{prop}

\begin{proof}
Since \(u\in C^\infty(B_2)\), it belongs to
\(W^{1,1}_{H,\mathrm{loc}}(B_2)\).  Proposition~\ref{prop:app-levelsets-finite-perimeter}
therefore gives local finite horizontal perimeter of \(E_t\) for a.e.
\(t>0\).  Sard's theorem gives another full-measure set of levels for which
\(t\) is a Euclidean regular value of \(u\).  Let \(\mathcal T\) be the
intersection of these two full-measure sets.  For \(t\in\mathcal T\), (i) and
(ii) are immediate.

For (iii), apply Proposition~\ref{prop:app-gauss-green-finite-perimeter} in any
open set compactly containing \(\operatorname{spt}Z\) and compactly contained in
\(B_2\).  Since \(Z\) vanishes near the boundary of that open set, no artificial
boundary term is created.  This is why the cutoff fields in Proposition 3.1
produce only the level-surface contribution.

For (iv), Proposition~\ref{prop:app-gauge-balls} gives finite horizontal
perimeter of \(B_s\) for a.e. \(s\).  For those \(s\), the intersection of the
two finite-perimeter sets \(E_t\) and \(B_s\) again has finite perimeter by the
BV product rule for characteristic functions.
\end{proof}

\begingroup
\ifJulinPaperTwo\color{black}\else\color{black}\fi

\ifJulinPaperTwo
\section{Extension from smooth to weak solutions}\label{sec:nonsmooth-extension}
\else
\section*{Reduction of the problem to one with smooth ingredients}
\fi

\begin{prop}[Passage to measurable coefficients and a continuous
nonlinearity]
Assume that the generalized Harnack estimate has been proved when the
coefficient matrix and the nonlinearity are smooth, with a constant depending
only on \(\lambda,\Lambda\), and \(Q\), and not on any derivatives or moduli of
continuity of the coefficients.  Then the same estimate holds when
\[
A\in L^\infty(B_2;\operatorname{Sym}_m),
\qquad
\lambda I\le A(x)\le\Lambda I\quad\text{a.e.},
\]
and \(f:[0,\infty)\to[0,\infty)\) is continuous and non-decreasing.
Consequently, under the additional hypothesis \({\bf(f-I)}\), the classical
Harnack inequality also holds in this nonsmooth setting, with exactly the same
structural dependence of the constant as in Theorem~\ref{gha-c}.
\end{prop}

\begin{proof}
Let \(u\in W^{1,2}(B_2,X)\) be a non-negative weak solution of
\begin{equation}\label{eq:final-approx-weak-equation}
\operatorname{div}_H(A(x)Xu)=f(u)
\qquad\text{in }B_2.
\end{equation}
We first pass the generalized estimate to \(u\), and then explain why the
classical estimate follows without any further smoothness approximation.

\medskip
\noindent\emph{1. Local regularity of the limiting weak solution.}
For every non-negative \(\phi\in C_c^1(B_2)\),
\[
-\int_{B_2}\langle AXu,X\phi\rangle
=\int_{B_2}f(u)\phi\ge0.
\]
Thus \(u\) is a weak subsolution of the homogeneous equation.  The local
De Giorgi--Moser estimate gives \(u\in L^\infty_{\rm loc}(B_2)\); see
\cite{MarchiDG} for the De Giorgi-class regularity theory associated with
H\"ormander vector fields.  Since
\(f\) is continuous and non-decreasing, \(f(u)\) is locally bounded.  The
interior oscillation estimate for linear equations with measurable uniformly
\(X\)-elliptic coefficients and bounded right-hand side therefore gives a
locally H\"older continuous representative of \(u\).  All infima and suprema
below refer to this representative.

\medskip
\noindent\emph{2. Approximation on a strictly interior ball.}
Fix
\[
0<r<s<1,
\qquad D_s:=B_{2s}.
\]
Extend \(A\) outside \(B_2\) by a fixed symmetric matrix whose spectrum is
contained in \([\lambda,\Lambda]\), and convolve the extension with a
non-negative smooth group approximate identity.  The details are as follows.

We describe this construction precisely.  Choose, for example,
\[
A_*:=\frac{\lambda+\Lambda}{2}I_m.
\]
After choosing a measurable representative of \(A\), define a matrix-valued
function on the whole group by
\[
\overline A(x):=
\begin{cases}
A(x),&x\in B_2,\\
A_*,&x\in\mathbb G\setminus B_2.
\end{cases}
\]
The values assigned on \(\partial B_2\), a null set for Haar measure, are
irrelevant.  The extension \(\overline A\) is measurable, essentially
bounded, and symmetric.  Moreover, for a.e. \(x\in\mathbb G\) and every
\(\xi\in\mathbb R^m\),
\begin{equation}\label{eq:final-extended-A-ellipticity}
\lambda|\xi|^2
\le
\langle\overline A(x)\xi,\xi\rangle
\le
\Lambda|\xi|^2.
\end{equation}

Let \(\rho\in C_c^\infty(\mathbb G)\) satisfy
\[
\rho\ge0,
\qquad
\int_{\mathbb G}\rho(y)\,dy=1,
\]
and, using the group dilations, set
\[
\rho_\epsilon(y)
:=
\epsilon^{-Q}\rho(\delta_{1/\epsilon}y),
\qquad \epsilon>0.
\]
The change of variables \(z=\delta_{1/\epsilon}y\), for which
\(dy=\epsilon^Qdz\), shows that
\(\int\rho_\epsilon=1\).  The supports of \(\rho_\epsilon\) shrink to the
identity as \(\epsilon\downarrow0\), so
\((\rho_\epsilon)_\epsilon\) is a non-negative group approximate identity.
Convolve each entry of \(\overline A\) and define
\begin{equation}\label{eq:final-A-mollification-definition}
A_\epsilon(x)
:=
(\rho_\epsilon*\overline A)(x)
=
\int_{\mathbb G}
\rho_\epsilon(y)\overline A(y^{-1}x)\,dy.
\end{equation}
Equivalently, after a change of variables, the nonsmooth factor can be kept
fixed and the smooth kernel differentiated.  It follows that every entry of
\(A_\epsilon\) is \(C^\infty\).  Since an integral of symmetric matrices is
symmetric, \(A_\epsilon(x)\) is symmetric.  Most importantly, positivity of
the kernel and \eqref{eq:final-extended-A-ellipticity} give, for every
\(x\in\mathbb G\) and \(\xi\in\mathbb R^m\),
\[
\begin{aligned}
\langle A_\epsilon(x)\xi,\xi\rangle
&=
\int_{\mathbb G}\rho_\epsilon(y)
\langle\overline A(y^{-1}x)\xi,\xi\rangle\,dy\\
&\ge
\lambda|\xi|^2\int_{\mathbb G}\rho_\epsilon(y)\,dy
=\lambda|\xi|^2,
\end{aligned}
\]
and the same calculation with the upper bound gives
\[
\langle A_\epsilon(x)\xi,\xi\rangle
\le\Lambda|\xi|^2.
\]
Thus convolution does not enlarge the ellipticity interval.

Finally, the approximate-identity theorem on homogeneous groups
\cite[Proposition~1.20(i),(iii)]{FollandSteinHardy} gives
\[
A_\epsilon\longrightarrow\overline A
\qquad\text{in }L^p_{\mathrm{loc}}(\mathbb G)
\quad\text{for every }1\le p<\infty.
\]
As explained in the appendix on approximate identities below, one may choose
a sequence \(\epsilon_j\downarrow0\) along which the convergence also holds
almost everywhere.  Consequently, the matrices
\[
A_j:=A_{\epsilon_j}
\]
are smooth and symmetric, satisfy the same ellipticity bounds,
and satisfy \(A_j(x)\to A(x)\) for a.e. \(x\in B_2\).  This is the sequence
used below.

Thus the resulting matrices satisfy
\begin{equation}\label{eq:final-approx-Aj}
\lambda I\le A_j(x)\le\Lambda I,
\qquad
A_j(x)\to A(x)\quad\text{for a.e. }x,
\end{equation}
and hence, by dominated convergence,
\begin{equation}\label{eq:final-approx-AjXu}
(A_j-A)Xu\longrightarrow0
\qquad\text{strongly in }L^2(D_s).
\end{equation}

Extend \(f\) continuously and monotonically to \(\mathbb R\) by
\[
\bar f(t):=
\begin{cases}
f(t),&t\ge0,\\
f(0),&t<0.
\end{cases}
\]
Convolving \(\bar f\) with a non-negative one-dimensional approximate
identity gives smooth, non-negative, non-decreasing functions
\(f_j:\mathbb R\to[0,\infty)\) satisfying
\begin{equation}\label{eq:final-approx-fj}
f_j\longrightarrow\bar f
\qquad\text{uniformly on compact intervals in }\mathbb R.
\end{equation}
If
\[
F_j(t):=\int_0^t f_j(\sigma)\,d\sigma,
\qquad
F(t):=\int_0^t f(\sigma)\,d\sigma,
\]
then
\begin{equation}\label{eq:final-approx-Fj}
F_j\longrightarrow\bar F
\qquad\text{uniformly on compact intervals in }\mathbb R,
\qquad
\bar F(t):=\int_0^t\bar f(\sigma)\,d\sigma,
\end{equation}
and \(\bar F=F\) on \([0,\infty)\).

For every \(j\), let \(u_j\) be the unique minimizer of
\[
\mathcal J_j(w)
=
\frac12\int_{D_s}\langle A_jXw,Xw\rangle
+\int_{D_s}F_j(w)
\]
over \(u+W_0^{1,2}(D_s,X)\).  Existence follows from the direct method:
the Dirichlet part is uniformly coercive on this affine space, \(F_j\) is
convex because \(f_j\) is non-decreasing, and the negative part of \(F_j\)
has at most linear growth because \(f_j\ge0\) is bounded on
\(( -\infty,0]\).  Strict convexity of the Dirichlet part and the horizontal
Poincar\'e inequality give uniqueness.

At this stage \(u_j\) is defined only as a member of the Sobolev affine space
\(u+W_0^{1,2}(D_s,X)\); no smoothness is assumed in the minimization problem.
This is important because the Sobolev affine space is closed under weak
convergence, whereas the corresponding class of smooth functions is not, and
the boundary datum \(u\) need not have a smooth trace.  The purpose of
regularizing \(A\) and \(f\) is instead to deduce, after the Euler equation has
been obtained, that the Sobolev minimizer has a smooth representative in the
interior of \(D_s\).

Let \(L_s:=\|u\|_{L^\infty(D_s)}\).  The competitor
\(\min\{u_j,L_s\}\) has the same Sobolev boundary data as \(u_j\).  Indeed,
the Sobolev boundary-order truncation theorem gives
\[
(u_j-L_s)_+\in W_0^{1,2}(D_s,X),
\]
because \(u_j-u\in W_0^{1,2}(D_s,X)\) and \(u\le L_s\).  Hence
\(\min\{u_j,L_s\}=u_j-(u_j-L_s)_+\) belongs to
\(u+W_0^{1,2}(D_s,X)\).  Neither its Dirichlet energy nor its potential
energy is larger.  Uniqueness gives
\begin{equation}\label{eq:final-approx-upper}
u_j\le L_s\qquad\text{a.e. in }D_s.
\end{equation}
In particular, \(0\le f_j(u_j)\le f_j(L_s)\), and the latter quantities are
uniformly bounded by \eqref{eq:final-approx-fj}.  First variation is therefore
legitimate.
More precisely, first take
\(\phi\in W_0^{1,2}(D_s,X)\cap L^\infty(D_s)\).  For every \(t\in\mathbb R\),
the perturbation \(u_j+t\phi\) belongs to the admissible affine space.  For
\(|t|\le1\),
\[
u_j+t\phi\le L_s+\|\phi\|_{L^\infty(D_s)},
\]
and the monotonicity of \(f_j\) gives an integrable bound for the difference
quotient of \(F_j(u_j+t\phi)\).  Dominated convergence therefore permits
differentiation at \(t=0\), and the minimality of \(u_j\) gives
\[
\left.\frac{d}{dt}\right|_{t=0}\mathcal J_j(u_j+t\phi)=0.
\]
Finally, truncating an arbitrary
\(\phi\in W_0^{1,2}(D_s,X)\) and passing to the limit, using
\(f_j(u_j)\in L^\infty(D_s)\), extends the identity to every Sobolev test
function.
Consequently,
\begin{equation}\label{eq:final-approx-uj-equation}
\int_{D_s}\langle A_jXu_j,X\phi\rangle
+\int_{D_s}f_j(u_j)\phi=0
\qquad
\text{for every }\phi\in W_0^{1,2}(D_s,X).
\end{equation}

\medskip
\noindent\emph{3. Strong and locally uniform convergence.}
Set \(z_j:=u_j-u\in W_0^{1,2}(D_s,X)\).  Subtracting the weak equation for
\(u\) from \eqref{eq:final-approx-uj-equation}, testing by \(z_j\), and
adding and subtracting \(f_j(u)\), gives
\[
\begin{aligned}
&\int_{D_s}\langle A_jXz_j,Xz_j\rangle
+\int_{D_s}\bigl(f_j(u_j)-f_j(u)\bigr)z_j\\
&\qquad=
\int_{D_s}\langle(A-A_j)Xu,Xz_j\rangle
+\int_{D_s}\bigl(f(u)-f_j(u)\bigr)z_j.
\end{aligned}
\]
The second integral on the left is non-negative.  Ellipticity, Cauchy--Schwarz,
and Poincar\'e's inequality consequently give
\[
\lambda\|Xz_j\|_{L^2(D_s)}
\le
\|(A-A_j)Xu\|_{L^2(D_s)}
+C_P\|f(u)-f_j(u)\|_{L^2(D_s)}.
\]
The right-hand side tends to zero by
\eqref{eq:final-approx-AjXu}, \eqref{eq:final-approx-fj}, and the boundedness
of \(u\) on \(D_s\).  Therefore
\begin{equation}\label{eq:final-approx-strong}
u_j\longrightarrow u
\qquad\text{strongly in }W^{1,2}(D_s,X).
\end{equation}

The right-hand sides \(f_j(u_j)\) are uniformly bounded on compact subsets by
\eqref{eq:final-approx-upper}.  The interior De Giorgi oscillation estimate,
together with the uniform \(L^2\) bound supplied by
\eqref{eq:final-approx-strong}, gives, for each \(K\Subset D_s\),
\[
\sup_j\|u_j\|_{C^\alpha(K)}<\infty
\]
for a structural \(\alpha\in(0,1)\).  Arzel\`a--Ascoli and the uniqueness of
the \(L^2\) limit in \eqref{eq:final-approx-strong} imply
\begin{equation}\label{eq:final-approx-uniform}
u_j\longrightarrow u
\qquad\text{locally uniformly in }D_s.
\end{equation}
We now justify the intended gain in regularity.  For fixed \(j\), the smooth
positive-definite matrix \(A_j\) has a smooth positive square root
\(B_j=A_j^{1/2}\).  If
\[
Y_{j,k}:=\sum_{\ell=1}^m (B_j)_{k\ell}X_\ell,
\qquad k=1,\dots,m,
\]
then the fields \(Y_{j,1},\dots,Y_{j,m}\) span the same horizontal distribution
as \(X_1,\dots,X_m\), and hence satisfy H\"ormander's bracket-generating
condition.  Moreover, the principal operator in
\eqref{eq:final-approx-uj-equation} is a smooth-coefficient H\"ormander
operator of sum-of-squares type, up to smooth first-order terms.  The local
subelliptic estimates of Rothschild--Stein \cite{RS} therefore apply.
See Proposition~\ref{prop:app-rothschild-stein-local-estimate} below.

We spell out the bootstrap because the right-hand side
\(-f_j(u_j)\) in \eqref{eq:final-approx-uj-equation} is initially only
H\"older continuous.  Write
\[
P_j:=X^*(A_jX).
\]
Equation~\eqref{eq:final-approx-uj-equation} is
\[
P_ju_j=-f_j(u_j)
\qquad\text{in }D_s.
\]
The preceding De Giorgi estimate gives
\(u_j\in C^\alpha_{\mathrm{loc}}(D_s)\).  Hence, on every compact subset and
for every \(0<s_0<\alpha\),
\[
u_j\in H^{s_0}_{\mathrm{loc}}(D_s),
\qquad
f_j(u_j)\in C^\alpha_{\mathrm{loc}}(D_s)
\subset H^{s_0}_{\mathrm{loc}}(D_s).
\]
Apply the localized estimate
\eqref{eq:app-rothschild-stein-higher-estimate} with \(P=P_j\).  The
subelliptic gain \(\delta>0\) is fixed on each compact subset, so
\[
u_j\in H^{s_0+\delta}_{\mathrm{loc}}(D_s).
\]
We now use the standard smooth-composition property of Sobolev spaces: if
\(t>0\), \(w\in H^t_{\mathrm{loc}}\cap L^\infty_{\mathrm{loc}}\), and
\(G\in C^\infty(\mathbb R)\), then
\(G(w)\in H^t_{\mathrm{loc}}\).  Since \(u_j\) is locally bounded and
\(f_j\) is smooth, this gives
\[
f_j(u_j)\in H^{s_0+\delta}_{\mathrm{loc}}(D_s).
\]
Reapplying the same localized estimate on a nested sequence of compact
subsets and iterating yields, for every integer \(n\ge1\),
\[
u_j,\ f_j(u_j)
\in H^{s_0+n\delta}_{\mathrm{loc}}(D_s).
\]
Thus \(u_j\in H^t_{\mathrm{loc}}(D_s)\) for every \(t>0\).  Euclidean Sobolev
embedding, applied in exponential coordinates, gives derivatives of every
order and therefore an ordinary smooth representative.  Consequently,
\begin{equation}\label{eq:final-approx-uj-smooth}
u_j\in C^\infty(D_s).
\end{equation}
This is an interior conclusion and does not assert smoothness up to
\(\partial D_s\); it is sufficient because \(\overline{B_{2r}}\Subset D_s\).

\medskip
\noindent\emph{4. Non-negative smooth approximants and the generalized
estimate.}
On \(\overline{B_{2r}}\), define
\[
\delta_j:=\max\{0,-\min_{\overline{B_{2r}}}u_j\},
\qquad
v_j:=u_j+\delta_j.
\]
By \eqref{eq:final-approx-uniform} and \(u\ge0\),
\(\delta_j\to0\).  Thus \(v_j\ge0\), \(v_j\in C^\infty(B_{2r})\), and
\[
\operatorname{div}_H(A_jXv_j)=\widehat f_j(v_j),
\qquad
\widehat f_j(t):=f_j(t-\delta_j).
\]
The functions \(\widehat f_j\) are smooth, non-negative, and non-decreasing.
Their primitives satisfy
\[
\widehat F_j(t)
:=\int_0^t\widehat f_j(\sigma)\,d\sigma
=F_j(t-\delta_j)-F_j(-\delta_j)
\longrightarrow F(t)
\]
uniformly for \(t\) in compact subsets of \([0,\infty)\), by
\eqref{eq:final-approx-Fj}, \(\delta_j\to0\), and \(\bar F(0)=0\).

Rescale \(B_{2r}\) to \(B_2\).  The rescaled nonlinearity is
\(r^2\widehat f_j\), whose primitive is \(r^2\widehat F_j\), while the
rescaled coefficient matrix retains the ellipticity constants
\(\lambda,\Lambda\).  The smooth generalized Harnack estimate therefore gives
\begin{equation}\label{eq:final-approx-smooth-estimate}
\int_{m_{j,r}}^{M_{j,r}}
\frac{dt}{r\sqrt{\widehat F_j(t)}+t}
\le C,
\end{equation}
where
\[
m_{j,r}:=\inf_{B_r}v_j,
\qquad
M_{j,r}:=\sup_{B_r}v_j.
\]
Because \(r<1\), \eqref{eq:final-approx-smooth-estimate} implies
\[
\int_{m_{j,r}}^{M_{j,r}}
\frac{dt}{\sqrt{\widehat F_j(t)}+t}
\le C.
\]
Local uniform convergence gives
\[
m_{j,r}\to m_r:=\inf_{B_r}u,
\qquad
M_{j,r}\to M_r:=\sup_{B_r}u.
\]
\begingroup

Here is the precise argument.  Since \(r<s\), one has
\[
\overline{B_r}\Subset D_s=B_{2s}.
\]
Therefore \eqref{eq:final-approx-uniform} gives
\[
\|u_j-u\|_{L^\infty(\overline{B_r})}\longrightarrow0.
\]
Moreover, \(v_j=u_j+\delta_j\) and \(\delta_j\to0\), so
\begin{equation}\label{eq:final-approx-vj-uniform-Br}
\|v_j-u\|_{L^\infty(\overline{B_r})}
\le
\|u_j-u\|_{L^\infty(\overline{B_r})}+\delta_j
\longrightarrow0.
\end{equation}
For any two bounded real-valued functions \(p,q\) on a set \(E\),
\[
\left|\sup_Ep-\sup_Eq\right|
\le\|p-q\|_{L^\infty(E)},
\qquad
\left|\inf_Ep-\inf_Eq\right|
\le\|p-q\|_{L^\infty(E)}.
\]
Indeed, \(p\le q+\|p-q\|_\infty\) and
\(q\le p+\|p-q\|_\infty\), and the corresponding infimum inequalities
follow in the same way.  Applying these estimates to \(p=v_j\), \(q=u\),
and \(E=\overline{B_r}\), and using
\eqref{eq:final-approx-vj-uniform-Br}, gives
\[
\left|M_{j,r}-M_r\right|
\le\|v_j-u\|_{L^\infty(\overline{B_r})}\to0,
\qquad
\left|m_{j,r}-m_r\right|
\le\|v_j-u\|_{L^\infty(\overline{B_r})}\to0.
\]
Because \(u\) and \(v_j\) are continuous and \(B_r\) is dense in
\(\overline{B_r}\), their suprema and infima over \(B_r\) agree with those
over \(\overline{B_r}\).  This proves the asserted convergence of the
endpoints.
\endgroup
For \(m_r<a<b<M_r\), the interval \([a,b]\) lies inside
\([m_{j,r},M_{j,r}]\) for all large \(j\).
\begingroup

Here is the precise endpoint argument.  The strict inequalities
\(m_r<a\) and \(b<M_r\) leave a positive distance between \([a,b]\) and
the endpoints of \([m_r,M_r]\).  Set
\[
\eta:=\frac12\min\{a-m_r,M_r-b\}>0.
\]
Since \(m_{j,r}\to m_r\) and \(M_{j,r}\to M_r\), there is an index
\(j_0\) such that, whenever \(j\ge j_0\),
\[
|m_{j,r}-m_r|<\eta,
\qquad
|M_{j,r}-M_r|<\eta.
\]
It follows that
\[
m_{j,r}<m_r+\eta<a,
\qquad
M_{j,r}>M_r-\eta>b.
\]
Therefore every \(t\in[a,b]\) satisfies
\(m_{j,r}\le t\le M_{j,r}\), which is exactly
\[
[a,b]\subset[m_{j,r},M_{j,r}]
\qquad\text{for all }j\ge j_0.
\]
Notice also that \(m_r\ge0\), because \(u\ge0\); hence
\(a>m_r\) implies \(a>0\).  Thus the denominator in the integrand is
uniformly separated from zero on \([a,b]\).  More explicitly, uniform
convergence \(\widehat F_j\to F\) on \([a,b]\), together with
\(t\ge a>0\), gives
\[
\frac{1}{\sqrt{\widehat F_j(t)}+t}
\longrightarrow
\frac{1}{\sqrt{F(t)}+t}
\qquad\text{uniformly for }t\in[a,b].
\]
Since the integrands are non-negative and
\([a,b]\subset[m_{j,r},M_{j,r}]\), estimate
\eqref{eq:final-approx-smooth-estimate} (with the subsequent removal of the
factor \(r\)) gives
\[
\int_a^b
\frac{dt}{\sqrt{\widehat F_j(t)}+t}
\le C.
\]
Passing to the limit proves the next displayed inequality.
\endgroup
Consequently,
\[
\int_a^b\frac{dt}{\sqrt{F(t)}+t}\le C.
\]
Letting \(a\downarrow m_r\) and \(b\uparrow M_r\), by monotone convergence,
and then letting \(r\uparrow1\), again by monotone convergence, proves
\[
\int_{\inf_{B_1}u}^{\sup_{B_1}u}
\frac{dt}{\sqrt{F(t)}+t}
\le C.
\]
At every stage \(C\) depends only on \(\lambda,\Lambda,Q\); in particular,
it is independent of \(j\), of derivatives of \(A_j,f_j\), and of any
modulus of continuity.  This proves the generalized Harnack inequality for
measurable \(A\) and continuous \(f\).

\medskip
\noindent\emph{5. Passage of the classical Harnack inequality.}
Assume now that \(f\) satisfies \({\bf(f-I)}\), with parameters
\(1<\varrho<\theta\).  No new approximation of \(f\) is required.  The
horizontal local Keller--Osserman estimate,
Theorem~\ref{thm:blue-horizontal-local-scale}, is already a weak-solution
statement for measurable uniformly elliptic coefficients.  Applied to \(u\)
in \(B_2\), it gives, with
\[
M:=\sup_{B_1}u,
\qquad m:=\inf_{B_1}u,
\]
\begin{equation}\label{eq:final-approx-M-scale}
\frac{\sqrt{F(M/2)}}{M}\le C
\end{equation}
when \(M>0\).  The constant now depends only on
\(\lambda,\Lambda,Q,\theta,\varrho\).

For \(0<t\le M/4\), monotonicity gives
\[
F(t)\le t f(t),
\qquad
\frac{\sqrt{F(t)}}{t}
\le \left(\frac{f(t)}{t}\right)^{1/2}
\le \left(\frac{f(M/4)}{M/4}\right)^{1/2}.
\]
Also
\[
F(M/2)-F(M/4)
\ge\frac{M}{4}f(M/4),
\]
and hence \eqref{eq:final-approx-M-scale} implies
\begin{equation}\label{eq:final-approx-integrand-control}
\frac{\sqrt{F(t)}}{t}\le C,
\qquad 0<t\le M/4.
\end{equation}
If \(m=0<M\), then \eqref{eq:final-approx-integrand-control} would give
\[
\int_0^{M/4}\frac{dt}{\sqrt{F(t)}+t}
\ge c\int_0^{M/4}\frac{dt}{t}=\infty,
\]
contradicting the generalized Harnack estimate just proved.  Thus either
\(M=0\), in which case the conclusion is immediate, or \(m>0\).  In the
latter case,
\[
c\log\frac{M}{m}
\le
\int_{m/4}^{M/4}\frac{dt}{\sqrt{F(t)}+t}
\le
\int_{m/4}^{m}\frac{dt}{t}
+\int_m^M\frac{dt}{\sqrt{F(t)}+t}
\le \log4+C.
\]
Therefore \(\log(M/m)\le C\), and exponentiation yields
\[
\sup_{B_1}u=M\le C m=C\inf_{B_1}u.
\]
This proves the classical Harnack inequality for measurable \(A\) and
continuous \(f\), with no dependence on their smoothness.
\end{proof}

\endgroup

\begingroup
\ifJulinPaperTwo\color{black}\else\color{black}\fi

\ifJulinPaperTwo
\section{Approximate identities on Carnot groups}
\label{sec:approximate-identities}
\else
\section*{Appendix: Approximate identities on Carnot groups}
\fi

We record the precise approximation result used in the passage from smooth to
measurable coefficients.  Let \(\mathbb G\) be a Carnot group, let \(dy\)
denote Haar measure, let \(Q\) be the homogeneous dimension, and let
\((\delta_r)_{r>0}\) be the group dilations.  We use the convolution
convention
\[
(k*g)(x)
:=
\int_{\mathbb G}k(y)g(y^{-1}x)\,dy.
\]

\begin{defn}[Approximate identity]
A family \((k_\epsilon)_{\epsilon>0}\subset L^1(\mathbb G)\) is an
approximate identity if
\[
\sup_{\epsilon>0}\|k_\epsilon\|_{L^1}<\infty,
\qquad
\int_{\mathbb G}k_\epsilon(y)\,dy=1,
\]
and, for every neighborhood \(U\) of the identity element,
\[
\int_{\mathbb G\setminus U}|k_\epsilon(y)|\,dy
\longrightarrow0
\qquad\text{as }\epsilon\downarrow0.
\]
\end{defn}

Choose \(\rho\in C_c^\infty(\mathbb G)\) such that
\[
\rho\ge0,
\qquad
\int_{\mathbb G}\rho(y)\,dy=1,
\]
and define
\begin{equation}\label{eq:app-approximate-identity-kernel}
\rho_\epsilon(y)
:=
\epsilon^{-Q}\rho(\delta_{1/\epsilon}y).
\end{equation}
Since Haar measure scales according to
\(d(\delta_\epsilon y)=\epsilon^Qdy\), a change of variables gives
\[
\int_{\mathbb G}\rho_\epsilon=1,
\qquad
\|\rho_\epsilon\|_{L^1}=1.
\]
If \(\operatorname{supp}\rho\subset B(e,R)\), then
\[
\operatorname{supp}\rho_\epsilon
\subset B(e,\epsilon R).
\]
Thus the mass of \(\rho_\epsilon\) is concentrated in arbitrarily small
neighborhoods of the identity, and \((\rho_\epsilon)_\epsilon\) is an
approximate identity.

\begin{prop}[Approximation by convolution]
\label{prop:app-carnot-approximate-identity}
Let \(1\le p<\infty\).
\begin{enumerate}[(i)]
\item If \(g\in L^p(\mathbb G)\), then
\[
\rho_\epsilon*g\longrightarrow g,
\qquad
g*\rho_\epsilon\longrightarrow g
\quad\text{in }L^p(\mathbb G).
\]
\item If \(g\in L^p_{\mathrm{loc}}(\mathbb G)\), then both convergences
hold in \(L^p_{\mathrm{loc}}(\mathbb G)\).
\item If \(g\in L^1_{\mathrm{loc}}(\mathbb G)\), then
\(\rho_\epsilon*g\) and \(g*\rho_\epsilon\) are smooth wherever the
convolutions are defined.
\end{enumerate}
\end{prop}

\begin{proof}
Part (i), as well as the smoothness assertion for a smooth compactly supported
kernel, is the standard approximate-identity theorem for homogeneous groups;
see \cite[Proposition~1.20(i),(iii)]{FollandSteinHardy}.

We explain the localization in (ii).  Let \(K\Subset\mathbb G\).  Choose an
open set \(V\) with \(K\Subset V\Subset\mathbb G\) and a cutoff
\(\chi\in C_c^\infty(V)\) satisfying \(\chi=1\) in a neighborhood of
\(K\).  Then \(\chi g\in L^p(\mathbb G)\).  For all sufficiently small
\(\epsilon\), the support condition in
\eqref{eq:app-approximate-identity-kernel} implies that the convolutions of
\(g\) and \(\chi g\) agree on \(K\).  Applying part (i) to \(\chi g\) gives
\[
\|\rho_\epsilon*g-g\|_{L^p(K)}  \longrightarrow0,
\]
and the right-convolution statement follows in the same way.

For completeness, smoothness can also be seen directly.  After the change of
variables \(z=y^{-1}x\), unimodularity of a Carnot group gives
\[
(\rho_\epsilon*g)(x)
=
\int_{\mathbb G}\rho_\epsilon(xz^{-1})g(z)\,dz.
\]
On compact subsets, derivatives in \(x\) may be passed under the integral,
because the differentiated kernel is smooth and compactly supported while
\(g\) is locally integrable.  Hence \(\rho_\epsilon*g\in C^\infty\).  The
argument for \(g*\rho_\epsilon\) is analogous.
\end{proof}

We next state carefully how the theorem yields almost-everywhere convergence.
The theorem itself supplies \(L^p_{\mathrm{loc}}\) convergence.  Fix a bounded
ball \(B_3\) containing \(\overline{B_2}\).  From
\[
\rho_\epsilon*g\longrightarrow g
\qquad\text{in }L^2(B_3),
\]
one can choose \(\epsilon_j\downarrow0\) so that
\[
\|\rho_{\epsilon_j}*g-g\|_{L^2(B_3)}^2\le2^{-j}.
\]
Consequently,
\[
\int_{B_3}
\sum_{j=1}^\infty
|\rho_{\epsilon_j}*g-g|^2\,dy
\le
\sum_{j=1}^\infty2^{-j}<\infty.
\]
Tonelli's theorem implies that the series inside the integral is finite for
a.e. \(y\in B_3\); in particular,
\begin{equation}\label{eq:app-approximate-identity-ae-subsequence}
(\rho_{\epsilon_j}*g)(y)\longrightarrow g(y)
\qquad\text{for a.e. }y\in B_3.
\end{equation}
Thus no pointwise convergence assertion for the entire family is needed.

Finally, apply the preceding facts componentwise to the bounded measurable
matrix \(\overline A\).  Define
\[
A_j:=\rho_{\epsilon_j}*\overline A.
\]
Then \(A_j\in C^\infty\), and
\eqref{eq:app-approximate-identity-ae-subsequence} gives
\[
A_j(x)\longrightarrow\overline A(x)
\qquad\text{for a.e. }x\in B_2.
\]
If \(\overline A\) is symmetric and satisfies
\[
\lambda I\le\overline A(x)\le\Lambda I,
\]
then positivity and unit mass of \(\rho_{\epsilon_j}\) show that the same
inequalities hold for every \(A_j\).  In particular,
\[
|A_j-A|\le2\Lambda
\qquad\text{a.e. in }B_2.
\]
Therefore, for \(Xu\in L^2(B_2)\),
\[
|(A_j-A)Xu|^2
\le4\Lambda^2|Xu|^2,
\]
and dominated convergence gives
\[
(A_j-A)Xu\longrightarrow0
\qquad\text{strongly in }L^2(B_2).
\]
This is precisely the convergence required in the stability argument.

\endgroup

\ifJulinPaperTwo
\color{black}
\section{Subelliptic regularity for the regularized operators}
\label{sec:regularized-subelliptic-regularity}
\else
{\color{black}
\fi
\begin{prop}[Rothschild--Stein local subelliptic estimate]
\label{prop:app-rothschild-stein-local-estimate}
Let \(U\subset\mathbb R^N\) be open, and let
\(Y_1,\dots,Y_m\) be real \(C^\infty\) vector fields on \(U\).  Assume that
there is an integer \(\kappa\ge1\) such that, at every point of \(U\), the
fields \(Y_1,\dots,Y_m\) together with their iterated commutators of length at
most \(\kappa\) span \(\mathbb R^N\).  Let
\[
P:=\sum_{k=1}^mY_k^*Y_k,
\]
where the adjoints are taken with respect to a fixed smooth positive density.
Then, for every open set \(V\Subset U\), there is a constant \(C_V\) such
that
\begin{equation}\label{eq:app-rothschild-stein-basic-estimate}
\|\varphi\|_{H^{1/\kappa}(\mathbb R^N)}^2
\le
C_V\left(
\sum_{k=1}^m\|Y_k\varphi\|_{L^2}^2
+\|\varphi\|_{L^2}^2
\right)
=
C_V\left(
(P\varphi,\varphi)_{L^2}
+\|\varphi\|_{L^2}^2
\right)
\end{equation}
for every \(\varphi\in C_c^\infty(V)\).  Equivalently, by duality,
\begin{equation}\label{eq:app-rothschild-stein-dual-estimate}
\|\varphi\|_{H^{1/\kappa}}
\le
C_V\left(
\|P\varphi\|_{H^{-1/\kappa}}
+\|\varphi\|_{L^2}
\right).
\end{equation}

More generally, let \(\chi,\widetilde\chi\in C_c^\infty(U)\), with
\(\widetilde\chi=1\) in a neighborhood of \(\operatorname{spt}\chi\).  There
is a fixed gain \(\delta>0\), depending on the vector fields and the chosen
compact sets but independent of \(s\), such that, for every
\(s\in\mathbb R\), there are constants \(C_s>0\) and an integer \(N_s\ge0\)
for which every
distribution \(v\in\mathcal D'(U)\) satisfies the localized estimate
\begin{equation}\label{eq:app-rothschild-stein-higher-estimate}
\|\chi v\|_{H^{s+\delta}}
\le
C_s\left(
\|\widetilde\chi Pv\|_{H^s}
+\|\widetilde\chi v\|_{H^{-N_s}}
\right),
\end{equation}
whenever the terms on the right are finite.  In particular,
\[
Pv\in C^\infty(U)
\quad\Longrightarrow\quad
v\in C^\infty(U).
\]
Thus \(P\) is locally \(C^\infty\)-hypoelliptic.  This is the local
subelliptic estimate and regularity theorem of Rothschild--Stein
\cite[Theorem~18]{RS}; the exponent \(1/\kappa\) in
\eqref{eq:app-rothschild-stein-basic-estimate} is the sharp gain associated
with commutators of length at most \(\kappa\).
\end{prop}

\begin{rem}[Application to the regularized operators]
For the matrices \(A_j\) used above, let \(B_j=A_j^{1/2}\) and
\[
Y_{j,k}=\sum_{\ell=1}^m(B_j)_{k\ell}X_\ell.
\]
Since \(B_j\) is smooth and invertible, the fields \(Y_{j,k}\) generate the
same horizontal distribution as the original Carnot-group fields \(X_\ell\).
They therefore satisfy the preceding bracket condition, and
\[
X^*(A_jX)=\sum_{k=1}^mY_{j,k}^*Y_{j,k}.
\]
Consequently Proposition~\ref{prop:app-rothschild-stein-local-estimate} is the
local estimate invoked in the bootstrap leading to
\eqref{eq:final-approx-uj-smooth}.
\end{rem}
\ifJulinPaperTwo\else
}
\fi

\end{document}